%% file: singpot-final.tex
\documentclass[10pt]{article}
\usepackage{amssymb,amsmath, amsthm, amsopn, amsfonts,bm,mathtools}
\usepackage[english]{babel}
\usepackage{color}\definecolor{rd}{rgb}{1,0.3,0.35}

\newtheorem{thm}{Theorem}[section]
\newtheorem{cor}[thm]{Corollary}
\newtheorem{lem}[thm]{Lemma}
\newtheorem{definition}[thm]{Definition}
\newtheorem{prop}[thm]{Proposition}
\def\proof{\noindent{\bf Proof. }}
\newtheorem{remark}[thm]{Remark}
\def\ds{\displaystyle}
\def\nm{\noalign{\medskip}}

\def\Z{{\mathcal{Z}}}
\def\O{{\mathcal{O}}}

\def\F{{\mathcal{F}}}
\def\P{{\mathcal{P}}}
\def\L{{\mathcal{L}}}
\def\H{{\mathcal{H}}}
\def\S{{\mathcal{S}}}
\def\D{{\mathcal{D}}}
\def\N{{\bf N}}
\newcommand{\field}[1]{\mathbb{#1}}
\newcommand{\rz}{\field{R}}
\newcommand{\cz}{\field{C}}
\newcommand{\nz}{\field{N}}

\def\11{{\rm 1~\hspace{-1.2ex}l} }
\def\d{{\rm{d}}}
\def\p{{\mathbb{P}}}
\def\fin{{$\hfill\square$\\}}
\def\Tr{{\rm{Tr}}}

\def\supp{{\rm supp~}}
\def\Real{{\mathrm{Re~}}}
\def\Imag{\mathrm{Im~}}
\begin{document}
\input{fig4tex}

\title{Mean field propagation of infinite dimensional Wigner measures with a singular
  two-body interaction potential}
\author{
Z.~Ammari\thanks{zied.ammari@univ-rennes1.fr,  IRMAR, Universit{\'e} de Rennes I, UMR-CNRS 6625,
campus de Beaulieu, 35042 Rennes Cedex, France.}\quad
F.~Nier\thanks{francis.nier@univ-rennes1.fr,  IRMAR, Universit{\'e} de Rennes I, UMR-CNRS 6625,
campus de Beaulieu, 35042 Rennes Cedex, France, and INRIA project team
MICMAC.}
}
\maketitle
\begin{abstract}
We consider the  quantum dynamics of many bosons systems in the mean field limit with
a singular pair-interaction potential, including the attractive or repulsive Coulombic case in
three dimensions. By using a measure transportation technique developed
in \cite{AGS}, we show that Wigner measures propagate along the nonlinear Hartree flow.
Such property was   previously proved only for bounded potentials in
our works \cite{AmNi2,AmNi3}
with a slightly different strategy.
\end{abstract}

\section{Introduction}
The evolution of a non relativistic system of many quantum  particles
is described by an $n$-body Schr{\"o}dinger equation.
The mean field limit consists in replacing this problem by a non
linear $1$-particle problem, by considering a one generic particle
interacting with the average field of all the others, when
the number of particles is large and the interaction potential is weak.
It is common knowledge that this
approximation starts to be very effective when the number of particles
exceeds a few tens.
In the last decades, many works have been devoted to justify this
limit. Most of them considered the mean field dynamics of well
prepared quantum states, coherent states or Hermite states, by
following and extending the phase-space approach, also known as the Hepp method (see
\cite{FGS,FKS,GiVe1,GiVe2,Hep,KnPi,RoSc}),
or by studying the BBGKY hierarchy of reduced density matrices
(see \cite{BGM,ChPa,ErYa,ESY1,KlMa,Spo}).
Some of these results deal with very singular pair
interaction potentials in \cite{BEGMY,ErYa,ESY1,KnPi} or considered
the rate of convergence
(see \cite{Ana, RoSc, KnPi}), sometimes motivated by the modelling
of Bose-Einstein condensates (see a.e.\cite{Aft,ESY2,LSSY}).
In this article, we continue our program, which consists in deriving
the mean field limit, for general initial data in the bosonic
framework.
Our strategy is inspired by older attempts to give substance to the
formal link between bosonic Quantum Field Theory and the finite
dimensional microlocal or phase-space analysis (see
\cite{Ber,Fol,Fol2,KrRa,Las}).
With this respect, the small parameter $\varepsilon=\frac{1}{n}$
asymptotics is the infinite dimensional version of semiclassical
analysis.
And it has been realized in the 90's, that the Wigner (or
semiclassical) measures provide a powerful tool in order to obtain the leading term
in the semiclassical limit (see
\cite{Ger,GMMP,HMR,LiPa}),  because they flexibly and
efficiently incorporate a priori estimates (see \cite{Bur,Bur2,FeGe,Mil,Nie}).

In \cite{AmNi1} Wigner measures were introduced in the infinite
dimensional setting and their main properties were studied. The above-mentioned work
exploited and clarified the intimate relationship between
pseudo-differential calculus, phase-space geometry and the probability
approach, inherent to bosonic QFT. In \cite{AmNi2}, the dynamics for
well prepared data and bounded interaction potentials was reconsidered
within this approach. The general propagation result was obtained in
\cite{AmNi3} for bounded interaction potentials. In particular, we
showed that the BBGKY hierarchy dynamics is a projected picture of
the evolution of the Wigner measure, for which there is a closed equation.
One difficulty which was solved in \cite{AmNi3} is concerned with the
integration of a weak Liouville equation valid after testing with
cylindrical or polynomial observables: Such classes of observables are
not preserved by the nonlinear Hamiltonian mean field flow. For
bounded interaction potentials, the number conservation allows
polynomial approximations of the nonlinear deformation in balls of the
phase-space.
This is done by adapting a truncated Dyson expansion approach presented
in \cite{FGS,FKP,FKS}.
In applications, an important case is the $2$-body Coulomb
interaction since it models the general non relativistic motion of
charged (or gravitational) particles. Again there are results about
the mean field problem for specific initial data (see \cite{BEGMY,
  KnPi}), but the approach we have followed in \cite{AmNi3} essentially
fails.
With a singular pair interaction potential, a solution to this problem
is provided by measure transportation techniques developed for
optimal transport theory (see \cite{AGS, Vil1}). Hence, the dynamical
mean field limit relies even more on the fact that Wigner measures are
\underline{probability measures} on the phase-space.

\bigskip
We now expose our main result.
The Hamiltonian of an $n$-body quantum system, with a pair interaction
potential, is given by the
Schr{\"o}dinger operator
$$
H_\varepsilon^{(n)}= \varepsilon\sum_{i=1}^n -\Delta_{x_i}+\varepsilon^2
\sum_{1\leq i<j\leq n} V(x_i-x_j)\,,
$$
where $\varepsilon$ is a positive parameter and $x_i,x_j\in\rz^{d}$\,.
We assume that the particles obey Bose statistics. So, we consider
$H_\varepsilon^{(n)}$ as an operator acting on the  space
$L_s^{2}(\mathbb{R}^{dn})$ of symmetric square integrable
functions. This means that
$$\Psi\in L_s^{2}(\mathbb{R}^{dn})\; \mbox{ iff } \; \Psi\in L^{2}(\mathbb{R}^{dn}) \;\mbox{ and }
\,\Psi(x_1,\cdots,x_n)=\Psi(x_{\sigma_1},\dots,x_{\sigma_n}) \mbox{ a.e}
$$
for any permutation $\sigma$ on the symmetric group $\mathfrak{S}_n$\,. The mean field asymptotics
is concerned with the limit as $\varepsilon\to 0$ and $n\varepsilon\to
1$\,, where $n=\left[\frac{1}{\varepsilon}\right]$
represents the number of particles of the system. \\
Let $\H$  be the direct sum of Hilbert spaces of the form
\[
\H=\bigoplus_{n=0}^\infty L^{2}_s(\rz^{dn})\,,
\]
and consider the Hamiltonian of the many-bosons system (with arbitrary number of particles) as
\begin{equation}
  \label{eq.hamq}
H_{\varepsilon}=  \bigoplus_{n=0}^\infty H_\varepsilon^{(n)}\,.
\end{equation}
An obvious feature of the operator $H_\varepsilon$ is the
conservation of the number of particles. Hence, it is
useful to define the number operator
$$
{\bf N}= \bigoplus_{n=0}^\infty \varepsilon n\,\11_{L^2_s(\rz^{dn})}.
$$
The free Hamiltonian, corresponding to $V=0$\,, will be denoted by
$H_{\varepsilon}^{0}$:
$$
H_{\varepsilon}^{0}=\bigoplus_{n=0}^{\infty}H_{\varepsilon}^{0,(n)}\,,\quad
H_{\varepsilon}^{0,(n)}=\varepsilon \sum_{i=1}^{n}-\Delta_{x_{i}}\,.
$$
Second quantization is a natural framework for the study of many-body
problems  and, even more, it helps to understand the mean field limit and the structures behind it.
However, the result can be presented without using the language of quantum field theory.
We just mention that
the operator $H_{\varepsilon}$ can be formally rewritten as
\begin{eqnarray*}
H_{\varepsilon}=\int_{\rz^{d}}\nabla a^{*}(x).\nabla a(x)~dx +
\frac{1}{2}\int_{\rz^{2d}}V(x-y)
a^{*}(x)a^{*}(y)a(x)a(y)~dx dy\,,
\end{eqnarray*}
with the $\varepsilon$-dependent canonical commutation relations
 $\left[a(x)\,,\, a^*(y)\right]=\varepsilon\delta(x-y)\,$\,.
It is interpreted  as  the Wick quantization of the classical Hamiltonian
\begin{eqnarray}
\label{eq.classical}
h(z,\bar z)= \int_{\rz^{d}}|\nabla z(x)|^2~dx
+\frac{1}{2}\int_{\rz^{2d}}|z(x)|^{2}|z(y)|^{2}V(x-y)~dx  dy\,.
\end{eqnarray}
In our analysis, an operator which violates the number of particles
conservation,
will play an important role,
namely the Weyl operator. Such operators are given for $f\in L^2(\rz^d)$  by
$$
W(f)=e^{\frac{i}{\sqrt{2}} [a^*(f)+a(f)]}\,,
$$
where $a^*(f), a(f)$ are the creation-annihilation operators on $\H$ satisfying
 the $\varepsilon$-canonical commutation relations (CCR):
\begin{eqnarray*}
[a(f_1),a^*(f_2)]=\varepsilon\langle f_1,f_2\rangle_{L^2(\rz^d)}\;\11, \;\;\;[a^*(f_1),a^*(f_2)]=0=[a(f_1),a(f_2)]\,.
\end{eqnarray*}
Accurate definitions on second quantized operators can be found in Appendix \ref{se.WWaW}.

\bigskip

Our approach is based on Wigner measures which are Borel probability
measures on the infinite dimensional phase-space $\Z_0:=L^2(\rz^d;\cz)$\,.
The states of the many-bosons system are positive trace-class operators
on $\H$ of normalized trace equal to $1$ (i.e.: normal states or density operators).
To every family of those states $(\varrho_\varepsilon )_{\varepsilon\in
(0,\bar \varepsilon)}$\,, we  asymptotically assign, when $\varepsilon \to 0$\,, at least one Borel probability
measure $\mu$ on $\Z_0:=L^2(\rz^d;\cz)$\,, called Wigner measure, such that there exists a
sequence $(\varepsilon_k)_{k\in\nz}$\,, such that $\lim_{k\to\infty} \varepsilon_k=0$ and
$$
\lim_{k\to 0}\Tr[\varrho_{\varepsilon_k} \,W(\sqrt{2}\pi \xi)]=\mathcal{F}^{-1}(\mu)(\xi)\,,
$$
under the sole uniform estimate
$\Tr\left[\varrho_{\varepsilon}{\bf N}^{\delta}\right]\leq C_{\delta}$
for some $\delta>0$\,. Here $\mathcal{F}^{-1}(\mu)$ is the inverse Fourier transform of $\mu$\,.\\
The problem of the mean field dynamics questions whether the
asymptotic quantities, namely Wigner measures,  as $\varepsilon\to 0$ associated with
$$
\varrho_{\varepsilon}(t)=
e^{-i\frac{t}{\varepsilon}H_{\varepsilon}}\varrho_{\varepsilon}
e^{i\frac{t}{\varepsilon}H_{\varepsilon}}\,,\quad t\in\rz
$$
are transported by the flow $\Phi(t,s)=\Phi(t-s)$ generated by the classical
Hamiltonian
$h(z,\bar z)$ and given, after writing $z_{t}=\Phi(t,s)(z_{s})$\,, by
\begin{equation}
  \label{eq.hartintro}
i\partial_{t}z_{t}= (\partial_{\bar
  z}h)(z_{t},\bar{z}_{t})=-\Delta z_{t}+
V*|z_{t}|^{2} z_{t}\,.
\end{equation}
After checking that the Hamiltonian \eqref{eq.hamq} has a self-adjoint
realization so that the quantum dynamics are well defined on
$\mathcal{H}$ and after checking that the mean field flow is well
defined on $\Z_{1}=H^{1}(\rz^{d})$\,, our main result is stated below.\\
Throughout the paper, we assume that the real valued potential $V$  satisfies the assumptions
\newtagform{Anombre}{(A}{)}
\usetagform{Anombre}
\addtocounter{equation}{-3}
\begin{eqnarray}
\label{eq.hypsym}
 && V(-x)=V(x)\in \rz\,, \\
\label{eq.hypbd}
&&V(1-\Delta)^{-1/2}\in {\cal L}(\Z_{0})\,,\\
\label{eq.hypcomp}
\text{and}
&&
(1-\Delta)^{-1/2}V(1-\Delta)^{-1/2}\in {\cal L}^{\infty}(\Z_{0})\,.
\end{eqnarray}
\usetagform{default}
\addtocounter{equation}{0}
We use the notation ${\cal L}({\mathfrak h})$ for the space of bounded
operators on the Hilbert space ${\mathfrak h}$  and ${\cal L}^{p}({\mathfrak
  h})$\,, $1\leq p\leq +\infty$\,, for the Schatten classes, ${\cal
  L}^{\infty}({\mathfrak h})$ being the space of compact operators for $p=+\infty$\,.

\begin{thm}
\label{th.main}
Let $(\varrho_{\varepsilon})_{\varepsilon\in(0,\bar\varepsilon)}$
be a family of normal states on $\H$
with a single Wigner measure $\mu_{0}$
such that the bound
\begin{equation}
\label{eq.hypmom}
{\rm Tr}[(\mathbf{N}+H_{\varepsilon}^{0})^{\delta}\varrho_{\varepsilon}
] \leq C_{\delta}< +\infty\,,
\end{equation}
holds uniformly w.r.t $\varepsilon\in (0,\bar\varepsilon)$ for some
$\delta>0$\,.\\
Then for all $t\in \rz$\,, the family
$(e^{-i\frac{t}{\varepsilon}H_{\varepsilon}}\varrho_{\varepsilon}e^{i\frac{t}{\varepsilon}
H_{\varepsilon}})_{\varepsilon\in
(0,\bar \varepsilon)}$
has a unique Wigner measure $\mu_{t}$ which is a Borel measure on
$\Z_{1}=H^{1}(\rz^d)$\,.
This measure $\mu_{t}=\Phi(t,0)_{*}\mu_{0}$ is the push forward of
the initial measure $\mu_{0}$  by the flow associated
with \eqref{eq.hartintro}, well defined on $\Z_{1}$\,.
\end{thm}
In a formal level the proof of the above theorem is rather simple. Writing first the integral formula
$$
\Tr[\varrho_\varepsilon(t) W(\xi)]=
\Tr[\varrho_\varepsilon  W(\xi)]+i \int_0^t
\Tr[ \varrho_\varepsilon(s) W(\xi) \sum_{j=1}^4 \varepsilon^{j-1} \O_j] \,ds\,,
$$
where $\varrho_\varepsilon(t)=e^{-it/\varepsilon H_\varepsilon}\varrho_\varepsilon e^{it/\varepsilon H_\varepsilon}$ and
$\O_j$ are some Wick quantized observables. By taking
the limit as $\varepsilon\to 0$\,,  the only term $j=1$
is left in the r.h.s. So, we formally end up
with a transport equation on the Wigner measures
$$
\partial_{t}\mu + i\left\{h\,,\, \mu\right\}=0,\quad
\left\{h,\mu\right\}=\partial_{z}h\,\partial_{\bar z}\mu
-\partial_{z}\mu\,\partial_{\bar z}h
$$
which is then solved by appealing to the results in \cite{AGS}.

\bigskip
\noindent\textbf{Outline:}
The self-adjointness of the
Hamiltonian
$H_\varepsilon$ and the existence
of a global flow on $\Z_{1}=H^1(\rz^d)$ for the Hartree equation \eqref{eq.hartintro} are proved in Section~\ref{se.dyn}.
The derivation of the mean field dynamics is done in Section \ref{se.dermeanfield} where
 Theorem \ref{th.main} is proved. Some additional properties are stated in Section~\ref{se.addppties}:
in particular, we draw the link with former results on bounded potential and reduced
density matrices  and provide non trivial examples elucidated by the
Wigner measure approach.
 The article ends with several appendices dedicated to second quantization,
absolutely continuous curves in $\mathrm{Prob}_{2}(\Z)$ as well as some
weak $L^p$ conditions for the potential $V$  ensuring the fulfillment
of the assumptions (A\ref{eq.hypbd}) and (A\ref{eq.hypcomp}).

\section{Well defined dynamics}
\label{se.dyn}

In this section we shall prove~that:
\begin{itemize}
\item the quantum dynamics is well defined, namely $H_{\varepsilon}$
  has a natural self-adjoint realization;
\item the mean field dynamics is well defined on
  $\Z_{1}=H^{1}(\rz^d)$\,, with additional useful estimates.
\end{itemize}

\subsection{Self-adjoint realization of $H_\varepsilon$}
\label{se.sareal}
The Hamiltonian $H_\varepsilon$ has a particular structure explained
in a general framework in Appendix \ref{se.comsa}. \\
Let $V$ be a real-valued Lebesgue measurable function a.e. finite and
satisfying the assumptions
(A\ref{eq.hypsym}) and (A\ref{eq.hypbd}). The multiplication operator
$$
V_\varepsilon^{(n)}=\varepsilon^2 \sum_{1\leq i<j\leq n} V(x_i-x_j)\,
$$
with its natural domain $\D(V_\varepsilon^{(n)})=\{\Psi\in L^2_s(\rz^{dn}): V_\varepsilon^{(n)}
\Psi\in L^2_s(\rz^{dn})\}$ is self-adjoint on $L^2_s(\rz^{dn})$ as well as
the differential operator
$$
H_\varepsilon^{0,(n)}=\varepsilon \sum_{i=1}^n -\Delta_{x_i}\,,\quad \mbox{ with } \quad
\D(H_\varepsilon^{0,(n)})=L^2_s(\rz^{dn})\cap H^2(\rz^{dn})\,.
$$
Therefore,  according to Appendix \ref{se.comsa}
$$
V_\varepsilon=\sum_{n=0}^\infty V_\varepsilon^{(n)}\,, \quad \mbox{ and } \quad
H^0_\varepsilon=
\sum_{n=0}^\infty H_\varepsilon^{0,(n)}\,
$$
endowed with their natural domains are self-adjoint  on $\H$\,.

\begin{prop}
Under the assumptions (A\ref{eq.hypsym}) and (A\ref{eq.hypbd}):\\
(i) The operator
$$
H_\varepsilon^{(n)}:= H_\varepsilon^{0,(n)}+ V_\varepsilon^{(n)}
$$
is self-adjoint on $\D( H_\varepsilon^{0,(n)})\subset\D( V_\varepsilon^{(n)})$\,.\\
(ii) The operator
$$
H_\varepsilon:=\sum_{n=0}^\infty H_{\varepsilon}^{(n)}, \quad
\D(H_\varepsilon):=\{\Psi\in\H\,, \sum_{n=0}^{\infty} \|H_\varepsilon^{(n)}\Psi^{(n)}\|^2<\infty\},
$$
is self-adjoint and essentially self-adjoint on $\oplus_{n\in\nz}^{alg} \D_n$ where $\D_n$ is any core of
$H_\varepsilon^{0,(n)}$\,.
\end{prop}
\proof
(i)
By assumption (A\ref{eq.hypbd}),
$V_\varepsilon^{(n)}$ is  infinitesimally small with respect to $H_\varepsilon^{0,(n)}$\,. So that,
$\D( H_\varepsilon^{0,(n)})\subset\D( V_\varepsilon^{(n)})$ and the operator $H_{\varepsilon}^{(n)}=H_\varepsilon^{0,(n)}+V_\varepsilon^{(n)}$ is self-adjoint on
the domain of $H_\varepsilon^{0,(n)}$ by Kato-Rellich theorem. \\
(ii) Applying Proposition \ref{graded}, we see  that
$H_\varepsilon$ is self-adjoint and essentially self-adjoint on
$\oplus_{n\in \nz}^{alg} \D_n$\,. 

\fin
\bigskip
\noindent
Later, it will be useful to use the reference operator
\begin{eqnarray}
\label{eq.opref}
S_\varepsilon(\lambda)=
\sum_{n=0}^\infty H_\varepsilon^{0,(n)}+\varepsilon n+\lambda (\varepsilon n)^3 \,.
\end{eqnarray}
which is self-adjoint by Proposition \ref{graded}. Moreover, by  functional calculus
 of strongly commuting self-adjoint operators we observe that
$\D(S_\varepsilon(\lambda))$ is invariant with respect to the parameter
$\lambda>0$\,.
\begin{prop}
\label{pr.sa}
Under the assumptions (A\ref{eq.hypsym}) and (A\ref{eq.hypbd}),
for any $\lambda>0$\,, the operator $V_\varepsilon$
is  $S_\varepsilon(\lambda)$-bounded with
$$
\forall \,\Psi\in \D(S_\varepsilon(\lambda)), \quad
\|V_\varepsilon\Psi\|_\H\leq \lambda \| V(1-\Delta)^{-1/2}\|_{\L(L^2(\rz^d))}
 \;\|S_\varepsilon(\lambda^{-2}) \Psi\|_\H
\,.
$$
Therefore $H_{\varepsilon}$ is essentially self-adjoint on
$\D(S_\varepsilon(\lambda))$\,.
\end{prop}
\proof
 The multiplication operator by $V(x_{1}-x_{2})$\,, at
 least defined as a symmetric operator from $\mathcal{S}(\rz^{2d})$
 into $\mathcal{S'}(\rz^{2d})$\,, satisfies
\begin{eqnarray}
\nonumber
e^{ix_{2}D_{x_{1}}}V(x_1-x_2) (1-\Delta_{x_{1}})^{-1/2}e^{-ix_{2}D_{x_{1}}}
&=&
e^{ix_{2}D_{x_{1}}} V(x_1-x_2) e^{-ix_{2}D_{x_{1}}}(1-\Delta_{x_{1}})^{-1/2}\\
\label{eq.V1D}
&=&V(x_{1})(1-\Delta_{x_{1}})^{-1/2}\in \mathcal{L}(L^{2}(\rz^{2d}))\,.
\end{eqnarray}
For $\Psi\in\H,\Phi\in \D(S_\varepsilon(\lambda))$\,, taking advantage of the symmetry of those wave functions,
we compute
\begin{eqnarray*}
  \langle \Psi\,,\,V_\varepsilon \Phi\rangle
&=&
\sum_{n=2}^\infty
  \langle \Psi^{(n)}\,,\, \frac{n(n-1)}{2}\varepsilon^{2} V(x_1-x_2) \Phi^{(n)}\rangle_{L_s^{2}(\rz^{dn})}\\
 &=& \sum_{n=2}^\infty\langle \Psi^{(n)}\,,\, \frac{n(n-1)}{2}\varepsilon^{2} V(x_1-x_2) (1-\Delta_{x_1})^{-1/2}
  (1-\Delta_{x_1})^{1/2}\Phi^{(n)}\rangle_{L_s^{2}(\rz^{dn})}\,.
\end{eqnarray*}
By noticing that
\begin{eqnarray*}
(n^2\varepsilon^2)^2 \|(1-\Delta_{x_{1}})^{1/2}\Phi^{(n)}\|_{L^{2}(\rz^{nd})}^{2}&=&
(\varepsilon n)^3 \langle \Phi^{(n)},
\varepsilon \sum_{i=1}^n  (1-\Delta_{x_i})\Phi^{(n)}\rangle_{L^2_s(\rz^{dn})}\\
&=&\|{\bf N}^{3/2} ({\bf N}+H^0_\varepsilon)^{1/2} \Phi^{(n)}\|_{\H}^2\,,
\end{eqnarray*}
the Cauchy-Schwarz inequality leads to
\begin{eqnarray}
\label{veps}
  | \langle \Psi\,,\, \,V_\varepsilon\Phi\rangle|\leq
  \|V(1-\Delta)^{-1/2}\|
~\|\Psi\|_{\H}
~\|{\bf N}^{3/2} \, ({\bf N}+H_\varepsilon^0)^{1/2}\Phi\|_{\H}\,.
\end{eqnarray}
Now with the inequality $ab\leq (\lambda a)^2+(b/\lambda)^2$\,, we see that
\begin{eqnarray}
\nonumber
\|{\bf N}^{3/2} \, ({\bf N}+H^0_\varepsilon)^{1/2}\Phi\|_{\H}^2&=&
\langle \Phi, {\bf N}^3 (\N+H^0_\varepsilon) \Phi\rangle\\ \nonumber
&\leq&
\langle \Phi, \lambda^{-2}{\bf N}^6+\lambda^2 (\N+H^0_\varepsilon)^2 \Phi\rangle\\
\label{NKest}
&\leq& \lambda^2 \|S_\varepsilon(\lambda^{-2}) \Phi\|_\H^2\,.
\end{eqnarray}
Putting together \eqref{veps} and \eqref{NKest} yields the estimate.\\
To prove the last statement, observe that $\oplus_{n\in\nz}^{alg}\D(H_\varepsilon^{0,(n)})$ is a core for $H_\varepsilon$\,.
Owing to the inclusions
$$
\oplus_{n\in\nz}^{alg}\D(H_\varepsilon^{0,(n)})\subset \D(S_\varepsilon(\lambda))
\subset \D(V_\varepsilon)\cap\D(H_\varepsilon^0)\subset \D(H_\epsilon)\,,
$$
$H_\epsilon$ is essentially self-adjoint on  $\D(S_\varepsilon(\lambda))$.
\fin

We end this section with some invariance properties of the domain $\D(S_\varepsilon(\lambda))$
with respect to the Hamiltonian $H_\varepsilon$ and the Weyl operators.
\begin{prop}
\label{invH}
For any $\lambda>0$ and $t\in\rz$
$$
e^{-i\frac{t}{\varepsilon} H_\varepsilon}  \D(S_\varepsilon(\lambda))\subset \D(S_\varepsilon(\lambda))\,.
$$
Moreover there exists $C_{\lambda}>0$ such that
$$
\|S_\varepsilon(\lambda) e^{-i\frac{t}{\varepsilon} H_\varepsilon}
(S_\varepsilon(\lambda)+1)^{-1}\|_{\L(\H)}\leq C_{\lambda}\,,
\quad \mbox{ for all } t\in \rz.
$$
\end{prop}
\proof
For any $\Psi\in \D(S_\varepsilon(\lambda))\subset \D(H_\varepsilon)$\,,
 observe that $e^{-i\frac{t}{\varepsilon} H_\varepsilon}
\Psi$ belongs to $\D(\N^3)\cap \D(H_\varepsilon)$  since
$\D(S_\varepsilon(\lambda))$ is contained in $\D(\N^3)$  and
$H_\varepsilon$ strongly commutes with $\N$\,. Proposition \ref{pr.sa} implies
$\D(S_\varepsilon(\beta)+V_\varepsilon)=\mathcal{D}(S_{\varepsilon}(\beta))=\mathcal{D}(S_{\varepsilon}(\lambda))$
when $\beta>0$ is large enough, so that for any $\Phi\in
\D(S_\varepsilon(\lambda))$
\begin{eqnarray*}
\langle (S_\varepsilon(\beta)+V_\varepsilon)\,\Phi, e^{-i\frac{t}{\varepsilon} H_\varepsilon}\Psi\rangle_\H &=&
\langle  (H_\varepsilon+\N+\beta\N^3)\,\Phi, e^{-i\frac{t}{\varepsilon} H_\varepsilon}\Psi\rangle_\H\\ &=&
\langle\Phi, (H_\varepsilon+\N+\beta\N^3) e^{-i\frac{t}{\varepsilon} H_\varepsilon}\Psi\rangle_\H
\end{eqnarray*}
and hence $e^{-i\frac{t}{\varepsilon} H_\varepsilon}\Psi$ belongs to $\D(
(S_\varepsilon(\beta)+V_\varepsilon)^*)=\D(S_\varepsilon(\beta))=\D(S_\varepsilon(\lambda))$\,.\\
Again for $\beta$ large enough
\begin{eqnarray*}
1+S_\varepsilon(\beta)+V_\varepsilon= (1+V_\varepsilon (S_\varepsilon(\beta{})+1)^{-1}) (S_\varepsilon(\beta{})+1)\,,
\end{eqnarray*}
and
$$
(1+\N+\beta\N^3+H_\varepsilon)^{-1}= (S_\varepsilon(\beta)+1)^{-1} (1+V_\varepsilon (S_\varepsilon(\beta)+1)^{-1})^{-1}\,.
$$
Therefore, the operators $ (1+S_\varepsilon(\beta)) (1+\N+\beta\N^3+H_\varepsilon)^{-1}$ and
 $(1+\N+\beta\N^3+H_\varepsilon)(S_{\varepsilon}(\beta)+1)^{-1}$ are
 bounded. Thus, we conclude that
 \begin{multline*}
   S_{\varepsilon}(\lambda)e^{-i\frac{t}{\varepsilon}H_{\varepsilon}}(1+S_{\varepsilon}(\lambda))^{-1}
=
S_{\varepsilon}(\lambda)(1+S_{\varepsilon}(\beta))^{-1}(1+S_{\varepsilon}(\beta))
(1+\N+\beta\N^3+H_\varepsilon)^{-1}\\
\circ e^{-i\frac{t}{\varepsilon}H_{\varepsilon}}
(1+\N+\beta\N^3+H_\varepsilon)(1+S_{\varepsilon}(\beta))^{-1}
(1+S_{\varepsilon}(\beta))(1+S_{\varepsilon}(\lambda))^{-1}\,,
 \end{multline*}
 is bounded.
\fin

\begin{prop}
\label{weylinv}
For any $\xi\in H^2(\rz^d)$ and any $\lambda>0$\,, the domain $\D(S_\varepsilon(\lambda))$  is invariant
under the action of the Weyl operator $W(\xi)$ with
$$
\|(S_\varepsilon(\lambda)+1)^{-1} W(\xi) S_\varepsilon(\lambda) \|_{\L(\H)}\leq C_{\lambda,\xi} \,,
$$
uniformly w.r.t $\varepsilon\in(0,\bar\varepsilon)$ for some constant $C_{\lambda,\xi}>0$\,.
\end{prop}
\proof
For all $\Phi,\Psi\in \D(S_{\varepsilon}(\lambda))$\,, one can write
\begin{eqnarray*}
\langle \Phi, W(\xi)^* S_\varepsilon(\lambda) W(\xi) \Psi\rangle=
\langle \Phi, (S_\varepsilon(\lambda)+Q_{\varepsilon}^{Wick}) \Psi\rangle,
\end{eqnarray*}
where $Q_{\varepsilon}$ is the following polynomial
$$
Q_{\varepsilon}(z)=\langle z+\frac{i\varepsilon}{\sqrt{2}}\xi, -\Delta
(z+\frac{i\varepsilon}{\sqrt{2}}\xi)\rangle_{\Z_0}
-
\langle z, -\Delta
z\rangle_{\Z_0}
+P_{\varepsilon}(z+\frac{i\varepsilon}{\sqrt{2}}\xi)-P_{\varepsilon}(z)
$$
and
$P_{\varepsilon}(z)=|z|_{\Z_{0}}^{6}+3\varepsilon|z|_{\Z_{0}}^{4}+\varepsilon^{2}|z|_{\Z_{0}}^{2}$
is the complete  Wick symbol of  $\N^3$\,, according to
Proposition~\ref{symbcalc} or by direct computation\,.
The assumption $\xi\in H^{2}(\rz^{d})$ ensures that $Q_{\varepsilon}$ is uniformly
bounded in $\oplus_{p+q\leq 3}\P_{p,q}(\Z_{0})$ and the number estimate of
Proposition \ref{pr.wick-estimate2}
says that
$
Q_{\varepsilon}^{Wick}\langle \N\rangle^{-\frac{3}{2}}$ is a bounded operator and therefore
$$
Q_{\varepsilon}^{Wick}(S_{\varepsilon}(\lambda)+1)^{-1}\in \mathcal{L}(\H)\,.
$$
Hence for $\Psi\in \mathcal{D}(S_{\varepsilon}(\lambda))$\,,
$$
S_\varepsilon(\lambda) W(\xi)\Psi=W(\xi)
\left[\frac{S_{\varepsilon}(\lambda)}{S_{\varepsilon}(\lambda)+1}+ Q_{\varepsilon}^{Wick}(S_{\varepsilon}(\lambda)+1)^{-1}\right](S_{\varepsilon}(\lambda)+1) \Psi\,
$$
and $W(\xi)\Psi$ belongs to $\D(S_\varepsilon(\lambda))$\,, with
$\|S_{\varepsilon}(\lambda)W(\xi)\Psi\|\leq C_{\lambda,\xi}\|(S_{\varepsilon}(\lambda)+1)\Psi\|$\,.
\fin

\begin{prop}
\label{pr.estimHdG}
  For any function $\chi\in \mathcal{C}^{\infty}_{0}(\rz^{2})$ and $\lambda>0$\,, the
  operator $\chi(\mathbf{N}, H_{\varepsilon})$ satisfies
$$
\forall k\in \nz, \quad
\|\mathbf{N}^{k} S_\varepsilon(\lambda)\chi(\mathbf{N},H_{\varepsilon})\|_{\mathcal{L}(\H)}\leq C_{\lambda,\chi}^{k+1}
$$
for some $C_{\lambda,\chi}>0$\,.
\end{prop}
\proof
 The operators $\mathbf{N}$\,,  $H_{\varepsilon}$ (like $\mathbf{N}$ and
 $\N+H_\varepsilon^0$) are strongly commuting
 self-adjoint operators so that the functional calculus is well
 defined for the pair $(\mathbf{N}, H_{\varepsilon})$\,. With a
 cut-off function $\chi_{1}\in \mathcal{C}^{\infty}_{0}(\rz)$ such
 that $\chi_{1}(x)\equiv 1$ on a neighborhood of $\supp \chi$\,, the
 operator
 $\N^{k} (1+\mathbf{N}+\beta\mathbf{N}^{3}+H_\varepsilon)\chi_{1}(\N)\chi(H_\varepsilon,\N)$
 is bounded with
$$
\|(1+\mathbf{N}+\beta \mathbf{N}^{3}+H_{\varepsilon})\N^{k}\chi_{1}(\N)\chi(H,\N)\|_{\mathcal{L}(\H)}
\leq C_{\beta}C_{\chi}^{k}\,.
$$
For sufficiently large $\beta$\,, Proposition~\ref{pr.sa} says
$$
\|(1+S_\varepsilon(\beta))
(1+\mathbf{N}+\beta\mathbf{N}^{3}+H_{\varepsilon})^{-1}\|_{\mathcal{L}(\H)}\leq
C_{\beta}'\,.
$$
This is done with
\begin{multline*}
\N^{k}S_{\varepsilon}(\lambda)\chi(\N,H_{\varepsilon})
=S_{\varepsilon}(\lambda)(1+S_{\varepsilon}(\beta))^{-1}(1+S_{\varepsilon}(\beta)
(1+\mathbf{N}+\beta\mathbf{N}^{3}+H_{\varepsilon})^{-1}
\\
\circ (1+\mathbf{N}+\beta\mathbf{N}^{3}+H_{\varepsilon})\mathbf{N}^{k}\chi(\N,H_{\varepsilon})\,,
\end{multline*}
and $C_{\lambda,\chi}=\max\left\{ C_{\chi},\;
  C_{\beta}C_{\beta}'\|S_{\varepsilon}(\lambda)(1+S_{\varepsilon}(\beta))^{-1}\|
\right\}$\,.
\fin
\subsubsection{Mean field dynamics}
\label{se.dynmean}
We shall use another more convenient writing of the Cauchy problem
\begin{equation}
  \label{eq.cauchy1}
  \left\{
 \begin{array}[c]{l}
      i\partial_{t}z_{t}=-\Delta z_{t}+ V*|z_{t}|^2 z_{t}\\
     z_{t=0}=z_{0}\,.
 \end{array}
\right.
\end{equation}
After setting $\tilde{z}_{t}=e^{it(-\Delta)}z_{t}=e^{-it\Delta}z_{t}$ it becomes
\begin{equation}
  \label{eq.cauchy2}
  \left\{
 \begin{array}[c]{l}
      i\partial_{t}\tilde z_{t}=e^{-it\Delta}\left[ V*|e^{it
        \Delta}\tilde{z}_{t}|^{2} (e^{it\Delta}\tilde{z}_{t})\right]\\
     \tilde{z}_{t=0}=z_{0}\,.
 \end{array}
\right.
\end{equation}
\begin{prop}
\label{pr.hartree}
Assume (A\ref{eq.hypsym}) and (A\ref{eq.hypbd}).
  For any $z_{0}\in \Z_{1}=H^{1}(\rz^{d})$ the Cauchy problem
  \eqref{eq.cauchy1} admits a unique solution $(t\mapsto z_{t})\in
  \mathcal{C}^{0}(\rz;H^{1}(\rz^{d}))\cap \mathcal{C}^{1}(\rz;H^{-1}(\rz^{d}))$\,.
More precisely, the Cauchy problem \eqref{eq.cauchy2}, which is
equivalent to \eqref{eq.cauchy1}, admits a unique solution in
$\mathcal{C}^{1}(\rz;H^{1}(\rz^{d}))$\,.
Moreover these solutions verify
\begin{eqnarray}
\label{eq.consL2}
  &&|z_{t}|_{L^{2}}=|\tilde{z}_{t}|_{L^{2}}=|z_{0}|_{L^{2}}\\
\label{eq.consen}
\text{and}
&&
h(z_{t},\overline{z_{t}})=
h(z_{0},\overline{z_{0}})\,,\\
\nonumber
\text{for}
&&
h(z,\overline{z})=
\int_{\rz^{d}}|\nabla z|^{2}(x)~dx+
\frac{1}{2}\int_{\rz^{2d}}V(x-y)|z(x)|^{2}|z(y)|^2~dxdy\,.
\end{eqnarray}
Finally,  the time-dependent velocity field defined on $\rz\times
  \Z_{1}$ by
$$
v(t,z)= e^{-it\Delta} ([V*|e^{it\Delta}z|^{2}]e^{it\Delta}z)
$$
satisfies the estimates
\begin{eqnarray}
  \label{eq.vtL2}
&&  |v(t,z)|_{\Z_{0}}\leq \|V(1-\Delta)^{-1/2}\|\;
|z|_{\Z_{0}}^{2}|z|_{\Z_{1}}\\
\label{eq.vtH1}
\text{and}&&
|v(t,z)|_{\Z_{1}}\leq \|V(1-\Delta)^{-1/2}\|\; |z|_{\Z_{1}}^{2}|z|_{\Z_{0}}\,.
\end{eqnarray}
\end{prop}
\proof
The first results are standard (see e.g. \cite{Caz,Gin})
in the analysis on nonlinear evolution
equation. Nevertheless, we recall the details of the proof because it
also contains \eqref{eq.vtL2}\eqref{eq.vtH1}, which is crucial in our analysis.\\
By considering the second formulation \eqref{eq.cauchy2}, it suffices
to prove that the mapping $z\to (V*|z|^{2})z$ is locally Lipschitz in $H^{1}(\rz^{d})$\,.
After noticing that the distributional derivative of $(V*|z|^{2})z$
 or more generally of $(V*(\overline{z_{1}}z_{2})z_{3})$ is
\begin{equation}
  \label{eq.multlinV}
\partial_{x}[(V*(\overline{z_{1}}z_{2}))z_{3}]=(V*(\partial_{x}\overline{z}_{1}z_{2}+
\overline{z}_1\partial_{x}z_{2}))z_{3}
+ (V*(\overline{z_{1}}z_{2}))(\partial_{x}z_{3})\,,
\end{equation}
it is reduced to the estimate of $V*(\overline{z_{1}}z_{2})z_{3}$ in $L^{2}$ in
terms of the $L^{2}$ and $H^{1}$- norms of $z_{1},z_{2},z_{3}$\,.
For $\xi\in L^{2}(\rz^{d})$\,, write
$$
\langle\xi\,,\,  (V*(\overline{z_{1}}z_{2}))z_{3}\rangle_{L^{2}(\rz^{d})}
=\langle z_{1}\otimes \xi\,,\, V(x_{1}-x_{2})z_{2}\otimes z_{3}\rangle_{L^{2}(\rz^{2d})}\,.
$$
When $\tilde{b}$ is the multiplication operator by $V(x_{1}-x_{2})$\,,
the estimate \eqref{eq.V1D} says that
$\tilde{b}(1-\Delta_{x_{1}})^{-1/2}$ is bounded, with
$$
|\langle \xi\otimes z_{1}\,,\, V(x_{1}-x_{2})z_{2}\otimes z_{3}\rangle_{L^{2}(\rz^{2d})}|
\leq
\|V(1-\Delta)^{-1/2}\|~|\xi|_{L^{^{2}}}|z_{1}|_{L^{2}}|z_{2}|_{H^{1}}|z_{3}|_{L^{2}}\,.
$$
A symmetric  version of \eqref{eq.V1D} says
$\tilde{b}(1-\Delta_{x_{2}})^{-1/2}$ is bounded, with
$$
|\langle \xi\otimes z_{1}\,,\, V(x_{1}-x_{2})z_{2}\otimes z_{3}\rangle_{L^{2}(\rz^{2d})}|
\leq
\|V(1-\Delta)^{-1/2}\|~|\xi|_{L^{^{2}}}|z_{1}|_{L^{2}}|z_{2}|_{L^{2}}|z_{3}|_{H^{1}}\,.
$$
Finally the symmetry of the expression
$V*(\overline{z_{1}}z_{2})z_{3}$ w.r.t  the exchange of
$\overline{z_{1}}$ and $z_{2}$ gives
$$
|\langle \xi\otimes z_{1}\,,\, V(x_{1}-x_{2})z_{2}\otimes
z_{3}\rangle_{L^{2}(\rz^{2d})}|
\leq
\|V(1-\Delta)^{-1/2}\|~|\xi|_{L^{^{2}}}|z_{1}|_{H^{1}}|z_{2}|_{L^{2}}|z_{3}|_{L^{2}}\,.
$$
Thus we have proved, owing to \eqref{eq.multlinV},
\begin{equation}
  \label{eq.estimVL2}
|V*(\overline{z_{1}}z_{2})z_{3}|_{L^{2}}\leq \|V(1-\Delta)^{-1/2}\|
\min_{\sigma\in \mathfrak{S}_{3}}|z_{\sigma(1)}|_{H^{1}}|z_{\sigma(2)}|_{L^{2}}|z_{\sigma(3)}|_{L^{2}}
\end{equation}
which gives
\begin{equation}
  \label{eq.estimVH1}
| (V*(\overline{z_{1}}z_{2}))z_{3}|_{H^{1}}
\leq
\|V(1-\Delta)^{-1/2}\|
\min_{\sigma\in \mathfrak{S}_{3}}|z_{\sigma(1)}|_{H^{1}}|z_{\sigma(2)}|_{H^{1}}|z_{\sigma(3)}|_{L^{2}}\,.
\end{equation}
Since $z\mapsto e^{it\Delta}z$ preserves the $L^{2}$ and $H^{1}$ norms, the
velocity field estimates \eqref{eq.vtL2} and \eqref{eq.vtH1} are
consequences of \eqref{eq.estimVL2} and \eqref{eq.estimVH1}.\\
For the sake of completeness, let us finish the proof of the global
well-posedness of the Cauchy problem.
The estimate \eqref{eq.estimVH1} provides the Lipschitz property of $z\to V*|z|^{2}z$ in
$H^{1}(\rz^{d})$\,. This implies the local in time existence and
uniqueness of a solution to \eqref{eq.cauchy2} in
$\mathcal{C}^{1}((-T_{z_{0}},T_{z_{0}}); H^{1}(\rz^{d}))$\,, and
therefore the local in time existence and uniqueness of a solution to
\eqref{eq.cauchy1} in
$\mathcal{C}^{0}((-T_{z_{0}},T_{z_{0}});H^{1}(\rz^{d}))\cap
\mathcal{C}^{1}([-T_{z_{0}},T_{z_{0}}]; H^{-1}(\rz^{d}))$\,.
The global in time existence then comes as usual from the control of
$|z_{t}|_{H^{1}}=|\tilde{z}_{t}|_{H^{1}}$ deduced from the
conservations of \eqref{eq.consL2} and \eqref{eq.consen}.
For \eqref{eq.consL2}, take the real part of the scalar product of each member of
\eqref{eq.cauchy1} with $\overline{z}_{t}$\,. This implies
$\partial_{t}|z_{t}|^{2}_{L^{2}}=0$\,.\\
For \eqref{eq.consen} take the scalar product with
$\chi(-R^{-1}\Delta)\partial_{t} z_{t}$ where $\chi\in
\mathcal{C}^{\infty}_{0}(\rz)$ satisfies $0\leq\chi\leq 1$ and
$\chi\equiv 1$ in a neighborhood of $0$\,, with $R>0$\,:
\begin{eqnarray*}
0&=&2\Real\langle\partial_{t}z_{t}\,,\,
\chi(-R^{-1}\Delta)i\partial_{t}z_{t} \rangle
\\
&=&\partial_{t}\langle z_{t}\,,\, -\Delta\chi(-R^{-1}\Delta)z_{t}\rangle
+2\Real \langle \partial_{t}z_{t}\,,\, \chi(-R^{-1}\Delta)[(V*|z_{t}|^{2})z_{t}]\rangle
\end{eqnarray*}
Integrating this identity from $0$ to $t$ and taking the limit as
$R\to \infty$ with the help of \eqref{eq.estimVH1} gives
$$
\int_{\rz^d}|\nabla z_{t}|^{2} dx-\int_{\rz^d}|\nabla z_{0}|^{2} \,dx
+2\int_{0}^{t}\Real \langle \partial_{s}z_{s}\,,\, (V*|z_{s}|^{2})z_{s}\rangle~ds\,=0\,.
$$
Due to the symmetry of $V(x)=V(-x)$\,, the last integrand equals
\begin{eqnarray*}
\Real \langle \partial_{s}z_{s}\,,\, (V*|z_{s}|^{2})z_{s}\rangle
&=&
\int_{\rz^{2d}} \partial_{s}(|z_{s}(x)|^{2})V(x-y)|z_{s}(y)|^2~dxdy\\
&=&
\frac{1}{2}\partial_{s }\int_{\rz^{2d}}
|z_{s}(x)|^{2}V(x-y)|z_{s}(y)|^2~dxdy\,.
\end{eqnarray*}
The conserved quantities \eqref{eq.consL2} and \eqref{eq.consen}
combined with \eqref{eq.estimVL2} imply $|z_{t}|_{H^{1}}\leq
C|z_{0}|_{H^{1}}$ for some constant independent of $t\in
(-T_{z_{0}},T_{z_{0}})$\,, and hence $T_{z_{0}}=+\infty$\,.
\fin

\section{Derivation of the mean field dynamics}
\label{se.dermeanfield}
This section contains the proof of our main Theorem \ref{th.main}.
Below, we recall from our previous work
\cite{AmNi1} the notion of infinite dimensional Wigner measures and collect some of their
properties.
We will often make use
of Weyl and Wick quantization throughout this section. So, we suggest
first the reading of Appendix \ref{se.WWaW}.\\
Two phase-spaces will be necessary for this analysis:
$\Z_{0}=L^{2}(\rz^{d}; \cz)$  (resp. $\Z_{1}=H^{1}(\rz^{d};\cz)$)
endowed with its scalar product $\langle ~,~\rangle$
(resp. $\langle z_{1}, z_{2}\rangle_{\Z_{1}}= \langle z_{1}\,,\, (1-\Delta)z_{2}\rangle$), its norm
$|z|_{\Z_{0}}^{2}=\langle z\,,\, z\rangle=|z|_{L^{2}}^{2}$
 (resp. $|z|_{\Z_{1}}^{2}=|z|_{H^{1}}^{2}$), its real scalar product
${\rm Re}~\langle z\,,\, z\rangle$ (resp. ${\rm Re}~\langle z_{1}\,,\,
z_{2}\rangle_{\Z_{1}}$). Only on $\Z_{0}$\,, we will use the  symplectic
structure with $\sigma(z_{1},z_{2})=\Imag \langle z_{1}\,,\, z_{2}\rangle$\,.
Meanwhile, the real euclidean structure on $\Z_{1}$ is important
especially when the Liouville transport equation is written as a
gradient flow according to Appendix~\ref{se.abscont}.
\subsection{Wigner measures}
\label{se.wigmeas}
  The Wigner measures are defined after the next result proved in
\cite[Theorem 6.2]{AmNi1}.
\begin{thm}
\label{th.wig-measure}
Let $\left(\varrho_{\varepsilon}\right)_{\varepsilon\in (0,\bar\varepsilon)}$ be a family
of normal states on $\mathcal{H}$ parametrized by $\varepsilon$\,.
Assume   ${\rm Tr}[\varrho_{\varepsilon}\N^\delta]$ $\leq C_{\delta}$ uniformly
w.r.t. $\varepsilon\in (0,\overline{\varepsilon})$
 for some fixed $\delta>0$ and $C_{\delta}\in (0,+\infty)$\,.
Then for every sequence $(\varepsilon_{n})_{n\in\nz}$ with $\lim_{n\to\infty}\varepsilon_n= 0$\,,
there exist a subsequence $(\varepsilon_{n_k})_{k\in\nz}$ and a Borel
probability measure $\mu$
on $\Z_{0}$\,, such that
\begin{eqnarray*}
\lim_{k\to\infty} \Tr[\varrho_{\varepsilon_{n_k}} b^{Weyl}]=
\int_{\Z_0} b(z) \; d\mu(z)\,,
\end{eqnarray*}
for all $b$ in the cylindrical Schwartz space $\S_{cyl}(\Z_0)$\, defined in Subsection \ref{se.weylAwick}.\\
Moreover this probability measure $\mu$ satisfies $\ds \int_{\Z_0} |z|_{\Z_0}^{2\delta} \, d\mu(z) <\infty$\,.
\end{thm}
\begin{definition}\ \\
  \label{de.setwig}
The set of Wigner measures associated with a family
$(\varrho_{\varepsilon})_{\varepsilon\in (0,\bar\varepsilon)}$
(resp. a sequence $(\varrho_{\varepsilon_{n}})_{n\in\nz}$)  which
satisfies the assumptions of Theorem~\ref{th.wig-measure} is
denoted by
$$
\mathcal{M}(\varrho_{\varepsilon}, \varepsilon\in
(0,\bar\varepsilon))\,, \quad (\textrm{resp.}\ \mathcal{M}(\varrho_{\varepsilon_{n}}, n\in\nz))\,.
$$
Moreover this definition can be extended to any family
$(\varrho_{\varepsilon})_{\varepsilon\in (0,\bar\varepsilon)}$ such
that
$$
\|(1+\N)^{\delta}\varrho_{\varepsilon}(1+\N)^{\delta}\|_{\mathcal{L}^{1}(\H)}\leq
C_{\delta}
$$ for some $\delta>0$ with the decomposition $\varrho_{\varepsilon}=
\lambda_{\varepsilon}^{R,+}\varrho_{\varepsilon}^{R,+}-
\lambda_{\varepsilon}^{R,-}\varrho_{\varepsilon}^{R,-}+
i\lambda_{\varepsilon}^{I,+}\varrho_{\varepsilon}^{I,+}-i\lambda_{\varepsilon}^{I,-}\varrho_{\varepsilon}^{I,-}$\,.
\end{definition}
Wigner measures are in practice identified via their characteristic functions according to the relation
\begin{eqnarray*}
\mathcal{M}(\varrho_{\varepsilon}, \varepsilon\in
(0,\bar\varepsilon))=\{\mu\}
&\Leftrightarrow & \lim_{\varepsilon\to 0}\Tr[\varrho_\varepsilon
\,W(\sqrt{2}\pi \xi)]=\mathcal{F}^{-1}(\mu)(\xi)\\
&\Leftrightarrow & \lim_{\varepsilon\to 0}\Tr[\varrho_\varepsilon
\,W(\xi)]=\int_{\Z_{0}}e^{i\sqrt{2}{\rm Re}~\langle \xi,z\rangle}d\mu(z)\,.
\end{eqnarray*}
The expression $\mathcal{M}(\varrho_{\varepsilon}, \varepsilon\in
(0,\bar\varepsilon))=\left\{\mu\right\}$ simply means that the family
$(\varrho_{\varepsilon})_{\varepsilon\in (0,\bar\varepsilon)}$ is "pure" in the sense
$$
\lim_{\varepsilon\to
  0}\Tr\left[\varrho_{\varepsilon} b^{Weyl}\right]=\int_\Z b(z)~d\mu\,,
$$
for all cylindrical symbol $b$ without extracting a
subsequence. Actually the general case can be reduced to this one, after
reducing the range of parameters to $\varepsilon\in \left\{\varepsilon_{n_{k}},
  k\in\nz\right\}$\,. For checking properties of the elements of
$\mathcal{M}(\varrho_{\varepsilon}, \varepsilon\in (0,\bar
\varepsilon))$\,, extracting a subsequence in this way allows to suppose without
loss of generality $\mathcal{M}(\varrho_{\varepsilon}, \varepsilon\in
(0,\bar \varepsilon))=\left\{\mu\right\}$\,.\\
A simple a priori estimate argument allows to extend the convergence
to symbols which have a polynomial growth and to take Wick
quantized symbols, with compact kernels, belonging to
$\P_{alg}^{\infty}(\Z_{0})=\oplus_{p,q\in\nz}^{alg}\P_{p,q}^{\infty}(\Z_{0})$
 (see \cite[Corollary 6.14]{AmNi1}).
\begin{prop}
 \label{pr.polycomp}
Let $\left(\varrho_{\varepsilon}\right)_{\varepsilon\in (0,\bar\varepsilon)}$ be a family
of normal states on $\L(\mathcal{H})$ parametrized by $\varepsilon$
such that ${\rm Tr}[\varrho_{\varepsilon}\N^{\alpha} ]\leq C_{\alpha}$ holds
uniformly with respect to $\varepsilon\in (0,\bar \varepsilon)$\,, for
all $\alpha\in\nz$\,,and such that $\mathcal{M}(\varrho_{\varepsilon},
\varepsilon\in (0,\bar \varepsilon))=\left\{\mu\right\}$\,.
Then the convergence
\begin{equation}
  \label{eq.convpolWW}
\lim_{\varepsilon\to 0}
\Tr[\varrho_{\varepsilon} b^{Wick}]=
\int_{\Z_0} b(z) \; d\mu(z)\,,
\end{equation}
holds for any $b\in \P_{alg}^{\infty}(\Z_0)$\,.
\end{prop}
A variant of the above result was provided in \cite[Theorem 6.13]{AmNi1}.
\begin{prop}
\label{pr.wigwick}
Assume that the family of operators
$(\varrho^{\varepsilon})_{\varepsilon\in (0,\overline{\varepsilon})}$
satisfies
$$\left\|(1+\N)^{\alpha}\varrho^{\varepsilon}(1+\N)^{\alpha}\right\|_{\L^{1}(\H)}\leq
C_{\alpha}$$
 uniformly w.r.t $\varepsilon\in
(0,\overline{\varepsilon})$ for all $\alpha\in\nz$\,.
For any fixed $\beta$ belonging to $\P_{alg}^{\infty}(\Z_0)$
  the family $(\beta^{Wick}\varrho^{\varepsilon})_{\varepsilon\in
  (0,\overline{\varepsilon})}$ satisfies the assumptions of
Definition ~\ref{de.setwig} and
\begin{equation}
  \label{eq.wigpol2}
\mathcal{M}(\beta^{Wick}\varrho^{\varepsilon})
=
\left\{\beta\mu\,,\;\mu\in \mathcal{M}(\varrho^{\varepsilon})\right\}
\,.
\end{equation}
\end{prop}
A closely related question is whether Wigner measures are completely identified
via Wick-quantized observable.
Of course this is related with the Hamburger moment problem even in
finite dimension and we again refer to \cite{AmNi1} for further
discussions about this.
\subsection{Weak mean field limit of the dynamics in terms of the
  characteristic function}
\label{se.meancyl}

After some extraction process and for some specific initial
 data $(\varrho_\varepsilon)_{\varepsilon\in(0,\bar\varepsilon)}$\,, a family $(\mu_{t})_{t\in\rz}$ of
 measures can be defined and solves weakly a transport equation.
We consider on $L^2_s(\rz^{2d})$ the (unbounded) multiplication operators
$$
\tilde V= \frac{1}{2} V(x_1-x_2) \quad \mbox{ and } \quad
\tilde V_s= (e^{-is \Delta_{x_1}}\otimes e^{-is \Delta_{x_2}})\tilde V
(e^{is \Delta_{x_1}}\otimes e^{is \Delta_{x_2}})\,,
$$
and respectively associate with them  the polynomials, well defined on  $\Z_1=H^1(\rz^d)$\,,
$$
V(z)=\frac{1}{2} \langle z^{\otimes 2}, V(x-y)z^{\otimes 2}\rangle_{L^2_s(\rz^{2d})}
\quad \mbox{ and } \quad  V_s(z)=\langle z^{\otimes 2}, \tilde V_s\,
z^{\otimes 2}\rangle_{L^2_s(\rz^{2d})}\,,\quad z\in \Z_{1}\,.
$$
Instead of considering
$$
\varrho_{\varepsilon}(t)=e^{-i\frac{t}{\varepsilon}H_{\varepsilon}}\varrho_{\varepsilon}e^{i\frac{t}{\varepsilon}H_{\varepsilon}}\,,
$$
we will rather work with
\begin{equation}
  \label{eq.deftildrho}
\tilde{\varrho}_{\varepsilon}(t)=
e^{i\frac{t}{\varepsilon}H_{\varepsilon}^{0}}e^{-i\frac{t}{\varepsilon}H_{\varepsilon}}
\varrho_{\varepsilon}
e^{i\frac{t}{\varepsilon}H_{\varepsilon}}e^{-i\frac{t}{\varepsilon}H_{\varepsilon}^{0}}\,.
\end{equation}
Our assumptions will be made in terms of the operator
$S_{\varepsilon}(1)$ already introduced in \eqref{eq.opref} and
which can be rewritten as a Wick observable.
\begin{definition}
  \label{de.opref}
The operator $S_{\varepsilon}$ is defined by
$$
S_{\varepsilon}=\sum_{n=0}^\infty
H_\varepsilon^{0,(n)}+\varepsilon n+ (\varepsilon n)^3
= \d\Gamma(1-\Delta)+ \N^{3}\,,
$$
with domain
$\mathcal{D}(S_{\varepsilon})=\left\{ \Psi\in \mathcal{H}\,,\;
\sum_{n=0}^{\infty} \|(H_\varepsilon^{0,(n)}+\varepsilon n+
(\varepsilon n)^3)\Psi^{(n)}\|_{L^{2}_{s}(\rz^{dn})}^2<\infty\right\}$ and
$H_{\varepsilon}^{0,(n)}=\d\Gamma(1-\Delta)\big|_{\bigvee^{n} \Z_{0}}$\,.
\end{definition}
Remember that it is self-adjoint with this domain (see
\eqref{eq.opref}). Moreover it can be written
$S_{\varepsilon}=s_{\varepsilon}^{Wick}$ with
$$
s_{\varepsilon}(z)= \langle
z\,,\,(1-\Delta)z\rangle+\left[|z|_{\Z_{0}}^{6}+3\varepsilon
  |z|_{\Z_{0}}^{4}+\varepsilon^{2}|z|_{\Z_{0}}^{2} \right]\,.
$$
\begin{prop}
\label{conteq}
Let  $\left(\varrho_{\varepsilon}\right)_{\varepsilon\in(0,\bar
  \varepsilon)}$ be a family of normal states on $\H$ satisfying for some finite constant $ C>0$
 the estimate
$$
{\rm Tr}[(1+S_\varepsilon)\varrho_{\varepsilon} (1+S_\varepsilon) ]\leq C \quad \mbox{ uniformly w.r.t }
\,\varepsilon\in(0,\bar\varepsilon)\,\,.
$$
The operator $S_{\varepsilon}$ is the one given in
Definition~\ref{de.opref} and $\tilde{\varrho}_{\varepsilon}(t)$ is
the operator given by \eqref{eq.deftildrho}.
Then for any sequence $(\varepsilon_{n})_{n\in\nz}$ in $(0,\bar\varepsilon)$
such that $\lim_{n\to\infty} \varepsilon_n=0$
there exist a subsequence $(\varepsilon_{n_k})_{k\in\nz}$ and a family $(\tilde{\mu}_t)_{t\in\rz}$ of
 Borel probability measures on $\Z_0$  satisfying for any $ t\in\rz$
\begin{eqnarray*}
\mathcal{M}( \tilde{\varrho}_{\varepsilon_{n_{k}}}(t), k\in\nz)=\{\tilde{\mu}_t\}\,,
\end{eqnarray*}
with the Liouville equation
\begin{eqnarray}
\label{fouriermu}
\tilde{\mu}_t(e^{i \sqrt{2}{\rm Re~}\langle \xi, .\rangle})&=&\tilde{\mu}_0(e^{i \sqrt{2}{\rm
    Re}\langle \xi, .\rangle})-2\sqrt{2} i
\int_0^t \tilde{\mu}_s(e^{i \sqrt{2}{\rm Re~}\langle \xi, z\rangle} \,{\rm Im~}\langle z^{\otimes 2}, \tilde V_s \xi\otimes z\rangle) \,ds\,,
\nonumber\\
&=& \tilde{\mu}_0(e^{i \sqrt{2}{\rm Re~}\langle \xi, .\rangle})+i
\int_0^t \tilde{\mu}_s(\{V_s(.);e^{i \sqrt{2} {\rm Re~}\langle \xi, .\rangle}\}) \,ds\,,
\quad \mbox{for all } \,\xi\in \Z_1\,.
\end{eqnarray}
\end{prop}
\proof
The proof uses several preliminary lemmas stated below.
The first step is to prove the existence
of Wigner measures defined for all times $t\in\rz$\,.
 This is done in Proposition \ref{prop-1}. Let us now prove the Liouville equation.\\
By Lemma \ref{derive} we have
\begin{equation}
  \label{eq.traceWxi}
\Tr[\tilde\varrho_\varepsilon(t) W(\xi)]= \Tr[\varrho_\varepsilon W(\xi)]+
i\int_0^t \Tr[\tilde\varrho_\varepsilon(s) W(\xi) \sum_{j=1}^4 \varepsilon^{j-1} b_j(s,\xi)^{Wick} ]
\,ds\,,
\end{equation}
where $b_j$ are the following polynomials
\begin{eqnarray*}
b_1(s,\xi)= -2 \sqrt{2} \;{\rm Im~}\langle z^{\otimes 2}, \tilde V_s \;\xi\otimes z\rangle &&
b_2(s,\xi)=-{\rm Re~}\langle z^{\otimes 2},\tilde V_s\; \xi^{\otimes 2}\rangle+
2 \langle \xi\vee z, \tilde V_s \xi\vee z\rangle\\
b_3(s,\xi)=\sqrt{2} {\rm Im~}\langle \xi^{\otimes 2}, \tilde V_s\; \xi\otimes z\rangle &&
b_4(s,\xi)= \frac{1}{4} \langle \xi^{\otimes 2}, \tilde V_s \;\xi^{\otimes 2}\rangle\,.
\end{eqnarray*}
With the number estimate in Proposition \ref{wick-estimate},
Lemma~\ref{comweyl} below will ensure that
the sum in the r.h.s over $j=2,\cdots ,4$ converges
 to $0$ when $\varepsilon\to 0$\,. On the other hand, the term
with $j=1$ has a limit  according to Lemma \ref{comparg} applied with
$\tilde{\varrho}_{\varepsilon}(t)$ after noticing that
$\Tr\left[(1+S_{\varepsilon})\tilde{\varrho}_{\varepsilon}(t)(1+S_{\varepsilon})\right]\leq
C'$ owing to $\|(1+S_{\varepsilon})^{\pm
  1}e^{i\frac{t}{\varepsilon}H_{0}^{\varepsilon}}e^{-i\frac{t}{\varepsilon}H_{\varepsilon}}(1+S_{\varepsilon})^{\mp
1}\|\leq C''$ due to Proposition~\ref{invH}.
\fin

The above proof is completed in essentially three steps: 1) The
relation~\eqref{eq.traceWxi} is first established by extending
Wick-calculus arguments to the case when $V$ is unbounded, and rough
estimates for  $b_{j}(s,\xi)^{Wick}$\,,
$j=1,\ldots,4$\,, are given; 2) An Ascoli type argument,
relying on these rough estimates
 allows to make the subsequence
extraction $(\varepsilon_{n_{k}})_{k\in\nz}$ uniform for all
$t\in\rz$\,; 3) An additional compactness argument is given in order to
ensure the convergence of the term with $j=1$ in \eqref{eq.traceWxi}.

\subsubsection{Wick calculus with unbounded kernels}
The results presented in this paragraph would be direct applications of
the Wick calculus given in Proposition~\ref{symbcalc} for a bounded
potential $V\in L^{\infty}(\rz^{d})$\,. Although the algebra is the
same as in the bounded case, justifying the formulas for unbounded
potentials fulfilling
(A\ref{eq.hypsym}), (A\ref{eq.hypbd}) and (A\ref{eq.hypcomp})
 requires some analysis.
\label{se.wicksing}
\begin{lem}
\label{comweyl}
The identity
\begin{eqnarray}
\label{comfor}
\left(V_s^{Wick} W(\xi)-W(\xi) V_s^{Wick}\right)\Psi=W(\xi) \,
\left(\sum_{j=1}^4 \varepsilon^j b_j(s,\xi)^{Wick}\right)\Psi\,,
\end{eqnarray}
holds
for any $\xi\in H^2(\rz^d)$ and $\Psi\in\D(S_\varepsilon)$\,, with
$S_{\varepsilon}$ given by Definition~\ref{de.opref}. Additionally,
for all $\Psi\in \D(S_{\varepsilon})\subset
\D(S_{\varepsilon}^{1/2})\subset\D(\N^{\frac{3}{2}})$\,, the estimates
\begin{multline}
  \label{eq.estimbjwick}
\|b_{j}(s,\xi)^{Wick}~\Psi\|\leq C(1+|\xi|_{\Z_{1}}^{4})\|
(1+\N)^{\frac{3}{2}}\Psi\|\leq C'
(1+|\xi|_{\Z_{1}}^{4})\|(1+S_{\varepsilon})^{1/2}\Psi\|\,,
\end{multline}
hold uniformly w.r.t $j\in\left\{1,\ldots,4\right\}$\,, $s\in\rz$\,, when
$\xi\in H^{1}(\rz^{d})$ \,.
\end{lem}
\proof
We first remark that, owing to the assumption (A\ref{eq.hypbd}) and the
estimate \eqref{eq.estimVL2}, the polynomials $b_j(s,\xi)$\,, $j=1,\cdots, 4$ belong to the set
$\oplus_{p,q\leq 3} \P_{p,q}(\Z_0)$\,, with
$$
|b_{j}|_{\oplus_{p+q\leq 3}\P_{p,q}(\Z_{0})}\leq C(1+|\xi|_{\Z_{1}}^{4})\,.
$$
Hence, Proposition \ref{wick-estimate} and Proposition \ref{pr.sa}
 prove \eqref{eq.estimbjwick} with
\begin{eqnarray}
\label{2eq}
\D(S_\varepsilon)\subset\D(\N^{3/2})\subset \D(b_j(s,\xi)^{Wick})\,,~j=1,\ldots,4~,\quad \mbox{ and } \quad
\D(S_\varepsilon)\subset  \D(V_s^{Wick})\,.
\end{eqnarray}
By Proposition \ref{weylinv} the domain
$\D(S_\varepsilon)$ is invariant under the action of $W(\xi)$ for all $\xi\in H^{2}(\rz^{d})$\,.
A Taylor expansion yields, for all $z\in\Z_1$\,, the equality
$$
V_s(z+\frac{i\varepsilon}{\sqrt{2}} \xi)=V_s(z)+\sum_{j=1}^4 \varepsilon^j b_j(s,\xi)[z]\,.
$$
The formula \eqref{comfor} is standard for bounded $\tilde{V}$ due to
$W^{*}(\xi)b^{Wick}W(\xi)=b(.+\frac{i\varepsilon}{\sqrt{2}}\xi)^{Wick}$
when $b\in \P_{alg}(\Z_{0})=\oplus_{p,q\in\nz}^{alg}\P_{p,q}(\Z_{0})$\,. Let us
reconsider the proof of this result for our unbounded $\tilde{V}$\,.\\
With the previous estimates, the quantity $\mathcal{A}(t)=
\langle \Phi, W(t\xi) V_s(.+\frac{i\varepsilon}{\sqrt{2}} t\xi)^{Wick} W(t\xi)^*
\Psi\rangle$ is well defined for all $\Phi,\Psi\in\D(S_\varepsilon)$ with
\begin{eqnarray*}
\mathcal{A}(t)=\sum_{j=1}^4\varepsilon^j t^j  \langle\Phi, W(t\xi) b_j(s,\xi)^{Wick} W(t\xi)^*
\Psi\rangle +  \langle\Phi, W(t\xi) V_s^{Wick}  W(t\xi)^*
\Psi\rangle\,.
\end{eqnarray*}
We first establish   in a weak sense the equality \eqref{comfor}:
Differentiate $\mathcal{A}(t)$ for any $\Psi,\Phi\in \D(S_\varepsilon)$\,,
\begin{eqnarray*}
\frac{d}{dt} \mathcal{A}(t)
&=&\sum_{j=1}^4 \varepsilon^j   \langle\Phi, W(t\xi) \left\{[i\phi(\xi), b_j(s,\xi)^{Wick}] t^j +
j t^{j-1} b_j(s,\xi)^{Wick}\right\} W(t\xi)^*
\Psi\rangle \\
&& \hspace{.2in} +  \langle\Phi, W(t\xi) [i\phi(\xi),V_s^{Wick}]  W(t\xi)^*
\Psi\rangle\,\\
&=& \sum_{j=0}^3 \varepsilon^j t^j \langle\Phi, W(t\xi) \left\{[i\phi(\xi), b_j(s,\xi)^{Wick}]  +
(j+1) b_{j+1}(s,\xi)^{Wick}\right\} W(t\xi)^*
\Psi\rangle
\end{eqnarray*}
where $b_0(z)=V_s(z)$\,. Now, a direct calculation with
$\phi(\xi)=\frac{1}{\sqrt{2}}(a(\xi)+a^{*}(\xi))$ gives
$$
[i\phi(\xi), b_j(s,\xi)^{Wick}]=-(j+1) b_{j+1}(s,\xi)^{Wick}
$$
for $j=0,\cdots, 3$\,. Therefore $\mathcal{A}(1)=\mathcal{A}(0)$
and, knowing \eqref{2eq}, we conclude that
\begin{eqnarray}
\label{1eq}
W(\xi) \left(V_s^{Wick}+\sum_{j=1}^4 \varepsilon^j b_j(s,\xi)^{Wick}\right) W(\xi)^* \Psi=V_s^{Wick} \Psi\,.
\end{eqnarray}
for any $\Psi\in\D(S_\varepsilon)$\,. With $\Psi=W(\xi)\tilde \Psi$ in \eqref{1eq} for
any $\tilde\Psi\in\D(S_\varepsilon)$\,, while $\Psi\in
\mathcal{D}(S_{\varepsilon})$ owing to  $\xi\in H^{2}(\rz^{d})$\,, the claimed equality is obtained.\\

\fin
\begin{lem}
\label{derive}
Let $(\varrho_\varepsilon)_{\varepsilon\in(0,\bar\varepsilon)}$ a family of normal states on $\H$\,. Assume that
$\varrho_\varepsilon (S_\varepsilon +1)\in \L^1(\H)$ for all
$\varepsilon\in(0,\bar\varepsilon)$\,, with $S_{\varepsilon}$
given by Definition~\ref{de.opref} and
$\tilde{\varrho}_{\varepsilon}(t)$ by \eqref{eq.deftildrho}. Then for any $\xi\in
H^{2}(\rz^{d})$\,, the map
$s\mapsto \Tr[\tilde \varrho_\varepsilon(s) \;W(\xi)]$ belongs to $C^1(\rz)$
and the following  integral formula holds true
$$
\Tr[\tilde \varrho_\varepsilon(t) W(\xi)]=
\Tr[\varrho_\varepsilon  W(\xi)]+\frac{i}{\varepsilon} \int_0^t
\Tr[\tilde \varrho_\varepsilon(s) W(\xi) \sum_{j=1}^4 \varepsilon^j b_j(s,\xi)^{Wick}] \,ds\,.
$$
\end{lem}
\proof
Write
\begin{eqnarray*}
&&\hspace{-.2in} \Tr[(\tilde \varrho_\varepsilon(t)-\tilde \varrho_\varepsilon(s)) W(\xi)]\\
&&=\Tr\left[\varrho_\varepsilon (S_\varepsilon +1) (S_\varepsilon +1)^{-1}
\left(e^{i\frac{t}{\varepsilon} H_\varepsilon} e^{-i\frac{t}{\varepsilon} H_\varepsilon^0} - e^{i\frac{s}{\varepsilon} H_\varepsilon}
e^{-i\frac{s}{\varepsilon} H_\varepsilon^0}\right)
W(\xi) e^{i\frac{s}{\varepsilon} H_\varepsilon^0}
e^{-i\frac{s}{\varepsilon} H_\varepsilon}\right]\\
&&\hspace{.1in} +\Tr\left[ \varrho_\varepsilon
e^{i\frac{t}{\varepsilon} H_\varepsilon} e^{-i\frac{ t}{\varepsilon} H_\varepsilon^0}
W(\xi) (S_\varepsilon +1) (S_\varepsilon +1)^{-1}
\left(e^{i\frac{t}{\varepsilon} H_\varepsilon^0} e^{-i\frac{t}{\varepsilon} H_\varepsilon}- e^{i\frac{s}{\varepsilon}
H_\varepsilon^0} e^{-i\frac{s}{\varepsilon} H_\varepsilon^0}\right)
\right]\,.
\end{eqnarray*}
The following limits  hold true on $\D(S_\varepsilon )$
\begin{eqnarray*}
 \lim_{s\to t}  \frac{1}{t-s}(S_\varepsilon +1)^{-1}
\left(e^{i\frac{t}{\varepsilon} H_\varepsilon} e^{-i\frac{t}{\varepsilon} H_\varepsilon^0} - e^{i\frac{s}{\varepsilon} H_\varepsilon}
e^{-i\frac{s}{\varepsilon} H_\varepsilon^0}\right)=\frac{i}{\varepsilon}
(S_\varepsilon +1)^{-1} e^{i\frac{t}{\varepsilon} H_\varepsilon} (H_\varepsilon-H_\varepsilon^0)
e^{-i\frac{t}{\varepsilon} H_\varepsilon^0}\,\\
\lim_{s\to t}  \frac{1}{t-s}(S_\varepsilon +1)^{-1}
\left(e^{i\frac{t}{\varepsilon} H_\varepsilon^0} e^{-i\frac{t}{\varepsilon} H_\varepsilon}- e^{i\frac{s}{\varepsilon}
H_\varepsilon^0} e^{-i\frac{s}{\varepsilon} H_\varepsilon^0}\right)=
\frac{i}{\varepsilon} (S_\varepsilon +1)^{-1} e^{i\frac{t}{\varepsilon} H_\varepsilon^0}
(H_\varepsilon^0-H_\varepsilon) e^{-i\frac{t}{\varepsilon} H_\varepsilon}\,,
\end{eqnarray*}
by Stone's theorem and the invariance of $\D(S_\varepsilon )$ w.r.t
$e^{itH_\varepsilon^0}$ and $e^{it H_\varepsilon}$\,.
By using the estimate in Proposition \ref{invH},
the latter limits are limits in $\L(\H)$ w.r.t the strong convergence topology.
After noticing that $ \varrho_\varepsilon
e^{i\frac{ t}{\varepsilon} H_\varepsilon} e^{-i\frac{ t}{\varepsilon} H_\varepsilon^0}
W(\xi) (S_\varepsilon +1)$ is trace class when $\xi\in H^{2}(\rz)$, owing to
Proposition~\ref{weylinv} and
Proposition~\ref{invH}, we take the trace and let $s\to t$\,.\\
Now  integrating the derivative from $0$ to $t$ yields
$$
\Tr[\tilde \varrho_\varepsilon(t) W(\xi)]=
\Tr[\varrho_\varepsilon  W(\xi)]+\frac{i}{\varepsilon} \int_0^t
 \Tr \left[\tilde \varrho_\varepsilon(s) \, \left(V_s^{Wick} W(\xi)-W(\xi) V_s^{Wick}\right)\right] \,ds\,.
$$
When $\xi\in H^{2}(\rz^{d})$\,, the equality
$$
\Tr \left[(1+S_\varepsilon)\tilde\varrho_\varepsilon(s)
\left(V_s^{Wick} W(\xi)-W(\xi) V_s^{Wick}\right)(1+S_\varepsilon)^{-1}\right]=
\Tr[\tilde \varrho_\varepsilon(s) W(\xi) \sum_{j=1}^4 \varepsilon^j b_j(s,\xi)^{Wick}] \,,
$$
makes sense,
since $(1+S_\varepsilon)\tilde\varrho_{\varepsilon}(s)\in \L^1(\H)$ and by Lemma~\ref{comweyl}
\begin{equation*}
\left(V_s^{Wick} W(\xi)-W(\xi)
  V_s^{Wick}\right)(1+S_\varepsilon)^{-1}=W(\xi)
\left[
\sum_{j=1}^{4}\varepsilon^{j}b_{j}(s,\xi)^{Wick}\right](1+S_{\varepsilon})^{-1}\quad \text{in}~\L(\H)\,.
\end{equation*}
\fin

\bigskip

\subsubsection{Subsequence extraction for all times}
\label{se.extrac}

The  first step in the proof of Proposition \ref{conteq} is to show the existence
of Wigner measures for all times. This is
accomplished below  by following merely the same lines
as \cite[Proposition 3.3]{AmNi3}.

\begin{prop}
\label{prop-1}
Let  $\left(\varrho_{\varepsilon}\right)_{\varepsilon\in(0,\bar
  \varepsilon)}$ be a family of normal states on $\H$ satisfying for some finite constant $ C>0$
 the estimate
$$
{\rm Tr}[\varrho_{\varepsilon} (1+S_\varepsilon) ]\leq C \quad \mbox{ uniformly w.r.t }
\varepsilon\in(0,\bar\varepsilon)\,\,.
$$
The operator $S_{\varepsilon}$ and $\tilde{\varrho}_{\varepsilon}(t)$
are
respectively given by Definition~\ref{de.opref}
 and \eqref{eq.deftildrho}.
Then for any sequence $(\varepsilon_{n})_{n\in\nz}$ in $(0,\bar\varepsilon)$
 such that $\lim_{n\to\infty} \varepsilon_n=0$
there exists a subsequence $(\varepsilon_{n_k})_{k\in\nz}$ and a family
of Borel probability measures on $\Z_{0}$\,, $(\tilde{\mu}_t)_{t\in\rz}$\,, satisfying
\begin{eqnarray*}
\mathcal{M}( \tilde{\varrho}_{\varepsilon_{n_{k}}}(t) ,~ k\in\nz)=\{\tilde{\mu}_t\}\,,
\end{eqnarray*}
for any $ t\in\rz$\,.
\end{prop}
\proof
We only sketch the proof and essentially indicate the points which differ from
\cite[Proposition 3.3]{AmNi3}.
Let us write
\begin{eqnarray*}
G_\varepsilon(t,\xi)=\Tr[\tilde\varrho_{\varepsilon}(t) W(\xi)]\;.
\end{eqnarray*}
By using Proposition \ref{wick-estimate} and $(1+\N)\leq 2(1+\N^{3})\leq
2(1+S_{\varepsilon})$\,,  one can prove like in \cite{AmNi3} that
\begin{eqnarray}
\label{3eq}
|G_\varepsilon(s,\xi)-G_\varepsilon(s,\eta)|\leq C |\xi-\eta|_{\Z_0}^{\frac{1}{2}} \,
(|\xi|^2_{\Z_0}+|\eta|_{\Z_0}^2+1)^{\frac{1}{4}}
\end{eqnarray}
for some constant $C>0$\,. We have
\begin{eqnarray*}
|G_\varepsilon(t,\eta)-G_\varepsilon(s,\xi)|\leq  \left|{\rm Tr}\left
[\left(\tilde\varrho_\varepsilon(t)-\tilde \varrho_\varepsilon(s)\right)
W(\xi)\right]\right|+
\left|G_\varepsilon(s,\xi)-G_\varepsilon(s,\eta)\right|\,.
\end{eqnarray*}
On the other hand, by making use of Lemma \ref{derive} we get
\begin{eqnarray*}
\label{equi-eq1}
\left|{\rm Tr}[\tilde\varrho_\varepsilon(t)-\tilde \varrho_\varepsilon(s)]
W(\xi)]\right|&\leq&  \left|\int_s^t {\rm Tr}[ \tilde\varrho_\varepsilon(w)
\sum_{j=1}^{4} \varepsilon^{j-1} b_j(w,\xi)^{wick}] \,dw\right|\\
&\leq& C_0 |t-s| \|(1+S_{\varepsilon})^{1/2}\varrho_\varepsilon (1+S_{\varepsilon})^{1/2}\|_{\L^1(\H)} \, \\
&&\times\sup_{w\in[t,s]} \|(1+S_{\varepsilon})^{-1/2} \left[\sum_{j=1}^{4} \varepsilon^{j-1} b_j(w,\xi)^{wick}\right] (1+S_{\varepsilon})^{-1/2}\|_{\L(\H)}\\
&\leq& C_1 |t-s| (1+|\xi|_{\Z_1})^4\,,
\end{eqnarray*}
when $\xi\in H^{2}(\rz^{d})$\,. Taking an approximation $\xi_{n}\in
H^{2}(\rz)$\,, $n\in \nz$\,, such that
 $\lim_{n\to\infty}|\xi-\xi_{n}|_{\Z_{1}}=0$\,,
$\Z_{1}=H^{1}(\rz^{n})$\,, and taking the limit as $n\to \infty$ of
the left-hand side with the help of \eqref{3eq}, allows first to
extend the previous inequality to any $\xi\in\Z_{1}$\,.\\
Thus, we conclude that
\begin{equation}
 \label{eq.unifcont}
\left|G_{\varepsilon}(t,\eta)-G_{\varepsilon}(s,\xi)\right|
\leq \tilde{C} \left( |t-s| (|\xi|_{\Z_{1}}+1)^{4}+\, |\eta-\xi|_{\Z_{0}}\, \sqrt{|\eta|_{\Z_{0}}^2+|\xi|_{\Z_{0}}^2}\right),
\end{equation}
holds for all $(s,\xi),(t,\eta)\in \rz\times\Z_{1}$\,,
uniformly w.r.t. $\varepsilon\in (0,\overline{\varepsilon})$\,. Remember also
the uniform estimate $\left|G_{\varepsilon}(s,\xi)\right|\leq
1$\,.\\
\noindent
Now, we apply the same Ascoli type argument the one used in \cite[Proposition
3.3]{AmNi3} in order to prove the
existence of a subsequence $(\varepsilon_{n_k})_k$ and a continuous function
$G(.,.):\rz\times \Z_1\to \cz$ such that  $G_{\varepsilon_k}(t,\xi)$ converges to $G(t,\xi)$ for any
$t\in \rz$ and $\xi\in \Z_1$\,. Furthermore \eqref{3eq} allows to extend $G(.,.)$ to a continuous
function on $\rz\times \Z_0$\,.  An ``$\delta/3$''-argument
 shows that for any $(t,\xi)\in\rz\times\Z_0$\,,  $\lim_{n\to\infty}$ $ G_{\varepsilon_{n}} (t,\xi)$
exists and equals  $G(t,\xi)$\,, so that
$G(t,.)$ is a norm continuous normalized function of positive type.
Therefore, for any $t\in\rz$\,,  $G(t,.)$ is a characteristic
function of weak distribution (or projective family of probability measures)  $\tilde\mu_t$ on $\Z_0$\,. Finally the proof is ended as in \cite[Proposition 3.3]{AmNi3}.
\fin

\subsubsection{An additional compactness argument}
\label{se.addcomp}
Here, the compactness assumption~(A\ref{eq.hypcomp}) is converted into some
compactness property of the Wick symbol $b_{1}$\,. It allows to refer
indirectly to Proposition~\ref{pr.polycomp} and to take the limit
as $\varepsilon\to 0$ in the term with $j=1$ in
\eqref{eq.traceWxi}. With the rough estimates used in
Proposition~\ref{prop-1}, the terms in \eqref{eq.traceWxi}
 corresponding to  $j>1$  with a
factor $\varepsilon^{j-1}$ will vanish as $\varepsilon\to 0$\,. The
next Lemma applied with $\tilde{\varrho}_{\varepsilon}(s)$ in the
integral term of \eqref{eq.traceWxi},
will end the proof of Proposition~\ref{conteq}.
\begin{lem}
\label{comparg}
Let $\varrho_\varepsilon$ be a family of normal states on $\H$
satisfying  for some finite constant $C>0$ the estimate
$$
\Tr\left[(1+S_{\varepsilon})\varrho_{\varepsilon}(1+S_{\varepsilon})\right]\leq C \quad
\mbox{ uniformly w.r.t } \, \varepsilon\in(0,\bar\varepsilon)\,.
$$
Here $S_{\varepsilon}$ is given by Definition~\ref{de.opref}\,. Assume that
$\mathcal{M}(\varrho_\varepsilon{, \varepsilon\in (0,\bar\varepsilon)})=\{\mu\}$\,, then for any  $\xi\in \Z_1$\,,
\begin{eqnarray*}
\lim_{\varepsilon\to 0} \Tr[\varrho_\varepsilon W(\xi) b_1(s,\xi)^{Wick}] =
\int_{\Z_0} e^{\sqrt{2} i {\rm Re~}\langle\xi,z\rangle}
b_1(s,\xi)[z] \;d\mu(z)\,.
\end{eqnarray*}
\end{lem}
\proof
 The polynomial $b_1(s,\xi)\in \P_{1,2}+\P_{2,1}$ splits into two similar terms, namely
 $$
B_1(z)=\langle \xi\otimes z, \tilde V_s z^{\otimes 2}\rangle
\quad \mbox{ and } \quad B_2(z)=\langle z^{\otimes 2}, \tilde V_s (z\otimes \xi)\rangle\,
 $$
with their associated operators
$$
\tilde{B}_{1}=(\langle
 \xi|\otimes\11)\, \tilde{V_{s}}\in \L(L_s^2(\rz^{2d}),L^2(\rz^d)) \quad
 \mbox{ and }  \quad
 \tilde{B}_{2}=\S_2
\tilde{V_{s}}\, (\11\otimes |\xi\rangle) \in \L(L^2(\rz^d),L_s^2(\rz^{2d}))\,.
 $$
Let $\chi\in C_0^\infty(\rz)$ with $\chi(x)=1$ if $|x|\leq 1$\,, $\chi(x)=0$ if
$|x|\geq 2$ and $0\leq\chi\leq 1$\,. For $m\in\mathbb{N}^*$\,,
  set $\chi_m(x)=\chi(\frac{x}{m})$ and
define
$$
\tilde B_{1,m}=\chi_m(|D_x|) \tilde B_1 \,(\11\otimes \chi_m(|D_x|)) \S_2 \quad \mbox{ and } \quad
\tilde B_{2,m}=\S_2 (\11\otimes \chi_m(|D_x|))\,  \tilde B_2  \chi_m(|D_x|)
$$
as bounded operators in $\L(L^2_s(\rz^{2d}), L^2(\rz^d))$ and $\L(L^2(\rz^d),L^2_s(\rz^{2d}))$ respectively.
We claim that both operators $ \tilde B_{1,m}$ and $\tilde B_{2,m}$ are compact. Actually,
$\tilde B_{2,m}=\tilde  B_{1,m}^* $ and
\begin{eqnarray*}
\tilde  B_{1,m}=\frac{1}{2}
(\11\otimes e^{-is \Delta_{x_{2}}}) (\langle e^{is \Delta}\xi|\otimes \11)
(e^{-ix_{1}D_{x_{2}}} \chi_m(|D_{x_2}|) V(x_2) \chi_m(|D_{x_2}|) e^{ix_{1}D_{x_{2}}} )
(e^{is \Delta_{x_{1}}}\otimes e^{is \Delta_{x_{2}}})\mathcal{S}_{2}\,.
\end{eqnarray*}
Moreover, the linear norm continuous application
\begin{eqnarray*}
A\in \L(L^2(\rz^d))&\longmapsto &(\langle e^{is\Delta}\xi|\otimes \11)e^{-ix_{1}D_{x_{2}}}
(1\otimes A) \in \L(L^2_s(\rz^{2d}),L^2(\rz^d))
\end{eqnarray*}
preserves the class of  Hilbert-Schmidt operators since
$$
\|(\langle e^{is\Delta}\xi|\otimes \11)e^{-ix_{1}D_{x_{2}}}
(1\otimes A)\|_{\L^2(L^2_s(\rz^{2d}),L^2(\rz^d))}=|\xi|_{\Z_0}
\|A\|_{\L^2(\Z_0)}
$$
comes by computing the Schwartz kernel with
$\|K\|_{\mathcal{L}^{2}(L^{2}(\mu);L^{2}(\nu))}^{2}=\int |K(x,y)|^{2}~d\nu(x)d\mu(y)$\,.
Hence it maps compact operators  into compact operators, because the
space of compact operators,
$\mathcal{L}^{\infty}$\,, is the norm closure of $\mathcal{L}^{2}$ in
$\mathcal{L}$\,.
 Therefore, by taking $A=
\chi_m(|D_{x}|) \, V(x) \, \chi_m(|D_{x}|)$ which is compact by assumption (A\ref{eq.hypcomp}), we
conclude that $\tilde B_{1,m}$ and $\tilde B_{2,m}$ are compact. \\
Now, writing for $j=1,2$
\begin{eqnarray}
\label{cv1}
&&
|\Tr[\varrho_\varepsilon W(\xi) B_j^{Wick}]-\mu(e^{\sqrt{2} i {\rm Re~}\langle\xi,z\rangle}
B_j(z))|\leq |\Tr[\varrho_\varepsilon W(\xi)
(B_j^{Wick}-B_{j,m}^{Wick})]|\\\label{cv2}
&&\hspace{5cm}
+
|\Tr[\varrho_\varepsilon W(\xi) B_{j,m}^{Wick}]-\mu(e^{\sqrt{2} i {\rm Re~}\langle\xi,z\rangle}
B_{j,m}(z))|\\\label{cv3}
&&\hspace{5cm}
+|\mu(e^{\sqrt{2} i {\rm Re~}\langle\xi,z\rangle}
B_{j,m}(z))-\mu(e^{\sqrt{2} i {\rm Re~}\langle\xi,z\rangle}
B_j(z))|\,,
\end{eqnarray}
with $B_{j,m}\in \P^{\infty}_{alg}(\Z_0)$\,. The right hand side \eqref{cv2}
converges to $0$ owing to
Proposition \ref{pr.wigwick}. Since $s-\lim_{m\to \infty} \chi_m(|D_x|)=\11$\,, the
polynomials $B_{j,m}(z)$ converge to $B_j(z)$ for any $z\in\Z_0$\,, while the estimate
$$
|B_{j,m}(z)|\leq c |\xi|_{\Z_1}\, \|(1-\Delta)^{-1/2} V\| \,|z|_{\Z_0}^3
$$
holds true uniformly w.r.t $m$ for  some constant $c>0$\,. Additionally,
the estimate
 $\Tr[\varrho_\varepsilon \N^{3/2}]\leq C$
implies
$$
\int_{\Z_0} |z|_{\Z_0}^3 \, d\mu\leq C\,.
$$
Therefore the
dominated convergence theorem applies and
the right hand side \eqref{cv3} tends to $0$ as $m\to\infty$\,.
It remains to prove the convergence of the r.h.s of \eqref{cv1}. Writing
\begin{eqnarray*}
\Tr[\varrho_\varepsilon W(\xi) (B_j^{Wick}-B_{j,m}^{Wick})]
&=&\Tr[(S_\varepsilon+1) \varrho_\varepsilon (S_\varepsilon+1)
(S_\varepsilon+1)^{-1} W(\xi) (S_\varepsilon+1) \\
&&\hspace{1in} \times (S_\varepsilon+1)^{-1}  (B_j^{Wick}-B_{j,m}^{Wick}) (S_\varepsilon+1)^{-1}]
\end{eqnarray*}
and referring to Proposition \ref{weylinv} lead to the estimate
\begin{eqnarray*}
|\Tr[\varrho_\varepsilon W(\xi) (B_j^{Wick}-B_{j,m}^{Wick})] | \leq c
\|(S_\varepsilon+1)^{-1}  (B_j^{Wick}-B_{j,m}^{Wick}) (S_\varepsilon+1)^{-1}\|_{\L(\H)}\,.
\end{eqnarray*}
By functional calculus of strongly commuting self-adjoint operators we see that
$(S_\varepsilon+1)^{-\emph{}1} (\sqrt{\N}+\d\Gamma(1-\Delta)+1)$ is uniformly bounded  with
respect to $\varepsilon\in(0,\bar\varepsilon)$\,.
By applying Lemma \ref{le.compest} (with $A=1-\Delta$), we conclude that
\begin{eqnarray*}
|\Tr[\varrho_\varepsilon W(\xi) (B_j^{Wick}-B_{j,m}^{Wick})] | &\lesssim& \|(1-\Delta_{x_2})^{-1/2}
(\tilde B_j-\tilde B_{j,m})  (1-\Delta_{x_2})^{-1/2} \|\\
&\lesssim&  |\xi|_{\Z_0} \, \|(1-\Delta)^{-1/2} V\|_{\L(\Z_0)} \, \|(1-\Delta_{x})^{-1/2}(1-\chi_m(|D_x|))
\|_{\L(\Z_0)}
\end{eqnarray*}
Again by functional calculus
$\|(1-\Delta_{x})^{-1/2}(1-\chi_m(|D_x|))
\|$ is estimated by
$\frac{1}{m}$ and the r.h.s of \eqref{cv1} goes to $0$ as $m\to\infty$ uniformly
w.r.t $\varepsilon\in(0,\bar\varepsilon)$\,. Finally, a
"$\delta/3$-argument" with the established convergence of \eqref{cv1}, \eqref{cv2} and \eqref{cv3}
yields the result.
\fin
\subsection{Asymptotic a priori estimates}
\label{se.asapriori}

In this section, a priori information on Wigner measures are derived from a
priori estimates on the state $\varrho_{\varepsilon}$\,. In
particular, we shall
prove the next result.
\begin{prop}
\label{pr.asapriori}
Let $S_{\varepsilon}$ be the operator given by
Definition~\ref{de.opref} and assume that the family of normal states
$\left(\varrho_{\varepsilon}\right)_{\varepsilon\in
  (0,\bar\varepsilon)}$ satisfies
$$
\forall \alpha\in\nz, \exists C_{\alpha}>0,\forall
\varepsilon\in (0,\bar \varepsilon),\quad\Tr\left[(1+S_{\varepsilon})\varrho_{\varepsilon}(1+S_{\varepsilon})(1+\N)^{\alpha}\right]\leq C_{\alpha}\,,
$$
and $\mathcal{M}(\varrho_{\varepsilon},
\varepsilon\in (0,\bar \varepsilon))=\left\{ \mu\right\}$\,. Then the
measure $\mu$ is carried by $\Z_{1}$\,, its restriction to $\Z_{1}$ is
 a Borel probability measure on $\Z_{1}$ and
\begin{equation}
   \label{eq.estimz6}
\int_{\Z_{0}}|z|_{\Z_{1}}^{4}|z|_{\Z_{0}}^{2}~d\mu(z)
=
\int_{\Z_{1}}|z|_{\Z_{1}}^{4}|z|_{\Z_{0}}^{2}~d\mu(z)<+\infty\,.
\end{equation}
\end{prop}
The proof of the above proposition requires the two next Lemmas.
\begin{lem}
\label{le.poswick}
Let $\tilde{b}$ be a non negative (self-adjoint) operator on
$\bigvee^{p}\Z_{0}$ and assume that the family of normal states
$(\varrho_{\varepsilon})_{\varepsilon\in (0,\bar \varepsilon)}$\,,
with $\Tr\left[\varrho_{\varepsilon}\N^{\alpha}\right]\leq C_{\alpha}$ for
all $\alpha\in\nz$\,,
satisfies
$$
\Tr\left[\varrho_{\varepsilon}b^{Wick}\right]\leq C\quad
\text{and}\quad \mathcal{M}(\varrho_{\varepsilon},\;\varepsilon\in
(0,\bar \varepsilon))=\left\{\mu\right\}\,.
$$
Then $\Z_{0}\ni z\mapsto b(z)=\langle z^{\otimes p}\,,\, \tilde{b}
z^{\otimes p}\rangle\in [0,+\infty]$ is a Borel function on $\Z_{0}$ and
$\int_{\Z_{0}}b(z)~d\mu(z)\leq C$\,.
\end{lem}
\proof
When $b\in \P_{p,p}^{\infty}(\Z_{0})$ has a compact kernel $\tilde{b}$
we know after Proposition~\ref{pr.polycomp} (see  \cite[Corollary
6.14]{AmNi1} for a complete proof) that
$$
C\geq \lim_{\varepsilon\to
  0}\Tr\left[\varrho_{\varepsilon}b^{Wick}\right]=\int_{\Z_{0}} b(z)~d\mu(z)\,.
$$
We use the fact that $\tilde{b}\to b^{Wick}$ is operator monotone, in
the following sense: if the  (possibly unbounded)
  non negative operators $\tilde{b}_{1}, \tilde{b}_{2}$ in $\bigvee^{p}\Z_{0}$
satisfy $\tilde{b}_{2}\geq\tilde{b}_{1}\geq 0$\,, then the
 densely defined essentially
self-adjoint operators $b_{j}^{Wick}$\,, $j=1,2$ in  $\mathcal{H}$ satisfy
$b_{2}^{Wick}\geq b_{1}^{Wick}\geq 0$ \,.\\
By taking $\tilde{b}\in \mathcal{L}(\bigvee^{p}\Z)$\,, for $b\in
\P_{p,p}(\Z)$\,, as the supremum of $\tilde{b}_{n}$ with $\tilde{b}_{n}$
compact, we obtain firstly for all $n\in\nz$
$$
C\geq
\liminf_{\varepsilon\to
  0}~\Tr\left[\varrho_{\varepsilon}^{1/2}b^{Wick}\varrho_{\varepsilon}^{1/2}\right]
\geq \lim_{\varepsilon\to
  0}~\Tr\left[\varrho_{\varepsilon}^{1/2}b_{n}^{Wick}\varrho_{\varepsilon}^{1/2}\right]
= \int_{\Z_{0}}b_{n}(z)~d\mu(z)\,.
$$
Secondly,  the monotone convergence yields
$$
C\geq \sup_{n\in\nz}\int_{\Z_{0}}\langle z^{\otimes p}\,,\,
\tilde{b}_{n}z^{\otimes p}\rangle~d\mu(z)= \int_{\Z_{0}} b(z)~d\mu(z)\,.
$$
When $\tilde{b}$ is unbounded, it can be approximated by
$\tilde{b}_{n}=\frac{\tilde{b}}{1+\frac{\tilde{b}}{n}}\in
\mathcal{L}(\bigvee^{p}\Z_{0})$\,,
 for $n\geq 1$\,. Set $b_{n}(z)=\langle z^{\otimes p}\,,\,
 \tilde{b}_{n}z^{\otimes p}\rangle$\,. The function $b(z)=\langle z^{\otimes p}\,,\,
 \tilde{b}z^{\otimes p}\rangle=\sup_{n\in\nz}b_{n}(z)$ is a Borel
 function on $\Z_{0}$ as a supremum of a sequence of continuous
 functions. The uniform estimate
$$
\Tr\left[\varrho_{\varepsilon}b_{n}^{Wick}\right]\leq
\Tr\left[\varrho_{\varepsilon}b^{Wick}\right]\leq C
$$
with the result for $\tilde{b}_{n}\in \mathcal{L}(\bigvee^{p}\Z_{0})$
gives
$\int_{\Z_{0}}b_{n}(z)~d\mu(z)\leq C\,,$
for all $n\in\nz^{*}$\,. Again by monotone convergence, we get
$$
\int_{\Z_{0}}b(z)~d\mu(z)=\sup_{n\in\nz^{*}}\int_{\Z_{0}}b_{n}(z)~d\mu(z)\leq C\,.
$$
\fin
\begin{lem}
\label{le.estimA}
Let $A$ be a non negative, self-adjoint with domain $\D(A)$\,, operator
in $\Z_{0}$\,.
Assume that the family $(\varrho_{\varepsilon})_{\varepsilon\in
  (0,\bar \varepsilon)}$ satisfies the uniform estimate
$\Tr\left[\varrho_{\varepsilon}\N^{\alpha}\right]\leq C_{\alpha}$\,, for
all $\alpha\in\nz$\,, and $\mathcal{M}(\varrho_{\varepsilon},\,
\varepsilon\in (0,\bar\varepsilon))=\left\{\mu\right\}$\,.
Then the following implication hold:
\begin{eqnarray*}
   \left(\forall\varepsilon\in (0,\bar\varepsilon),\;
 \Tr\left[\varrho_{\varepsilon}\d\Gamma(A)\right]\leq C\right)
 &\Rightarrow&
 \left(\int_{\Z_{0}}\langle z\,,\,Az\rangle~d\mu(z)\leq C\right)\,,\\
   \left(\forall\varepsilon\in (0,\bar\varepsilon),\;
 \Tr\left[\varrho_{\varepsilon}\d\Gamma(A)^{2}\right]\leq C\right)
 &\Rightarrow&
 \left(\int_{\Z_{0}}\langle z\,,\,Az\rangle^{2}~d\mu(z)\leq C\right)\,,\\
  \left(\forall\varepsilon\in (0,\bar\varepsilon),\;
 \Tr\left[\varrho_{\varepsilon}\d\Gamma(A)^{2}\N\right]\leq C\right)
 &\Rightarrow&
 \left(\int_{\Z_{0}}\langle z\,,\,Az\rangle^{2}|z|_{\Z_{0}}^{2}~d\mu(z)\leq C\right)\,.
\end{eqnarray*}
In all the three cases, the measure $\mu$ is carried by the form
domain $Q(A)$ of $A$\,.
\end{lem}
\proof
The first implication is a direct application of
Lemma~\ref{le.poswick} applied with
$$
b(z)=\langle z\,,\, A z\rangle\quad,\quad \tilde{b}=A\quad,\quad b^{Wick}=\d\Gamma(A)\,.
$$
The second one is the consequence of
$$
\d\Gamma(A)^{2}= \left(\langle z^{\otimes 2}\,,\, (A\otimes A) z^{\otimes
  2}\rangle\right)^{Wick}+\varepsilon \d\Gamma(A^{2})\geq
 \left(\langle z^{\otimes 2}\,,\, (A\otimes A) z^{\otimes
  2}\rangle\right)^{Wick}\,
$$
and Lemma~\ref{le.poswick} with
$$
b(z)=\langle z\,,\, A z\rangle^{2}\quad\text{and}\quad
\tilde{b}=A\otimes A\,.
$$
For the last one, notice that $\N=\d\Gamma(1)$ and $\d\Gamma(A)$
commute so that
$$
\d\Gamma(A)^{2}\N\geq  \N\left(\langle z^{\otimes 2}\,,\, (A\otimes A) z^{\otimes
  2}\rangle\right)^{Wick}\,.
$$
With $\N=(|z|_{\Z_{0}}^{2})^{Wick}$\,, the composition formula of
Proposition~\ref{symbcalc} (extended to an unbounded $A$)
says that $\N\left(\langle z^{\otimes 2}\,,\, (A\otimes A) z^{\otimes
  2}\rangle\right)^{Wick}=b_{\varepsilon}^{Wick}$ with
$$
b_{\varepsilon}(z)= |z|_{\Z_{0}}^{2}\langle z\,,\, A z\rangle^{2}+
2\varepsilon \langle z\,,\, Az\rangle^{2}\,.
$$
Hence we get
$$
\d\Gamma(A)^{2}\N\geq \left(|z|_{\Z_{0}}^{2}\langle z\,,\, A z\rangle^{2}\right)^{Wick}\,.
$$
So, the result is again a consequence of Lemma~\ref{le.poswick} with
$$
b(z)=|z|_{\Z_{0}}^{2}\langle z\,,\, A z\rangle^{2}\quad, \quad
\tilde{b}=\frac{1}{3}(\11\otimes A\otimes A+ A\otimes \11\otimes A +
A\otimes A\otimes \11)\,.
$$
For the last statement it suffices to notice that the integrand is
infinite in the Borel subset of $\Z_{0}$\,,  $\Z_{0}\setminus Q(A)=\left\{z\in \Z_{0}, \langle
  z\,,Az\rangle=+\infty\right\}$\,.
\fin

\bigskip

\noindent\textbf{Proof of Proposition~\ref{pr.asapriori}:} With
$S_{\varepsilon}=\d\Gamma(1-\Delta)+\N^{3}$\,, while $\d\Gamma(1-\Delta)$ and
$\N$ commute, we know
$$
(1+S_{\varepsilon})(1+\N)\geq \d\Gamma(1-\Delta)\,.
$$
Hence Lemma~\ref{le.estimA} says that the measure $\mu$ is carried by
$Q(1-\Delta)=\Z_{1}$ with
\begin{equation}
  \label{eq.estimP2}
\int_{\Z_{0}}|z|_{\Z_{1}}^{2}~d\mu(z)=\int_{\Z_{1}}|z|_{\Z_{1}}^{2}~d\mu(z)\leq C\,.
\end{equation}
Let us check that $\mu$ is a Borel measure on
$(\Z_{1},|~|_{\Z_{1}})$\,. The tightness property is given by the
above inequality. According to \cite{AmNi1,Par,Sko}, it
suffices
 to check that
 \begin{eqnarray*}
&&G_{1}(\xi)=\int_{\Z_{1}}e^{-2i\pi\Real\langle \xi\,,\,
  z\rangle_{\Z_{1}}}~d\mu(z)\\
\text{with}
&&
\langle u\,, v\rangle_{\Z_{1}}= \langle u\,,\, (1-\Delta) v\rangle_{\Z_{0}}\,,
\end{eqnarray*}
is a positive type function which is continuous w.r.t $\xi$ restricted
to any finite dimensional subspace of $\Z_{1}$\,.\\
Consider the regularized version
$$
G_{1,n}(\xi)=\int_{\Z_{1}}e^{-2i\pi\Real\langle \frac{A}{1+\frac{A}{n}}\xi\,,\,
  z\rangle_{\Z_{0}}}~d\mu(z)
=
\int_{\Z_{0}}e^{-2i\pi\Real\langle \frac{A}{1+\frac{A}{n}}\xi\,,\,
  z\rangle_{\Z_{0}}}~d\mu(z)
$$
with $A=(1-\Delta)$\,. For  all $\xi \in \Z_{1}$
the pointwise convergence
$$
\forall z\in \Z_{1}\,,\quad
\lim_{n\to\infty}\langle \frac{A}{1+\frac{A}{n}}\xi\,,\,
  z\rangle_{\Z_{0}}=\langle \xi\,,\, z\rangle_{\Z_{1}}
$$
and the uniform bound
$$
|e^{-2i\pi\Real\langle \frac{A}{1+\frac{A}{n}}\xi\,,\,
  z\rangle_{\Z_{0}}}|\leq 1
$$
imply the pointwise convergence of the integrals
$$
\forall \xi\in \Z_{1},\quad \lim_{n\to \infty}G_{1,n}(\xi)=G_{1}(\xi)\,.
$$
But $G_{1,n}(\xi)$ equals $G((1+\frac{A}{n})^{-1}A\xi)$\,, where $G$ is
the characteristic function of $\mu$ in $\Z_{0}$:
$$
G(\eta)=\int_{\Z_{0}}e^{-2i\pi\Real\langle \eta\,,\,
  z\rangle_{\Z_{0}}}~d\mu(z)\,.
$$
Hence for every $n\in\nz$\,, the function $G_{1,n}(\xi)$ is a positive
type function. As a pointwise limit of $G_{1,n}$\,, the function $G_{1}$
is also a positive type function.\\
For the continuity, the equality
$$
G_{1}(\xi)-G_{1}(\xi')= \int_{\Z_{1}}
\left(
e^{-i\pi\Real\langle \xi-\xi',z\rangle_{\Z_{1}}}
- e^{i\pi\Real\langle
  \xi-\xi',z\rangle_{\Z_{1}}}\right)e^{-i\pi\Real\langle
\xi+\xi'\,,\,z\rangle_{\Z_{1}}}~d\mu(z)\,
$$
implies
\begin{equation}
  \label{eq.estimcontG1}
|G_{1}(\xi)-G_{1}(\xi')|\leq
2\pi|\xi-\xi'|_{\Z_{1}}\int_{\Z_{1}}|z|_{\Z_{1}}~d\mu(z)
\leq \pi \left(\int_{\Z_{1}}1+|z|^{2}_{\Z_{1}}~d\mu(z)\right)|\xi-\xi'|_{\Z_{1}}\,,
\end{equation}
and the function $G_{1}$ is a Lipschitz function on $\Z_{1}$\,.
This finishes the proof that $\mu$ is a Borel probability measure on
$\Z_{1}$\,.\\
For the inequality \eqref{eq.estimz6}, it suffices to notice the
inequality of (commuting) operators
$$
(1+S_{\varepsilon})^{2}(1+\N)\geq (\d\Gamma(1-\Delta)^{2})\N\,.
$$
Applying Lemma~\ref{le.estimA} yields
$$
\int_{\Z_{1}}|z|_{\Z_{1}}^{4}|z|_{\Z_{0}}^{2}~d\mu(z)\leq C\,.
$$
\fin
\subsection{Uniqueness of the mean field dynamics via measure transportation technique }
\label{se.meastransp}
Now we are in position to prove Theorem \ref{th.main}. This will be
done in three steps: 1) Writing a transport equation, in a weak sense
in $\Z_{1}$ for $\mu_{t}$\,; 2) Solving this equation as
$\mu_{t}=\Phi(t,0)_{*}\mu_{0}$ when the initial state
$\varrho_{\varepsilon}$ fulfills strong decay estimates; 3) Relaxing
the strong decay estimates.

\subsection{The transport equation on $\Z_1$}
\label{se.tranZ1}
We shall need similar notions about cylindrical functions, as the one
used in $\Z_{0}$ and recalled in Appendix~\ref{se.weylAwick}.
  Let  $\p_{1}$ denote the set of all finite rank orthogonal projections
on $\Z_{1}$ and  for a given $\wp\in\p_{1}$ let $L_{\wp,1}(dz)$ denote the
Lebesgue measure on the finite dimensional subspace $\wp\Z_{1}$\,, with
volume $1$ for a $\Z_{1}$-orthonormal hypercube. A
function $f:\Z_{1}\to\cz$ is said  cylindrical if there
exists $\wp\in\p_{1}$ and a function $g$ on $\wp\Z_{1}$ such that $
f(z)=g(\wp z),$
for all $z\in\Z_{1}$\,.
 In this case we say that $f$ is based on the
subspace $\wp\Z_{1}$\,. The set of $\mathcal{C}^{\infty}_{0}$ (resp. $\mathcal{S}$) cylindrical
functions on $\Z_{1}$\,, is denoted by
$\mathcal{C}^{\infty}_{0,cyl}(\Z_{1})$
(resp. $\mathcal{S}_{cyl}(\Z_{1})$). We shall also need
$\mathcal{C}^{\infty}_{0,cyl}(\Z_{1}\times \rz)$\,, in which the
algebraic tensor product
$\mathcal{C}^{\infty}_{0,cyl}(\Z_{1})\stackrel{alg}{\otimes}\mathcal{C}^{\infty}_{0}(\rz)$
is dense. Finally the Fourier transform of elements of
$\mathcal{S}_{cyl}(\Z_{1})$ is given by
\begin{eqnarray*}
\F_{1}[f](\xi)=\int_{\wp\Z_{1}} f(\xi) \;\;e^{-2\pi i
\,{\rm Re~}\langle z,\xi\rangle_{\Z_{1}}}~L_{\wp,1}(dz)\,,
\\
f(z)=\int_{\wp\Z_{1}} \F_{1}[f](\xi) \;\;e^{2\pi i
\,{\rm Re~}\langle z,\xi\rangle_{\Z_{1}}}~L_{\wp,1}(d\xi)\,.
\end{eqnarray*}
\begin{prop}
\label{pr.traZ1}
Let $S_{\varepsilon}$ and  $\tilde{\varrho}_{\varepsilon}(t)$  be the operators given by
Definition~\ref{de.opref} and \eqref{eq.deftildrho}.
Assume that the family of normal states
$\left(\varrho_{\varepsilon}\right)_{\varepsilon\in
  (0,\bar\varepsilon)}$ satisfies
$$
\forall \alpha\in\nz,\, \exists C_{\alpha}>0,\,
\forall
\varepsilon\in (0,\bar \varepsilon),\quad\Tr\left[(1+S_{\varepsilon})\varrho_{\varepsilon}(1+S_{\varepsilon})(1+\N)^{\alpha}\right]\leq C_{\alpha}\,,
$$
 and consider a subsequence $(\varepsilon_{k})_{k\in\nz}$\,,
 $\varepsilon_{k}\stackrel{k\to \infty}{\to}0$ such that
$$
\mathcal{M}(\tilde{\varrho}_{\varepsilon_{k}}(t),~k\in\nz)=\left\{\tilde{\mu}_{t}\right\}
$$
according to Proposition~\ref{conteq}. Then the measure
$\tilde{\mu}_{t}$ is a Borel probability measure on $\Z_{1}$ which
satisfies the following properties:
\begin{itemize}
\item $\int_{\Z_{1}}|z|_{\Z_{1}}^{2}~d\tilde{\mu}_{t}(z)+\int_{\Z_{1}}|z|_{\Z_{1}}^{4}|z|_{\Z_{0}}^{2}~d\tilde{\mu}_{t}(z)\leq
  C'$ for some $C'$ independent of $t\in\rz$\,.
\item When $(e_{n})_{n\in\nz^*}$ is a Hilbert basis of $\Z_{1}$ and
  $\Z_{1}$ is endowed with the distance $d_{\omega}(z_{1},z_{2})=\sqrt{\sum_{n\in\nz^*}\frac{|\langle
  z_{1}-z_{2},e_{n}\rangle|^{2}}{(1+n)^{2}}}$\,, $\tilde{\mu}_{t}$ is
narrowly continuous with respect to $t\in\rz$\,.
\item The measure $\mu_{t}$ is a solution to the Liouville equation
$$
\partial_{t}\tilde{\mu}_t+i\left\{V_{t}, \tilde{\mu}_{t}\right\}=0\,,
$$
in the weak sense,
\begin{equation}
  \label{eq.weakcont}
\forall f\in \mathcal{C}^{\infty}_{cyl,0}(\Z_{1}\times
\rz)\,,\quad
\int_{\rz}\int_{\Z_1}
\left(\partial_{t}f+i \left\{V_{t}, f\right\}\right)~d\tilde{\mu}_{t}(x)dt=0\,.
\end{equation}
\end{itemize}
\end{prop}
\proof
The Proposition~\ref{invH} as well as the commutations
$[e^{itH_{\varepsilon}},\N]=[e^{itH_{\varepsilon}^{0}},\N]=0$ ensure
$$
\forall \alpha\in\nz,\exists C_{\alpha}'>0,
\forall t\in\rz,\forall \varepsilon\in (0,\bar\varepsilon)\,,\quad
\Tr\left[(1+S_{\varepsilon})\tilde{\varrho}_{\varepsilon}(t)(1+S_{\varepsilon})(1+\N)^{\alpha}\right]\leq C_{\alpha}'\,.
$$
The Proposition~\ref{pr.asapriori} and \eqref{eq.estimP2} applied for
any $t\in\rz$\,, provides the first
results.\\
It remains to check  the narrow continuity and the Liouville equation.\\
\textbf{a)}
Take the $\Z_{1}$-characteristic function
$$
G_{1}(\eta,t)=\tilde{\mu}_{t}(e^{-2i\pi \Real\langle \eta\,,\, z\rangle_{\Z_{1}}})\,.
$$
The inequality \eqref{eq.estimcontG1}  and the uniform estimate
$\int_{\Z_{1}}(1+|z|_{\Z_{1}}^{2})~d\tilde{\mu}_{t}(z)\leq 1+C'$ ensures that the inequality
\begin{equation}
  \label{eq.contG1eta}
|G_{1}(\eta,t)-G_{1}(\eta',t)|\leq \pi(1+C')|\eta-\eta'|_{\Z_1}
\end{equation}
holds uniformly for all $\eta,\eta'\in\Z_{1}$ and all $t\in\rz$\,.
From the identity \eqref{fouriermu}, we deduce
$$
\tilde{\mu}_{t'}(e^{i \sqrt{2}{\rm Re~}\langle \xi, .\rangle})-\tilde{\mu}_{t}(e^{i \sqrt{2}{\rm
    Re}\langle \xi, .\rangle})
=-2\sqrt{2} i
\int_{t}^{t'}
\tilde{\mu}_s(e^{i \sqrt{2}{\rm Re~}\langle \xi, z\rangle} \,{\rm
  Im~}\langle z^{\otimes 2},
\tilde V_s \xi\otimes z\rangle) \,ds\,,
$$
The estimate \eqref{eq.estimVH1} implies
$$
|{\rm Im~}\langle z^{\otimes{2}}\,,\, \tilde{V}_{s}(\xi\otimes
z)\rangle|\leq C|z|_{\Z_{1}}^{2}|z|_{\Z_{0}}|\xi|_{H^{-1}(\rz^{d})}\,.
$$
Taking $\xi=\sqrt{2}\pi(1-\Delta)\eta$ with $\eta\in
H^{1}(\rz)=\Z_{1}$ leads to
$$
|G_{1}(\eta,t)-G_{1}(\eta,t')|\leq
4\pi C|\eta|_{\Z_{1}}|t-t'|\sup_{s\in[t,t']}\int_{\Z_{1}}|z|_{\Z_{1}}^{2}|z|_{\Z_{0}}~d\tilde{\mu}_{s}(z)
$$
and, with the uniform estimate $\int_{\Z_{1}}|z|_{\Z_{1}}^{4}|z|_{\Z_{0}}^{2}~d\tilde{\mu}_{t}(z)\leq
  C'$\,, to
\begin{equation}
    \label{eq.contG1t}
\forall \eta\in\Z_{1}, \forall t,t'\in\rz,\quad
|G_{1}(\eta,t)-G_{1}(\eta,t')|\leq 4\pi C(1+C')|\eta|_{\Z_{1}}|t-t'|\,.
\end{equation}
When $g\in \mathcal{S}_{cyl}(\Z_{1})$\,, based on $\wp\Z_{1}$\,, the
relation
$$
\int_{\Z_{1}}g(z)~d\tilde{\mu}_{t}(z)=\int_{\wp\Z_{1}}\mathcal{F}_{1}[g](\eta)G_{1}(\eta,t)~dL_{\wp,1}(z)\,,
$$
combined with the continuity properties
\eqref{eq.contG1eta}\eqref{eq.contG1t}, implies that $t\to
\int_{\Z_{1}}g(z)~d\tilde{\mu}_{t}(z)$ is continuous. This continuity
holds for all $g\in \mathcal{S}_{cyl}(\Z_{1})$\,. The
uniform weak tightness property
$\int_{\Z_{1}}|z|_{\Z_{1}}^{2}~d\tilde{\mu}_{t}(z)\leq C'$ and
Lemma~5.12-f) in \cite{AGS} ensure that $t\mapsto \tilde{\mu}_{t}$ is
narrowly continous when $\Z_{1}$ is endowed with the distance $d_{\omega}$\,.\\
\textbf{b)} Integrating \eqref{fouriermu}  with
$\mathcal{F}_{1}[g](\eta)~L_{\wp,1}(dz)$ also provides
$$
\int_{\Z_{1}}g(z)~d\tilde{\mu}_{t}(z)=\int_{\Z_{1}}g(z)~d\tilde{\mu}_{0}(z)
+i\int_{0}^{t}\int_{\Z_{1}}\left\{V_{s}\,,\, g\right\}(z)~d\tilde{\mu}_{s}(z)~ds\,.
$$
Hence for any $g\in \mathcal{S}_{cyl}(\Z_{1})$\,, the function $I_{g}:t\mapsto
\int_{\Z_{1}}g(z)~d\tilde{\mu}_{t}(z)$ belongs to
$\mathcal{C}^{1}(\rz)$ with
$$
\partial_{t}I_{g}(t)=i\int_{\Z_{1}}\left\{V_{t}\,,\, g\right\}(z)~d\tilde{\mu}_{t}(z)\,.
$$
By multiplying the above relation by $\varphi(t)$\,, with $\varphi\in
\mathcal{C}^{\infty}_{0}(\rz)$\,, and integrating by part proves
\eqref{eq.weakcont} when $f(t,z)=\varphi(t)g(z)$\,.
We conclude by the density of
$\mathcal{C}^{\infty}_{0,cyl}(\Z_{1})\stackrel{alg}{\otimes}\mathcal{C}^{\infty}_{0}(\rz)$
in $\mathcal{C}^{\infty}_{0,cyl}(\Z_{1}\times\rz)$\,.
\fin

\subsection{Uniqueness of the measure for regular initial data}
\label{se.reguniq}
According to the notations of \cite{AGS} and Appendix~\ref{se.abscont},
we consider the space $\textrm{Prob}_{2}(\Z_{1})$ of Borel probability
measures $\mu$ such that
$$
\int_{\Z_{1}}|z|_{\Z_{1}}^{2}~d\mu(z)< +\infty\,.
$$
On this space, we introduce the Wasserstein distance
\begin{equation}
  \label{eq.was2}
W_2(\mu_{1},\mu_{2})=\left[\min_{\mu\in
  \Gamma(\mu_{1},\mu_{2})}\int_{\Z_{1}^{2}}|z_{1}-z_{2}|_{\Z_{1}}^{2}~d\mu(z_{1},z_{2})
\right]^{1/2}
\end{equation}
where $\Gamma(\mu_{1},\mu_{2})$ is the set of Borel probability
measures $\mu$ on $\Z_{1}\times \Z_{1}$ with the marginals
$(\Pi_{1})_{*}\mu=\mu_{1}$ and $(\Pi_{2})_{*}\mu=\mu_{2}$\,.
\begin{prop}
  \label{pr.regdata}
Let $S_{\varepsilon}$ and  $\tilde{\varrho}_{\varepsilon}(t)$  be the operators given by
Definition~\ref{de.opref} and \eqref{eq.deftildrho}.
Assume that the family of normal states
$\left(\varrho_{\varepsilon}\right)_{\varepsilon\in
  (0,\bar\varepsilon)}$ satisfies
\begin{eqnarray*}
  &&
\forall \alpha\in\nz,\exists C_{\alpha}>0,
\forall
\varepsilon\in (0,\bar \varepsilon),\quad\Tr\left[(1+S_{\varepsilon})\varrho_{\varepsilon}(1+S_{\varepsilon})(1+\N)^{\alpha}\right]\leq C_{\alpha}
\\
\text{and}
&&
\mathcal{M}(\varrho_{\varepsilon}, \varepsilon\in (0,\bar
\varepsilon))=\left\{\mu_{0}\right\}\,.
\end{eqnarray*}
Then for any time $t\in \rz$\,, the family
$(\varrho_{\varepsilon}(t)=e^{-i\frac{t}{\varepsilon}H_{\varepsilon}}\varrho_{\varepsilon}e^{i\frac{t}{\varepsilon}H_{\varepsilon}})_{\varepsilon\in
(0,\bar \varepsilon)}$
admits a unique Wigner measure $\mu_{t}=\Phi(t,0)_{*}\mu_{0}$\,, where
$\Phi$ is the Hartree flow defined by \eqref{eq.hartintro} on
$\Z_{1}$\,. It is a Borel probability measure on $\Z_{1}$ with $t\mapsto \mu_t$
being  an absolutely continuous curve in $\textrm{Prob}_{2}(\Z_{1})$\, w.r.t
the Wasserstein distance $W_2$ and which satisfies
$$
\forall t\in\rz\,,\quad\int_{\Z_{1}}|z|_{\Z_{1}}^{4}|z|_{\Z_{0}}^{2}~d\mu_{t}(z)\leq C\,.
$$
\end{prop}
\proof
We still start with the state $\tilde{\varrho}_{\varepsilon}(t)$
defined in \eqref{eq.deftildrho}. Proposition~\ref{pr.hartree} says
that the group $\Phi(t,s)$  associated with \eqref{eq.hartintro} and the
dynamical system
$\tilde{\Phi}(t,s)$ associated with
$$
i\partial_{t}z=v(t,z)\quad,\quad v(t,z)=e^{-it\Delta} ([V*|e^{it\Delta}z|^{2}]e^{it\Delta}z)
$$
are well defined on $\Z_{1}$\,. Further it gives the estimate
for the velocity field
$$
|v(t,z)|_{\Z_{1}}\leq \|V(1-\Delta)^{-1/2}\|\;|z|_{\Z_{1}}^{2}|z|_{\Z_{0}}\,.
$$
When $\tilde{\mu}_{t}$ is the Wigner measure defined for all times and
associated with a subsequence $(\varepsilon_{n_{k}})_{k\in\nz}$\,,
we obtain
$$
\forall t\in\rz,\quad
\int_{\Z_{1}}|v(t,z)|_{\Z_{1}}^{2}~d\tilde{\mu}_{t}(z)\leq C
\int_{\Z_{1}}|z|_{\Z_{1}}^{4}|z|_{\Z_{0}}^{2}~d\tilde{\mu}_{t}(z)\leq C'\,.
$$
With Proposition~\ref{pr.traZ1}, $t\mapsto \tilde{\mu}_{t}$ is narrowly
continuous (on $(\Z_{1},d_{\omega})$) with respect to time
$t\in\rz$\,.
 According to Lemma~\ref{le.Poidiv},  the Liouville equation
 \eqref{eq.weakcont} is
 nothing but the weak form of
$$
\partial_{t}\mu+\nabla^{T}(v(t,z)\mu)=0\,.
$$
According to Proposition~\ref{pr.conteq}, the curve $\rz\ni
t\mapsto\tilde{\mu}_{t}\in \textrm{Prob}_{2}(\Z_{1})$ is
absolutely continuous for the Wasserstein distance $W_{2}$\,.\\
Therefore all the conditions of Proposition~\ref{pr.comptrmeas} are fulfilled and hence we
deduce that $\tilde{\mu}_{t}=\tilde{\Phi}(t,0)_{*}\mu_{0}$\,. 
Moreover this uniqueness implies
$\mathcal{M}(\tilde{\varrho}_{\varepsilon}(t), \varepsilon\in
(0,\bar{\varepsilon}))=\left\{\tilde{\mu}_{t}\right\}$ for the whole family
$(\tilde{\varrho}_{\varepsilon}(t))_{\varepsilon\in (0,\bar
  \varepsilon)}$ and all times $t\in\rz$\,.\\
Going back to
$\varrho_{\varepsilon}(t)=e^{-i\frac{t}{\varepsilon}H_{\varepsilon}^{0}}\tilde{\varrho}_{\varepsilon}(t)
e^{i\frac{t}{\varepsilon}H_{\varepsilon}^{0}}$\,, it gives
$\mu_{t}=\Phi(t,0)_{*}\mu_{0}$\,.\\
The last uniform estimate is given by Proposition~\ref{pr.traZ1}.
\fin

\subsection{Evolution of the Wigner measure for general data}
\label{se.gendata}
 We follow the truncation scheme used in \cite{AmNi3}. When the
 initial data satisfies only
$$
\|(\N+H_{\varepsilon}^{0})^{\delta/2}\varrho_{\varepsilon}(\N+H_{\varepsilon}^{0})^{\delta/2}\|\leq C_{\delta}
$$
for some $\delta>0$\,, we approximate $\varrho_{\varepsilon}$ by
$$
\varrho_{\varepsilon,R}=\frac{1}{\Tr\left[\chi_{R}(\N, H_{\varepsilon}^{0})\varrho_{\varepsilon}\chi_{R}(\N, H_{\varepsilon}^{0})\right]}
\chi_{R}(\N, H_{\varepsilon}^{0})\varrho_{\varepsilon}\chi_{R}(\N, H_{\varepsilon}^{0})
$$
as $R\to +\infty$ where $\chi_{R}(n,h)=\chi(\frac{n}{R}\,, \frac{h}{R})$,  with $0\leq \chi\leq 1$\,, $\chi\in
\mathcal{C}^{\infty}_{0}(\rz^{2})$ and $\chi\equiv 1$ in a
neighborhood of $0$\,. The time evolved state is defined by
$$
\varrho_{\varepsilon,R}(t)=e^{-i\frac{t}{\varepsilon}H_{\varepsilon}}\varrho_{\varepsilon,R}\,
e^{i\frac{t}{\varepsilon}H_{\varepsilon}}\,.
$$
The assumptions ensure that for all times
$$
\|\varrho_{\varepsilon}(t)-\varrho_{\varepsilon,R}(t)\|_{\mathcal{L}^{1}(\H)}\leq \nu(R)
$$
with $\nu$ independent of $(t,\varepsilon)$ and $\lim_{R\to
  \infty}\nu(R)=0$\,.
We recall the Proposition~2.10 of \cite{AmNi3}.
 \begin{prop}
\label{pr.compwig}
   Let $(\varrho_{\varepsilon}^{j})_{\varepsilon\in
     (0,\bar\varepsilon)}$\,, $j=1,2$\,, be two families (or sequences) of normal states
   on $\H$ such that
   $\Tr\left[\varrho_{\varepsilon}^{j}\,\N^{\delta}\right]\leq
   C_{\delta}$ uniformly w.r.t $\varepsilon\in (0,\bar \varepsilon)$
   for some $\delta>0$ and $C_{\delta}\in (0,+\infty)$\,. Assume further
   $\mathcal{M}(\varrho_{\varepsilon}^{j}, \varepsilon\in (0,\bar
   \varepsilon))=\left\{\mu_{j}\right\}$ for $j=1,2$\,. Then
$$
\int|\mu_{1}-\mu_{2}|\leq \liminf_{\varepsilon\to 0}
\|\varrho_{\varepsilon}^{1}-\varrho_{\varepsilon}^{2}\|_{\mathcal{L}^{1}(\H)}\,.
$$
 \end{prop}
\noindent\textbf{End of the proof of Theorem~\ref{th.main}:}
For $R\in (0,+\infty)$\,, the state $\varrho_{\varepsilon,R}$ fulfills
the conditions of Proposition~\ref{pr.regdata} except the uniqueness
of the Wigner measure at time $t=0$\,.
Out of any sequence $(\varepsilon_{n})_{n\in\nz}$\,, a
subsequence $(\varepsilon_{n_{k}})_{k\in\nz}$ can be extracted in order to ensure
$$
\mathcal{M}(\varrho_{\varepsilon_{n_{k}},R},
k\in\nz)=\left\{\mu_{0,R}\right\}\,.
$$
Thus after this extraction we obtain
$$
\forall t\in\rz\,,\quad \mathcal{M}(\varrho_{\varepsilon_{n_{k}},R}(t),
k\in\nz)=\left\{\Phi(t,0)_{*}\mu_{0,R}\right\}\,.
$$
Take $t\in \rz$
 and let $\mu$ belong to $\mathcal{M}(\varrho_{\varepsilon}(t),
\varepsilon\in (0,\bar \varepsilon))$\,. There exists a sequence
$(\varepsilon_{n})_{n\in\nz}$ such that
$$
\mathcal{M}(\varrho_{\varepsilon_{n}}(t),
n\in\nz)=\left\{\mu\right\}\,.
$$
After extracting a subsequence like above and by using
Proposition~\ref{pr.compwig}, we obtain
$$
\int|\mu-\Phi(t,0)_{*}\mu_{0}|\leq
\int|\mu-\Phi(t,0)_{*}\mu_{0,R}|+\int|\mu_{0,R}-\mu_{0}|
\leq 2\nu(R)\,,
$$
since the total variation of $\Phi(t,0)_{*}\mu_{0,R}-\Phi(t,0)_{*}\mu_{0}$ and 
$\mu_{0,R}-\mu_{0}$ are equal. 
Taking the limit as $R\to \infty$ implies $\mu=\Phi(t,0)_{*}\mu_{0}$ and
therefore
$$
\mathcal{M}(\varrho_{\varepsilon}, \varepsilon\in
(0,\bar\varepsilon))=
\left\{\Phi(t,0)_{*}\mu_{0}\right\}\,.
$$
This also proves that
$\lim_{R\to\infty}\int_{\Z_{0}}|\mu_{t}-\Phi(t,0)_{*}\mu_{0,R}|=0$\,, while all
the measures $\Phi(t,0)_{*}\mu_{0,R}$ are Borel probability measures
 carried by, and on, $\Z_{1}$\,.  This implies that
$\mu_{t}$ is carried by $\Z_{1}$ and is also a Borel measure on $\Z_{1}$\,.
This ends the proof of Theorem~\ref{th.main}
\fin
\bigskip
\noindent

\section{Complements}
\label{se.addppties}
Additional results are given in the three first paragraphs, concerned
with the BBGKY hierarchy or the propagation of energy. The fourth one
shows some examples and the last one is an informal discussion about
the classical mean field problem.

\subsection{BBGKY hierarchy}
\label{se.BBGKY}
Although the analysis here is different from our previous work \cite{AmNi3}
it is possible to combine  them, in order to strengthen the result of Theorem \ref{th.main}.
It is also interesting to reformulate our result in terms of reduced density matrices
since, in the literature, several mathematical results on mean field limit use the
 BBGKY hierarchy method (see for example \cite{BGM, BEGMY, KlMa}).
For a family of normal states
$(\varrho_\varepsilon)_{\varepsilon\in(0,\bar\varepsilon)}$ on $\H$
and $p\in\nz$\,,
the reduced density matrices $\gamma^{(p)}_\varepsilon\in
\L^1(L_s^2(\rz^{dp}))$ is defined according to
\begin{equation}
  \label{eq.relwickbbgky}
\Tr\left[\gamma^{(p)}_{\varepsilon}\tilde{b}\right]
=
\frac{\Tr\left[\varrho_{\varepsilon}\right]}{\Tr\left[\varrho_{\varepsilon}(|z|^{2p})^{Wick}\right]}
\Tr\left[\varrho_{\varepsilon}b^{Wick}\right]\,,\quad \forall \tilde b\in \L(L_s^2(\rz^{dp}))\,,
\end{equation}
with the convention that the right-hand side is $0$
when ${\Tr\left[\varrho_{\varepsilon}(|z|^{2p})^{Wick}\right]}=0$ and $p>0$\,.

\begin{thm}
\label{Piwig}
Let $(\varrho_{\varepsilon})_{\varepsilon\in(0,\bar\varepsilon)}$
be a family of normal states on $\H$\,, satisfying the hypothesis
of Theorem \ref{th.main}, with a single Wigner measure $\mu_{0}$
such that
\begin{equation}
\label{eq.hypconv}
\forall \alpha\in\nz,
\quad
\lim_{\varepsilon\to 0}
{\rm Tr}[\varrho_{\varepsilon}
\mathbf{N}^\alpha]=\int_{\Z_0}|z|^{2\alpha}~d\mu_{0}(z) < +\infty\,.
\end{equation}
Then for all $t\in \rz$\,, the convergence
$$
\lim_{\varepsilon\to
  0}\Tr\left[\varrho_{\varepsilon}(t)b^{Wick}\right]=\int_{\Z_0} b(\Phi(t,0)z)~d\mu_{0}(z)
=\int_{\Z_0} b(z)~d\mu_{t}(z)
$$
holds for any
$b\in\P_{alg}(\Z_0)=\oplus_{p,q\in\nz}^{alg}\P_{p,q}(\Z_0)$\,, with
$\mu_{t}=\Phi(t,0)_{*}\mu_{0}$\,.\\
Finally, the convergence of the reduced density matrices
$$
\lim_{\varepsilon\to
  0}\gamma^{(p)}_{\varepsilon}(t)=\frac{1}{\int_{\Z_0}|z|^{2p}~d\mu_{t}(z)}
\int_{\Z_0}|z^{\otimes p}\rangle\langle z^{\otimes p}|~d\mu_{t}(z)\,,
$$
holds in the $\mathcal{L}^{1}(L_s^2(\rz^{dp}))$-norm for all $p\in\nz$\,.
\end{thm}
\proof
By Theorem \ref{th.main} the
the family of normal states $(\varrho_\varepsilon(t))_{\varepsilon\in(0,\bar\varepsilon)}$
admits a single Wigner measure $\mu_t$ equal to
${\Phi(t,0)}_*\mu_0$\,. Since the quantum and classical flows preserve
the total number, the state $\varrho_\varepsilon(t)$ satisfies as well the
condition
\eqref{eq.hypconv}  for any time $t\in\rz$\,.
Then  \cite[Proposition 2.11, 2.13]{AmNi3} provide the claimed results.
\fin

\subsection{Moment upper bounds}
\label{se.momupper}

In \cite{AmNi1}, it was proved that the sole a priori estimate
$\Tr\left[\varrho_{\varepsilon}\N^{\delta}\right]\leq C_{\delta}$ for
a given $\delta>0$ (possibly small), with
$\mathcal{M}(\varrho_{\varepsilon}, \varepsilon\in
(0,\bar\varepsilon))=\left\{\mu\right\}$ leads to
$$
\int_{\Z_{0}}|z|^{2\delta}~d\mu(z) <+\infty\,.
$$
The a priori estimate, assumed in Theorem~\ref{th.main} at time $t=0$\,,
leads to
$$
\int_{\Z_{1}}|z|_{\Z_{1}}^{2\delta}~d\mu(z)< +\infty\,,
$$
according to the following result which is a variation of Lemma~\ref{le.poswick}.
\begin{prop}
\label{pr.momest}
Let $(A,\D(A))$ be a self-adjoint operator on $\Z_{0}$ such that $A\geq
\11$\,. If the family of normal states
$(\varrho_{\varepsilon})_{\varepsilon\in (0,\bar\varepsilon)}$
satisfies
$\Tr\left[\varrho_{\varepsilon}(\d\Gamma(A))^{\delta}\right]\leq
C_{\delta}$ for some $\delta>0$ and
$\mathcal{M}\left\{\varrho_{\varepsilon},
\varepsilon\in (0,\bar\varepsilon)\right\}=\left\{\mu\right\}$\,, then
$$
\int_{\Z_{0}}\langle z\,,\,Az\rangle^{\delta}~d\mu(z)< +\infty\,.
$$
\end{prop}
\proof
By Wick calculus (see Proposition~\ref{symbcalc} when $A$ is bounded), one gets
$$
\d\Gamma(A)^{k}\geq \left(\langle z\,,\,A
  z\rangle^{k}\right)^{Wick}\,, \forall k\in\nz\,.
$$
Let $(e_{j})_{j\in \nz}$ be an orthonormal basis of $\Z_{0}$ such that
$e_{j}\in \D(A)$ for all $j\in\nz$\,,  and set
$$
A_{J}=\sum_{j=0}^{J}A^{1/2}|e_{j}\rangle \langle
e_{j}|A^{1/2}=\sum_{j=0}^{J}|A^{1/2}e_{j}\rangle \langle A^{1/2}e_{j}|\,.
$$
 The inequality
$A^{\otimes k}\geq A_{J}^{\otimes k}$ holds for all $J\in\nz$\,, while $\tilde{b}\to b^{Wick}$ is
operator monotone when
restricted to operators $\tilde{b}$ acting in $\bigvee^{2k} \Z_0$\,.
Therefore, we obtain
\begin{eqnarray}
\nonumber  (1+\d\Gamma(A))^{n}
&=&
\sum_{k=0}^{n}C_{n}^{k}\d\Gamma(A)^{k}\geq
\left(\sum_{k=0}^{n}C_{n}^{k}\langle z\,,\,
  Az\rangle^{k}\right)^{Wick}
\\
\label{eq.AAJWick}
&\geq & \left(\sum_{k=0}^{n}C_{n}^{k}\langle z\,,\,
  A_{J}z\rangle^{k}\right)^{Wick}
=
\left[(1+\langle z\,,\, A_{J}z\rangle)^{n}\right]^{Wick}\,.
\end{eqnarray}
We shall use the same argument as the one in \cite[Theorem 6.2]{AmNi1} when
$A=\mathrm{Id}$\,, relying on the semiclassical calculus in finite dimension
(see \cite{BoLe,Hor,Mar,NaNi,Rob}).\\
Let $\wp_{J}$ be the orthogonal projection from $\Z_{0}$ onto
$\oplus_{j=0}^{J}A^{1/2}e_{j}$\,. The symbol $1+\langle z\,,\,
A_{J}z\rangle$ is a cylindrical symbol based on $\wp_{J}\Z_{0}$\,.
Since $\ker{A^{1/2}}=\left\{0\right\}$ and $e_{0},\ldots, e_{J}$
are linearly independent, the symbol
$$
(1+\langle z\,,\, A_{J}z\rangle)=1+\sum_{j=0}^{J} |\langle A^{1/2}e_{j},z\rangle|^{2}
$$
is an elliptic symbol on $\wp_{J}\Z_{0}\sim \cz^{J+1}$ in the
H{\"o}rmander class
$\mathcal{S}(1+|z|_{\cz^{J+1}}^{2},
\frac{|dz|_{\cz^{J+1}}^{2}}{1+|z|_{\cz^{J+1}}^{2}})$\,.
The functional calculus of Weyl $\varepsilon$-quantized elliptic operators  in
finite dimensions gives
\begin{equation}
  \label{eq.AJWeyl}
\forall s\in\rz,  \left[(1+\langle z\,,\,
  A_{J}z\rangle)^{Weyl}\right]^{s}
\geq (1-C_{J,s}\varepsilon)\left[(1+\langle z\,,\,
  A_{J}z\rangle)^{s}\right]^{Weyl}\,.
\end{equation}
The finite dimensional comparison of Wick and Weyl quantization, also
gives
\begin{equation}
  \label{eq.AJWW}
\forall n\in\nz, \left[(1+\langle z\,,\,
  A_{J}z\rangle)^{n}\right]^{Wick}\geq (1-C_{J,n}\varepsilon)\left[(1+\langle z\,, A_{J}z\rangle)^{Weyl}\right]^{n}\,.
\end{equation}
From \eqref{eq.AAJWick}\eqref{eq.AJWW} and the operator monotonicity
of $B\to B^{t}$ for $t\in (0,1]$\,, we deduce
$$
\forall s\in\rz\,,\quad (1+\d\Gamma(A))^{s}\geq
(1-C_{J,s}'\varepsilon)\left[(1+\langle z\,,\,
  A_{J}z\rangle)^{Weyl}\right]^{s}\,,
$$
and \eqref{eq.AJWeyl} gives
\begin{equation}
  \label{eq.minoJs}
\forall s\in\rz\,,\quad (1+\d\Gamma(A))^{s}\geq
(1-C_{J,s}''\varepsilon)\left[(1+\langle z\,,\,
  A_{J}z\rangle)^{s}\right]^{Weyl}\,.
\end{equation}
The definition of Wigner measures, recalled in
Theorem~\ref{th.wig-measure}, says
$$
\lim_{\varepsilon\to 0}\Tr\left[\varrho_{\varepsilon}b^{Weyl}\right]=\int_{\Z_{0}}b(z)~d\mu(z)\,,
$$
for all $b\in \mathcal{S}_{cyl}(\Z_{0})$\,, in particular the $b$'s
based on $\wp_{J}\Z_{0}$\,.
Take now $s=\delta$ in \eqref{eq.minoJs}. The a priori estimate
$$
\Tr\left[\varrho_{\varepsilon}\left[(1+\langle z\,,\,
    A_{J}z\rangle)^{\delta}\right]^{Weyl}\right]\leq (1+C_{J,\delta}\varepsilon)
\Tr\left[\varrho_{\varepsilon}(1+\d\Gamma(A))^{\delta}\right]\leq C_{\delta}(1+C_{J,\delta}\varepsilon)\,,
$$
and the ellipticity of $(1+\langle z\,,\, A_{J}z\rangle)^{\delta}$
allows to extend the above convergence to any cylindrical $b=f\circ
\wp_{J}$ with $f\in \mathcal{S}((1+|z|_{\cz^{J+1}}^{2})^{\delta},
\frac{|dz|_{\cz^{J+1}}^{2}}{1+|z|_{\cz^{J+1}}^{2}})$\,. In particular,
this leads to
$$
C_{\delta}'\geq \limsup_{\varepsilon\to
  0}\Tr\left[\varrho_{\varepsilon}(1+\d\Gamma(A))^{\delta}\right]
\geq \lim_{\varepsilon\to 0}\Tr\left[\varrho_{\varepsilon}\left[(1+\langle z\,,\,
    A_{J}z\rangle)^{\delta}\right]^{Weyl}\right]
=\int_{\Z_{0}}(1+\langle z\,, A_{J}z\rangle)^{\delta}~d\mu(z)\,.
$$
Since $A=\sup_{J}A_{J}$ with  $A_{J'}\geq A_{J}$ for $J'\geq J$\,, the
monotone convergence implies
$$
C_{\delta}'\geq \int_{\Z_{0}}(1+\langle z\,,\,
Az\rangle)^{\delta}~d\mu(z)\geq \int_{\Z_{0}}\langle z\,,\,Az\rangle^{\delta}~d\mu(z)\,.
$$
\fin
\begin{prop}
Within the framework of Theorem~\ref{th.main} with the assumption
${\rm Tr}[(\mathbf{N}+H_{\varepsilon}^{0})^{\delta}\varrho_{\varepsilon}
] \leq C_{\delta}$ for $\delta\leq 6$\,,
 the measure
  $\mu_{t}$ satisfies the additional estimate
$$
\int_{\Z_{1}}|z|_{\Z_{1}}^{\frac{2\delta}{3}}~d\mu_{t}(z)\leq C_{\delta}
$$
for all times $t\in\rz$\,.
\end{prop}
\proof
The functional calculus of commuting operators implies
$$
c_{\delta}(1+S_{\varepsilon})^{\delta/3}\leq
(1+\N+H_{\varepsilon}^{0})^{\delta}\leq (1+S_{\varepsilon})^{\delta}\,.
$$
Thus the initial state $\varrho_{\varepsilon}$\,, satisfies
$$
\Tr\left[\varrho_{\varepsilon}(1+S_{\varepsilon})^{\delta/3}\right]\leq
C_{\delta}'\,.
$$
From Proposition~\ref{invH}, we deduce
$$
e^{i\frac{t}{\varepsilon}H_{\varepsilon}}(1+S_{\varepsilon})^{2}e^{-i\frac{t}{\varepsilon}H_{\varepsilon}}
\leq C(1+S_{\varepsilon})^{2}\,.
$$
Since $B\to B^{s}$ is operator monotone for $s\in (0,1]$\,, this implies
$$
e^{i\frac{t}{\varepsilon}H_{\varepsilon}}(1+S_{\varepsilon})^{\delta/3}e^{-i\frac{t}{\varepsilon}H_{\varepsilon}}
\leq C^{\delta/6}(1+S_{\varepsilon})^{\delta/3}
$$
as soon as $\frac{\delta}{3}\leq 2$\,.
The inequality
$$
\Tr\left[\varrho_{\varepsilon}(t)(\d\Gamma(1-\Delta))^{\delta/3}\right]
\leq
\Tr\left[\varrho_{\varepsilon}(t)(1+S_{\varepsilon})^{\delta/3}\right]
\leq C_{\delta}'
$$
and the previous Proposition~\ref{pr.momest} applied with
$A=(1-\Delta)$\,, yields the result.
\fin
A more accurate version of this last result is given below by making
use of the conservation of energy.
\subsection{Convergence of moments and energy conservation}
\label{se.proen}
For a family $(\varrho_{\varepsilon})_{\varepsilon\in
  (0,\bar\varepsilon)}$ of normal states with a single Wigner measure $\mu_0$
the condition \eqref{eq.hypconv} is an important and non trivial assumption.
Indeed, we proved in \cite{AmNi3} the following equivalence
\begin{eqnarray}
\label{cond-PI}
\left(\forall \alpha\in\nz,~
  \lim_{\varepsilon\to 0}\Tr\left[\varrho_{\varepsilon}{\bf N}^{\alpha}\right]
=\int_\Z |z|^{2\alpha}~d\mu_0(z)\right)
\Leftrightarrow
\left(
\forall b\in \P_{alg}(\Z_{0}),~
\lim_{\varepsilon\to
  0}\Tr\left[\varrho_{\varepsilon}b^{Wick}\right]=\int_\Z b(z)~d\mu_0\right)\,.
\end{eqnarray}
Hence the condition \eqref{eq.hypconv}, although it involves only the
number operator, is exactly the one which leads to a good asymptotic
behaviour of the reduced density matrices.
\begin{prop}
Let $(\varrho_{\varepsilon})_{\varepsilon\in(0,\bar\varepsilon)}$
be a family of normal states on $\H$\,, satisfying the hypothesis
of Theorem \ref{th.main}, with a single Wigner measure $\mu_{0}$\,.
Assume $\Tr[\varrho_\varepsilon\; \N^\alpha]\leq C_\alpha$ uniformly w.r.t
$\varepsilon\in(0,\bar \varepsilon)$\,, for all $\alpha\in\nz$\,. \\
Then for every $\alpha\in\nz$\,, the quantity
$$
 \liminf_{\varepsilon\to 0} \;\Tr[\varrho_\varepsilon(t) \,\N^\alpha]-\int_{\Z_0} |z|_{\Z_0}^{2\alpha} \;d\mu_t(z),
$$
does not depend on time when $\varrho_\varepsilon(t)=e^{-i\frac{t}{\varepsilon}H_{\varepsilon}}
\varrho_{\varepsilon}e^{i\frac{t}{\varepsilon}H_{\varepsilon}}$ and
$\mathcal{M}(\varrho_\varepsilon(t), \varepsilon\in(0,\bar\varepsilon))=\left\{\mu_{t}\right\}$\,.\\
The condition \eqref{eq.hypconv} is satisfied
by $(\varrho_\varepsilon(t))_{\varepsilon\in
(0,\bar \varepsilon)}$ and $\mu_t$\,, for all times $t\in\rz$\,, as soon
as it is true for one $t_{0}\in\rz$\,.
\end{prop}
\proof
According to Theorem \ref{th.main}, we know that $\mathcal{M}(\varrho_\varepsilon(t), \varepsilon\in(0,\bar\varepsilon))=\{\mu_t\}$ with $\mu_{t}=\Phi(t,0)_{*}\mu_{0}$\,. The conservation
of the $|.|_{\Z_0}$-norm  by the nonlinear flow  $\Phi(t,0)$ yields
$$
\int_{\Z_0} |z|_{\Z_0}^{2\alpha} \;d\mu_t(z)=\int_{\Z_0} |z|_{\Z_0}^{2\alpha} \;d\mu_0(z)
$$
for any $t\in\rz$\,. On the other hand, $H_\varepsilon$ and $\N$ are strongly commuting self-adjoint
operators and therefore $\Tr[\varrho_\varepsilon(t)
\,\N^\alpha]=\Tr[\varrho_\varepsilon \N^\alpha]$ for every
 $\alpha\in \nz$\,.
\fin

\begin{prop}
Let $(\varrho_{\varepsilon})_{\varepsilon\in(0,\bar\varepsilon)}$
be a family of normal states on $\H$ with a single Wigner measure $\mu_{0}$ satisfying the hypothesis
of Theorem \ref{th.main} with $\delta=2$ and condition \eqref{eq.hypconv}. Then for any $t\in\rz$
\begin{equation}
\label{eq.energy}
\lim_{\varepsilon\to 0}
{\rm Tr}[\varrho_{\varepsilon}(t)
\, H_\varepsilon]=\int_{\Z_1} h(z,\bar z)~d\mu_{t}(z)\quad\in (-\infty,\infty)\,
\end{equation}
where $\varrho_\varepsilon(t)=e^{-i\frac{t}{\varepsilon}H_{\varepsilon}}
\varrho_{\varepsilon}e^{i\frac{t}{\varepsilon}H_{\varepsilon}}$\,,
$\mu_{t}=\Phi(t,0)_{*}\mu_{0}$ and $h(z,\bar z)$ is the classical
energy given in
\eqref{eq.classical}, and both sides of the identity do not depend on time.
\end{prop}
\proof
With the energy conservation, it suffices to prove \eqref{eq.energy} for $t=0$\,.  Let
$\chi\in C_0^\infty(\rz)$ such that $0\leq \chi\leq 1$\,, $\chi(s)=1$ if $|s|\leq 1$ and
$\chi(s)=0$ if $|s|\geq 2$\,. For $m\in\nz^*$\,, set
 $\chi_m(x)=\chi(\frac{x}{m})$\,. Let $B_1(z)$ and $B_2(z)$ be respectively the polynomial $\langle z, -\Delta z\rangle$ and $B_2(z)=V(z)=\frac{1}{2}\langle z^{\otimes 2}, V(x-y) z^{\otimes 2}\rangle$ well defined for $z\in \Z_1$\,.
Remember that although the kernels of  $B_1$ and $B_2$ are unbounded operators their Wick
quantization still have a meaning as densely
 defined operators on $\H$ (see Appendix \ref{se.WWaW}). \\
Write for $j=1,2$
\begin{eqnarray}
\label{eq.estgama1}
|\Tr[\varrho_\varepsilon \, B_j^{Wick}]- \int_{\Z_1} B_j(z) \,d\mu_0(z)|&\leq&
|\Tr[\varrho_\varepsilon \,(B_j^{Wick}-B_{j,m}^{Wick})]|\\ \label{eq.estgama2}&&+
|\Tr[\varrho_\varepsilon \,B_{j,m}^{Wick}]-\int_{\Z_0} B_{j,m}(z) \,d\mu_0(z)|\\ \label{eq.estgama3}
&&+ |\int_{\Z_1} B_{j,m}(z) \,d\mu_0(z)- \int_{\Z_1} B_{j}(z)\,d\mu_0(z)|\,,
\end{eqnarray}
where $B_{1,m}(z)=\langle z, -\Delta  [\chi_m(-\Delta)] z\rangle$ and $B_{2,m}(z)=
\frac{1}{2}\langle z^{\otimes 2}, V(x-y)  [\chi_m(-\Delta_x)] z^{\otimes 2}\rangle$\,.
Observe that Lemma \ref{le.boundedest} leads to
\begin{multline*}
\| (\d\Gamma(-\Delta)+\N+1)^{-1} \,\d\Gamma(\Delta 
  [(1-\chi_{m})(-\Delta)]) \, (\d\Gamma(-\Delta)+\N+1)^{-1}\|\\ \leq
\|\frac{-\Delta}{(1-\Delta)^2} [(1-\chi_{m})(-\Delta)]\|\stackrel{m\to\infty}{\to} 0\,
\end{multline*}
and
\begin{multline*}
\| (\d\Gamma(-\Delta)+\N^2+1)^{-1/2} \, (B_2^{Wick}-B_{2,m}^{Wick}) \,
(\d\Gamma(-\Delta)+\N^2+1)^{-1/2}\| \\
\leq C_{V} \|(1-\Delta)^{-1/2}
(1-\chi_{m})(-\Delta)\|
\stackrel{m\to\infty}{\to} 0\,.
\end{multline*}
Therefore the r.h.s \eqref{eq.estgama1} tends to $0$ when $m\to\infty$ thanks to the regularity of
$\varrho_\varepsilon$ and by noticing that
$(\d\Gamma(-\Delta)+\N^2+1)^{1/2}(\d\Gamma(-\Delta)+\N+1)^{-1}$ is bounded.
Now, since $B_{j,m}, j=1,2$ belong to $\P_{alg}(\Z_0)$ then by the statement \eqref{cond-PI}, proved in
\cite[Proposition 2.12]{AmNi3}, the r.h.s \eqref{eq.estgama2}
converges to 0 when $\varepsilon\to 0$\,.
Further, by the dominated convergence theorem and with the help of Lemma \ref{le.estimA},
the r.h.s \eqref{eq.estgama3} vanishes  as $m\to\infty$\,.
Hence an $\delta/3$-argument gives
$$
\lim_{\varepsilon\to 0} \Tr[\varrho_\varepsilon \, B_j^{Wick}]= \int_{\Z_1} B_j(z) \,d\mu_0(z)\, \quad
\mbox{ for } \quad j=1,2\,.
$$
Thus \eqref{eq.energy} is proved.
\fin

\subsection{Examples}
\label{se.examp}
We give here two examples, other can be found in our previous articles
\cite{AmNi1,AmNi2,AmNi3}. The first one recalls that the
transport of the Wigner measure takes into account some
correlations. The second one is about the mean field dynamics of
states, which do not satisfy \eqref{eq.hypconv} and makes a connection
with Bose-Einstein condensation.

\subsection{Deformed tori}
\label{se.deftorii}
For two elements $\psi_{1},\psi_{2}\in \Z_{1}\subset \Z_{0}$ such that
$\|\psi_{1}\|=\|\psi_{2}\|=1$ and $\langle \psi_{1}\,,\,
\psi_{2}\rangle=0$\,, the space $\Z_{0}$ can be decomposed into
$$
\Z_{0}=\cz \psi_{1}\stackrel{\perp}{\oplus}\cz \psi_{1}\stackrel{\perp}{\oplus}\psi^{\perp}\,.
$$
This decomposition is second-quantized into the Hilbert tensor product
$$
\mathcal{H}=\Gamma_{s}(\Z_{0})=\Gamma_{s}(\cz \psi_{1})\otimes
\Gamma_{s}(\cz \psi_{2})\otimes \Gamma_{s}(\psi^{\perp})\,,
$$
which allows an analysis by separating the variables. The number
observable is now
$$
\N=(\N_{1}\otimes \11\otimes \11) \oplus
(\11\otimes \N_{2}\otimes \11)
\oplus (\11\otimes \11\otimes \N')\,,
$$
simply written as $\N=\N_{1}+\N_{2}+\N'$ and
where $\N_{1}$\,, $\N_{2}$ and $\N'$ are respectively the number
operators on $\Gamma_{s}(\cz \psi_{1})$\,, $\Gamma_{s}(\cz \psi_{2})$
and $\Gamma_{s}(\psi^{\perp})$\,.
Consider in this decomposition, the state
$$
\varrho_{\varepsilon}= \varrho_{\varepsilon}^{1}\otimes
\varrho_{\varepsilon}^{2}\otimes (|\Omega'\rangle\langle \Omega'|)
$$
where $|\Omega'\rangle$ is the vacuum state of
$\Gamma_{s}(\psi^{\perp})$ and
\begin{eqnarray*}
  && \varrho_{\varepsilon}^{1}=|\psi_{1}^{\otimes n_{1}}\rangle
  \langle \psi_{1}^{\otimes n_{1}}|\quad,\quad
\varrho_{\varepsilon}^{2}=|\psi_{2}^{\otimes n_{2}}\rangle
  \langle \psi_{2}^{\otimes n_{2}}|\,,\\
\text{with}&& \lim_{\varepsilon\to 0}\varepsilon
n_{1}=\lim_{\varepsilon\to 0}\varepsilon n_{2}=\frac{1}{2}\,.
\end{eqnarray*}
In $\mathcal{H}=\Gamma_{s}(\Z_{0})$\,, this state is explicitly written
(see \cite{AmNi3}) as
\begin{eqnarray}
  \label{eq.vpower1}
&&  \varrho_{\varepsilon}=|\psi^{\vee (n_{1},n_{2})}\rangle \langle
  \psi^{\vee (n_{1},n_{2})}|\\
\label{eq.vpower2}
\text{with}&&
\psi^{\vee  (n_{1},n_{2})}
=\frac{1}{\sqrt{\varepsilon^{n_{1}+n_{2}}n_{1}!n_{2}!}}
\overbrace{a^{*}(\psi_{1})\ldots a^{*}(\psi_{1})}^{n_{1}~\text{times}}
\overbrace{a^{*}(\psi_{2})\ldots a^{*}(\psi_{2})}^{n_{2}~\text{times}}|\Omega\rangle\,.
\end{eqnarray}
The state satisfies
$$
\lim_{\varepsilon\to
    0}\Tr\left[\N^{k}\varrho_{\varepsilon}\right]=(\frac{1}{2}+\frac{1}{2})^{k}=1\,,
$$
owing to $\N=\N_{1}+\N_{2}+\N'$\,. Moreover, with \eqref{eq.vpower1}\eqref{eq.vpower2},
$\N+H^{0}_{\varepsilon}=\d\Gamma(1-\Delta)$ and the help of Wick
calculus, it also fulfills
$$
\lim_{\varepsilon\to
    0}\Tr\left[(\N+H^{0}_{\varepsilon})\varrho_{\varepsilon}\right]
=\frac{|\psi_{1}|_{\Z_{1}}^{2}+|\psi_{2}|_{\Z_{1}}^{2}}{2}
\,.
$$
Meanwhile the separation of variables allows to compute explicitly
the (it is unique) Wigner measure of
$(\varrho_{\varepsilon})_{\varepsilon\in (0,\bar\varepsilon)}$
\begin{eqnarray*}
  &&\mu_{0}=\delta_{\frac{\sqrt{2}}{2}\psi_{1}}^{S^{1}}\otimes
\delta_{\frac{\sqrt{2}}{2}\psi_{1}}^{S^{1}}\otimes \delta_{0}\quad \text{on}\quad
\Z_{1}= (\cz \psi_{1})\times(\cz\psi_{2})\times\psi^{\perp}\,,\\
\text{with}
&&
\delta_{u}^{S^{1}}=\frac{1}{2\pi}\int_{0}^{2\pi} \delta_{e^{i\theta}u}~d\theta\,.
\end{eqnarray*}
We get
$$
\int_{\Z_{1}}|z|^{2k}~d\mu_{0}(z)=\int_{\Z_{1}}\left(|z_{1}|^{2}+|z_{2}|^{2}+|z'|^{2}\right)^{k}~d\mu_{0}(z)=1
=\lim_{\varepsilon\to \infty}\Tr\left[\N^{k}\varrho_{\varepsilon}\right]\,.
$$
Hence all the assumptions of Theorem~\ref{th.main} and
Theorem~\ref{Piwig} are fulfilled.\\
This measure is carried by a torus in  $\Z_{1}$ better described by
using an other orthonormal basis of $\cz \psi_{1}\oplus \cz \psi_{2}$:
\begin{eqnarray*}
  && \psi_{0}=\frac{\sqrt{2}}{2}(\psi_{1}+\psi_{2})\quad,\quad
  \psi_{\frac{\pi}{2}}=
i\frac{\sqrt{2}}{2}(\psi_{1}-\psi_{2})\,,\\
&&
\psi_{\varphi}=\cos(\varphi)\psi_{0}+\sin(\varphi)\psi_{\frac{\pi}{2}}\,,\\
&&
\frac{\sqrt{2}}{2}(e^{i\theta}\psi_{1}+e^{i\theta'}\psi_{2})=e^{i\frac{\theta+\theta'}{2}}\psi_{\frac{\theta-\theta'}{2}}\,,\\
&&
\mu_{0}=\frac{1}{2\pi}\int_{0}^{2\pi}\delta_{\psi_{\varphi}}^{S^{1}}~d\varphi\,.
\end{eqnarray*}
Two elements $e^{i\theta}\psi_{\varphi}$ and
$e^{i\theta'}\psi_{\varphi'}$ in the support of $\mu_{0}$ are equal
when
$$
\left(\theta'=\theta\;\text{and}\; \varphi'=\varphi\right)
\quad\text{or} \quad\left(\theta'=\theta+\pi\;\text{and}\; \varphi'=\varphi+\pi\right)\,.
$$
Hence a one to one parametrization of the torus can be done by
$\varphi\in [0,2\pi)$ and $\theta\in [\varphi, \varphi+\pi)$\,.

Let $\psi_{\varphi}(t)=\Phi(t,0)\psi_{\varphi}$\,, be the solution to the Hartree equation
$$
\left\{
  \begin{array}[c]{l}
    i\partial_{t}\psi_{\varphi}(t)= -\Delta \psi_{\varphi}(t)+
    (V*|\psi_{\varphi}(t)|^{2})\psi_{\varphi}(t)
\\
\psi_{\varphi}(t=0)=
\psi_{\varphi}=e^{i\frac{\pi}{4}}\cos(\varphi)\psi_{1}+
e^{-i\frac{\pi}{4}}\sin(\varphi)\psi_{2}
  \end{array}
\right.\,,
$$
The gauge invariance of the equation says that for any $\theta\in
[0,2\pi]$\,, $e^{i\theta}\psi_{\varphi}(t)=\Phi(t,0)\left[e^{i\theta}\psi_{\varphi}\right]$\,.
By applying the result of Theorem~\ref{th.main} and
Theorem~\ref{Piwig} we get
\begin{eqnarray*}
&&
\mu_{t}=\frac{1}{2\pi}\int_{0}^{2\pi}\delta_{\psi_{\varphi}(t)}^{S^{1}}~d\varphi
=
\frac{1}{4\pi^{2}}\int_{0}^{2\pi}\int_{0}^{2\pi}\delta_{e^{i\theta}\psi_{\varphi}(t)}~d\varphi
d\theta\\
&&
\forall p\in\nz,\quad
\lim_{\varepsilon\to
  0}\gamma_{\varepsilon}^{(p)}(t)=\frac{1}{2\pi}\int_{0}^{2\pi}
|[\psi_{\varphi}(t)]^{\otimes p}\rangle\langle
[\psi_{\varphi}(t)]^{\otimes p}|~d\varphi\,.
\end{eqnarray*}
Since the Hartree flow is nonlinear, the complete hierarchy of
reduced density matrices have to be taken into account if one wants to
write evolution equation for them. More simply,
they can be computed after solving an autonomous equation
for the Wigner measure.
Due to the nonlinear term the dynamics of correlations is by far nontrivial.
This can also be thought geometrically: The initial measure is
initially supported by a torus which lies in a $2$-dimensional complex
vector space (think of the circle in the plane $\rz\psi_{0}\oplus
\rz\psi_{\frac{\pi}{2}}$); along the time evolution, the measure
$\mu_{t}$ is
still carried by a torus in $\Z_{1}$\,, which  nevertheless,
is a priori not embedded in any finite dimensional subspace\,.
\figinit{0.8pt}

\figpt 0:$0$ (0,0)
\figpt 10: $\rz\psi_{0}$ (100,0)
\figpt 11: (-100,0)
\figpt 20: (0,100)
\figpt 21: (0,-100)
\figptrot 22: $\rz\psi_{\frac{\pi}{2}}$= 20 /0, -38/
\figptrot 23:=21 /0,-38/

\figpttraC 100:$0$=0/250,0/
\figpttraC 110:$\rz\psi_{0}$=10/250,0/
\figpttraC 111:=11/250,0/
\figpttraC 120:$\psi^{\perp}$=20/250,0/
\figpttraC 121:=21/250,0/
\figpttraC 122:$\rz\psi_{\frac{\pi}{2}}$=22/250,0/
\figpttraC 123:=23/250,0/

\def\haxis{80}\def\vaxis{30}
\def\startangle{-20}\def\stopangle{90}\def\inclin{10}
\figptell 4:: 0 ; \haxis,\vaxis( 0, \inclin) 
\figptell 5:: 0 ; \haxis,\vaxis(90, \inclin) 

\figptrot 30:= 21/0,20/
\figptsinterlinell 31 [0, \haxis,\vaxis,\inclin; 0,30]
\figpttraC 35:$\psi_{\varphi}$=31/0,0/
\figpttraC 32:$\times$=31/-50,0/
\figpttraC 34:=32/-22,0/

\figpt 150: (270,-30)
\figpt 151: (275, -10)
\figpt 152:$\psi_{\varphi}{(t)}$ (310, 20)
\figpt 153: (280, 30)
\figpt 154: (250, 60)
\figpt 155: (230, 50)
\figpt 156: (230, 30)

\figpthom 160:= 150/100, -1/
\figpthom 161:= 151/100, -1/
\figpthom 162:= 152/100, -1/
\figpthom 163:= 153/100, -1/
\figpthom 164:= 154/100, -1/
\figpthom 165:= 155/100, -1/
\figpthom 166:= 156/100, -1/

\figpttraC 157:$\times$=152/30,0/
\figpttraC 158:=157/22,0/
\psbeginfig{}
\psarrow[23,22]
\psarrow[11,10]
\psarrow[32,31]
\pscirc 34(10)

\psarrow[123,122]
\psarrow[111,110]
\psarrow[121,120]
\psset(width=1.5)
\psarcellPA 0,4,5(-180,180)
\pscurve[150,150,151,152,153,154,155,156,156]
\psset(dash=2)
\pscurve[160,160,161,162,163,164,165,166,166]
\psset(width=\defaultwidth)
\psset(dash=\defaultdash)
\psarrow[157,152]
\pscirc 158 (10)

\psendfig
\figvisu{\figBoxA}{\vbox{
    \begin{center}
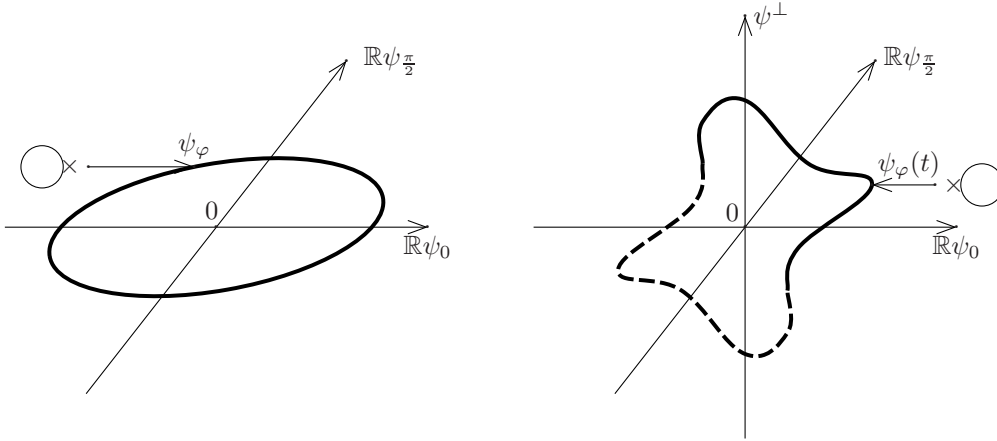
\bf Fig.1: Evolution of the measure initially carried by a
torus in $\cz\psi_{0}\oplus \cz \psi_{\frac{\pi}{2}}$\,.\\
The complex gauge parameter $e^{i\theta}$ is represented by the small circle.
\end{center}
}}{
\figsetmark{.}
\figwriten 0:(4)
\figwrites 10,110:(4)
\figwritenw 100:(4)
\figwritee 22,122:(4)
\figwritew32:(4)
\figwritee 157:(4)
\figwritene 152:(4)
\figwriten 35:(4)
\figwritee 120:(4)
}
\centerline{\box\figBoxA}
\noindent In Figure~1, the deformed torus for time $t\neq 0$\,, has to be imagined in
the infinite-dimensional phase-space $\Z_{1}\subset
\Z_{0}$\,. Contrary to the picture, there might be no intersection
with the real plane $\rz \psi_{0}\oplus \rz \psi_{\frac{\pi}{2}}$\,.\\
This discussion can also be extended to higher dimensional tori after
taking a finite (or countable) orthornormal family $(\psi_{n})_{1\leq
  n\leq N}$ for
building the initial states $\varrho_{\varepsilon}$ with a measure
$\prod_{j=1}^{N}\delta_{\lambda_{j}\psi_{j}}^{S^{1}}$ (see
\cite{AmNi3})\,.\\

\subsection{Propagation without the convergence of moments}
\label{se.bec}

In \cite{AmNi1} we considered the thermodynamic limit of a free Bose
gas on a torus with the one particle energy given by $-\Delta$\,. We
showed that in the regime which may exhibit
a Bose condensation, the condition
\eqref{eq.hypconv} fails and illustrates what we called a dimensional
defect of compactness, in opposition to the phase space or microlocal
defect of compactness (see \cite{PGe,Tar}). Others examples were given. In \cite{AmNi3} the
propagation result for bounded interactions but without any
compactness condition, cannot be applied for such initial states. With
Theorem~\ref{th.main} the propagation holds for this kind of initial
states. Since our analysis is valid on $\rz^{d}$ the analysis for the
torus does not apply directly and we adapt the presentation of the
Bose-Einstein condensation.\\
Moreover the dimensional defect of compactness which plays with all
the directions of the phase-space $\Z_{0}=L^{2}(\rz^{d})$\,, can be
geometrically thought in the one particle phase-space
$T^{*}\rz^{d}$\,.
The condition~\eqref{eq.hypmom}, which leads to estimates of
$\int_{\Z_{0}}|z|_{H^{1}}^{2\delta}~d\mu$\,, suggests that the
dimensional defect of compactness is due to mass going to $\infty$ in
the position variable rather than in the momentum variable, in
$T^{*}\rz^{d}$\,.
The mean field limit that we consider here, can be tested by using the
harmonic
oscillator Hamiltonian
$A=-\partial_{x}^{2}+\frac{x^{2}}{4}-\frac{d}{4}$\,. The motivated
reader will then see that the dimensional defect of compactness
$\varrho_{\varepsilon}$ is incompatible with the
condition~\eqref{eq.hypmom}.\\

We work in dimension $d\geq 2$\,. Let $e_{0}$ be an $L^{2}$-normalized
 $\mathcal{C}^{\infty}$ function supported in the hypercube
$(-\frac{1}{2},\frac{1}{2})^{d}$ and set
$$
\forall k\in\nz^{d}\,,\quad e_{k}(x)=e_{0}(x-k)\,.
$$
The family $(e_{k})_{k\in\nz}$ is orthonormal in
$\Z_{0}=L^{2}(\rz^{d})$\,. The spanned Hilbert
subspace and the corresponding orthogonal projection are respectively
denoted by $\Z_{e}$ and  $\Pi_{e}$\,, $\mathrm{Ran}\,\Pi_{e}=\Z_{e}$\,.
Note that
$$
\Pi_{e}(-\Delta)\Pi_{e}= \left(\int_{\rz}|\nabla e_{0}|^{2}\right)\Pi_{e}\,.
$$
Consider now the self-adjoint operator defined on $\Z_{0}=L^{2}(\rz^{d})$ by
$$
A=\sum_{k\in\nz^{d}}|k| |e_{k}\rangle \langle
e_{k}|\,,\quad |k|=\sum_{j=1}^{d}k_{j}\,,
$$
which restricted to $\Z_{e}$ is unitarily equivalent to the
harmonic
oscillator Hamiltonian
$A=-\partial_{x}^{2}+\frac{x^{2}}{4}-\frac{d}{4}$ on $\rz^{d}$\,.
We use the tensor decomposition
\begin{eqnarray*}
  && \Z_{0}=\Z_{e}\stackrel{\perp}{\oplus}\Z_{e}^{\perp}\quad,\quad
  \mathcal{H}=\Gamma_{s}(\Z_{0})= \Gamma_{s}(\Z_{e})\otimes
  \Gamma_{s}(\Z_{e}^{\perp})\,,\\
&& \N= \N_{e}\otimes \11+
\11\otimes\N_{e}^{\perp}=\N_{e}+\N_{e}^{\perp}\quad,\quad
|\Omega\rangle=|\Omega_{e}\rangle\otimes |\Omega_{e}^{\perp}\rangle\,,\\
&& \forall B_{e}\in \mathcal{L}(\Z_{e})\,, \|B_{e}\|\leq 1, \quad
\Gamma(B_{e})\otimes (|\Omega_{e}^{\perp}\rangle\langle
\Omega_{e}^{\perp}|)= \Gamma(\Pi_{e}B_{e}\Pi_{e})\,,\\
&& \forall B\in \mathcal{L}(\Z_{0})\,,\, \|B\|\leq 1,\quad
\Gamma(\Pi_{e}B\Pi_{e})= \Gamma(\Pi_{e}\circ B|_{\Z_{e}})\otimes (|\Omega_{e}^{\perp}\rangle\langle
\Omega_{e}^{\perp}|)\,.
\end{eqnarray*}
In particular the last relation with $B=e^{i\varepsilon t\Delta}$
differentiated at time $t=0$ gives
$$
\Gamma(\Pi_{e})\d\Gamma(-\Delta)\Gamma(\Pi_{e})=\left(\int_{\rz^{d}}|\nabla
  e_{0}|^{2}\right)\N_{e}\,.
$$
Consider on $\mathcal{H}$\,,
the $\varepsilon$-dependent gauge invariant (tensorized) quasi-free state
\begin{eqnarray}
  \label{eq.rhoBECe}
\varrho_{\varepsilon}&=&
\frac{1}{\Tr\left[\Gamma(\Pi_{e}
e^{-\beta_{\varepsilon}(A-\mu_{\varepsilon})}\Pi_{e})\right]}
\Gamma(\Pi_{e}e^{-\beta_{\varepsilon}(A-\mu_{\varepsilon})}\Pi_{e})
\\
\nonumber
&=&\frac{1}{\Tr\left[\Gamma(\Pi_{e} Z_{\varepsilon}e^{-\beta_{\varepsilon} A}\Pi_{e})\right]}
\Gamma(\Pi_{e}Z_{\varepsilon}e^{-\beta_{\varepsilon} A}\Pi_{e})\,,\\
\nonumber
&=&
\frac{1}{\Tr\left[\Gamma(Z_{\varepsilon}e^{-\beta_{\varepsilon} A|_{\Z_{e}\cap \D(A)}})\right]}
\Gamma(\Pi_{e})\Gamma(Z_{\varepsilon}e^{-\beta_{\varepsilon} A|_{\Z_{e}\cap \D(A)}})\Gamma(\Pi_{e})\,.
\end{eqnarray}
The chemical potential $\mu_{\varepsilon}$ is negative of order $\varepsilon^{1-1/d}$ and the
temperature is large according to
$$
Z_{\varepsilon}=e^{\beta_{\varepsilon}
  \mu_{\varepsilon}}=1-\frac{\varepsilon}{\nu_{C}}
\quad,\quad \beta_{\varepsilon}=\varepsilon^{1/d}\,.
$$
With the $\varepsilon$-dependent definition of $a(f)$\,, $a^{*}(f)$\,,
$\left[a(g),a^{*}(f)\right]=\varepsilon \langle g,f\rangle$\,, and
$W(f)$\,, this
quasi-free state is characterized by the two-point function
\begin{eqnarray}
\label{eq.twopoint}
&&  \Tr\left[\varrho_{\varepsilon}a^{*}(f)a(g))\right]
=\varepsilon
\left\langle \Pi_{e}g\,,
  Z_{\varepsilon}e^{-\beta_{\varepsilon} A_{\varepsilon}}(1-Z_{\varepsilon}e^{-\beta_{\varepsilon} A_{\varepsilon}})^{-1} \Pi_{e}f
\right\rangle\,.
\\
\label{eq.charfunc}
\text{or}&&
\Tr\left[\varrho_{\varepsilon}W(f)\right]
=\exp\left[-\varepsilon\langle \Pi_{e}f\,,\,(1+Z_{\varepsilon}e^{-\beta_{\varepsilon}
    A_{\varepsilon}})(1-Z_{\varepsilon}e^{-\beta_{\varepsilon} A_{\varepsilon}})^{-1}\Pi_{e}f\rangle/4\right]\,.
\end{eqnarray}
In particular the total number (multiplied by $\varepsilon$) is given by
\begin{eqnarray}
\nonumber
\Tr\left[\varrho_{\varepsilon}\N\right]&=&\Tr\left[\varrho_{\varepsilon}\N_{e}\right]
\sum_{k\in\nz^{d}}\Tr\left[\varrho_{\varepsilon}a^{*}(e_{k})a(e_{k})\right]
\\
  \label{eq.densbos}
&=&\nu_{C}+ \nu+
r(\varepsilon)\quad\text{with}~\lim_{\varepsilon\to
  0}r(\varepsilon)=0\,,\\
\text{and}
&&\nu=\int_{\rz^{d}}\frac{e^{-|u|}}{1-e^{-|u|}}~du=
|S^{d-1}|\int_{0}^{+\infty}\frac{e^{- t}}{1-e^{-t}}t^{d-1}~dt\,,
\quad (d\geq 2)\,.
\end{eqnarray}
We deduce
\begin{eqnarray*}
 &&\lim_{\varepsilon\to  0}\Tr\left[\varrho_{\varepsilon}\N\right]=\nu_{C}+\nu\quad,\quad
\\
&&\lim_{\varepsilon\to
  0}\Tr\left[\varrho_{\varepsilon}H_{\varepsilon}^{0}\right]=
\lim_{\varepsilon\to
  0}\Tr\left[\varrho_{\varepsilon}\Gamma(\Pi_{e})\d\Gamma(-\Delta)\Gamma(\Pi_{e})\right]=
\left(\int_{\rz^{d}}|\nabla e_{0}|^{2}\right)(\nu_{C}+\nu)\,,
\end{eqnarray*}
and the condition \eqref{eq.hypmom} of Theorem~\ref{th.main} is
satisfied.\\
Actually $\nu_{C}>0$ corresponds, in the analysis of the free Bose gas
(see \cite{AmNi1}), to the density associated with the condensate
phase. In the scaling that we consider, it is the other part  which
produces the dimensional defect of compactness.
Let us compute  the
Wigner measure, by considering the limit of
$\Tr\left[\varrho_{\varepsilon}W(\sqrt{2}\pi f)\right]$ as
$\varepsilon\to 0$\,. With
$$
f=\sum_{k\in\nz^{d}}f_{k}e_{k}+f^{\perp}\,,\quad
|f|^{2}=\sum_{k\in\nz^{d}}|f_{k}|^{2}+|f^{\perp}|^{2}=|\Pi_{e}f|^{2}+|f^{\perp}|^{2}\,,
$$ the expression
\eqref{eq.charfunc} gives
\begin{equation}
  \label{eq.explicit}
\Tr\left[\varrho^{\varepsilon}W(\sqrt{2}\pi f)\right]
=
e^{-\varepsilon \pi^{2}\left|\Pi_{e} f\right|^{2}/2}
\times
\exp\left[-\varepsilon \pi^{2}
\sum_{k\in\nz^{d}}\left|f_{k}\right|^{2}\frac{Z_{\varepsilon}
e^{-\varepsilon^{1/d}|k|}}{(1-Z_{\varepsilon}e^{-\varepsilon^{1/d}|k|})}\right]
\stackrel{\varepsilon\to 0}{\to} e^{-\pi^{2}\nu_{C}|f_{0}|^{2}}\,.
\end{equation}
The family $(\varrho_{\varepsilon})_{\varepsilon\in
  (0,\bar\varepsilon)}$ admits the unique Wigner measure
$$
\mu_{0}=
\frac{e^{-\frac{|z_{0}|^{2}}{\nu_{C}}}}{
\pi\nu_{C}}\otimes\delta_{0}
\quad \text{on}~\Z_{0}= (\cz e_{0})\times e_{0}^{\perp}\,,
$$
which is carried by $\cz e_{0}\subset \Z_{1}$ and which can also be written
$$
\mu_{0}=\int_{\cz e_{0}} \frac{e^{-\frac{|z|^{2}}{\nu_{C}}}}{
\pi\nu_{C}} \delta_{ze_{0}}~L_{\cz e_{0}}(dz)
=\int_{0}^{+\infty}\frac{e^{-\frac{u}{\nu_{C}}}}{\nu_{C}}
\delta_{\sqrt{u}e_{0}}^{S^{1}}~du\,.
$$
In particular, we get
\begin{eqnarray*}
  &&\int_{\Z_{1}}|z|_{\Z_{0}}^{2}d\mu_{0}(z)=\nu_{C} <
  \nu_{C}+\nu=\lim_{\varepsilon\to 0}\Tr\left[\varrho_{\varepsilon}\N\right]\\
\text{and}
&&\int_{\Z_{1}}|z|_{\Z_{1}}^{2}d\mu_{0}(z)=\nu_{C}\left(\int_{\rz^{d}}|\nabla
  e_{0}|^{2}\right) <
  (\nu_{C}+\nu) \left(\int_{\rz^{d}}|\nabla
  e_{0}|^{2}\right)=\lim_{\varepsilon\to 0}\Tr\left[\varrho_{\varepsilon}\d\Gamma(1-\Delta)\right]\,,
\end{eqnarray*}
and the condition~\eqref{eq.hypconv} does not hold.
Even at time $t=0$\,, no formula  is available  for the reduced density
matrices in terms of the Wigner measure.
Nevertheless the time-dependent Wigner measure of
$\varrho_{\varepsilon}(t)=e^{-i\frac{t}{\varepsilon}H_{\varepsilon}}\varrho_{\varepsilon}e^{i\frac{t}{\varepsilon}H_{\varepsilon}}$
is given
by Theorem~\ref{th.main}, since the condition
\eqref{eq.hypmom} is verified.
Consider the solutions to the Hartree initial value problems
$$
\left\{
  \begin{array}[c]{l}
    i\partial_{t}\psi_{u}=-\Delta
    \psi_{u}+(V*|\psi_{u}|^{2})\psi_{u}\\
  \psi_{u}(t=0)=\sqrt{u}e_{0}\,,\quad u\in (0,+\infty)\,.
  \end{array}
\right.
$$
Then the Wigner measure of
$\varrho_{\varepsilon}(t)=e^{-i\frac{t}{\varepsilon}H_{\varepsilon}}\varrho_{\varepsilon}
e^{i\frac{t}{\varepsilon}H_{\varepsilon}}$ is given by
$$
\mu_{t}=\int_{0}^{+\infty}\frac{e^{-\frac{u}{\nu_{C}}}}{\nu_{C}}
\delta_{\psi_{u}(t)}^{S^{1}}~du\,.
$$
Again like in the example of the previous section, the measure
$\mu_{t}$ is carried by surface containing $0$ and topologically equivalent to
$\cz$\,, but this $2$-dimensional surface does not remain a priori in any finite
dimensional subspace of $\Z_{1}$ for $t\neq 0$\,.
\subsection{About the classical mean field problem}
\label{se.clas}
The classical analogue of our analysis is the derivation of the Vlasov
equation
$$
\left\{
  \begin{array}[c]{l}
\partial_{t}f+v.\partial_{x}f-\frac{1}{m}(\partial_{x}V_{f}(x,t)).\partial_{v}f=0\\
f(t,x,v)=f_{0}(x,v)\\
V_{f}(x,t)=V*\varrho_{f}(x,t)\quad,\quad \varrho_{f}(x,t)=\int_{\rz^{d}}f(x,v,t)~dv
\end{array}
\right.
$$
where $f(x,v,t)$  represents the particle density in the
$1$-particle phase space $\rz^{2d}_{x,v}$\,, from the classical Hamilton many body
system
$$
\left\{
  \begin{array}[c]{l}
    \dot{x}_{i}=v_{i}\,,\\
   \dot{v}_{i}=-\frac{1}{mN}(\sum_{j=1}^{N}\partial_{x_{i}}V(x_{i}-x_{j}))
  \end{array}
\right.\,,\quad  i=1,\ldots, N\,,
$$
in the limit $N\to \infty$\,. This problems is still open for singular
potential and C.~Villani, in a recent survey article about the Landau
damping \cite{Vil}
 quotes the work of Hauray-Jabin \cite{HaJa} as the most advanced one
 in this direction. It works for a potential such that $|\nabla
 V|=\mathcal{O}(|x|^{-s})$\,, $s\in (0,1)$\,, and does not include the
 Coulomb interaction.\\
Indirectly our result, justifies the mean field model up to Coulomb
interaction in dimension $d=3$\,. In \cite{LiPa} and more recently
\cite{AFFGP}, the Vlasov equation is proved to be the semiclassical
limit of the semiclassical Hartree equation. This means that there are
two ``semiclassical'' limits, one in the phase-space $L^{2}(\rz^{d};\cz)$
with the small parameter $1/N$\,, another one on the phase-space
$T^{*}\rz^{d}\sim\rz^{2d}$ for the one particle nonlinear
problem. This double asymptotic regime
 is well presented in \cite{FGS,FKS,GMP}.\\

A possible strategy, for deriving directly the classical mean field
limit from the classical many body problem, consists in adapting our
approach by, as usual, replacing traces by integrals.
For information, we refer the reader to the presentation \cite{Der} by
J. Derezinski of the classical analogue of second quantization.
Of course
classical mechanics, although living in the commutative world, is
often more singular than quantum mechanics, from the analysis point of
view. With the Coulomb interaction, the Kustaanheimo-Stiefel
desingularization of the hamiltonian flow may be useful (see
a.e. \cite{CJK, HeSi, Ker, Kna, KuSt}).

\bigskip

\appendix
\begin{center}
{\bf\Large Appendix}
\end{center}

\section{Commuting self-adjoint operators on a graded Hilbert space}
\label{se.comsa}
We briefly study the general structure of self-adjoint operators on a graded Hilbert space.
Properties collected in this section are useful for the analysis of the
quantum Hamiltonian \eqref{eq.hamq}. In this appendix, the small
parameter is not required and we work with $\varepsilon =1$\,.\\
 Remember that a graded Hilbert space
$\H$ is a direct sum of Hilbert spaces $\H_n, n\in\nz,$ of the form
$$
\H=\bigoplus_{n=0}^\infty \H_n\,.
$$
Let $(A_n)_{n\in\nz}$ be a sequence of self-adjoint operators where each $A_n$ acts on $\H_n$\,.
We define the operator
\begin{eqnarray}
\label{co1}
\D(A)=\{\Psi\in\H : \sum_{n=0}^\infty \|A_n \Psi^{(n)}\|_{\H_n}^2<\infty\},\quad
A\,\Psi=\sum_{n=0}^\infty A_n \,\Psi^{(n)}\,, \mbox{ for all } \Psi\in \D(A)\,.
\end{eqnarray}
Taking in particular  $A_n=n \11_{\H_n}$  for $n\in\nz$\,, we obtain the {\it number } operator
\begin{eqnarray}
\label{co2}
N=\sum_{n=0}^\infty n \11_{\H_n}\,.
\end{eqnarray}
We say that two self-adjoint operators $B$ and $C$ on a Hilbert space {\it strongly commute} if their spectral projections mutually commute. This is equivalent to the commutation of their resolvents
for some $z\in \cz\setminus \rz$ and also to the commutation of their associated unitary groups. More precisely, $B$ and $C$ strongly commute if and only if for all $t,s\in\rz$
$$
e^{i t C} e^{is B}=e^{is B} e^{i t C}\,.
$$
 \begin{prop}
 \label{graded}
Let $A$ and $N$ be the operators given by (\ref{co1})-(\ref{co2}). The following assertions hold:\\
(i) $A$ and $N$ are self-adjoint.\\
(ii) For any bounded Borel function on $\rz$\,,
$$
f(A)=\sum_{n=0}^\infty f(A_n)\,.
$$
(iii) The operators $A$ and $N$ strongly commute.\\
(iv) If $\D_n$ is a core for $A_n$ for each $n\in\nz$ then $
\oplus_{n\in\nz}^{alg} \D_n$ is a core for $A$\,.\\
(v) For any real polynomial $p$ the operator $
A+p(N)_{|\D(A)\cap\D(p(N))}$ is essentially self-adjoint and
$$
\overline{A+p(N)_{|\D(A)\cap\D(p(N))}}=\sum_{n=0}^\infty A_n+p(n)\11_{\H_n}\,.
$$
\end{prop}
\proof
(i) Clearly, $A$ is a densely defined operator. It is also symmetric, since for any $\Psi,\Phi\in\D(A)$
$$
\langle \Phi, A \Psi\rangle_\H=\sum_{n=0}^\infty \langle \Phi^{(n)}, A_n \Psi^{(n)}\rangle_{\H_n}=
\sum_{n=0}^\infty \langle A_n \Phi^{(n)}, \Psi^{(n)}\rangle_{\H_n}=\langle A\Phi,  \Psi\rangle_\H\,.
$$
For any $\Psi\in \D(A)$ and $\Phi\in \D(A^*)$\,,
$$
\langle A^* \Phi, \Psi\rangle=\sum_{n=0}^\infty \langle \Phi^{(n)},A_n \Psi^{(n)}\rangle_{\H_n}
$$
Hence the inequality holds
$$
\left|\sum_{n=0}^\infty \langle \Phi^{(n)}, A_n
  \Psi^{(n)}\rangle_{\H_n}\right|\leq \|A^*\Phi\|_\H \,
\|\Psi\|_{\H}\,.
$$
By taking any $\Psi^{(n)}\in \D(A_{n})$\,, this means  $\Phi^{(n)}\in
\D(A_n^*)=\D(A_n^{*})=\D(A_{n})$\,. The extension to any $\Psi\in
\mathcal{H}$ gives $\Phi\in \D(A)$\,. This proves that $A$ and $N$ are self-adjoint.\\
(ii) For each $n\in\nz$\,, the map $t\mapsto e^{i t A} e^{-i t A_n} \Psi^{(n)}$ is of class $C^1$ for any
$\Psi^{(n)}\in\D(A_n)$ by  Stone's theorem with the derivative
$$
\frac{d}{dt} \,e^{i t A} e^{-i t A_n} \Psi^{(n)}=i e^{i t A} (A-A_n) e^{-i t A_n} \Psi^{(n)}=0\,.
$$
Hence, for any $\Psi\in \oplus_{n\in \nz}^{alg} \D(A_n)$ (and then for any $\Psi\in\H$ since
$\oplus_{n\in\nz}^{alg} \D(A_n)$ is dense in $\H$)
we see that for all $t\in\rz$
\begin{eqnarray}
\label{co3}
e^{itA} \,\Psi=\sum_{n=0}^\infty e^{itA_n} \Psi^{(n)}\,.
\end{eqnarray}
By functional calculus we extend the identity \eqref{co3} to any bounded Borel function $f$ on $\rz$\,.\\
(iii) By using (ii), we get for all $s,t\in\rz$ and $\Psi\in\H$
$$
e^{i t N} e^{i s A} \Psi= e^{i t N} \sum_{n=0}^\infty e^{i s A_n } \Psi^{(n)}=
\sum_{n=0}^\infty  e^{i t n} e^{i s A_n } \Psi^{(n)}=e^{i s A} e^{i t N}  \Psi\,.
$$
(iv) The algebraic direct sum
$\D_{fin}=\oplus^{alg}_{n\in\nz}\D(A_n)\subset \D(A)$ is dense in $\H$ and
invariant with respect to the group $(e^{i t
  A})_{t\in\rz}$\,. Therefore, $\D_{fin}$ is a core for $A$\,.  On the other hand,
the subspace $\D_{fin}^0=\oplus^{alg}_{n\in\nz}\D_n$ satisfies
$$
A_{|\D_{fin}}\subset \overline{A_{|\D^0_{fin}}}\subset A_{|\D(A)}\,.
$$
Hence $\D^0_{fin}$ is also a core for $A$ since $\overline{A_{|\D_{fin}}}=A_{|\D(A)}$\,.\\
(v) The operator $B=\sum_{n=0}^\infty A_n+p(n)\11_{\H_n}$ (with its natural domain) is
self-adjoint by assertion (i). It is clear that
$$\D_{fin}=\oplus^{alg}_{n\in\nz}\D(A_n)\subset \D(A)\cap\D(p(N))\subset \D(B)\,,$$
and furthermore
$$
B_{|\D_{fin}} =A+p(N)_{|\D_{fin}}\subset A+p(N)_{|\D(A)\cap\D(p(N))}\subset B_{|\D(B)}\,.
$$
Therefore, the operator $ A+p(N)_{|\D(A)\cap\D(p(N))}$ is essentially self-adjoint since $\overline{B_{|\D_{fin}}}=B_{|\D(B)}$\,.
\fin

\section{Second quantization}
\label{se.WWaW}
 For the reader's convenience; the general framework of second
  quantization and some related notations are recalled.\\
The phase-space, a complex separable Hilbert space, is denoted by $\mathcal{Z}$
with the scalar product $\langle .,.\rangle$\,. The symmetric Fock space over $\mathcal{Z}$ is defined as the following
direct Hilbert sum
\begin{eqnarray*}
\Gamma_s(\Z)=\bigoplus_{n=0}^\infty\bigvee^n
\mathcal{Z}\,,
\end{eqnarray*}
where $\bigvee^n \mathcal{Z}$ is the $n$-fold symmetric tensor product.
The orthogonal projection from $\Z^{\otimes n}$ onto the closed subspace $\bigvee^n\Z$ is given by
$$
\S_{n}(\xi_{1}\otimes\xi_{2}\cdots\otimes \xi_{n})
=\frac{1}{n!}\sum_{\sigma\in \mathfrak{S}_n}\xi_{\sigma(1)}\otimes\xi_{\sigma(2)}\cdots\otimes \xi_{\sigma(n)}\,.
$$
Algebraic direct sums or tensor products are denoted with a $~^{alg}~$
superscript. Hence
$$
\H_{fin}=\mathop{\bigoplus}_{n\in   \mathbb{N}}^{alg}\bigvee^n \mathcal{Z}
$$
denotes the subspace of vectors with a finite number of particles.
The creation and annihilation operators  $a^*(z)$ and $a(z)$\,, parameterized by $\varepsilon>0$\,,
are then defined by :
\begin{eqnarray*}
a(z) \varphi^{\otimes n}&=&\sqrt{\varepsilon n} \; \;\langle z,\varphi\rangle \varphi^{\otimes (n-1)}\\
a^*(z)\varphi^{\otimes n}&=&\sqrt{\varepsilon (n+1)}  \;\;\S_{n+1}
(\;z\otimes \varphi^{\otimes n})\,,\;\;\, \forall \varphi, z\in\Z.
\end{eqnarray*}
They extend to closed operators and they are adjoint of one another.
They also satisfy the $\varepsilon$-canonical commutation relations (CCR):
\begin{eqnarray}
\label{ccr}
[a(z_1),a^*(z_2)]=\varepsilon\langle z_1,z_2\rangle \11, \;\;\;[a^*(z_1),a^*(z_2)]=0=[a(z_1),a(z_2)]\,.
\end{eqnarray}
The Weyl operators are given for $z\in\Z$  by
$$
W(z)=e^{\frac{i}{\sqrt{2}} [a^*(z)+a(z)]}\,,
$$
and they satisfy Weyl commutation relations in the Fock space
\begin{eqnarray}
\label{eq.Weylcomm}
W(z_1) W(z_2)=e^{-\frac{i\varepsilon}{2} {\rm Im~}\langle z_1, z_2\rangle} \;W(z_1+z_2), \,\,z_1,z_2\in\Z\,.
\end{eqnarray}
The   number operator is also parametrized by $\varepsilon>0$\,,
\begin{eqnarray*}
\mathbf{N}_{|\bigvee^n \Z}=\varepsilon{n} \11_{|\bigvee^n \Z}\,.
\end{eqnarray*}
 For any self-adjoint operator $A:\Z\supset \D(A)\to \Z,$
the operator $\d\Gamma(A)$ is the self-adjoint operator given by
\begin{eqnarray*}
 \d\Gamma(A)_{|\bigvee^{n,\textrm{alg}}\D(A)}=\varepsilon\left[
\sum_{k=1}^{n}\11\otimes\cdots\otimes\underbrace{A}_{k}\otimes
\cdots\otimes \11\right]\,.
\end{eqnarray*}

\subsection{Weyl, Anti-Wick quantized operators}
\label{se.weylAwick}
Let  $\p$ denote the set of all finite rank orthogonal projections
on $\Z$ and  for a given $\wp\in\p$ let $L_{\wp}(dz)$ denote the
Lebesgue measure on the finite dimensional subspace $\wp\Z$\,, with
volume $1$ for an orthonormal hypercube in $\wp\Z$\,. A
function $f:\Z\to\cz$ is said  {\it cylindrical} if there
exists $\wp\in\p$ and a function $g$ on $\wp\Z$ such that $ f(z)=g(\wp z),$
for all $z\in\Z$\,. In this case we say that $f$ is based on the
subspace $\wp\Z$\,. We set $\S_{cyl}(\Z)$ to be the   cylindrical
Schwartz space:
$$
(f\in \S_{cyl}(\Z))\Leftrightarrow
\left(\exists \wp\in \p,\exists g\in \S(\wp\Z), \quad f(z)= g(\wp z)\right)\,.
$$
The Fourier transform of a function $f\in\S_{cyl}(\Z)$ based on the
subspace $\wp\Z$ is defined as
\begin{eqnarray*}
\F[f](\xi)=\int_{\wp\Z} f(z) \;\;e^{-2\pi i
\,{\rm Re~}\langle z,\xi\rangle}~L_{\wp}(d z)
\end{eqnarray*}
and its inverse Fourier transform as
\begin{eqnarray*}
f(z)=\int_{\wp\Z} \F[f](\xi) \;\;e^{2\pi i
\,{\rm Re~}\langle z,\xi\rangle}~L_{\wp}(d\xi)\,.
\end{eqnarray*}
With any  symbol $b\in\S_{cyl}(\Z)$ based on $\wp\Z$\,,
a {\it Weyl observable} can be associated  according to
\begin{eqnarray}
\label{weyl-obs} b^{Weyl}=\int_{\wp\Z} \F[b](z) \;\;\; W(\sqrt{2}\pi
z)~L_{\wp}(dz)\,.
\end{eqnarray}
Notice that $b^{Weyl}$ is a well defined bounded operator on $\H$ for all
$b\in\S_{cyl}(\Z)$  and that this quantization of cylindrical symbols depends on the parameter
$\varepsilon$\,.

We also recall the  Anti-Wick quantization  through its usual finite dimensional relation to Weyl
operators :
\begin{eqnarray}
\label{eq.AWW1}
b^{A-Wick}&= &
\left(b\mathop{*}_{\wp\Z}
\frac{e^{-\frac{|z|_{\wp\Z}^{2}}{\varepsilon/2}}}{(\pi\varepsilon/2)^{\textrm{dim}\wp\Z}}\right)^{Weyl}
\\
\label{eq.AWW2}
&=&
\int_{\wp\Z}
\F[b](\xi)\;\;W(\sqrt{2}\pi\xi)\;\;e^{-\frac{\varepsilon \pi^2}{2} |\xi|^2_{\wp\Z}}
\;L_{\wp}(d\xi)\,,
\end{eqnarray}
for any $b\in \mathcal{S}(\wp\Z)$ by setting
$b\mathop{*}_{\wp\Z}\gamma(z)=\int_{\wp\Z}b(z)\gamma(z-z')~L_{\wp}(dz')\,.$

\subsection{Wick quantized operators}
For any $p,q\in \mathbb{N}$\,, the space
$\P_{p,q}(\Z)$ of complex-valued polynomials on $\Z$  is defined with
the following continuity condition:
\begin{eqnarray*}
 b\in\P_{p,q}(\Z) \mbox{ iff there
 exists }  \tilde b\in\L(\bigvee^p\Z,\bigvee^q\Z)
\end{eqnarray*}
such that:
$$
b(z)=\langle z^{\otimes q}, \tilde{b} z^{\otimes p}\rangle\,.
$$
 On these spaces the
norms are given by
$|b|_{p,q}=\|\tilde{b}\|_{\mathcal{L}(\bigvee^{p}\Z;\bigvee^{q}\Z)}$\,.\\
The subspace of $\P_{p,q}(\Z)$ made of polynomials $b$ such that
$\tilde{b}$ is a compact operator  is
denoted by $\mathcal{P}^{\infty}_{p,q}(\Z)$\,.\\
The {\it Wick monomial} of a 'symbol' $b\in \P_{p,q}(\Z)$ is the linear
operator
 $$
b^{Wick}:\H_{fin}\to\H_{fin}
$$
defined as :
\begin{eqnarray}
\label{wick}
b^{Wick}_{|\bigvee^n \Z}=1_{[p,+\infty)}(n)\frac{\sqrt{n!
(n+q-p)!}}{(n-p)!} \;\varepsilon^{\frac{p+q}{2}} \;\S_{n-p+q}\left(\tilde{b}\otimes \11^{(n-p)}\right)\,.
\end{eqnarray}
where $\tilde{b}\otimes \11^{(n-p)}$ is the operator with the action $(\tilde{b}\otimes \11^{(n-p)} \varphi^{\otimes n}= (\tilde b \varphi^{\otimes p})
\otimes \varphi^{ \otimes (n-p)}$\,. Notice that $b^{Wick}$ depends on
the scaling parameter $\varepsilon$\,. When $\tilde{b}$ is an
unbounded operator with domain $\D(\tilde{b})$ containing $\bigvee^{p,
alg}\D$\,, the formula
\eqref{wick} makes sense when applied to $\Psi\in
\bigvee^{n,alg}\D$\,.

\begin{prop}
\label{wick-estimate}
For $b\in\P_{p,q}(\Z)$\,, the following number estimate holds
\begin{equation}
  \label{eq.2bis}
\left|\left\langle {\bf N}\right\rangle^{-\frac{q}{2}}b^{Wick}
\left\langle
 {\bf N}\right\rangle^{-\frac{p}{2}}\right|_{\mathcal{L}(\mathcal{H})}\leq
 \left|b\right|_{\P_{p,q}}\,.
\end{equation}
\end{prop}
\noindent An important  property of our class of Wick polynomials is
that a composition of
$b_{1}^{Wick}\circ b_{2}^{Wick}$ with
$b_1,b_2\in\oplus_{p,q}^{alg}\P_{p,q}(\Z)$
is a Wick polynomial with symbol in $\oplus_{p,q}^{\rm
  alg}\P_{p,q}(\Z)$\,.
This was checked with a convenient writing
in \cite{AmNi1} and widely used also in \cite{AmNi2,AmNi3}.\\
We need some notations:
  For $b\in\P_{p,q}(\Z)$\,, the $k$-th differential is well defined
  according to
$$
\partial_{z}^{k}b(z)\in (\bigvee^{k}\Z)^{*}\quad\text{and}\quad
\partial_{\overline{z}}^{k}b(z)\in \bigvee^{k}\Z\,,
$$
for any fixed $z\in \Z$\,. Actually $(\bigvee^{k}\Z)^{*}$ is the dual
of $(\bigvee^{k}\Z)$ with a $\cz$-bilinear duality bracket.
For two symbols $b_{i}\in \P_{p_{i},q_{i}}(\Z)$\,,
$i=1,2$\,, and any $k\in \nz$\,, the new symbol
$\partial_{z}^{k}b_{1}.\partial_{\bar
  z}^{k}b_{2}$ is now defined by
\begin{equation}
\label{eq.dotprod}
\partial_z^k b_1 \; .\;\partial_{\bar z}^k b_2 (z) =\langle \partial_z^k b_1(z),
\partial_{\bar z}^k b_2(z)\rangle_{(\bigvee^k \Z)^{*},\bigvee^{k}\Z}\quad.
\end{equation}
We also use the following notation for multiple Poisson brackets:
\begin{eqnarray*}
\{b_1,b_2\}^{(k)}&=&\partial^k_z b_1
.\partial^k_{\bar z} b_2 -\; \partial^k_z b_2 .\partial^k_{\bar z} b_1,\\
\{b_1,b_2\}&=&\{b_1,b_2\}^{(1)}.
\end{eqnarray*}
With these notations,
the composition formula of Wick symbols has  a
very familiar form.
\begin{prop}
\label{symbcalc}
Let $b_{1}\in \P_{p_{1},q_{1}}(\Z)$ and  $b_{2}\in \P_{p_{2},q_{2}}(\Z)$\,.\\
For any $k\in \left\{0,\ldots, \min\left\{p_{1},q_{2}\right\}\right\}$\,,
$\partial_{z}^{k}b_{1}.\partial_{\bar z}^{k}b_{2}$ belongs to
$\P_{ p_{1}+p_{2}-k,q_{1}+q_{2}-k}(\Z)$ with the estimate
$$
|\partial_{z}^{k}b_{1}.\partial_{\bar z}^{k}b_{2}|_{\P_{ p_{1}+p_{2}-k,q_{1}+q_{2}-k}}\leq
\frac{p_{1}!}{(p_{1}-k)!}\frac{q_{2}!}{(q_{2}-k)!}
|b_{1}|_{\P_{p_{1},q_{1}}}|b_{2}|_{\P_{p_{2},q_{2}}}\,.
$$
The formulas
\[
\begin{array}{rlcl}
(i)&
b_1^{Wick} \circ b_2^{Wick}&=&\left(\ds\sum_{k=0}^{\min\{p_1,q_2\}}
\;\;\frac{\varepsilon^k}{k!}  \;\;\;\partial^{k}_{z} b_1 .\partial^{k}_{\bar z}
b_2 \right)^{Wick}\\ \nm\ds & &=& \left(\ds e^{\varepsilon \langle \partial_z,\partial_{\bar
\omega}\rangle}
b_1(z) b_2(\omega)\left|_{z=\omega}\right. \right)^{Wick}\,,\\ \nm\ds
(ii)&
[b_1^{Wick},b_2^{Wick}]&=&\left(\ds\sum_{k=1}^{\max\{\min\{p_1,q_2\}\,,\,
\min\{p_{2},q_{1}\}\}} \;\;\frac{\varepsilon^k}{k!}  \;\;\{b_1
,b_2\}^{(k)} \right)^{Wick}\,,
\end{array}
\]
hold as identities on $\H_{fin}$\,.
\end{prop}
\noindent Combined with Proposition~\ref{wick-estimate} and
$(b^{Wick})^{*}=\left(\overline{b(z)}\right)^{Wick}$
this also gives the
\begin{prop}
\label{pr.wick-estimate2}
For $b\in\P_{p,q}(\Z)$\,, $\langle \N\rangle^{-\frac{(p+q)}{2}}b^{Wick}$  and
$b^{Wick}\langle \N\rangle^{-\frac{(p+q)}{2}}$ extend as bounded operators on
$\H$ with norm smaller that
$C_{p,q}\|\tilde{b}\|_{\mathcal{L}(\bigvee^{p}\Z;\bigvee^{q}\Z)}$\,,
for all $\varepsilon\in (0,\bar\varepsilon)$\,.
\end{prop}
\noindent We will also need some more particular estimates stated in the following two lemmata.
\begin{lem}
\label{le.compest}
Let $A$ be a self-adjoint operator on $\Z$ with $A\geq \11$\,. For any polynomials $b_1\in
\P_{1,2}(\Z)$ and $b_2\in\P_{2,1}(\Z)$ the estimates below hold true:
\begin{eqnarray*}
(i) &&\|(\d\Gamma(A)+\sqrt{\N}+1)^{-1}\; b_1^{Wick} \;(\d\Gamma(A)+\sqrt{\N}+1)^{-1}\|\leq \| A^{-1/2}\, \tilde b_1 \,(\11\otimes A^{-1/2})\|_{\L(\bigvee^2\Z,\Z)}\,, \\
(ii) &&\|(\d\Gamma(A)+\sqrt{\N}+1)^{-1}\; b_2^{Wick}\; (\d\Gamma(A)+\sqrt{\N}+1)^{-1}\|\leq \|(\11\otimes A^{-1/2})\, \tilde b_2 \,A^{-1/2}\|_{\L(\Z,\bigvee^2\Z)}\,.
\end{eqnarray*}
\end{lem}
\begin{remark}
  The term $\sqrt{\N}$ can be absorbed in $\d\Gamma(A)+1$\,, if one accepts
  constants larger than $1$ as factors of the right-hand sides of (i)
  and (ii).
\end{remark}
\proof
The estimate (ii) follows from (i) by taking the adjoint. Let us prove (i). \\
For $\Phi,\Psi\in\oplus_{n\in\nz}^{alg} \bigvee^{n, {alg}}\D(A)$\,, we write
\begin{eqnarray*}
\langle \Psi, b_1^{Wick} \Phi\rangle
&=&\sum_{n=2}^\infty \varepsilon^{3/2} \sqrt{n (n-1)^2} \,
\langle \Psi^{(n-1)}, (\tilde b_1\otimes \11^{(n-2)}) \Phi^{(n)}\rangle\\
&=& \sum_{n=2}^\infty \varepsilon^{3/2} \sqrt{n (n-1)^2}
\langle
(A^{1/2}\otimes\11^{(n-2)})\Psi^{(n-1)}\,,\,
B_A
(\11\otimes A^{1/2}\otimes\11^{(n-2)}) \; \Phi^{(n)}\rangle\,,
\end{eqnarray*}
with
$$
B_{A}= [(A^{-1/2}\tilde
  b_1)(\11 \otimes A^{-1/2})]
\otimes \11^{(n-2)}\,.
$$
Hence, by  the Cauchy-Schwarz inequality, we get
\begin{eqnarray*}
|\langle \Psi, b_1^{Wick} \Phi\rangle |&\leq&
\|A^{-1/2}\tilde b_1 \, (\11\otimes A^{-1/2})\|_{\L(\bigvee^2\Z,\Z)}
\\
&&\hspace{.1in}\times
 \left(\sum_{n=2}^\infty \varepsilon^{3/2} \sqrt{n (n-1)^2} \|(A^{1/2}\otimes \11^{(n-2)})
\Psi^{(n-1)}\|^2\right)^{1/2}  \\
 &&\hspace{.1in} \times \left(\sum_{n=2}^\infty
\varepsilon^{3/2} \sqrt{n
     (n-1)^2} \|(\11\otimes A^{1/2}\otimes \11^{(n-2)})
\Phi^{(n)}\|^2\right)^{1/2}\,.
\end{eqnarray*}
Now, observe that
\begin{eqnarray*}
 \varepsilon^{3/2} \sqrt{n (n-1)^2} \|(A^{1/2}\otimes \11^{(n-2)})
\Psi^{(n-1)}\|^2 &\leq& 2\sqrt{\varepsilon (n-1)} \langle \Psi^{(n-1)}, (\varepsilon (n-1))
(A\otimes\11^{(n-2)})\Psi^{(n-1)}\rangle \\
&\leq&2\langle \Psi^{(n-1)}, \sqrt{\N} \d\Gamma(A)\Psi^{(n-1)}\rangle
\end{eqnarray*}
and
\begin{eqnarray*}
\varepsilon^{3/2} \sqrt{n (n-1)^2} \|(\11\otimes A^{1/2}\otimes \11^{(n-2)})
\Phi^{(n)}\|^2&\leq& \sqrt{n\varepsilon} \langle \Phi^{(n)}, n\varepsilon
\, (\11\otimes A\otimes \11^{(n-2)})\Phi^{(n)}\rangle
\\ &\leq & \sqrt{n\varepsilon} \langle \Phi^{(n)}, n\varepsilon \, (A\otimes
\11^{(n-1)})\Phi^{(n)}\rangle
\\
&\leq&
\langle \Phi^{(n)}, \sqrt{\N} \d\Gamma(A)\Phi^{(n)}\rangle\,.
\end{eqnarray*}
On the other hand, with the inequality $2ab\leq a^2+b^2$\,, we see that
\begin{eqnarray*}
  &&2\langle \Psi^{(n-1)}, \sqrt{\N} \d\Gamma(A)\Psi^{(n-1)}\rangle\leq \langle \Psi^{(n-1)}, (\N+\d\Gamma(A)^2)
\Psi^{(n-1)}\rangle \leq \|(\sqrt{\N}+\d\Gamma(A))\Psi^{(n-1)}\|^2\\
\text{and}&&
2\langle \Phi^{(n)}, \sqrt{\N} \d\Gamma(A)\Phi^{(n)}\rangle\leq \langle \Phi^{(n)}, (\N+\d\Gamma(A)^2)
\Psi^{(n)}\rangle \leq \|(\sqrt{\N}+\d\Gamma(A))\Phi^{(n)}\|^2\,,
\end{eqnarray*}
where the last inequalities come from $2\sqrt{\N}\d\Gamma(A)\geq 0$\,.
Therefore, we obtain
$$
|\langle \Psi, b_1^{Wick} \Phi\rangle |\leq \|A^{-1/2}\tilde b_1 \, (\11\otimes  A^{-1/2})\| \;\;
\|(\sqrt{\N}+\d\Gamma(A))\Psi\| \;\; \|(\sqrt{\N}+\d\Gamma(A))\Phi\|
$$
and hence the estimate extends to $\Phi,\Psi\in \D(\sqrt{\N}+\d\Gamma(A))\cap\H_{fin}$\,. This means
 that the operator $(\d\Gamma(A)+\sqrt{\N}+1)^{-1} b_1^{Wick} (\d\Gamma(A)+\sqrt{\N}+1)^{-1}_{|\D(\sqrt{\N}+\d\Gamma(A))\cap\H_{fin}}$ extends to a bounded operator
 satisfying (i).
\fin
\begin{lem}
\label{le.boundedest}
Let $A, B$ two self-adjoint operators on $\Z$ with $\D(A)\subset\D(B)$ and $B\geq 0$\,. Let
$C$ be a self-adjoint operator on $\bigvee^{2} \Z$ such that
$\D(C) \subset\D(B_2)$ where $B_2=B\otimes \11+\11\otimes B$\,. Then the estimates below hold true:
\begin{eqnarray*}
(i) &&\|(\d\Gamma(B)+\N+1)^{-1}\; \d\Gamma(A) \;(\d\Gamma(B)+\N+1)^{-1}\|\leq
\| (1+B)^{-1}\, A \,(1+B)^{-1}\|_{\L(\Z)}\,, \\
(ii) &&\|(\d\Gamma(B)+\N^2+1)^{-1/2}\; C^{Wick}\; (\d\Gamma(B)+\N^2+1)^{-1/2}\|\leq \|(1+B_2)^{-1/2}\,
C \,(1+B_2)^{-1/2}\|_{\L(\bigvee^2\Z)}\,.
\end{eqnarray*}
\end{lem}
\proof
We follow a similar argument as in the proof of Lemma
\ref{le.compest}. Indeed, the Cauchy-Schwarz inequality
gives for every $\Psi,\Phi\in \oplus_{n\in\nz}^{alg} \bigvee^{n,alg}\D(B)$
\begin{eqnarray*}
|\langle \Psi, \d\Gamma(A) \,\Phi\rangle |&\leq& \|(1+B)^{-1} A \, (1+B)^{-1}\|_{\L(\Z)}
\;  \left(\sum_{n=1}^\infty \varepsilon n \| ((1+B)\otimes \11^{(n-1)})
\Psi^{(n)}\|^2\right)^{1/2}  \\
 &&\hspace{.1in} \times \left(\sum_{n=1}^\infty \varepsilon n \|((1+B)\otimes \11^{(n-1)})
\Phi^{(n)}\|^2\right)^{1/2}\,.
\end{eqnarray*}
Now,  observe that
\begin{eqnarray*}
 \varepsilon n \|((1+B)\otimes \11^{(n-1)})
\Psi^{(n)}\|^2 &=&  \langle \Psi^{(n)}, \d\Gamma((1+B)^2)\Psi^{(n)}\rangle\\
&\leq & \|\d\Gamma(1+B)\Psi^{(n)}\|^2_{\bigvee^n\Z}\,,
\end{eqnarray*}
since in the sense of quadratic forms $\d\Gamma((1+B)^2)\leq \d\Gamma(1+B)^2$\,.
Hence we obtain
$$
|\langle \Psi, \d\Gamma(A) \,\Phi\rangle | \leq \|(1+B)^{-1} A \, (1+B)^{-1}\|_{\L(\Z)} \,
\|\d\Gamma(1+B)\Psi\| \,\|\d\Gamma(1+B)\Phi\|\,.
$$
This proves (i). \\
Expressing $C^{Wick}$ as a quadratic form for $\Psi,\Phi \in \oplus_{n\in\nz}^{alg }
\bigvee^{n, alg}\D(B)$ and then applying the Cauchy-Schwarz inequality
yield
\begin{eqnarray*}
|\langle \Psi, C^{Wick} \,\Phi\rangle |&\leq&
\|(\11+B\otimes \11+\11\otimes B)^{-1/2} C \, (\11+B\otimes \11+\11\otimes B)^{-1/2}\|_{\L(\bigvee^2\Z)}
\\ &&\hspace{.1in}
\times
\;  \left(\sum_{n=2}^\infty \varepsilon^2 n (n-1) \|
[(\11+B\otimes \11+\11\otimes B)^{1/2}\otimes \11^{(n-2)}]
\Psi^{(n)}\|^2\right)^{1/2}  \\
 &&\hspace{.1in} \times
\left(\sum_{n=2}^\infty \varepsilon^2 n (n-1) \|[(\11+B\otimes \11+\11\otimes B)^{1/2}
 \otimes \11^{(n-2)}]\Phi^{(n)}\|^2\right)^{1/2}\,.
\end{eqnarray*}
Due to the symmetry of the vector $\Phi^{(n)}$ we remark that
\begin{eqnarray*}
\varepsilon^2 n^2 \|[(\11+B\otimes \11+\11\otimes B)^{1/2}\otimes \11^{(n-2)}]\Phi^{(n)}\|^2
&=&\varepsilon^2 n^2 \langle \Phi^{(n)}, [
(\11+B\otimes \11+\11\otimes B)\otimes \11^{(n-2)}]
\Phi^{(n)}\rangle\\&= &
{
\langle \Phi^{(n)}, (1+\d\Gamma(B))\, \N^2 \,\Phi^{(n)}\rangle}\\
&\leq& \| (1+\d\Gamma(B)+\N^2) \,\Phi^{(n)}\|^2\,
\end{eqnarray*}
So we obtain
\begin{eqnarray*}
|\langle \Psi, C^{Wick} \,\Phi\rangle |&\leq& \|(1+B_2)^{-1/2} C \, (1+B_2)^{-1/2}\|_{\L(\bigvee^2\Z)}
\; \| (1+\d\Gamma(B)+\N^2) \,\Psi\| \\
 &&\hspace{.1in} \times \| (1+\d\Gamma(B)+\N^2) \,\Phi\|\,.
\end{eqnarray*}
This proves (ii).\fin

\section{Absolutely continuous curves in $\mathrm{Prob}_{2}(\Z_{1,\rz})$}
\label{se.abscont}
This section firstly gathers  results presented in \cite{AGS} about Borel probability
measures on a separable \underline{real} Hilbert space which are weak solutions to
continuity equations.
In a second step, we shall adapt it to a complex Hilbert space
$\Z_{1}$ endowed with its real euclidean structure.

\subsection{Absolutely continuous curves in $\mathrm{Prob}_{2}(E)$}
\label{se.abscontE}
Let $E$ be a real Hilbert space, with scalar product
$\langle\,,\,\rangle$ and norm $|~|$\,.
The symbol $\mathrm{Prob}_{p}(E)$ (resp. $\mathrm{Prob}(E)$) refers to the set of Borel
probability measures $\mu$ on $E$ such that
$\int_{E}|x|^{p}~d\mu(x)<+\infty$ (resp. with no momentum condition), and
we simply work with $p=2$\,. On $\mathrm{Prob}_{2}(E)$\,,
 the $2$-Wasserstein distance, $W_{2}$\,, is defined by
$$
W_2^{2}(\mu^{1},\mu^{2}):=\min\left\{\int_{E^{2}}|x_{1}-x_{2}|_{E}^{2}~d\mu(x_{1},x_{2})\,;\;
\Pi_{j,*}\mu=\mu_{j}\right\}\,,
$$
where $\Pi_{j}:E^{2}\to E$ is the natural projection, $j=1,2$\,.
The narrow convergence of a sequence $(\mu_{n})_{n\in\nz}$ of
$\mathrm{Prob}_{2}(E)$\,,  with a uniform control of
$\int_{E}|x|^{2}~d\mu_{n}$
is equivalent to the $W_{2}$ convergence on $\mathrm{Prob}_{2}(E)$
(see Proposition~7.1.5 in \cite{AGS})\,.
 Remember also that the tightness property
of subsets of $\mathrm{Prob}_{2}(E)$ can be checked in the infinite
dimensional case with the weak topology,
or after introducing a Hilbert basis $(e_{n})_{n\in\nz^{*}}$\,, with the
distance $d_{\omega}(x_{1},x_{2})=\sqrt{\sum_{n\in\nz^{*}}\frac{|\langle
  x_{1}-x_{2},e_{n}\rangle|^{2}}{(1+n)^{2}}}$\,.
This use of weak or $d_{\omega}$ topology, is done also when
considering probability measures on the set of absolutely continuous
curves in $E$\,.\\
This tightness property is called the weak tightness property in
\cite{AGS} since it refers to the weak topology on $E$\,.
Especially when one considers the narrow convergence in
$\textrm{Prob}_{2}(E)$\,, there is a weak narrow convergence and a strong
narrow convergence (see the discussions about this in Chapter~5 and 7
of \cite{AGS}). The terms ``narrow convergence'' or ``narrow continuity'' refer  \textbf{to
the strong ones} and we shall specify ``weak narrow convergence'' and
``weak narrow continuity'' when necessary. \\
We recall two results of \cite{AGS} and give a complete proof in the
infinite dimensional case of the second one, for the sake of
completeness (it is left as an exercise to the reader in \cite{AGS}).

The following result is the second part of Theorem~8.3.1 in \cite{AGS}
with $p=2$\,. Although it is not clearly stated in Theorem~8.3.1 in
\cite{AGS}, the proof contains a ``weak$\Rightarrow$ strong'' result
about the narrow continuity w.r.t time.
\begin{prop}
\label{pr.conteq}
Let $I$ be an open interval in $\rz$\,.
If a weakly narrowly continuous curve $\mu_{t}:I\to \mathrm{Prob}_{2}(E)$
satisfies the continuity equation
\begin{equation}
  \label{eq.conteq}
\partial_{t}\mu_{t}+\nabla^{T}(v_{t}\mu_{t})=0
\end{equation}
in the weak sense
\begin{equation}
  \label{eq.conteqw}
\int_{I}\int_{E}
\left(\partial_{t}\varphi(x,t)+\langle v_{t}(x)\,,\,
  \nabla_{x}\varphi(x,t)\rangle_{E}\right)~d\mu_{t}(x)dt=0\,,\quad
\forall \varphi \in \mathcal{C}^{\infty}_{0, \textrm{cyl}}(E\times I)\,,
\end{equation}
for some Borel velocity field $v_{t}$\,, with
$|v_{t}|_{L^{2}(E,\mu_{t})}\in L^{1}(I)$\,, then $\mu_{t}:I\to
\mathrm{Prob}_{2}(E)$ is absolutely continuous with $W_{2}(\mu_{t'},\mu_{t})\leq
\int_{t}^{t'}|v_{s}|_{L^{2}(E,\mu_{s})}~ds$\,. Moreover for
Lebesgue almost every $t\in I$\,, $v_{t}$ belongs to the closure in
$L^{2}(E,\mu_{t})$ of the subspace spanned by $\left\{\nabla\varphi,
  \varphi\in \mathcal{C}^{\infty}_{0,\textrm{cyl}}(E)\right\}$\,.
\end{prop}
\proof
  The proof is given in \cite{AGS}. We simply insist here on the
  ``weak$\Rightarrow$ strong'' narrow continuity argument. The proof
  of Theorem~8.3.1 ends with the following statements. For any time
$t\in I$\,, $\mu_{t}$ is the weak narrow limit of a sequence
$(\hat{\mu}_{t}^{d})_{d\in\nz}$ (of which the definition is recalled
below) which satisfies
$$
W_{2}(\hat{\mu}_{t_{2}}^{d},\hat{\mu}_{t_{1}}^{d})\leq
\int_{t_{1}}^{t_{2}}|v_{t}|_{L^{2}(E,\mu_{t})}~dt\,.
$$
Then the authors refer to the  weak
narrow lower semicontinuity of $W_{2}$\,,
$$
W_{2}(\mu_{t_{2}},\mu_{t_{1}})\leq \liminf_{d\to\infty}W_{2}(\hat{\mu}_{t_{2}}^{d},\hat{\mu}_{t_{1}}^{d})
$$
stated in their Lemma~7.1.4 and relation (7.1.11).\\
This implies absolute continuity in terms of $W_{2}$ and the narrow
continuity w.r.t time:
$$
W_{2}(\mu_{t_{1}},\mu_{t_{2}})\leq \int_{t_{1}}^{t_{2}}|v_{t}|_{L^{2}(E,\mu_{t})}~dt\,.
$$
Finally notice that the weak narrow continuity of $\mu_{t}$ suffices
for the strong narrow continuity of $\hat{\mu}_{t}^{d}$ w.r.t $t\in I$
because $\hat{\mu}_{t}^{d}$ is constructed after taking the image of
$\mu_{t}$ via a finite rank projection.
\fin

The previous result concerns non regular (non Lipschitz) vector fields
for which there is no uniqueness result for the Cauchy
problem. Remember that the infinite dimensional case, which relies on
the cylindrical integration of $v_{t}$ and cylindrical disintegration of the
measure $\mu_{t}$\,, requires the introduction of such singular vector
fields (see the proof of Theorem~8.3.1 in \cite{AGS}). Nevertheless an
interpretation of the continuity equation
\eqref{eq.conteq}\eqref{eq.conteqw} in terms of characteristic curves can be
done via a probabilistic representation. For the sake of completeness,
we adapt the proof of Theorem~8.2.1  stated  in \cite{AGS} for the
finite dimensional case, to our infinite dimensional case.
For $T\in (0,+\infty)$\,, consider the set $\Gamma_{T}=
\mathcal{C}^{0}([-T,T];E)$ endowed with the norm
$|\gamma|_{\infty,T}=\max_{t\in [-T,T]}|\gamma(t)|_{E}$ or for weak
topology argument with the distance $\max_{t\in
  [-T,T]}d_{\omega}(\gamma_{1}(t),\gamma_{2}(t))$\,. For  a Borel
probability measure $\bm{\eta}$ defined on $E\times \Gamma_{T}$\,,
consider the time dependent Borel probability measure
$\mu_{t}^{\bm{\eta}}$ defined by
\begin{equation}
  \label{eq.defmueta}
  \int_{E}\varphi~d\mu_{t}^{\bm{\eta}}=\int_{E\times
    \Gamma_{T}}\varphi(\gamma(t))~d\bm{\eta}(x,\gamma)\,,\quad
\forall \varphi\in \mathcal{C}^{0}_{b,cyl}(E),\, t\in [-T,T]\,.
\end{equation}
The measure $\mu_{t}^{\bm{\eta}}$ is the push-forward of
$\bm{\eta}$ by the evaluation map
$$
e_{t}: (x,\gamma)\in E\times \Gamma_{T}\to \gamma(t)\in E\,,\quad
\text{for}~t\in [-T,T]\,.
$$
\begin{prop}
  \label{pr.probrep}
Let $\mu_{t}:[-T,T]\to \mathrm{Prob}_{2}(E)$ be a  $W_{2}$-continuous
solution to the continuity equation
\eqref{eq.conteq}\eqref{eq.conteqw}, with $I=(-T,T)$\,, for a suitable Borel vector field
$v(t,x)=v_{t}(x)$ such that $|v_{t}|_{L^{2}(E,\mu_{t})}\in
L^{1}([-T,T])$\,.
Then there exists a Borel probability measure $\bm{\eta}$ in $E\times
\Gamma_{T}$ such that
\begin{description}
\item[(i)] $\bm{\eta}$ is concentrated on the set of pairs
  $(x,\gamma)$ such that $\gamma\in AC^{2}([-T,T];E)$ is a solution to
  the ODE $\dot\gamma(t)=v_{t}(\gamma(t))$ for Lebesgue
  almost every $t\in (-T,T)$ with $\gamma(0)=x$\,;
\item[(ii)] $\mu_{t}=\mu_{t}^{\bm{\eta}}$ for any $t\in [-T,T]$\,,
  with $\mu_{t}^{\bm{\eta}}$ defined as in \eqref{eq.defmueta}.
\end{description}
Conversely, any $\bm{\eta}$ satisfying (i) and
$$
\int_{0}^{T}\int_{E\times
  \Gamma_{T}}|v_{t}(\gamma(t))|_{E}~d\bm{\eta}(x,\gamma)dt < +\infty\,,
$$
induces via \eqref{eq.defmueta} a solution to the continuity equation,
with $\mu_{0}=\gamma(0)_{*}\bm{\eta}$\,.
\end{prop}
\begin{remark}
  The notation $AC^{2}([-T,T];E)$ refers to the set of absolutely
  continuous curves in $E$ with $L^{2}([-T,T];E)$ derivative.
We keep the notation $\Phi_{*}\mu$ of differential geometry, for the
push-forward or direct image of a measure $\mu$\,, by the Borel map
$\Phi$\,.
\end{remark}
\proof
The result is proved in \cite{AGS} when $E$ is finite dimensional.
The proof of the second (converse) part of the statement, is exactly
the same as in finite dimension, after replacing regular (Lipschitz) test
functions by cylindrical ones. We now show, for the first part, how the infinite
dimensional case is deduced from the finite dimensional result,
following an approximation scheme like in the proof of
\cite[Theorem~8.3.1]{AGS}. After introducing an Hilbert
basis $(e_{n})_{n\in\nz^{*}}$ of $E$\,, the maps $\pi^{d}:E\to \rz^{d}$\,,
$\pi^{d,T}:\rz^{d}\to E$ and $\hat{\pi}^{d}: E\to E$ are
defined according to
\begin{eqnarray*}
  \pi^{d}(x)&=& (\langle e_{1},x\rangle,\ldots, \langle e_{d},
  x\rangle)\,,\\
\pi^{d,T}(y_{1},\ldots,y_{d})&=&
\sum_{j=1}^{d}y_{j}e_{j}\,,\\
\hat{\pi}^{d}&=&\pi^{d,T}\circ \pi^{d}\,.
\end{eqnarray*}
With the measure $\mu_{t}\in \mathrm{Prob}_{2}(E)$\,, the measure
$\mu_{t}^{d}\in \mathrm{Prob}_{2}(\rz^{d})$ is defined by
$\mu_{t}^{d}=\pi^{d}_{*}\mu_{t}$ and $\left\{\mu_{t,y}, y\in
  \rz^{d}\right\}$ denotes the disintegration of $\mu_{t}$ w.r.t
$\mu_{t}^{d}$\,. Within the space $E$ endowed with the basis
$(e_{n})_{n\in\nz^{*}}$\,, $\hat{\mu}_{t}^{d}$  is nothing but
$\mu_{t}^{d}\otimes \delta_{0}$ in the decomposition $Z=F_{d}\times
F_{d}^{\perp}$ with $F_{d}=\textrm{span}(e_{1},\ldots,e_{d})$\,.
 The vector field $v_{t}^{d}$
(resp. $\hat{v}_{t}^{d}$) is defined on $\rz^{d}$ (resp. on $E$)
by
\begin{eqnarray*}
&&v_{t}^{d}(y)=\int_{(\pi^{d})^{-1}(y)}\pi^{d}v_{t}(x)~d\mu_{t,y}(x)\,,
\quad y\in\rz^{d}
\\
\text resp.&&
\hat{v}_{t}^{d}(y) =
\int_{(\hat\pi^{d})^{-1}(\hat{\pi}^{d}y)}\hat\pi^{d}v_{t}(x)~d\mu_{t,{
  \pi^{d}y}}(x)\,,\quad y\in E\,.
\end{eqnarray*}
Within the proof of Theorem~8.3.1 in \cite{AGS}, it was checked that
$\mu_{t}^{d}$ (resp. $\hat\mu_{t}^{d}$) is a weak solution to the continuity equation
\begin{eqnarray*}
&&\partial_{t}\mu_{t}^{d}+ \nabla^{T}(v_{t}^{d}\mu_{t}^{d})=0\,,
\\\text{resp.}&&
\partial_{t}\hat\mu_{t}^{d}+ \nabla^{T}(\hat v_{t}^{d}\hat\mu_{t}^{d})=0\,.
\end{eqnarray*}
with the following properties:
\begin{description}
\item[1)] $|\hat{v}_{t}^{d}|_{L^{2}(E, \hat\mu_{t}^{d})}=|v_{t}^{d}|_{L^{2}(\rz^{d}, \mu_{t}^{d})}\leq
  |v_{t}|_{L^{2}(E,d\mu_{t})}$\,;
\item[2)] $W_2(\mu_{t_{1}}^{d},\mu_{t_{2}}^{d})\leq
  \int_{t_{1}}^{t_{2}}|v_{t}^{d}|_{L^{2}(\rz^{d},\mu_{t}^{d})}~dt\leq
 \int_{t_{1}}^{t_{2}}|v_{t}|_{L^{2}(E,\mu_{t})}~dt$\,, for
 $-T<t_{1}\leq t_{2}< T$\,;
\item[3)] the
  sequence $({\hat\mu}_{t}^{d})_{d\in\nz^{*}}$ converges weakly narrowly to
  $\mu_{t}$ with the estimate
\begin{equation}
\label{est.chavit}
W_2(\mu_{t_{1}},\mu_{t_{2}})\leq \liminf_{d\to
  \infty}W_2(\mu_{t_{1}}^{d}, \mu_{t_{2}}^{d})\leq
\int_{t_{1}}^{t_{2}}|v_{t}|_{L^{2}(E,\mu_{t})}~dt\,,\quad
-T< t_{1}\leq t_{2}<T\,.
\end{equation}
\end{description}
Additionally a time rescaling argument (see Lemma~1.1.4 and
Theorem~8.1.3 in \cite{AGS}) allows to assume without restriction
$$
|v_{t}|_{L^{2}(E,\mu_{t})}\in L^{\infty}((-T,T))\,.
$$
The set of continuous maps from $[-T,T]$ to $\rz^{d}$ is denoted by
$\Gamma_{T}^{d}$\,.
The mapping from $\Gamma_{T}$ to $\Gamma_{T}^{d}$\,,
still denoted by $\pi^{d}$ is defined
by $[\pi^{d}\gamma](t)=\pi^{d}(\gamma(t))$\,.
By using the finite dimensional result, stated in Theorem~8.2.1 of
\cite{AGS}, there exists for any $d\in\nz$ a probability measure,
${\bm \eta}^{d}$\,,
 on $\rz^{d}\times \Gamma_{T}^{d}$ such that the
properties (i) and (ii) hold when $(\mu_{t},v_{t},E)$ is replaced by
$(\mu_{t}^{d},v_{t}^{d}, \rz^{d})$\,. Equivalently the result can be
formulated in $E$ after using
$({\hat\mu}_{t}^{d}, {\hat v}_{t}^{d},E)$ instead of
$(\mu_{t}^{d},v_{t}^{d}, \rz^{d})$  and using
${\hat{\bm \eta}}^{d}=(\pi^{d,T}\times \pi^{d,T})_{*}{\bm \eta}^{d}$\,.
Hence we have a sequence $({\hat{\bm\eta}^{d}})_{d\in\nz}$ of
probability measures on $E\times \Gamma_{T}$ which satisfy
\begin{eqnarray}
\nonumber
&&\int_{E}\varphi\circ \pi^{d}~d\hat{\mu}_{t}^{d}
=
\int_{\rz^{d}}\varphi~d\mu_{t}^{d}= \int_{\rz^{d}\times
  \Gamma_{T}^{d}}\varphi(\gamma(t))~d{\bm \eta}^{d}(x,\gamma)
\\
\label{eq.relmueta}
&&\hspace{1cm}=
\int_{E\times
  \Gamma_{T}}\varphi\circ\pi^{d}(\gamma(t))~d{\hat{\bm\eta}}^{d}(x,\gamma)\,,\;
\forall \varphi\in \mathcal{C}^{0}_{b}(\rz^{d}),\, t\in[-T,T]\,,
\end{eqnarray}
where $\varphi\circ\pi^{d}$ can be replaced by $\varphi\circ
\hat{\pi}^{d}$ with $\varphi\in \mathcal{C}^{0}_{b}(F_{d})$\,.\\
After some regularization done \cite{AGS}~pp179-180, it is proved that any measure
${\hat{\bm \eta}}^{d}$ satisfies
$$
\int_{\Gamma_{T}}\int_{-T}^{T}|\dot\gamma(t)|^{2}~dt
d{\hat{\bm\eta}^{d}}\leq
\int_{-T}^{T}\int_{E}|v_{t}^{d}(x)|^{2}~d{\hat\mu}_{t}^{d}dt
\leq \int_{-T}^{T}\int_{E}|v_{t}(x)|^{2}~d{\hat\mu}_{t}dt\,.
$$
Since the functional $g\to \int_{-T}^{T}|\dot g(t)|^{2}~dt$\,,
defined on $\left\{g\in \Gamma_{T}\,,\quad g(0)=0\right\}$ and
set to $+\infty$ if $g\not\in AC^{2}([-T,T];E)$\,,  has compact sublevel
sets in $\Gamma_{T}$\,,
 the two
mappings
$$
r^{1}:(x,\gamma)\in E\times \Gamma_{T}\to x\in E\quad, \quad
r^{2}:(x,\gamma)\in E\times \Gamma_{T}\to g_{\gamma,x}=\gamma-x\in \Gamma_{T}
$$
give rise to (weakly) tight families of marginals
$(r^{1}_{*}{\hat{\bm\eta}^{d}})_{d\in\nz}=({\hat\mu}_{0}^{d})_{d\in\nz}$
and
$(r^{2}_{*}{\hat{\bm\eta}^{d}})_{d\in\nz}$\,. Remember that the
compactness of subsets of $E$ or $\Gamma_{T}$ is considered with the weak
topology on $E$ or the distance $d_{\omega}$\,.
Hence the family $({\hat{\bm \eta}}^{d})_{d\in\nz}$ is (weakly) tight in
$\mathrm{Prob}(E\times \Gamma_{T})$ and we take for $\bm\eta$ a weak narrow
limit point of ${\hat{\bm \eta}}^{d}$\,. By assuming the test function
$\varphi$ in \eqref{eq.relmueta} to depend only on $d'$ coordinates
with $d'\leq d$\,, and by taking the limit $d\to +\infty$ while $d'$ and
$\varphi$ are fixed, we get
$$
\int_{E}(\varphi\circ \pi^{d'})~d\mu_{t}=
\int_{E\times
  \Gamma}(\varphi\circ\pi^{d'})(\gamma(t))~d{\bm\eta}(x,\gamma)\,,
$$
for all $\varphi\in \mathcal{C}^{0}_{b}(\rz^{d'})$ and $t\in[-T,T]$\,,
where $\varphi\circ\pi^{d'}$ can then be replaced by any cylindrical
function or Borel bounded function on $E$\,.
It remains to prove the condition {\textbf(i)} for $\bm \eta$\,,  namely that
this measure is concentrated on curves verifying
$\dot\gamma=v(\gamma(t))$ for Lebesgue almost every $t\in (-T,T)$ ($\gamma(0)=x$ is
already known)\,.\\
The estimate (8.2.6) used in \cite{AGS} for the finite dimensional
case, provides the inequality
$$
\int_{\rz^{d}\times
  \Gamma_{T}^{d}}|\gamma(t)-x-\int_{0}^{t}w_{s}(\gamma(s))~ds|^{2}~d{\bm\eta^{d}}(x,\gamma)
\leq (2T)\int_{-T}^{T}\int_{\rz^{d}}|v_{t}^{d}-w_{t}|^{2}d\mu_{t}^{d}~dt\,,
$$
for any family $w_{s}(x)=w(s,x)$
uniformly  bounded continuous function from $[-T,T]\times\rz^{d}$ to
$\rz^{d}$\,. After assuming that $w$ actually belongs to
$\mathcal{C}^{0}_{b}([-T,T]\times \rz^{d'};\rz^{d'})$ with $d'\leq d$
fixed and, by using $\hat{w}_{t}=\pi^{d',T}\circ w_{t}\circ \pi_{d'}\in
\mathcal{C}^{0}_{b}([-T,T]\times E;E)$\,, taking the limit as $d\to
\infty$ gives
 $$
\int_{E\times
  \Gamma_{T}}|\gamma(t)-x-\int_{0}^{t}\hat{w}_{s}(\gamma(s))~ds|^{2}~d{\bm\eta}(x,\gamma)
\leq
(2T)\limsup_{d\to\infty}\int_{-T}^{T}\int_{\rz^{d}}|\hat{v}_{t}^{d}-\hat{w}_{t}|^{2}d\hat\mu_{t}^{d}~dt\,.
$$
But the condition 1) for
$|\hat{v}_{t}^{d}|_{L^{2}(E,\hat\mu_{t}^{d})}$ is easily extended  to
$$
|\hat v_{t}^{d}-\hat w_{t}|_{L^{2}(E,\hat\mu_{t}^{d})}
\leq
| v_{t}-\hat w_{t}|_{L^{2}(E, \mu_{t})}\,,
$$
by the same argument, relying on
\begin{eqnarray*}
  \left|\int_{E} \langle \hat v_{t}^{d}-\hat w_{t}\,,\,
    \chi\rangle~d\hat\mu_{t}^{d}\right|
&\leq& \left|\int_{E}\langle \hat\pi_{d}(v_{t}(x)-\hat w_{t}(x))\,,\,
  \chi(\hat{\pi}^{d}(x))\rangle~d\mu_{t}(x)\right| \\
&\leq& |v_{t}-\hat{w}_{t}|_{L^{2}(E;\mu_{t})}|\chi|_{L^{2}(E,
  \hat{\mu}_{t}^{d})}\,,\quad
\forall \chi\in L^{2}(E,\hat\mu_{t}^{d})\,.
\end{eqnarray*}
This uniform upper
bound leads to
$$
\int_{E\times
  \Gamma_{T}}|\gamma(t)-x-\int_{0}^{t}\hat{w}_{s}(\gamma(s))~ds|^{2}~d{\bm\eta}(x,\gamma)
\leq
(2T) \int_{0}^{T}\int_{E}|v_{t}(x)-\hat w_{t}(x)|^{2}~d\mu_{t}dt\,.
$$
According to the last statement of Proposition~\ref{pr.conteq},
$v_{t}$ can be approximated in $L^{2}(E,\mu_{t})$ by a sequence of
bounded regular cylindrical functions, $(\hat{w}_{t,n})_{n\in\nz}$\,. By
possibly truncating with respect to times $t\to \hat{w}_{t,n}$ so that
$|v_{t}-\hat{w}_{t,n}|_{L^{2}(E,\mu_{t})}\leq~1$ for a.e. $t\in (-T,T)$ and all
$n\in\nz^{*}$\,,  Lebesgue's
Theorem implies
$$
\int_{E\times
  \Gamma_{T}}|\gamma(t)-x-\int_{0}^{t}v_{s}(\gamma(s))~ds|^{2}~d\hat{\bm\eta}(x,\gamma) =0\,,
$$
which ends the proof.
\fin
Below is a consequence of the above probabilistic interpretation when
the Cauchy problem $\dot\gamma(t)=v_{t}(\gamma(t))$\,, $\gamma(0)=x$ admits a unique
solution for all $x\in E$\,. The fact that we have to pass by the
probabilistic representation is a real question. Contrary to the
finite dimensional case, the well-posedness of the Cauchy problem, even
with the standard Picard's contraction argument,
defining a flow on the whole space $E$\,, does not give a
representation formula for observables. The point is that the natural
observables, or test functions, are cylindrical functions, a property
which is not generally preserved by the nonlinear flow.
\begin{prop}
\label{pr.flowmeas}
Let $\mu_{t}:\rz\to \mathrm{Prob}_{2}(E)$ be a $W_{2}$-continuous
solution to the continuity equation
\eqref{eq.conteq}\eqref{eq.conteqw} for a suitable Borel velocity
field $v(t,x)=v_{t}(x)$ such that $|v_{t}|_{L^{2}(E,\mu_{t})}\in
L^{1}([-T,T])$ for all $T>0$\,. Assume additionally that the Cauchy
problem
$$
\dot\gamma(t)=v_{t}(\gamma(t))\,,\quad \gamma(s)=x
$$
or
$$
\gamma(t)=x+\int_{s}^{t}v_{s}(\gamma(s))~ds\,,
$$
admits a unique global continuous solution on $\rz$ for all $s\in\rz$
and all $x\in E$\,, such that $\gamma(t)=\Phi(t,s)\gamma(s)$ defines a
Borel flow on $E$ (i.e. $\Phi(t,s):E\to E$ is a Borel function for
all $t,s\in\rz$). Then the measure $\mu_{t}$ satisfies
$$
\forall t,s\in\rz,\quad \mu_{t}=\Phi(t,s)_{*}\mu_{s}\,.
$$
\end{prop}
\proof
It suffices to work with  $t\in [-T,T]$ as in Proposition~\ref{pr.probrep}.
  Since the evaluation map $e_{t}:E\times (x,\gamma)\Gamma_{T}\to
  \gamma(t)\in E$ is a continuous, thus Borel, map. The relation
  $\mu_{t}=\mu_{t}^{\bm\eta}$ defined according to \eqref{eq.defmueta}
  extends to any bounded Borel function $\varphi$ on $E$:
$$
\int_{E}\varphi~d\mu_{t}=\int_{E\times \Gamma_{T}}\varphi(\gamma(t))~d{\bm\eta}(x,\gamma)\,.
$$
By using $\gamma(t)=\Phi(t,s)\gamma(s)$\,, with $\Phi(t,s)$ Borel, we
deduce
$$
\int_{E}\varphi~d\mu_{t}=\int_{E\times \Gamma_{T}}[\varphi\circ
\Phi(t,s)](\gamma(s))~d{\bm\eta}(x,\gamma)=
\int_{E}[\varphi\circ \Phi(t,s)]~d\mu_{s}
$$
which is nothing but $\mu_{t}=\Phi(t,s)_{*}\mu_{s}$\,.
\fin

\subsection{Application to Hamiltonian fields}
\label{se.abscontZ1}
We finally specify how  these results apply to
our case, when the phase-space $\Z_{0}$ is a complex Hilbert space and the
velocity field is associated with a (singular) Hamiltonian vector
field, only defined on $\Z_{1}\subset \Z_{0}$\,.\\
Consider a complex Hilbert triple $\Z_{1}\subset \Z_{0}\subset
\Z_{-1}$\,, with
$\Z_{1}$ densely continuously embedded in $\Z_{0}$ and $\Z_{-1}$ being
the dual of $\Z_{1}$ for the duality bracket extending
$\langle z_{1}\,,\,z_{2}\rangle_{\Z_{0}}$\,. The dual of
a complex Hilbert space $\Z$ while keeping the $\cz$-bilinear duality
bracket, written $u\,.\,v$ in \eqref{eq.dotprod}, is still denoted by $\Z^{*}$\,.  In the case
treated in the article $\Z_{0}=L^{2}(\rz^{d},dx)$\,,
$\Z_{1}=H^{1}(\rz^{d})$ and $\Z_{-1}=H^{-1}(\rz^{d})$\,.
The space $\Z_{0}$ is endowed with its scalar product
$\langle z_{1}\,,\,z_{2}\rangle_{\Z_{0}}$\,, real euclidean structure with
$\Real \langle z_{1}\,,\,z_{2}\rangle_{\Z_{0}}$ and its symplectic structure
$\sigma(z_{1},z_{2})=\Imag \langle z_{1}\,,\, z_{2}\rangle_{\Z_{0}}$\,.
On $\Z_{1}$ we will use the hermitian $\langle
z_{1}\,,\,z_{2}\rangle_{\Z_{1}}$ and euclidean  scalar product
$$
\langle z_{1}\,,\,z_{2}\rangle_{\Z_{1},\rz}=\Real\langle
z_{1}\,,\, z_{2}\rangle_{\Z_{1}}\,.
$$
For a cylindrical function $f\in \mathcal{S}_{cyl}(\Z_{0})$\,, based on
$\wp \Z_{0}$\,,  the
differentials $\partial_{z}f(z)$ and $\partial_{\overline{z}}f(z)$ are defined
\begin{eqnarray*}
  &&\partial_{z}f(z)=\int_{\wp \Z_{0}} i\pi\langle \xi|
 e^{2i\pi \Real    \langle z\,,\,\xi\rangle}\mathcal{F}[f](\xi)~L_{\wp}(d\xi)\,\\
&& \partial_{\overline{z}}f(z)=\int_{\wp \Z_{0}} i\pi |\xi\rangle e^{2i\pi \Real
    \langle z\,,\,\xi\rangle}\mathcal{F}[f](\xi)~L_{\wp}(d\xi)\,.
\end{eqnarray*}
Hence $\partial_{z}f(z)$ is a continuous $\cz$-linear form on $\Z_{0}$
while $\partial_{\overline{z}}f(z)\in \Z_{0}$\,. This notation is
coherent with
the definition of $\partial_{z}b(z)$ and $\partial_{\overline{z}}b(z)$
when $b$ is a Wick symbol in $\oplus_{p,q}^{alg}\P_{p,q}(\Z)$\,.\\
A function $f\in
\mathcal{S}_{cyl}(\Z_{1})$ is given by
$$
f(z)=\varphi(\langle\xi_{1}\,,\, z\rangle_{\Z_{1}},\ldots ,\langle
\xi_{N}\,,\, z\rangle_{\Z_{1}})=\varphi(\langle \eta_{1}\,,\,
z\rangle,\ldots, \langle \eta_{N}\,,\, z\rangle)
$$
with $\varphi(w_{1},\ldots w_{N})\in \mathcal{S}(\rz^{N})$ and
 $\xi_{1},\ldots,\xi_{N}\in \Z_{1}$ and
$\eta_{1},\ldots,\eta_{N}\in \Z_{-1}$\,, such that $\langle \xi_{j}\,,\,
z\rangle_{\Z_{1}}=\langle\eta_{j}\,,\, z\rangle$ for all
$z\in\Z_{1}$\,. The derivatives $\partial_{z}f$ and
$\partial_{\overline{z}}f$ are thus given by
\begin{eqnarray*}
\forall z\in \Z_{1}\,,&&  \partial_{z}f(z)=\sum_{j=1}^{N}\partial_{w_{j}}\varphi(\langle \eta_{1}\,,\,
z\rangle,\ldots, \langle \eta_{N}\,,\, z\rangle) \langle
  \eta_{j}|\quad\in \Z_{1}^{*}\,,\\
 \forall z\in \Z_{1},&&
\partial_{\overline{z}}f(z)=\sum_{j=1}^{N}\partial_{\overline{w_{j}}}\varphi(\langle \eta_{1}\,,\,
z\rangle,\ldots, \langle \eta_{N}\,,\, z\rangle)
  |\eta_{j}\rangle\quad\in \Z_{-1}\,.
\end{eqnarray*}
When $h(z)$ is an unbounded polynomial on $\Z_{0}$ but which happens
to be a real-valued Fr{\'e}chet $\mathcal{C}^{1}$-function on $\Z_{1}$\,, the derivatives
 $\partial_{z}h(z)$ and $\partial_{\overline{z}}h(z)$
 are defined only for $z\in \Z_{1}$ and we have
$$
\forall z\in \Z_{1}\,,\quad \partial_{\overline{z}}h(z)\in
\Z_{1}\quad,\quad
\partial_{z}h(z)\in \Z_{-1}^{*}\,.
$$
When $f\in \mathcal{S}_{cyl}(\Z_{0})$ (resp. $g\in
\mathcal{S}_{cyl}(\Z_{1})$ or $h$) is real valued
differentiating $f(z+te)$ at $t=0$\,, $t\in\rz$\,, for any $e\in\Z_{0}$
(resp. any $e\in \Z_{1}$) leads to
\begin{eqnarray}
\label{eq.realZ0}
  && \overline{\partial_{z}f(z).u}=\langle u\,,\,
  \partial_{\overline{z}}f(z)\rangle\,,\quad z\in \Z_{0}\,,\, u\in
  \Z_{0}\\
\label{eq.realZ1}
&& \overline{\partial_{z}g(z).u}=\langle
u\,,\,\partial_{\overline{z}}g(z)\rangle\,,
\quad z\in \Z_{1}\,,\, u\in \Z_{1}\,.\\
\label{eq.realh}&& \overline{\partial_{z}h(z).u}=\langle
u\,,\, \partial_{\overline{z}}h(z)\rangle\,,\quad z\in \Z_{1}, u\in \Z_{-1}\,.
\end{eqnarray}
Note that the Poisson bracket
$$
i\left\{h,b\right\}(z)=i\left(\partial_{z}h.\partial_{\overline{z}}b-
\partial_{z}b.\partial_{\overline{z}}h\right)(z)\,,
\quad z\in\Z_{1}
$$
is well defined for $b\in \mathcal{S}_{cyl}(\Z_{1};\rz)$  and our aim
is to write it  as the real scalar product
$$
\langle v(z)\,,\, (\nabla b)(z)\rangle_{\Z_{1},\rz}\,,\quad z\in \Z_{1}\,.
$$
\begin{definition}
  For a cylindrical function on $\Z_{1}$\,, $f\in
  \mathcal{S}_{cyl}(\Z_{1})$\,,
 the gradients $\nabla_{\overline{z}}$
 and $\nabla$ are  defined by
 \begin{eqnarray*}
&&\forall z\in \Z_{1}, u\in \Z_{1},\quad \langle u\,,\,
\nabla_{\overline{z}}f(z)\rangle_{\Z_{1}}=
\langle u\,,\, \partial_{\overline{z}}f(z)\rangle\,,\\
&&
\nabla=2\nabla_{\overline{z}}\,.
\end{eqnarray*}
\end{definition}
\begin{remark}
  \begin{itemize}
  \item  Although it is not necessary, these definitions can be
    justified by introducing a complex conjugation $u\to  \bar{u}$ on
    $\Z_{0}$\,, which remains a conjugation on $\Z_{1}$\,,
 that is an isometric $\cz$-antilinear application such
    that $
\overline{\langle u,v\rangle_{\Z_{0,1}}}=\langle \overline{u}\,,\,
\overline{v}\rangle_{\Z_{0,1}}$\,. When $\Z_{0}=L^{2}(\rz^{d};\cz)$ and
$\Z_{1}=H^{1}(\rz^{d})$ this is the usual pointwise complex conjugation.\\
For real valued functions, set
$$
\nabla_{R}f=\nabla_{\overline{z}}f+\nabla_{z}f\quad\text{and}\quad
\nabla_{I}f=\frac{1}{i}\left(\nabla_{\overline{z}}f-\nabla_{z}f\right)\,.
$$
so that
$$
\nabla_{\overline{z}}f=\frac{1}{2}\left(\nabla_{R}f+i\nabla_{I}f\right)\,.
$$
Similarly, an element $X$  of $\Z_{1}$ can be written
$X=X_{R}+iX_{I}$ with $\overline{X_{R,I}}=X_{R,I}$ or
$X=\begin{pmatrix}
  X_{R}\\X_{I}
\end{pmatrix}
$ and the real scalar product
$$
\langle X\,,\, Y\rangle_{\Z_{1},\rz}
=\Real\langle X\,,\, Y\rangle_{\Z_{1}}=\langle
X_{R}\,,Y_{R}\rangle_{\Z_{1}}+\langle X_{I}\,,\, Y_{I}\rangle_{\Z_{1}}\,.
$$
Then the definition of the gradient of a real cylindrical function $f$
becomes
$$
\nabla f=
\begin{pmatrix}
  \nabla_{R}f\\
\nabla_{I}f
\end{pmatrix}\,.
$$
\item It is important to notice that we do not use the
  $\Z_{1}$-gradient for the real valued function $h(z)$\,, but keep
  the derivative, $\partial_{\overline{z}}h(z)$ modeled on the duality
  bracket $\langle~,~\rangle$\,. With a complex conjugation and since
  $h$ is real valued, it can be decomposed into
  $\partial_{\overline{z}}h=
\frac{1}{2}\left(\partial_{R}h+i\partial_{I}h\right)$ and
$$
-i\partial_{\overline{z}}h=\frac{1}{2}\partial_{I}h -\frac{i}{2}\partial_{R}h\,.
$$
\end{itemize}
\end{remark}
\begin{lem}
\label{le.Poidiv}
  With the above notations and assumptions the equality
$$
\forall z\in\Z_{1}\,,\quad
i\left\{h,b\right\}(z)=2\Real \langle
-i \partial_{\overline{z}}h(z)\,,\,
\nabla_{\overline{z}}b(z)\rangle_{\Z_{1}}
=
\langle v(z)\,,\,\nabla b(z)\rangle_{\Z_{1},\rz}\,,
$$
holds for any $b\in \mathcal{S}_{cyl}(\Z_{1}; \rz)$ with
$v(z)=-i\partial_{\overline{z}}h(z)$\,.
\end{lem}
\proof
It suffices to compute
\begin{eqnarray*}
  i\left\{h\,,\, b\right\}
&=&
i\left[\partial_{z}h.\partial_{\overline{z}}b
  -\partial_{z}b.\partial_{\overline{z}}h\right]\\
&=& i\left[ \overline{\langle \partial_{\overline{z}}b\,,\, \partial_{\overline{z}}h\rangle_{\Z_{0}}}
  -\overline{\langle \partial_{\overline{z}}h\,,\, \partial_{\overline{z}}b\rangle_{\Z_{0}}}\right]
\\
&=&
-2\Imag\langle \partial_{\overline{z}}h\,,\, \partial_{\overline{z}}b\rangle=
2\Real\langle-i\partial_{\overline{z}}h\,,\, \partial_{\overline{z}}{b}\rangle
=2\Real \langle -i\partial_{\overline{z}}h\,,\,
\nabla_{\overline{z}}b\rangle_{\Z_{1}}
\\
&=& \langle-i\partial_{\overline{z}}h\,,\, \nabla b\rangle_{\Z_{1},\rz}\,.
\end{eqnarray*}
\fin

\begin{prop}
\label{pr.comptrmeas}
Let $\Z_{1}\subset \Z_{0}\subset \Z_{-1}$  be a Hilbert triple of separable complex Hilbert spaces.
 Consider a time dependent real sesquilinear form $z\to h(z,t)$ on
 $\Z_{1}$ which is Fr{\'e}chet-$\mathcal{C}^{1}$ and such that
 $\Z_{1}\times \rz\ni (z,t)\to
 \left(\partial_{\overline{z}}h(z,t)\,,\, \partial_{z}h(z,t)\right)\in
 \Z_{1}\times \Z_{-1}^{*}$ is strongly continuous. Assume also that
the time-dependent Hamilton equation
$$
i\partial_{t}z_{t}=\partial_{\bar z}h(z_{t}, \bar z_{t}, t)\quad, \quad z_{t=s}=z
$$
admits a unique continuous solution $z_{t}=\Phi(t,s)z$ for all
$t,s\in\rz$ and all $z\in\Z_{1}$\,, with $\Phi(t,s):\Z_{1}\to\Z_{1}$ Borel. \\
Consider a time dependent measure $\mu(t)\in
\mathrm{Prob}_{2}(\Z_{1})$ which satisfies
\begin{itemize}
\item $t\to \mu_{t}\in \mathrm{Prob}_{2}(\Z_{1})$ is
  $W_{2}$-continuous.
\item For all $T>0$\,, $|\partial_{\bar z}h(t)|_{L^{2}(\Z,\mu_{t})}\in
  L^{1}([-T,T])$\,.
\item The time-dependent probability measure $\mu_{t}$ is a weak
  solution to
$$
\partial_{t}\mu+i\left\{h(t)\,,\,\mu\right\}=0\,,
$$
namely for all $\varphi\in \mathcal{C}^{\infty}_{0,cyl}(\Z_{1}\times \rz;\rz)$\,,
$$
\int_{\rz}\int_{\Z_{1}}
\left(\partial_{t}\varphi(z,t)+i\left\{h,\varphi\right\}(z,t)\right)~d\mu_{t}(z)dt=0\,.
$$
\end{itemize}
Then the measure $\mu_{t}$ satisfies
$$
\forall t,s\in\rz,\quad \mu_{t}=\Phi(t,s)_{*}\mu_{s}
$$
and it is unique when $\mu_{0}$ is fixed.
\end{prop}
\proof
We  apply Proposition~\ref{pr.flowmeas} while $E=\Z_{1}$ is
endowed with its euclidean structure $\langle z_{1}\,,\,
z_{2}\rangle_{\Z_{1},\rz} $ $=\Real \langle z_{1}\,,\, z_{2}
\rangle_{\Z_{1}}$\,.
Lemma~\ref{le.Poidiv} says that the weak Liouville equation is
$$
\forall \varphi\in \mathcal{C}^{\infty}_{0,cyl}(\Z_{1}\times \rz)\,,\quad
\int_{\rz}\int_{\Z_{1}}\left(\partial_{t}\varphi(z,t) + \langle
  v\,,\, \nabla \varphi\rangle(t,z)\right)~d\mu_{t}(z)dt=0\,,
$$
with $v(z,t)=-i\partial_{\overline{z}}h (z,t)$\,. The measure
$\mu_{t}$ is a weak solution to
$$
\partial_{t}\mu + \nabla^{T}(v\mu)=0
$$
where $\nabla$ and $\nabla^{T}$ are defined according to the real
structure on $\Z_{1}$\,. Our hypothesis on $\mu$ and $h$ cover all the
assumptions
of Proposition~\ref{pr.flowmeas}
\fin

\section{Weak $L^p$ conditions for the potential $V$}
\label{se.Lp}
Let $L_{\rz^d}$ be the Lebesgue measure on $\rz^d$\,. Let $0<p<\infty$\,, a Lebesgue measurable function
$f:\rz^d\to \cz$ is said to  belong to weak-$L^p(\rz^{d})$, or shortly in
$L^{p,\infty}(\rz^d)$\,, if there exists a constant $c>0$ such that for all $t>0$
$$
L_{\rz^d}\{x:|f(x)|>t\} \leq c^p/t^p\,.
$$
Two functions in $L^{p,\infty}(\rz^d)$ are equal  if they are  equal $L_{\rz^d}$-almost everywhere.
The quantity
\begin{eqnarray*}
\|f\|_{p,\infty}&=&\inf\{c: L_{\rz^d}\{x:|f(x)|>t\} \leq c^p/t^p, \forall t>0\}\\
&=&\sup_{t>0}\{t L_{\rz^d}\{x:|f(x)|>t\}^{1/p}\}
\end{eqnarray*}
defines a complete quasi-norm on $L^{p,\infty}(\rz^d)$ with $\|f\|_{p,\infty}\leq \|f\|_p$\,.

\noindent By combining Hunt and Marcinkiewicz interpolation theorems according to \cite{Gr,RS,BL}),
the Young and H{\"o}lder inequalities  can be extended to weak $L^p$ spaces.
\begin{prop}[generalized Young's inequality]
Let $1<p,q,r<\infty$ such that $\frac{1}{p}+\frac{1}{q}=1+\frac{1}{r}$\,. There exists a constant $c_{p,q}>0$
such that for all $f\in L^p(\rz^d)$ and  $g\in L^{q,\infty}(\rz^d)$
$$
\|f*g\|_r\leq c_{p,q} \, \|f\|_p \,\|g\|_{q,\infty}\,.
$$
\end{prop}
\begin{prop}[generalized H{\"o}lder inequality]
Let $1<p,q,r<\infty$ satisfying $\frac{1}{p}+\frac{1}{q}=\frac{1}{r}$\,. There exists a constant
$c_{p,q}$ such that for all $f\in L^{p,\infty}(\rz^d)$
and $g\in L^q(\rz^d)$
$$
\|f . g\|_r\leq c_{p,q} \, \|f\|_{p,\infty} \,\|g\|_q\,.
$$
\end{prop}
\begin{prop}[Hardy inequality]
\label{hardy}
Suppose that $d\geq 3$ and $V\in L^{d,\infty}(\rz^d)$\,. There exists a constant $c>0$ such that
for all $u\in H^1(\rz^d)$
$$
\|V u\|_2\leq c \|u\|_{H^1(\rz^d)}\,.
$$
\end{prop}
\proof
For $u\in L^2(\rz^d)$\,, we can write $(1-\Delta)^{-1/2} u(x)=G*u(x)$ with $G$ the inverse Fourier
transform of $(1+|x|^2)^{-1/2}$\,. It is not difficult to prove that
$G\in L^{\frac{d}{d-1},\infty}$ (see \cite[Exercice 50]{RS}). Hence, we conclude that
\begin{eqnarray*}
\ds \|V\,( 1-\Delta)^{-1/2} u\|_2=\| V \, G*u\|_2& \overset{\textrm{ H{\"o}lder}}{\leq} & C_1\,
\|V\|_{d,\infty} \, \|G*u\|_{\frac{2d}{d-2}} \\ \ds
&\overset{{\rm Young}}{\leq}& C_2\, \|V\|_{d,\infty} \,\|G\|_{\frac{d}{d-1},\infty} \,\|u\|_2
\end{eqnarray*}
\fin\\
The above proposition provides a class of potentials which are bounded multiplication operators
from $H^1(\rz^d)$ into $L^2(\rz)$ when the dimension $d\geq 3$\,. For
lower dimension, the  Sobolev embeddings give at once:
\begin{description}
  \item[$\bullet$ ] if $d=1$\,, $V\in L^2(\rz)+L^\infty(\rz)$ then $V\in \L(H^1(\rz), L^2(\rz))$\,.
  \item[$\bullet$ ] if $d=2$\,, $V\in L^p(\rz^2)+L^\infty(\rz^2)$ for $p>2$\,, then $V\in \L(H^1(\rz^2), L^2(\rz^2))$\,.
\end{description}

\smallskip

We denote by $L^p(\rz^d)+L_0^\infty(\rz^d)$ the space of Lebesgue measurable functions
$f$ such that there exists $(f_n)_{n\in\nz}\in L^p(\rz^d)^\nz$ satisfying $\lim_{n\to\infty}
\|f-f_n\|_\infty=0$\,.
\begin{lem}
\label{le.ld}
For $0<p<q$\,,
$$
L^{q,\infty}(\rz^d)\subset L^p(\rz^d)+L_0^{\infty}(\rz^d)
$$
\end{lem}
\proof
For $\epsilon>0$\,,  decompose each $f\in L^{q,\infty}(\rz^d)$ into a sum $f=f_\epsilon+f^\epsilon$ such that
$f_\epsilon=f 1_{|f|>\epsilon}$ and $f^\epsilon=f 1_{|f|\leq \epsilon}$\,.
 Observe that for any $\epsilon>0$
\begin{eqnarray*}
\|f_\epsilon\|_p^p&=& p \int_0^\infty t^{p-1} L_{\rz^d}\{x: |f_\epsilon(x)|>t\} \,dt \\
&=& p \int_\epsilon^\infty t^{p-1} L_{\rz^d}\{x: |f(x)|>t\}
\,dt+\epsilon^p L_{\rz^d}
\{x: |f(x)|>\epsilon\}\\
&\leq& c \int_\epsilon^\infty \frac{t^{p-1}}{t^q} \,dt+\epsilon^p L_{\rz^d}\{x: |f(x)|>\epsilon\}<\infty
\end{eqnarray*}
Moreover, when $\epsilon\to 0$
$$
\|f^\epsilon\|_\infty=\|f 1_{|f|\leq \epsilon}\|_\infty\leq \epsilon\to 0\,.
$$
Therefore, each $f\in L^{q,\infty}(\rz^d)$ belongs to the space $L^p(\rz^d)+L_0^{\infty}(\rz^d)$\,.
\fin
\begin{prop}
\label{prop-comp}
For any $V\in L^{2}(\rz^d)+L_0^\infty(\rz^d)$ such that $V (1-\Delta)^{-1/2}\in \L(L^2(\rz^d))$ the operator $(1-\Delta)^{-1/2} V (1-\Delta)^{-1/2}$ is compact.
\end{prop}
\proof
Let
$g(\xi)=(1+|\xi|^2)^{-1/2}$ and $g_m(\xi)=1_{[0,m]}(|\xi|) g(\xi)$\,. The following norm convergence holds
$$
\lim_{m\to\infty} g_m(D)\, V\, g_m(D)= g(D) \, V\, g(D)\,,
$$
using the fact that $\lim_{m\to\infty}\|g_m(D)-g(D)\|_{\L(L^2(\rz^d))}=0$ and $\|V \, g(D)\|_{\L(L^2(\rz^d))}<\infty$\,.  \\
By Lemma \ref{le.ld}, there exist $V_n\in L^2(\rz^d)$ such that $\lim_{n\to\infty}\|V-V_n\|_\infty=0$\,. We
observe now that the Hilbert-Schmidt norm of $g_m(D) V_n(x) g_m(D)$ is
$$
\|g_m(D)\, V_n \, g_m(D)\|_{\L^2(L^2(\rz^d))}\leq\|V_n\|_2\; \|g_m\|_2^2<\infty.
$$
Therefore by norm convergence, the operator  $g(D) \,V \,g(D)$ is compact.\fin

\begin{cor}
The potential $V$ satisfies the assumptions (A\ref{eq.hypbd})-(A\ref{eq.hypcomp}) in the following cases:
\begin{description}
  \item[$\bullet$ ] if $d=1$ and  $V\in L^2(\rz)+L_0^\infty(\rz)$\,,
  \item[$\bullet$ ] if $d=2$ and  $V\in L^p(\rz^2)+L_0^\infty(\rz^2)$
    with $p>2$\,,
  \item[$\bullet$ ] if $d\geq 3$ and $V\in L^{d,\infty}(\rz^d)$\,.
\end{description}
\end{cor}
\proof
Combine Proposition~\ref{hardy}, Lemma \ref{le.ld} and
Proposition~\ref{prop-comp}, with in dimension $d=2$ the
observation
$$
L^p(\rz^2)+L^\infty_0(\rz^2)\subset L^2(\rz^2)+L^\infty_0(\rz^2)\quad\text{for}~p>2\,.
$$
\fin

\bigskip
\noindent
In particular, in dimension $d=3$ the Coulomb potential $V(x)=\pm\frac{1}{|x|}$
satisfies the assumptions (A\ref{eq.hypsym}), (A\ref{eq.hypbd}) and
(A\ref{eq.hypcomp})  because  $\frac{1}{|x|}\in L^{3,\infty}(\rz^3)$\,.


\bigskip
\bibliographystyle{amsalpha}

\end{document}

%% file: fig4tex.tex
%
%
%
%
\ifx\figforTeXisloaded\relax \else\global\let\figforTeXisloaded=\relax\fi
\message{version 1.8.4}
\catcode`\@=11
\ifx\ctr@ln@m\undefined\else%
    \immediate\write16{*** Fig4TeX WARNING : \string\ctr@ln@m\space already defined.}\fi
\def\ctr@ln@m#1{\ifx#1\undefined\else%
    \immediate\write16{*** Fig4TeX WARNING : \string#1 already defined.}\fi}
\ctr@ln@m\ctr@ld@f
\def\ctr@ld@f#1#2{\ctr@ln@m#2#1#2}
\ctr@ld@f\def\ctr@ln@w#1#2{\ctr@ln@m#2\csname#1\endcsname#2}
{\catcode`\/=0 \catcode`/\=12 /ctr@ld@f/gdef/BS@{\}}
\ctr@ld@f\def\ctr@lcsn@m#1{\expandafter\ifx\csname#1\endcsname\relax\else%
    \immediate\write16{*** Fig4TeX WARNING : \BS@\expandafter\string#1\space already defined.}\fi}
\ctr@ld@f\edef\colonc@tcode{\the\catcode`\:}
\ctr@ld@f\edef\semicolonc@tcode{\the\catcode`\;}
\ctr@ld@f\def\t@stc@tcodech@nge{{\let\c@tcodech@nged=\z@%
    \ifnum\colonc@tcode=\the\catcode`\:\else\let\c@tcodech@nged=\@ne\fi%
    \ifnum\semicolonc@tcode=\the\catcode`\;\else\let\c@tcodech@nged=\@ne\fi%
    \ifx\c@tcodech@nged\@ne%
    \immediate\write16{}
    \immediate\write16{!!!=============================================================!!!}
    \immediate\write16{ Fig4TeX WARNING :}
    \immediate\write16{ The category code of some characters has been changed, which will}
    \immediate\write16{ result in an error (message "Runaway argument?").}
    \immediate\write16{ This probably comes from another package that changed the category}
    \immediate\write16{ code after Fig4TeX was loaded. If that proves to be exact, the}
    \immediate\write16{ solution is to exchange the loading commands on top of your file}
    \immediate\write16{ so that Fig4TeX is loaded last. For example, in LaTeX, we should}
    \immediate\write16{ say :}
    \immediate\write16{\BS@ usepackage[french]{babel}}
    \immediate\write16{\BS@ usepackage{fig4tex}}
    \immediate\write16{!!!=============================================================!!!}
    \immediate\write16{}
    \fi}}
\ctr@ld@f\def\FigforTeX{F\kern-.05em i\kern-.05em g\kern-.1em\raise-.14em\hbox{4}\kern-.19em\TeX}
\ctr@ln@w{newdimen}\epsil@n\epsil@n=0.00005pt
\ctr@ln@w{newdimen}\Cepsil@n\Cepsil@n=0.005pt
\ctr@ln@w{newdimen}\dcq@\dcq@=254pt
\ctr@ln@w{newdimen}\PI@\PI@=3.141592pt
\ctr@ln@w{newdimen}\DemiPI@deg\DemiPI@deg=90pt
\ctr@ln@w{newdimen}\PI@deg\PI@deg=180pt
\ctr@ln@w{newdimen}\DePI@deg\DePI@deg=360pt
\ctr@ld@f\chardef\t@n=10
\ctr@ld@f\chardef\c@nt=100
\ctr@ld@f\chardef\@lxxiv=74
\ctr@ld@f\chardef\@xci=91
\ctr@ld@f\mathchardef\@nMnCQn=9949
\ctr@ld@f\chardef\@vi=6
\ctr@ld@f\chardef\@xxx=30
\ctr@ld@f\chardef\@lvi=56
\ctr@ld@f\chardef\@@lxxi=71
\ctr@ld@f\chardef\@lxxxv=85
\ctr@ld@f\mathchardef\@@mmmmlxviii=4068
\ctr@ld@f\mathchardef\@ccclx=360
\ctr@ld@f\mathchardef\@dccxx=720
\ctr@ln@w{newcount}\p@rtent \ctr@ln@w{newcount}\f@ctech \ctr@ln@w{newcount}\result@tent
\ctr@ln@w{newdimen}\v@lmin \ctr@ln@w{newdimen}\v@lmax \ctr@ln@w{newdimen}\v@leur
\ctr@ln@w{newdimen}\result@t\ctr@ln@w{newdimen}\result@@t
\ctr@ln@w{newdimen}\mili@u \ctr@ln@w{newdimen}\c@rre \ctr@ln@w{newdimen}\delt@
\ctr@ld@f\def\degT@rd{0.017453 }  
\ctr@ld@f\def\rdT@deg{57.295779 } 
\ctr@ln@m\v@leurseule
{\catcode`p=12 \catcode`t=12 \gdef\v@leurseule#1pt{#1}}
\ctr@ld@f\def\repdecn@mb#1{\expandafter\v@leurseule\the#1\space}
\ctr@ld@f\def\arct@n#1(#2,#3){{\v@lmin=#2\v@lmax=#3%
    \maxim@m{\mili@u}{-\v@lmin}{\v@lmin}\maxim@m{\c@rre}{-\v@lmax}{\v@lmax}%
    \delt@=\mili@u\m@ech\mili@u%
    \ifdim\c@rre>\@nMnCQn\mili@u\divide\v@lmax\tw@\c@lATAN\v@leur(\z@,\v@lmax)
    \else%
    \maxim@m{\mili@u}{-\v@lmin}{\v@lmin}\maxim@m{\c@rre}{-\v@lmax}{\v@lmax}%
    \m@ech\c@rre%
    \ifdim\mili@u>\@nMnCQn\c@rre\divide\v@lmin\tw@
    \maxim@m{\mili@u}{-\v@lmin}{\v@lmin}\c@lATAN\v@leur(\mili@u,\z@)%
    \else\c@lATAN\v@leur(\delt@,\v@lmax)\fi\fi%
    \ifdim\v@lmin<\z@\v@leur=-\v@leur\ifdim\v@lmax<\z@\advance\v@leur-\PI@%
    \else\advance\v@leur\PI@\fi\fi%
    \global\result@t=\v@leur}#1=\result@t}
\ctr@ld@f\def\m@ech#1{\ifdim#1>1.646pt\divide\mili@u\t@n\divide\c@rre\t@n\m@ech#1\fi}
\ctr@ld@f\def\c@lATAN#1(#2,#3){{\v@lmin=#2\v@lmax=#3\v@leur=\z@\delt@=\tw@ pt%
    \un@iter{0.785398}{\v@lmax<}%
    \un@iter{0.463648}{\v@lmax<}%
    \un@iter{0.244979}{\v@lmax<}%
    \un@iter{0.124355}{\v@lmax<}%
    \un@iter{0.062419}{\v@lmax<}%
    \un@iter{0.031240}{\v@lmax<}%
    \un@iter{0.015624}{\v@lmax<}%
    \un@iter{0.007812}{\v@lmax<}%
    \un@iter{0.003906}{\v@lmax<}%
    \un@iter{0.001953}{\v@lmax<}%
    \un@iter{0.000976}{\v@lmax<}%
    \un@iter{0.000488}{\v@lmax<}%
    \un@iter{0.000244}{\v@lmax<}%
    \un@iter{0.000122}{\v@lmax<}%
    \un@iter{0.000061}{\v@lmax<}%
    \un@iter{0.000030}{\v@lmax<}%
    \un@iter{0.000015}{\v@lmax<}%
    \global\result@t=\v@leur}#1=\result@t}
\ctr@ld@f\def\un@iter#1#2{%
    \divide\delt@\tw@\edef\dpmn@{\repdecn@mb{\delt@}}%
    \mili@u=\v@lmin%
    \ifdim#2\z@%
      \advance\v@lmin-\dpmn@\v@lmax\advance\v@lmax\dpmn@\mili@u%
      \advance\v@leur-#1pt%
    \else%
      \advance\v@lmin\dpmn@\v@lmax\advance\v@lmax-\dpmn@\mili@u%
      \advance\v@leur#1pt%
    \fi}
\ctr@ld@f\def\c@ssin#1#2#3{\expandafter\ifx\csname COS@\number#3\endcsname\relax\c@lCS{#3pt}%
    \expandafter\xdef\csname COS@\number#3\endcsname{\repdecn@mb\result@t}%
    \expandafter\xdef\csname SIN@\number#3\endcsname{\repdecn@mb\result@@t}\fi%
    \edef#1{\csname COS@\number#3\endcsname}\edef#2{\csname SIN@\number#3\endcsname}}
\ctr@ld@f\def\c@lCS#1{{\mili@u=#1\p@rtent=\@ne%
    \relax\ifdim\mili@u<\z@\red@ng<-\else\red@ng>+\fi\f@ctech=\p@rtent%
    \relax\ifdim\mili@u<\z@\mili@u=-\mili@u\f@ctech=-\f@ctech\fi\c@@lCS}}
\ctr@ld@f\def\c@@lCS{\v@lmin=\mili@u\c@rre=-\mili@u\advance\c@rre\DemiPI@deg\v@lmax=\c@rre%
    \mili@u\@@lxxi\mili@u\divide\mili@u\@@mmmmlxviii%
    \edef\v@larg{\repdecn@mb{\mili@u}}\mili@u=-\v@larg\mili@u%
    \edef\v@lmxde{\repdecn@mb{\mili@u}}%
    \c@rre\@@lxxi\c@rre\divide\c@rre\@@mmmmlxviii%
    \edef\v@largC{\repdecn@mb{\c@rre}}\c@rre=-\v@largC\c@rre%
    \edef\v@lmxdeC{\repdecn@mb{\c@rre}}%
    \fctc@s\mili@u\v@lmin\global\result@t\p@rtent\v@leur%
    \let\t@mp=\v@larg\let\v@larg=\v@largC\let\v@largC=\t@mp%
    \let\t@mp=\v@lmxde\let\v@lmxde=\v@lmxdeC\let\v@lmxdeC=\t@mp%
    \fctc@s\c@rre\v@lmax\global\result@@t\f@ctech\v@leur}
\ctr@ld@f\def\fctc@s#1#2{\v@leur=#1\relax\ifdim#2<\@lxxxv\p@\cosser@h\else\sinser@t\fi}
\ctr@ld@f\def\cosser@h{\advance\v@leur\@lvi\p@\divide\v@leur\@lvi%
    \v@leur=\v@lmxde\v@leur\advance\v@leur\@xxx\p@%
    \v@leur=\v@lmxde\v@leur\advance\v@leur\@ccclx\p@%
    \v@leur=\v@lmxde\v@leur\advance\v@leur\@dccxx\p@\divide\v@leur\@dccxx}
\ctr@ld@f\def\sinser@t{\v@leur=\v@lmxdeC\p@\advance\v@leur\@vi\p@%
    \v@leur=\v@largC\v@leur\divide\v@leur\@vi}
\ctr@ld@f\def\red@ng#1#2{\relax\ifdim\mili@u#1#2\DemiPI@deg\advance\mili@u#2-\PI@deg%
    \p@rtent=-\p@rtent\red@ng#1#2\fi}
\ctr@ld@f\def\pr@c@lCS#1#2#3{\ctr@lcsn@m{COS@\number#3 }%
    \expandafter\xdef\csname COS@\number#3\endcsname{#1}%
    \expandafter\xdef\csname SIN@\number#3\endcsname{#2}}
\pr@c@lCS{1}{0}{0}
\pr@c@lCS{0.7071}{0.7071}{45}\pr@c@lCS{0.7071}{-0.7071}{-45}
\pr@c@lCS{0}{1}{90}          \pr@c@lCS{0}{-1}{-90}
\pr@c@lCS{-1}{0}{180}        \pr@c@lCS{-1}{0}{-180}
\pr@c@lCS{0}{-1}{270}        \pr@c@lCS{0}{1}{-270}
\ctr@ld@f\def\invers@#1#2{{\v@leur=#2\maxim@m{\v@lmax}{-\v@leur}{\v@leur}%
    \f@ctech=\@ne\m@inv@rs%
    \multiply\v@leur\f@ctech\edef\v@lv@leur{\repdecn@mb{\v@leur}}%
    \p@rtentiere{\p@rtent}{\v@leur}\v@lmin=\p@\divide\v@lmin\p@rtent%
    \inv@rs@\multiply\v@lmax\f@ctech\global\result@t=\v@lmax}#1=\result@t}
\ctr@ld@f\def\m@inv@rs{\ifdim\v@lmax<\p@\multiply\v@lmax\t@n\multiply\f@ctech\t@n\m@inv@rs\fi}
\ctr@ld@f\def\inv@rs@{\v@lmax=-\v@lmin\v@lmax=\v@lv@leur\v@lmax%
    \advance\v@lmax\tw@ pt\v@lmax=\repdecn@mb{\v@lmin}\v@lmax%
    \delt@=\v@lmax\advance\delt@-\v@lmin\ifdim\delt@<\z@\delt@=-\delt@\fi%
    \ifdim\delt@>\epsil@n\v@lmin=\v@lmax\inv@rs@\fi}
\ctr@ld@f\def\minim@m#1#2#3{\relax\ifdim#2<#3#1=#2\else#1=#3\fi}
\ctr@ld@f\def\maxim@m#1#2#3{\relax\ifdim#2>#3#1=#2\else#1=#3\fi}
\ctr@ld@f\def\p@rtentiere#1#2{#1=#2\divide#1by65536 }
\ctr@ld@f\def\r@undint#1#2{{\v@leur=#2\divide\v@leur\t@n\p@rtentiere{\p@rtent}{\v@leur}%
    \v@leur=\p@rtent pt\global\result@t=\t@n\v@leur}#1=\result@t}
\ctr@ld@f\def\sqrt@#1#2{{\v@leur=#2%
    \minim@m{\v@lmin}{\p@}{\v@leur}\maxim@m{\v@lmax}{\p@}{\v@leur}%
    \f@ctech=\@ne\m@sqrt@\sqrt@@%
    \mili@u=\v@lmin\advance\mili@u\v@lmax\divide\mili@u\tw@\multiply\mili@u\f@ctech%
    \global\result@t=\mili@u}#1=\result@t}
\ctr@ld@f\def\m@sqrt@{\ifdim\v@leur>\dcq@\divide\v@leur\c@nt\v@lmax=\v@leur%
    \multiply\f@ctech\t@n\m@sqrt@\fi}
\ctr@ld@f\def\sqrt@@{\mili@u=\v@lmin\advance\mili@u\v@lmax\divide\mili@u\tw@%
    \c@rre=\repdecn@mb{\mili@u}\mili@u%
    \ifdim\c@rre<\v@leur\v@lmin=\mili@u\else\v@lmax=\mili@u\fi%
    \delt@=\v@lmax\advance\delt@-\v@lmin\ifdim\delt@>\epsil@n\sqrt@@\fi}
\ctr@ld@f\def\extrairelepremi@r#1\de#2{\expandafter\lepremi@r#2@#1#2}
\ctr@ld@f\def\lepremi@r#1,#2@#3#4{\def#3{#1}\def#4{#2}\ignorespaces}
\ctr@ld@f\def\@cfor#1:=#2\do#3{%
  \edef\@fortemp{#2}%
  \ifx\@fortemp\empty\else\@cforloop#2,\@nil,\@nil\@@#1{#3}\fi}
\ctr@ln@m\@nextwhile
\ctr@ld@f\def\@cforloop#1,#2\@@#3#4{%
  \def#3{#1}%
  \ifx#3\Fig@nnil\let\@nextwhile=\Fig@fornoop\else#4\relax\let\@nextwhile=\@cforloop\fi%
  \@nextwhile#2\@@#3{#4}}

\ctr@ld@f\def\@ecfor#1:=#2\do#3{%
  \def\@@cfor{\@cfor#1:=}%
  \edef\@@@cfor{#2}%
  \expandafter\@@cfor\@@@cfor\do{#3}}
\ctr@ld@f\def\Fig@nnil{\@nil}
\ctr@ld@f\def\Fig@fornoop#1\@@#2#3{}
\ctr@ln@m\list@@rg
\ctr@ld@f\def\trtlis@rg#1#2{\def\list@@rg{#1}%
    \@ecfor\p@rv@l:=\list@@rg\do{\expandafter#2\p@rv@l|}}
\ctr@ln@w{newbox}\b@xvisu
\ctr@ln@w{newtoks}\let@xte
\ctr@ln@w{newif}\ifitis@K
\ctr@ln@w{newcount}\s@mme
\ctr@ln@w{newcount}\l@mbd@un \ctr@ln@w{newcount}\l@mbd@de
\ctr@ln@w{newcount}\superc@ntr@l\superc@ntr@l=\@ne        
\ctr@ln@w{newcount}\typec@ntr@l\typec@ntr@l=\superc@ntr@l 
\ctr@ln@w{newdimen}\v@lX  \ctr@ln@w{newdimen}\v@lY  \ctr@ln@w{newdimen}\v@lZ
\ctr@ln@w{newdimen}\v@lXa \ctr@ln@w{newdimen}\v@lYa \ctr@ln@w{newdimen}\v@lZa
\ctr@ln@w{newdimen}\unit@\unit@=\p@ 
\ctr@ld@f\def\unit@util{pt}
\ctr@ld@f\def\ptT@ptps{0.996264 }
\ctr@ld@f\def\ptpsT@pt{1.00375 }
\ctr@ld@f\def\ptT@unit@{1} 
\ctr@ld@f\def\setunit@#1{\def\unit@util{#1}\setunit@@#1:\invers@{\result@t}{\unit@}%
    \edef\ptT@unit@{\repdecn@mb\result@t}}
\ctr@ld@f\def\setunit@@#1#2:{\ifcat#1a\unit@=\@ne#1#2\else\unit@=#1#2\fi}
\ctr@ld@f\def\d@fm@cdim#1#2{{\v@leur=#2\v@leur=\ptT@unit@\v@leur\xdef#1{\repdecn@mb\v@leur}}}
\ctr@ln@w{newif}\ifBdingB@x\BdingB@xtrue
\ctr@ln@w{newdimen}\c@@rdXmin \ctr@ln@w{newdimen}\c@@rdYmin  
\ctr@ln@w{newdimen}\c@@rdXmax \ctr@ln@w{newdimen}\c@@rdYmax
\ctr@ld@f\def\b@undb@x#1#2{\ifBdingB@x%
    \relax\ifdim#1<\c@@rdXmin\global\c@@rdXmin=#1\fi%
    \relax\ifdim#2<\c@@rdYmin\global\c@@rdYmin=#2\fi%
    \relax\ifdim#1>\c@@rdXmax\global\c@@rdXmax=#1\fi%
    \relax\ifdim#2>\c@@rdYmax\global\c@@rdYmax=#2\fi\fi}
\ctr@ld@f\def\b@undb@xP#1{{\Figg@tXY{#1}\b@undb@x{\v@lX}{\v@lY}}}
\ctr@ld@f\def\ellBB@x#1;#2,#3(#4,#5,#6){{\s@uvc@ntr@l\et@tellBB@x%
    \setc@ntr@l{2}\figptell-2::#1;#2,#3(#4,#6)\b@undb@xP{-2}%
    \figptell-2::#1;#2,#3(#5,#6)\b@undb@xP{-2}%
    \c@ssin{\C@}{\S@}{#6}\v@lmin=\C@ pt\v@lmax=\S@ pt%
    \mili@u=#3\v@lmin\delt@=#2\v@lmax\arct@n\v@leur(\delt@,\mili@u)%
    \mili@u=-#3\v@lmax\delt@=#2\v@lmin\arct@n\c@rre(\delt@,\mili@u)%
    \v@leur=\rdT@deg\v@leur\advance\v@leur-\DePI@deg%
    \c@rre=\rdT@deg\c@rre\advance\c@rre-\DePI@deg%
    \v@lmin=#4pt\v@lmax=#5pt%
    \loop\ifdim\v@leur<\v@lmax\ifdim\v@leur>\v@lmin%
    \edef\@ngle{\repdecn@mb\v@leur}\figptell-2::#1;#2,#3(\@ngle,#6)%
    \b@undb@xP{-2}\fi\advance\v@leur\PI@deg\repeat%
    \loop\ifdim\c@rre<\v@lmax\ifdim\c@rre>\v@lmin%
    \edef\@ngle{\repdecn@mb\c@rre}\figptell-2::#1;#2,#3(\@ngle,#6)%
    \b@undb@xP{-2}\fi\advance\c@rre\PI@deg\repeat%
    \resetc@ntr@l\et@tellBB@x}\ignorespaces}
\ctr@ld@f\def\initb@undb@x{\c@@rdXmin=\maxdimen\c@@rdYmin=\maxdimen%
    \c@@rdXmax=-\maxdimen\c@@rdYmax=-\maxdimen}
\ctr@ld@f\def\c@ntr@lnum#1{%
    \relax\ifnum\typec@ntr@l=\@ne%
    \ifnum#1<\z@%
    \immediate\write16{*** Forbidden point number (#1). Abort.}\end\fi\fi%
    \set@bjc@de{#1}}
\ctr@ln@m\objc@de
\ctr@ld@f\def\set@bjc@de#1{\edef\objc@de{@BJ\ifnum#1<\z@ M\romannumeral-#1\else\romannumeral#1\fi}}
\s@mme=\m@ne\loop\ifnum\s@mme>-19
  \set@bjc@de{\s@mme}\ctr@lcsn@m\objc@de\ctr@lcsn@m{\objc@de T}
\advance\s@mme\m@ne\repeat
\s@mme=\@ne\loop\ifnum\s@mme<6
  \set@bjc@de{\s@mme}\ctr@lcsn@m\objc@de\ctr@lcsn@m{\objc@de T}
\advance\s@mme\@ne\repeat
\ctr@ld@f\def\setc@ntr@l#1{\ifnum\superc@ntr@l>#1\typec@ntr@l=\superc@ntr@l%
    \else\typec@ntr@l=#1\fi}
\ctr@ld@f\def\resetc@ntr@l#1{\global\superc@ntr@l=#1\setc@ntr@l{#1}}
\ctr@ld@f\def\s@uvc@ntr@l#1{\edef#1{\the\superc@ntr@l}}
\ctr@ln@m\c@lproscal
\ctr@ld@f\def\c@lproscalDD#1[#2,#3]{{\Figg@tXY{#2}%
    \edef\Xu@{\repdecn@mb{\v@lX}}\edef\Yu@{\repdecn@mb{\v@lY}}\Figg@tXY{#3}%
    \global\result@t=\Xu@\v@lX\global\advance\result@t\Yu@\v@lY}#1=\result@t}
\ctr@ld@f\def\c@lproscalTD#1[#2,#3]{{\Figg@tXY{#2}\edef\Xu@{\repdecn@mb{\v@lX}}%
    \edef\Yu@{\repdecn@mb{\v@lY}}\edef\Zu@{\repdecn@mb{\v@lZ}}%
    \Figg@tXY{#3}\global\result@t=\Xu@\v@lX\global\advance\result@t\Yu@\v@lY%
    \global\advance\result@t\Zu@\v@lZ}#1=\result@t}
\ctr@ld@f\def\c@lprovec#1{%
    \det@rmC\v@lZa(\v@lX,\v@lY,\v@lmin,\v@lmax)%
    \det@rmC\v@lXa(\v@lY,\v@lZ,\v@lmax,\v@leur)%
    \det@rmC\v@lYa(\v@lZ,\v@lX,\v@leur,\v@lmin)%
    \Figv@ctCreg#1(\v@lXa,\v@lYa,\v@lZa)}
\ctr@ld@f\def\det@rm#1[#2,#3]{{\Figg@tXY{#2}\Figg@tXYa{#3}%
    \delt@=\repdecn@mb{\v@lX}\v@lYa\advance\delt@-\repdecn@mb{\v@lY}\v@lXa%
    \global\result@t=\delt@}#1=\result@t}
\ctr@ld@f\def\det@rmC#1(#2,#3,#4,#5){{\global\result@t=\repdecn@mb{#2}#5%
    \global\advance\result@t-\repdecn@mb{#3}#4}#1=\result@t}
\ctr@ld@f\def\getredf@ctDD#1(#2,#3){{\maxim@m{\v@lXa}{-#2}{#2}\maxim@m{\v@lYa}{-#3}{#3}%
    \maxim@m{\v@lXa}{\v@lXa}{\v@lYa}
    \ifdim\v@lXa>\@xci pt\divide\v@lXa\@xci%
    \p@rtentiere{\p@rtent}{\v@lXa}\advance\p@rtent\@ne\else\p@rtent=\@ne\fi%
    \global\result@tent=\p@rtent}#1=\result@tent\ignorespaces}
\ctr@ld@f\def\getredf@ctTD#1(#2,#3,#4){{\maxim@m{\v@lXa}{-#2}{#2}\maxim@m{\v@lYa}{-#3}{#3}%
    \maxim@m{\v@lZa}{-#4}{#4}\maxim@m{\v@lXa}{\v@lXa}{\v@lYa}%
    \maxim@m{\v@lXa}{\v@lXa}{\v@lZa}
    \ifdim\v@lXa>\@lxxiv pt\divide\v@lXa\@lxxiv%
    \p@rtentiere{\p@rtent}{\v@lXa}\advance\p@rtent\@ne\else\p@rtent=\@ne\fi%
    \global\result@tent=\p@rtent}#1=\result@tent\ignorespaces}
\ctr@ld@f\def\FigptintercircB@zDD#1:#2:#3,#4[#5,#6,#7,#8]{{\s@uvc@ntr@l\et@tfigptintercircB@zDD%
    \setc@ntr@l{2}\figvectPDD-1[#5,#8]\Figg@tXY{-1}\getredf@ctDD\f@ctech(\v@lX,\v@lY)%
    \mili@u=#4\unit@\divide\mili@u\f@ctech\c@rre=\repdecn@mb{\mili@u}\mili@u%
    \figptBezierDD-5::#3[#5,#6,#7,#8]%
    \v@lmin=#3\p@\v@lmax=\v@lmin\advance\v@lmax0.1\p@%
    \loop\edef\T@{\repdecn@mb{\v@lmax}}\figptBezierDD-2::\T@[#5,#6,#7,#8]%
    \figvectPDD-1[-5,-2]\n@rmeucCDD{\delt@}{-1}\ifdim\delt@<\c@rre\v@lmin=\v@lmax%
    \advance\v@lmax0.1\p@\repeat%
    \loop\mili@u=\v@lmin\advance\mili@u\v@lmax%
    \divide\mili@u\tw@\edef\T@{\repdecn@mb{\mili@u}}\figptBezierDD-2::\T@[#5,#6,#7,#8]%
    \figvectPDD-1[-5,-2]\n@rmeucCDD{\delt@}{-1}\ifdim\delt@>\c@rre\v@lmax=\mili@u%
    \else\v@lmin=\mili@u\fi\v@leur=\v@lmax\advance\v@leur-\v@lmin%
    \ifdim\v@leur>\epsil@n\repeat\figptcopyDD#1:#2/-2/%
    \resetc@ntr@l\et@tfigptintercircB@zDD}\ignorespaces}
\ctr@ln@m\figptinterlines
\ctr@ld@f\def\inters@cDD#1:#2[#3,#4;#5,#6]{{\s@uvc@ntr@l\et@tinters@cDD%
    \setc@ntr@l{2}\vecunit@{-1}{#4}\vecunit@{-2}{#6}%
    \Figg@tXY{-1}\setc@ntr@l{1}\Figg@tXYa{#3}%
    \edef\A@{\repdecn@mb{\v@lX}}\edef\B@{\repdecn@mb{\v@lY}}%
    \v@lmin=\B@\v@lXa\advance\v@lmin-\A@\v@lYa%
    \Figg@tXYa{#5}\setc@ntr@l{2}\Figg@tXY{-2}%
    \edef\C@{\repdecn@mb{\v@lX}}\edef\D@{\repdecn@mb{\v@lY}}%
    \v@lmax=\D@\v@lXa\advance\v@lmax-\C@\v@lYa%
    \delt@=\A@\v@lY\advance\delt@-\B@\v@lX%
    \invers@{\v@leur}{\delt@}\edef\v@ldelta{\repdecn@mb{\v@leur}}%
    \v@lXa=\A@\v@lmax\advance\v@lXa-\C@\v@lmin%
    \v@lYa=\B@\v@lmax\advance\v@lYa-\D@\v@lmin%
    \v@lXa=\v@ldelta\v@lXa\v@lYa=\v@ldelta\v@lYa%
    \setc@ntr@l{1}\Figp@intregDD#1:{#2}(\v@lXa,\v@lYa)%
    \resetc@ntr@l\et@tinters@cDD}\ignorespaces}
\ctr@ld@f\def\inters@cTD#1:#2[#3,#4;#5,#6]{{\s@uvc@ntr@l\et@tinters@cTD%
    \setc@ntr@l{2}\figvectNVTD-1[#4,#6]\figvectNVTD-2[#6,-1]\figvectPTD-1[#3,#5]%
    \r@pPSTD\v@leur[-2,-1,#4]\edef\v@lcoef{\repdecn@mb{\v@leur}}%
    \figpttraTD#1:{#2}=#3/\v@lcoef,#4/\resetc@ntr@l\et@tinters@cTD}\ignorespaces}
\ctr@ld@f\def\r@pPSTD#1[#2,#3,#4]{{\Figg@tXY{#2}\edef\Xu@{\repdecn@mb{\v@lX}}%
    \edef\Yu@{\repdecn@mb{\v@lY}}\edef\Zu@{\repdecn@mb{\v@lZ}}%
    \Figg@tXY{#3}\v@lmin=\Xu@\v@lX\advance\v@lmin\Yu@\v@lY\advance\v@lmin\Zu@\v@lZ%
    \Figg@tXY{#4}\v@lmax=\Xu@\v@lX\advance\v@lmax\Yu@\v@lY\advance\v@lmax\Zu@\v@lZ%
    \invers@{\v@leur}{\v@lmax}\global\result@t=\repdecn@mb{\v@leur}\v@lmin}%
    #1=\result@t}
\ctr@ln@m\n@rminf
\ctr@ld@f\def\n@rminfDD#1#2{{\Figg@tXY{#2}\maxim@m{\v@lX}{\v@lX}{-\v@lX}%
    \maxim@m{\v@lY}{\v@lY}{-\v@lY}\maxim@m{\global\result@t}{\v@lX}{\v@lY}}%
    #1=\result@t}
\ctr@ld@f\def\n@rminfTD#1#2{{\Figg@tXY{#2}\maxim@m{\v@lX}{\v@lX}{-\v@lX}%
    \maxim@m{\v@lY}{\v@lY}{-\v@lY}\maxim@m{\v@lZ}{\v@lZ}{-\v@lZ}%
    \maxim@m{\v@lX}{\v@lX}{\v@lY}\maxim@m{\global\result@t}{\v@lX}{\v@lZ}}%
    #1=\result@t}
\ctr@ld@f\def\n@rmeucCDD#1#2{\Figg@tXY{#2}\divide\v@lX\f@ctech\divide\v@lY\f@ctech%
    #1=\repdecn@mb{\v@lX}\v@lX\v@lX=\repdecn@mb{\v@lY}\v@lY\advance#1\v@lX}
\ctr@ld@f\def\n@rmeucCTD#1#2{\Figg@tXY{#2}%
    \divide\v@lX\f@ctech\divide\v@lY\f@ctech\divide\v@lZ\f@ctech%
    #1=\repdecn@mb{\v@lX}\v@lX\v@lX=\repdecn@mb{\v@lY}\v@lY\advance#1\v@lX%
    \v@lX=\repdecn@mb{\v@lZ}\v@lZ\advance#1\v@lX}
\ctr@ln@m\n@rmeucSV
\ctr@ld@f\def\n@rmeucSVDD#1#2{{\Figg@tXY{#2}%
    \v@lXa=\repdecn@mb{\v@lX}\v@lX\v@lYa=\repdecn@mb{\v@lY}\v@lY%
    \advance\v@lXa\v@lYa\sqrt@{\global\result@t}{\v@lXa}}#1=\result@t}
\ctr@ld@f\def\n@rmeucSVTD#1#2{{\Figg@tXY{#2}\v@lXa=\repdecn@mb{\v@lX}\v@lX%
    \v@lYa=\repdecn@mb{\v@lY}\v@lY\v@lZa=\repdecn@mb{\v@lZ}\v@lZ%
    \advance\v@lXa\v@lYa\advance\v@lXa\v@lZa\sqrt@{\global\result@t}{\v@lXa}}#1=\result@t}
\ctr@ln@m\n@rmeuc
\ctr@ld@f\def\n@rmeucDD#1#2{{\Figg@tXY{#2}\getredf@ctDD\f@ctech(\v@lX,\v@lY)%
    \divide\v@lX\f@ctech\divide\v@lY\f@ctech%
    \v@lXa=\repdecn@mb{\v@lX}\v@lX\v@lYa=\repdecn@mb{\v@lY}\v@lY%
    \advance\v@lXa\v@lYa\sqrt@{\global\result@t}{\v@lXa}%
    \global\multiply\result@t\f@ctech}#1=\result@t}
\ctr@ld@f\def\n@rmeucTD#1#2{{\Figg@tXY{#2}\getredf@ctTD\f@ctech(\v@lX,\v@lY,\v@lZ)%
    \divide\v@lX\f@ctech\divide\v@lY\f@ctech\divide\v@lZ\f@ctech%
    \v@lXa=\repdecn@mb{\v@lX}\v@lX%
    \v@lYa=\repdecn@mb{\v@lY}\v@lY\v@lZa=\repdecn@mb{\v@lZ}\v@lZ%
    \advance\v@lXa\v@lYa\advance\v@lXa\v@lZa\sqrt@{\global\result@t}{\v@lXa}%
    \global\multiply\result@t\f@ctech}#1=\result@t}
\ctr@ln@m\vecunit@
\ctr@ld@f\def\vecunit@DD#1#2{{\Figg@tXY{#2}\getredf@ctDD\f@ctech(\v@lX,\v@lY)%
    \divide\v@lX\f@ctech\divide\v@lY\f@ctech%
    \Figv@ctCreg#1(\v@lX,\v@lY)\n@rmeucSV{\v@lYa}{#1}%
    \invers@{\v@lXa}{\v@lYa}\edef\v@lv@lXa{\repdecn@mb{\v@lXa}}%
    \v@lX=\v@lv@lXa\v@lX\v@lY=\v@lv@lXa\v@lY%
    \Figv@ctCreg#1(\v@lX,\v@lY)\multiply\v@lYa\f@ctech\global\result@t=\v@lYa}}
\ctr@ld@f\def\vecunit@TD#1#2{{\Figg@tXY{#2}\getredf@ctTD\f@ctech(\v@lX,\v@lY,\v@lZ)%
    \divide\v@lX\f@ctech\divide\v@lY\f@ctech\divide\v@lZ\f@ctech%
    \Figv@ctCreg#1(\v@lX,\v@lY,\v@lZ)\n@rmeucSV{\v@lYa}{#1}%
    \invers@{\v@lXa}{\v@lYa}\edef\v@lv@lXa{\repdecn@mb{\v@lXa}}%
    \v@lX=\v@lv@lXa\v@lX\v@lY=\v@lv@lXa\v@lY\v@lZ=\v@lv@lXa\v@lZ%
    \Figv@ctCreg#1(\v@lX,\v@lY,\v@lZ)\multiply\v@lYa\f@ctech\global\result@t=\v@lYa}}
\ctr@ld@f\def\vecunitC@TD[#1,#2]{\Figg@tXYa{#1}\Figg@tXY{#2}%
    \advance\v@lX-\v@lXa\advance\v@lY-\v@lYa\advance\v@lZ-\v@lZa\c@lvecunitTD}
\ctr@ld@f\def\vecunitCV@TD#1{\Figg@tXY{#1}\c@lvecunitTD}
\ctr@ld@f\def\c@lvecunitTD{\getredf@ctTD\f@ctech(\v@lX,\v@lY,\v@lZ)%
    \divide\v@lX\f@ctech\divide\v@lY\f@ctech\divide\v@lZ\f@ctech%
    \v@lXa=\repdecn@mb{\v@lX}\v@lX%
    \v@lYa=\repdecn@mb{\v@lY}\v@lY\v@lZa=\repdecn@mb{\v@lZ}\v@lZ%
    \advance\v@lXa\v@lYa\advance\v@lXa\v@lZa\sqrt@{\v@lYa}{\v@lXa}%
    \invers@{\v@lXa}{\v@lYa}\edef\v@lv@lXa{\repdecn@mb{\v@lXa}}%
    \v@lX=\v@lv@lXa\v@lX\v@lY=\v@lv@lXa\v@lY\v@lZ=\v@lv@lXa\v@lZ}
\ctr@ln@m\figgetangle
\ctr@ld@f\def\figgetangleDD#1[#2,#3,#4]{\ifps@cri{\s@uvc@ntr@l\et@tfiggetangleDD\setc@ntr@l{2}%
    \figvectPDD-1[#2,#3]\figvectPDD-2[#2,#4]\vecunit@{-1}{-1}%
    \c@lproscalDD\delt@[-2,-1]\figvectNVDD-1[-1]\c@lproscalDD\v@leur[-2,-1]%
    \arct@n\v@lmax(\delt@,\v@leur)\v@lmax=\rdT@deg\v@lmax%
    \ifdim\v@lmax<\z@\advance\v@lmax\DePI@deg\fi\xdef#1{\repdecn@mb{\v@lmax}}%
    \resetc@ntr@l\et@tfiggetangleDD}\ignorespaces\fi}
\ctr@ld@f\def\figgetangleTD#1[#2,#3,#4,#5]{\ifps@cri{\s@uvc@ntr@l\et@tfiggetangleTD\setc@ntr@l{2}%
    \figvectPTD-1[#2,#3]\figvectPTD-2[#2,#5]\figvectNVTD-3[-1,-2]%
    \figvectPTD-2[#2,#4]\figvectNVTD-4[-3,-1]%
    \vecunit@{-1}{-1}\c@lproscalTD\delt@[-2,-1]\c@lproscalTD\v@leur[-2,-4]%
    \arct@n\v@lmax(\delt@,\v@leur)\v@lmax=\rdT@deg\v@lmax%
    \ifdim\v@lmax<\z@\advance\v@lmax\DePI@deg\fi\xdef#1{\repdecn@mb{\v@lmax}}%
    \resetc@ntr@l\et@tfiggetangleTD}\ignorespaces\fi}
\ctr@ld@f\def\figgetdist#1[#2,#3]{\ifps@cri{\s@uvc@ntr@l\et@tfiggetdist\setc@ntr@l{2}%
    \figvectP-1[#2,#3]\n@rmeuc{\v@lX}{-1}\v@lX=\ptT@unit@\v@lX\xdef#1{\repdecn@mb{\v@lX}}%
    \resetc@ntr@l\et@tfiggetdist}\ignorespaces\fi}
\ctr@ld@f\def\Figg@tT#1{\c@ntr@lnum{#1}%
    {\expandafter\expandafter\expandafter\extr@ctT\csname\objc@de\endcsname:%
     \ifnum\B@@ltxt=\z@\ptn@me{#1}\else\csname\objc@de T\endcsname\fi}}
\ctr@ld@f\def\extr@ctT#1,#2,#3/#4:{\def\B@@ltxt{#3}}
\ctr@ld@f\def\Figg@tXY#1{\c@ntr@lnum{#1}%
    \expandafter\expandafter\expandafter\extr@ctC\csname\objc@de\endcsname:}
\ctr@ln@m\extr@ctC
\ctr@ld@f\def\extr@ctCDD#1/#2,#3,#4:{\v@lX=#2\v@lY=#3}
\ctr@ld@f\def\extr@ctCTD#1/#2,#3,#4:{\v@lX=#2\v@lY=#3\v@lZ=#4}
\ctr@ld@f\def\Figg@tXYa#1{\c@ntr@lnum{#1}%
    \expandafter\expandafter\expandafter\extr@ctCa\csname\objc@de\endcsname:}
\ctr@ln@m\extr@ctCa
\ctr@ld@f\def\extr@ctCaDD#1/#2,#3,#4:{\v@lXa=#2\v@lYa=#3}
\ctr@ld@f\def\extr@ctCaTD#1/#2,#3,#4:{\v@lXa=#2\v@lYa=#3\v@lZa=#4}
\ctr@ln@m\t@xt@
\ctr@ld@f\def\figinit#1{\t@stc@tcodech@nge\initpr@lim\Figinit@#1,:\initpss@ttings\ignorespaces}
\ctr@ld@f\def\Figinit@#1,#2:{\setunit@{#1}\def\t@xt@{#2}\ifx\t@xt@\empty\else\Figinit@@#2:\fi}
\ctr@ld@f\def\Figinit@@#1#2:{\if#12 \else\Figs@tproj{#1}\initTD@\fi}
\ctr@ln@w{newif}\ifTr@isDim
\ctr@ld@f\def\UnD@fined{UNDEFINED}
\ctr@ld@f\def\ifundefined#1{\expandafter\ifx\csname#1\endcsname\relax}
\ctr@ln@m\@utoFN
\ctr@ln@m\@utoFInDone
\ctr@ln@m\disob@unit
\ctr@ld@f\def\initpr@lim{\initb@undb@x\figsetmark{}\figsetptname{$A_{##1}$}\def\Sc@leFact{1}%
    \initDD@\figsetroundcoord{yes}\ps@critrue\expandafter\setupd@te\defaultupdate:%
    \edef\disob@unit{\UnD@fined}\edef\t@rgetpt{\UnD@fined}\gdef\@utoFInDone{1}\gdef\@utoFN{0}}
\ctr@ld@f\def\initDD@{\Tr@isDimfalse%
    \ifPDFm@ke%
     \let\Ps@rcerc=\Ps@rcercBz%
     \let\Ps@rell=\Ps@rellBz%
    \fi
    \let\c@lDCUn=\c@lDCUnDD%
    \let\c@lDCDeux=\c@lDCDeuxDD%
    \let\c@ldefproj=\relax%
    \let\c@lproscal=\c@lproscalDD%
    \let\c@lprojSP=\relax%
    \let\extr@ctC=\extr@ctCDD%
    \let\extr@ctCa=\extr@ctCaDD%
    \let\extr@ctCF=\extr@ctCFDD%
    \let\Figp@intreg=\Figp@intregDD%
    \let\Figpts@xes=\Figpts@xesDD%
    \let\n@rmeucSV=\n@rmeucSVDD\let\n@rmeuc=\n@rmeucDD\let\n@rminf=\n@rminfDD%
    \let\pr@dMatV=\pr@dMatVDD%
    \let\ps@xes=\ps@xesDD%
    \let\vecunit@=\vecunit@DD%
    \let\figcoord=\figcoordDD%
    \let\figgetangle=\figgetangleDD%
    \let\figpt=\figptDD%
    \let\figptBezier=\figptBezierDD%
    \let\figptbary=\figptbaryDD%
    \let\figptcirc=\figptcircDD%
    \let\figptcircumcenter=\figptcircumcenterDD%
    \let\figptcopy=\figptcopyDD%
    \let\figptcurvcenter=\figptcurvcenterDD%
    \let\figptell=\figptellDD%
    \let\figptendnormal=\figptendnormalDD%
    \let\figptinterlineplane=\figptinterlineplaneDD%
    \let\figptinterlines=\inters@cDD%
    \let\figptorthocenter=\figptorthocenterDD%
    \let\figptorthoprojline=\figptorthoprojlineDD%
    \let\figptorthoprojplane=\figptorthoprojplaneDD%
    \let\figptrot=\figptrotDD%
    \let\figptscontrol=\figptscontrolDD%
    \let\figptsintercirc=\figptsintercircDD%
    \let\figptsinterlinell=\figptsinterlinellDD%
    \let\figptsorthoprojline=\figptsorthoprojlineDD%
    \let\figptorthoprojplane=\figptorthoprojplaneDD%
    \let\figptsrot=\figptsrotDD%
    \let\figptssym=\figptssymDD%
    \let\figptstra=\figptstraDD%
    \let\figptsym=\figptsymDD%
    \let\figpttraC=\figpttraCDD%
    \let\figpttra=\figpttraDD%
    \let\figptvisilimSL=\figptvisilimSLDD%
    \let\figsetobdist=\figsetobdistDD%
    \let\figsettarget=\figsettargetDD%
    \let\figsetview=\figsetviewDD%
    \let\figvectDBezier=\figvectDBezierDD%
    \let\figvectN=\figvectNDD%
    \let\figvectNV=\figvectNVDD%
    \let\figvectP=\figvectPDD%
    \let\figvectU=\figvectUDD%
    \let\psarccircP=\psarccircPDD%
    \let\psarccirc=\psarccircDD%
    \let\psarcell=\psarcellDD%
    \let\psarcellPA=\psarcellPADD%
    \let\psarrowBezier=\psarrowBezierDD%
    \let\psarrowcircP=\psarrowcircPDD%
    \let\psarrowcirc=\psarrowcircDD%
    \let\psarrowhead=\psarrowheadDD%
    \let\psarrow=\psarrowDD%
    \let\psBezier=\psBezierDD%
    \let\pscirc=\pscircDD%
    \let\pscurve=\pscurveDD%
    \let\psnormal=\psnormalDD%
    }
\ctr@ld@f\def\initTD@{\Tr@isDimtrue\initb@undb@xTD\newt@rgetptfalse\newdis@bfalse%
    \let\c@lDCUn=\c@lDCUnTD%
    \let\c@lDCDeux=\c@lDCDeuxTD%
    \let\c@ldefproj=\c@ldefprojTD%
    \let\c@lproscal=\c@lproscalTD%
    \let\extr@ctC=\extr@ctCTD%
    \let\extr@ctCa=\extr@ctCaTD%
    \let\extr@ctCF=\extr@ctCFTD%
    \let\Figp@intreg=\Figp@intregTD%
    \let\Figpts@xes=\Figpts@xesTD%
    \let\n@rmeucSV=\n@rmeucSVTD\let\n@rmeuc=\n@rmeucTD\let\n@rminf=\n@rminfTD%
    \let\pr@dMatV=\pr@dMatVTD%
    \let\ps@xes=\ps@xesTD%
    \let\vecunit@=\vecunit@TD%
    \let\figcoord=\figcoordTD%
    \let\figgetangle=\figgetangleTD%
    \let\figpt=\figptTD%
    \let\figptBezier=\figptBezierTD%
    \let\figptbary=\figptbaryTD%
    \let\figptcirc=\figptcircTD%
    \let\figptcircumcenter=\figptcircumcenterTD%
    \let\figptcopy=\figptcopyTD%
    \let\figptcurvcenter=\figptcurvcenterTD%
    \let\figptinterlineplane=\figptinterlineplaneTD%
    \let\figptinterlines=\inters@cTD%
    \let\figptorthocenter=\figptorthocenterTD%
    \let\figptorthoprojline=\figptorthoprojlineTD%
    \let\figptorthoprojplane=\figptorthoprojplaneTD%
    \let\figptrot=\figptrotTD%
    \let\figptscontrol=\figptscontrolTD%
    \let\figptsintercirc=\figptsintercircTD%
    \let\figptsorthoprojline=\figptsorthoprojlineTD%
    \let\figptsorthoprojplane=\figptsorthoprojplaneTD%
    \let\figptsrot=\figptsrotTD%
    \let\figptssym=\figptssymTD%
    \let\figptstra=\figptstraTD%
    \let\figptsym=\figptsymTD%
    \let\figpttraC=\figpttraCTD%
    \let\figpttra=\figpttraTD%
    \let\figptvisilimSL=\figptvisilimSLTD%
    \let\figsetobdist=\figsetobdistTD%
    \let\figsettarget=\figsettargetTD%
    \let\figsetview=\figsetviewTD%
    \let\figvectDBezier=\figvectDBezierTD%
    \let\figvectN=\figvectNTD%
    \let\figvectNV=\figvectNVTD%
    \let\figvectP=\figvectPTD%
    \let\figvectU=\figvectUTD%
    \let\psarccircP=\psarccircPTD%
    \let\psarccirc=\psarccircTD%
    \let\psarcell=\psarcellTD%
    \let\psarcellPA=\psarcellPATD%
    \let\psarrowBezier=\psarrowBezierTD%
    \let\psarrowcircP=\psarrowcircPTD%
    \let\psarrowcirc=\psarrowcircTD%
    \let\psarrowhead=\psarrowheadTD%
    \let\psarrow=\psarrowTD%
    \let\psBezier=\psBezierTD%
    \let\pscirc=\pscircTD%
    \let\pscurve=\pscurveTD%
    }
\ctr@ld@f\def\un@v@ilable#1{\immediate\write16{*** The macro #1 is not available in the current context.}}
\ctr@ld@f\def\figinsert#1{{\def\t@xt@{#1}\relax%
    \ifx\t@xt@\empty\ifnum\@utoFInDone>\z@\Figinsert@\DefGIfilen@me,:\fi%
    \else\expandafter\FiginsertNu@#1 :\fi}\ignorespaces}
\ctr@ld@f\def\FiginsertNu@#1 #2:{\def\t@xt@{#1}\relax\ifx\t@xt@\empty\def\t@xt@{#2}%
    \ifx\t@xt@\empty\ifnum\@utoFInDone>\z@\Figinsert@\DefGIfilen@me,:\fi%
    \else\FiginsertNu@#2:\fi\else\expandafter\FiginsertNd@#1 #2:\fi}
\ctr@ld@f\def\FiginsertNd@#1#2:{\ifcat#1a\Figinsert@#1#2,:\else%
    \ifnum\@utoFInDone>\z@\Figinsert@\DefGIfilen@me,#1#2,:\fi\fi}
\ctr@ln@m\Sc@leFact
\ctr@ld@f\def\Figinsert@#1,#2:{\def\t@xt@{#2}\ifx\t@xt@\empty\xdef\Sc@leFact{1}\else%
    \X@rgdeux@#2\xdef\Sc@leFact{\@rgdeux}\fi%
    \Figdisc@rdLTS{#1}{\t@xt@}\@psfgetbb{\t@xt@}%
    \v@lX=\@psfllx\p@\v@lX=\ptpsT@pt\v@lX\v@lX=\Sc@leFact\v@lX%
    \v@lY=\@psflly\p@\v@lY=\ptpsT@pt\v@lY\v@lY=\Sc@leFact\v@lY%
    \b@undb@x{\v@lX}{\v@lY}%
    \v@lX=\@psfurx\p@\v@lX=\ptpsT@pt\v@lX\v@lX=\Sc@leFact\v@lX%
    \v@lY=\@psfury\p@\v@lY=\ptpsT@pt\v@lY\v@lY=\Sc@leFact\v@lY%
    \b@undb@x{\v@lX}{\v@lY}%
    \ifPDFm@ke\Figinclud@PDF{\t@xt@}{\Sc@leFact}\else%
    \v@lX=\c@nt pt\v@lX=\Sc@leFact\v@lX\edef\F@ct{\repdecn@mb{\v@lX}}%
    \ifx\TeXturesonMacOSltX\special{postscriptfile #1 vscale=\F@ct\space hscale=\F@ct}%
    \else\includegraphics{#1}\fi\fi%
    \message{[\t@xt@]}\ignorespaces}
\ctr@ld@f\def\Figdisc@rdLTS#1#2{\expandafter\Figdisc@rdLTS@#1 :#2}
\ctr@ld@f\def\Figdisc@rdLTS@#1 #2:#3{\def#3{#1}\relax\ifx#3\empty\expandafter\Figdisc@rdLTS@#2:#3\fi}
\ctr@ld@f\def\figinsertE#1{\FiginsertE@#1,:\ignorespaces}
\ctr@ld@f\def\FiginsertE@#1,#2:{{\def\t@xt@{#2}\ifx\t@xt@\empty\xdef\Sc@leFact{1}\else%
    \X@rgdeux@#2\xdef\Sc@leFact{\@rgdeux}\fi%
    \Figdisc@rdLTS{#1}{\t@xt@}\pdfximage{\t@xt@}%
    \setbox\Gb@x=\hbox{\pdfrefximage\pdflastximage}%
    \v@lX=\z@\v@lY=-\Sc@leFact\dp\Gb@x\b@undb@x{\v@lX}{\v@lY}%
    \advance\v@lX\Sc@leFact\wd\Gb@x\advance\v@lY\Sc@leFact\dp\Gb@x%
    \advance\v@lY\Sc@leFact\ht\Gb@x\b@undb@x{\v@lX}{\v@lY}%
    \v@lX=\Sc@leFact\wd\Gb@x\pdfximage width \v@lX {\t@xt@}%
    \rlap{\pdfrefximage\pdflastximage}\message{[\t@xt@]}}\ignorespaces}
\ctr@ld@f\def\X@rgdeux@#1,{\edef\@rgdeux{#1}}
\ctr@ln@m\figpt
\ctr@ld@f\def\figptDD#1:#2(#3,#4){\ifps@cri\c@ntr@lnum{#1}%
    {\v@lX=#3\unit@\v@lY=#4\unit@\Fig@dmpt{#2}{\z@}}\ignorespaces\fi}
\ctr@ld@f\def\Fig@dmpt#1#2{\def\t@xt@{#1}\ifx\t@xt@\empty\def\B@@ltxt{\z@}%
    \else\expandafter\gdef\csname\objc@de T\endcsname{#1}\def\B@@ltxt{\@ne}\fi%
    \expandafter\xdef\csname\objc@de\endcsname{\ifitis@vect@r\C@dCl@svect%
    \else\C@dCl@spt\fi,\z@,\B@@ltxt/\the\v@lX,\the\v@lY,#2}}
\ctr@ld@f\def\C@dCl@spt{P}
\ctr@ld@f\def\C@dCl@svect{V}
\ctr@ln@m\c@@rdYZ
\ctr@ln@m\c@@rdY
\ctr@ld@f\def\figptTD#1:#2(#3,#4){\ifps@cri\c@ntr@lnum{#1}%
    \def\c@@rdYZ{#4,0,0}\extrairelepremi@r\c@@rdY\de\c@@rdYZ%
    \extrairelepremi@r\c@@rdZ\de\c@@rdYZ%
    {\v@lX=#3\unit@\v@lY=\c@@rdY\unit@\v@lZ=\c@@rdZ\unit@\Fig@dmpt{#2}{\the\v@lZ}%
    \b@undb@xTD{\v@lX}{\v@lY}{\v@lZ}}\ignorespaces\fi}
\ctr@ln@m\Figp@intreg
\ctr@ld@f\def\Figp@intregDD#1:#2(#3,#4){\c@ntr@lnum{#1}%
    {\result@t=#4\v@lX=#3\v@lY=\result@t\Fig@dmpt{#2}{\z@}}\ignorespaces}
\ctr@ld@f\def\Figp@intregTD#1:#2(#3,#4){\c@ntr@lnum{#1}%
    \def\c@@rdYZ{#4,\z@,\z@}\extrairelepremi@r\c@@rdY\de\c@@rdYZ%
    \extrairelepremi@r\c@@rdZ\de\c@@rdYZ%
    {\v@lX=#3\v@lY=\c@@rdY\v@lZ=\c@@rdZ\Fig@dmpt{#2}{\the\v@lZ}%
    \b@undb@xTD{\v@lX}{\v@lY}{\v@lZ}}\ignorespaces}
\ctr@ln@m\figptBezier
\ctr@ld@f\def\figptBezierDD#1:#2:#3[#4,#5,#6,#7]{\ifps@cri{\s@uvc@ntr@l\et@tfigptBezierDD%
    \FigptBezier@#3[#4,#5,#6,#7]\Figp@intregDD#1:{#2}(\v@lX,\v@lY)%
    \resetc@ntr@l\et@tfigptBezierDD}\ignorespaces\fi}
\ctr@ld@f\def\figptBezierTD#1:#2:#3[#4,#5,#6,#7]{\ifps@cri{\s@uvc@ntr@l\et@tfigptBezierTD%
    \FigptBezier@#3[#4,#5,#6,#7]\Figp@intregTD#1:{#2}(\v@lX,\v@lY,\v@lZ)%
    \resetc@ntr@l\et@tfigptBezierTD}\ignorespaces\fi}
\ctr@ld@f\def\FigptBezier@#1[#2,#3,#4,#5]{\setc@ntr@l{2}%
    \edef\T@{#1}\v@leur=\p@\advance\v@leur-#1pt\edef\UNmT@{\repdecn@mb{\v@leur}}%
    \figptcopy-4:/#2/\figptcopy-3:/#3/\figptcopy-2:/#4/\figptcopy-1:/#5/%
    \l@mbd@un=-4 \l@mbd@de=-\thr@@\p@rtent=\m@ne\c@lDecast%
    \l@mbd@un=-4 \l@mbd@de=-\thr@@\p@rtent=-\tw@\c@lDecast%
    \l@mbd@un=-4 \l@mbd@de=-\thr@@\p@rtent=-\thr@@\c@lDecast\Figg@tXY{-4}}
\ctr@ln@m\c@lDCUn
\ctr@ld@f\def\c@lDCUnDD#1#2{\Figg@tXY{#1}\v@lX=\UNmT@\v@lX\v@lY=\UNmT@\v@lY%
    \Figg@tXYa{#2}\advance\v@lX\T@\v@lXa\advance\v@lY\T@\v@lYa%
    \Figp@intregDD#1:(\v@lX,\v@lY)}
\ctr@ld@f\def\c@lDCUnTD#1#2{\Figg@tXY{#1}\v@lX=\UNmT@\v@lX\v@lY=\UNmT@\v@lY\v@lZ=\UNmT@\v@lZ%
    \Figg@tXYa{#2}\advance\v@lX\T@\v@lXa\advance\v@lY\T@\v@lYa\advance\v@lZ\T@\v@lZa%
    \Figp@intregTD#1:(\v@lX,\v@lY,\v@lZ)}
\ctr@ld@f\def\c@lDecast{\relax\ifnum\l@mbd@un<\p@rtent\c@lDCUn{\l@mbd@un}{\l@mbd@de}%
    \advance\l@mbd@un\@ne\advance\l@mbd@de\@ne\c@lDecast\fi}
\ctr@ld@f\def\figptmap#1:#2=#3/#4/#5/{\ifps@cri{\s@uvc@ntr@l\et@tfigptmap%
    \setc@ntr@l{2}\figvectP-1[#4,#3]\Figg@tXY{-1}%
    \pr@dMatV/#5/\figpttra#1:{#2}=#4/1,-1/%
    \resetc@ntr@l\et@tfigptmap}\ignorespaces\fi}
\ctr@ln@m\pr@dMatV
\ctr@ld@f\def\pr@dMatVDD/#1,#2;#3,#4/{\v@lXa=#1\v@lX\advance\v@lXa#2\v@lY%
    \v@lYa=#3\v@lX\advance\v@lYa#4\v@lY\Figv@ctCreg-1(\v@lXa,\v@lYa)}
\ctr@ld@f\def\pr@dMatVTD/#1,#2,#3;#4,#5,#6;#7,#8,#9/{%
    \v@lXa=#1\v@lX\advance\v@lXa#2\v@lY\advance\v@lXa#3\v@lZ%
    \v@lYa=#4\v@lX\advance\v@lYa#5\v@lY\advance\v@lYa#6\v@lZ%
    \v@lZa=#7\v@lX\advance\v@lZa#8\v@lY\advance\v@lZa#9\v@lZ%
    \Figv@ctCreg-1(\v@lXa,\v@lYa,\v@lZa)}
\ctr@ln@m\figptbary
\ctr@ld@f\def\figptbaryDD#1:#2[#3;#4]{\ifps@cri{\edef\list@num{#3}\extrairelepremi@r\p@int\de\list@num%
    \s@mme=\z@\@ecfor\c@ef:=#4\do{\advance\s@mme\c@ef}%
    \edef\listec@ef{#4,0}\extrairelepremi@r\c@ef\de\listec@ef%
    \Figg@tXY{\p@int}\divide\v@lX\s@mme\divide\v@lY\s@mme%
    \multiply\v@lX\c@ef\multiply\v@lY\c@ef%
    \@ecfor\p@int:=\list@num\do{\extrairelepremi@r\c@ef\de\listec@ef%
           \Figg@tXYa{\p@int}\divide\v@lXa\s@mme\divide\v@lYa\s@mme%
           \multiply\v@lXa\c@ef\multiply\v@lYa\c@ef%
           \advance\v@lX\v@lXa\advance\v@lY\v@lYa}%
    \Figp@intregDD#1:{#2}(\v@lX,\v@lY)}\ignorespaces\fi}
\ctr@ld@f\def\figptbaryTD#1:#2[#3;#4]{\ifps@cri{\edef\list@num{#3}\extrairelepremi@r\p@int\de\list@num%
    \s@mme=\z@\@ecfor\c@ef:=#4\do{\advance\s@mme\c@ef}%
    \edef\listec@ef{#4,0}\extrairelepremi@r\c@ef\de\listec@ef%
    \Figg@tXY{\p@int}\divide\v@lX\s@mme\divide\v@lY\s@mme\divide\v@lZ\s@mme%
    \multiply\v@lX\c@ef\multiply\v@lY\c@ef\multiply\v@lZ\c@ef%
    \@ecfor\p@int:=\list@num\do{\extrairelepremi@r\c@ef\de\listec@ef%
           \Figg@tXYa{\p@int}\divide\v@lXa\s@mme\divide\v@lYa\s@mme\divide\v@lZa\s@mme%
           \multiply\v@lXa\c@ef\multiply\v@lYa\c@ef\multiply\v@lZa\c@ef%
           \advance\v@lX\v@lXa\advance\v@lY\v@lYa\advance\v@lZ\v@lZa}%
    \Figp@intregTD#1:{#2}(\v@lX,\v@lY,\v@lZ)}\ignorespaces\fi}
\ctr@ld@f\def\figptbaryR#1:#2[#3;#4]{\ifps@cri{%
    \v@leur=\z@\@ecfor\c@ef:=#4\do{\maxim@m{\v@lmax}{\c@ef pt}{-\c@ef pt}%
    \ifdim\v@lmax>\v@leur\v@leur=\v@lmax\fi}%
    \ifdim\v@leur<\p@\f@ctech=\@M\else\ifdim\v@leur<\t@n\p@\f@ctech=\@m\else%
    \ifdim\v@leur<\c@nt\p@\f@ctech=\c@nt\else\ifdim\v@leur<\@m\p@\f@ctech=\t@n\else%
    \f@ctech=\@ne\fi\fi\fi\fi%
    \def\listec@ef{0}%
    \@ecfor\c@ef:=#4\do{\sc@lec@nvRI{\c@ef pt}\edef\listec@ef{\listec@ef,\the\s@mme}}%
    \extrairelepremi@r\c@ef\de\listec@ef\figptbary#1:#2[#3;\listec@ef]}\ignorespaces\fi}
\ctr@ld@f\def\sc@lec@nvRI#1{\v@leur=#1\p@rtentiere{\s@mme}{\v@leur}\advance\v@leur-\s@mme\p@%
    \multiply\v@leur\f@ctech\p@rtentiere{\p@rtent}{\v@leur}%
    \multiply\s@mme\f@ctech\advance\s@mme\p@rtent}
\ctr@ln@m\figptcirc
\ctr@ld@f\def\figptcircDD#1:#2:#3;#4(#5){\ifps@cri{\s@uvc@ntr@l\et@tfigptcircDD%
    \c@lptellDD#1:{#2}:#3;#4,#4(#5)\resetc@ntr@l\et@tfigptcircDD}\ignorespaces\fi}
\ctr@ld@f\def\figptcircTD#1:#2:#3,#4,#5;#6(#7){\ifps@cri{\s@uvc@ntr@l\et@tfigptcircTD%
    \setc@ntr@l{2}\c@lExtAxes#3,#4,#5(#6)\figptellP#1:{#2}:#3,-4,-5(#7)%
    \resetc@ntr@l\et@tfigptcircTD}\ignorespaces\fi}
\ctr@ln@m\figptcircumcenter
\ctr@ld@f\def\figptcircumcenterDD#1:#2[#3,#4,#5]{\ifps@cri{\s@uvc@ntr@l\et@tfigptcircumcenterDD%
    \setc@ntr@l{2}\figvectNDD-5[#3,#4]\figptbaryDD-3:[#3,#4;1,1]%
                  \figvectNDD-6[#4,#5]\figptbaryDD-4:[#4,#5;1,1]%
    \resetc@ntr@l{2}\inters@cDD#1:{#2}[-3,-5;-4,-6]%
    \resetc@ntr@l\et@tfigptcircumcenterDD}\ignorespaces\fi}
\ctr@ld@f\def\figptcircumcenterTD#1:#2[#3,#4,#5]{\ifps@cri{\s@uvc@ntr@l\et@tfigptcircumcenterTD%
    \setc@ntr@l{2}\figvectNTD-1[#3,#4,#5]%
    \figvectPTD-3[#3,#4]\figvectNVTD-5[-1,-3]\figptbaryTD-3:[#3,#4;1,1]%
    \figvectPTD-4[#4,#5]\figvectNVTD-6[-1,-4]\figptbaryTD-4:[#4,#5;1,1]%
    \resetc@ntr@l{2}\inters@cTD#1:{#2}[-3,-5;-4,-6]%
    \resetc@ntr@l\et@tfigptcircumcenterTD}\ignorespaces\fi}
\ctr@ln@m\figptcopy
\ctr@ld@f\def\figptcopyDD#1:#2/#3/{\ifps@cri{\Figg@tXY{#3}%
    \Figp@intregDD#1:{#2}(\v@lX,\v@lY)}\ignorespaces\fi}
\ctr@ld@f\def\figptcopyTD#1:#2/#3/{\ifps@cri{\Figg@tXY{#3}%
    \Figp@intregTD#1:{#2}(\v@lX,\v@lY,\v@lZ)}\ignorespaces\fi}
\ctr@ln@m\figptcurvcenter
\ctr@ld@f\def\figptcurvcenterDD#1:#2:#3[#4,#5,#6,#7]{\ifps@cri{\s@uvc@ntr@l\et@tfigptcurvcenterDD%
    \setc@ntr@l{2}\c@lcurvradDD#3[#4,#5,#6,#7]\edef\Sprim@{\repdecn@mb{\result@t}}%
    \figptBezierDD-1::#3[#4,#5,#6,#7]\figpttraDD#1:{#2}=-1/\Sprim@,-5/%
    \resetc@ntr@l\et@tfigptcurvcenterDD}\ignorespaces\fi}
\ctr@ld@f\def\figptcurvcenterTD#1:#2:#3[#4,#5,#6,#7]{\ifps@cri{\s@uvc@ntr@l\et@tfigptcurvcenterTD%
    \setc@ntr@l{2}\figvectDBezierTD -5:1,#3[#4,#5,#6,#7]%
    \figvectDBezierTD -6:2,#3[#4,#5,#6,#7]\vecunit@TD{-5}{-5}%
    \edef\Sprim@{\repdecn@mb{\result@t}}\figvectNVTD-1[-6,-5]%
    \figvectNVTD-5[-5,-1]\c@lproscalTD\v@leur[-6,-5]%
    \invers@{\v@leur}{\v@leur}\v@leur=\Sprim@\v@leur\v@leur=\Sprim@\v@leur%
    \figptBezierTD-1::#3[#4,#5,#6,#7]\edef\Sprim@{\repdecn@mb{\v@leur}}%
    \figpttraTD#1:{#2}=-1/\Sprim@,-5/\resetc@ntr@l\et@tfigptcurvcenterTD}\ignorespaces\fi}
\ctr@ld@f\def\c@lcurvradDD#1[#2,#3,#4,#5]{{\figvectDBezierDD -5:1,#1[#2,#3,#4,#5]%
    \figvectDBezierDD -6:2,#1[#2,#3,#4,#5]\vecunit@DD{-5}{-5}%
    \edef\Sprim@{\repdecn@mb{\result@t}}\figvectNVDD-5[-5]\c@lproscalDD\v@leur[-6,-5]%
    \invers@{\v@leur}{\v@leur}\v@leur=\Sprim@\v@leur\v@leur=\Sprim@\v@leur%
    \global\result@t=\v@leur}}
\ctr@ln@m\figptell
\ctr@ld@f\def\figptellDD#1:#2:#3;#4,#5(#6,#7){\ifps@cri{\s@uvc@ntr@l\et@tfigptell%
    \c@lptellDD#1::#3;#4,#5(#6)\figptrotDD#1:{#2}=#1/#3,#7/%
    \resetc@ntr@l\et@tfigptell}\ignorespaces\fi}
\ctr@ld@f\def\c@lptellDD#1:#2:#3;#4,#5(#6){\c@ssin{\C@}{\S@}{#6}\v@lmin=\C@ pt\v@lmax=\S@ pt%
    \v@lmin=#4\v@lmin\v@lmax=#5\v@lmax%
    \edef\Xc@mp{\repdecn@mb{\v@lmin}}\edef\Yc@mp{\repdecn@mb{\v@lmax}}%
    \setc@ntr@l{2}\figvectC-1(\Xc@mp,\Yc@mp)\figpttraDD#1:{#2}=#3/1,-1/}
\ctr@ld@f\def\figptellP#1:#2:#3,#4,#5(#6){\ifps@cri{\s@uvc@ntr@l\et@tfigptellP%
    \setc@ntr@l{2}\figvectP-1[#3,#4]\figvectP-2[#3,#5]%
    \v@leur=#6pt\c@lptellP{#3}{-1}{-2}\figptcopy#1:{#2}/-3/%
    \resetc@ntr@l\et@tfigptellP}\ignorespaces\fi}
\ctr@ln@m\@ngle
\ctr@ld@f\def\c@lptellP#1#2#3{\edef\@ngle{\repdecn@mb\v@leur}\c@ssin{\C@}{\S@}{\@ngle}%
    \figpttra-3:=#1/\C@,#2/\figpttra-3:=-3/\S@,#3/}
\ctr@ln@m\figptendnormal
\ctr@ld@f\def\figptendnormalDD#1:#2:#3,#4[#5,#6]{\ifps@cri{\s@uvc@ntr@l\et@tfigptendnormal%
    \Figg@tXYa{#5}\Figg@tXY{#6}%
    \advance\v@lX-\v@lXa\advance\v@lY-\v@lYa%
    \setc@ntr@l{2}\Figv@ctCreg-1(\v@lX,\v@lY)\vecunit@{-1}{-1}\Figg@tXY{-1}%
    \delt@=#3\unit@\maxim@m{\delt@}{\delt@}{-\delt@}\edef\l@ngueur{\repdecn@mb{\delt@}}%
    \v@lX=\l@ngueur\v@lX\v@lY=\l@ngueur\v@lY%
    \delt@=\p@\advance\delt@-#4pt\edef\l@ngueur{\repdecn@mb{\delt@}}%
    \figptbaryR-1:[#5,#6;#4,\l@ngueur]\Figg@tXYa{-1}%
    \advance\v@lXa\v@lY\advance\v@lYa-\v@lX%
    \setc@ntr@l{1}\Figp@intregDD#1:{#2}(\v@lXa,\v@lYa)\resetc@ntr@l\et@tfigptendnormal}%
    \ignorespaces\fi}
\ctr@ld@f\def\figptexcenter#1:#2[#3,#4,#5]{\ifps@cri{\let@xte={-}%
    \Figptexinsc@nter#1:#2[#3,#4,#5]}\ignorespaces\fi}
\ctr@ld@f\def\figptincenter#1:#2[#3,#4,#5]{\ifps@cri{\let@xte={}%
    \Figptexinsc@nter#1:#2[#3,#4,#5]}\ignorespaces\fi}
\ctr@ld@f\let\figptinscribedcenter=\figptincenter
\ctr@ld@f\def\Figptexinsc@nter#1:#2[#3,#4,#5]{%
    \figgetdist\LA@[#4,#5]\figgetdist\LB@[#3,#5]\figgetdist\LC@[#3,#4]%
    \figptbaryR#1:{#2}[#3,#4,#5;\the\let@xte\LA@,\LB@,\LC@]}
\ctr@ln@m\figptinterlineplane
\ctr@ld@f\def\figptinterlineplaneDD{\un@v@ilable{figptinterlineplane}}
\ctr@ld@f\def\figptinterlineplaneTD#1:#2[#3,#4;#5,#6]{\ifps@cri{\s@uvc@ntr@l\et@tfigptinterlineplane%
    \setc@ntr@l{2}\figvectPTD-1[#3,#5]\vecunit@TD{-2}{#6}%
    \r@pPSTD\v@leur[-2,-1,#4]\edef\v@lcoef{\repdecn@mb{\v@leur}}%
    \figpttraTD#1:{#2}=#3/\v@lcoef,#4/\resetc@ntr@l\et@tfigptinterlineplane}\ignorespaces\fi}
\ctr@ln@m\figptorthocenter
\ctr@ld@f\def\figptorthocenterDD#1:#2[#3,#4,#5]{\ifps@cri{\s@uvc@ntr@l\et@tfigptorthocenterDD%
    \setc@ntr@l{2}\figvectNDD-3[#3,#4]\figvectNDD-4[#4,#5]%
    \resetc@ntr@l{2}\inters@cDD#1:{#2}[#5,-3;#3,-4]%
    \resetc@ntr@l\et@tfigptorthocenterDD}\ignorespaces\fi}
\ctr@ld@f\def\figptorthocenterTD#1:#2[#3,#4,#5]{\ifps@cri{\s@uvc@ntr@l\et@tfigptorthocenterTD%
    \setc@ntr@l{2}\figvectNTD-1[#3,#4,#5]%
    \figvectPTD-2[#3,#4]\figvectNVTD-3[-1,-2]%
    \figvectPTD-2[#4,#5]\figvectNVTD-4[-1,-2]%
    \resetc@ntr@l{2}\inters@cTD#1:{#2}[#5,-3;#3,-4]%
    \resetc@ntr@l\et@tfigptorthocenterTD}\ignorespaces\fi}
\ctr@ln@m\figptorthoprojline
\ctr@ld@f\def\figptorthoprojlineDD#1:#2=#3/#4,#5/{\ifps@cri{\s@uvc@ntr@l\et@tfigptorthoprojlineDD%
    \setc@ntr@l{2}\figvectPDD-3[#4,#5]\figvectNVDD-4[-3]\resetc@ntr@l{2}%
    \inters@cDD#1:{#2}[#3,-4;#4,-3]\resetc@ntr@l\et@tfigptorthoprojlineDD}\ignorespaces\fi}
\ctr@ld@f\def\figptorthoprojlineTD#1:#2=#3/#4,#5/{\ifps@cri{\s@uvc@ntr@l\et@tfigptorthoprojlineTD%
    \setc@ntr@l{2}\figvectPTD-1[#4,#3]\figvectPTD-2[#4,#5]\vecunit@TD{-2}{-2}%
    \c@lproscalTD\v@leur[-1,-2]\edef\v@lcoef{\repdecn@mb{\v@leur}}%
    \figpttraTD#1:{#2}=#4/\v@lcoef,-2/\resetc@ntr@l\et@tfigptorthoprojlineTD}\ignorespaces\fi}
\ctr@ln@m\figptorthoprojplane
\ctr@ld@f\def\figptorthoprojplaneDD{\un@v@ilable{figptorthoprojplane}}
\ctr@ld@f\def\figptorthoprojplaneTD#1:#2=#3/#4,#5/{\ifps@cri{\s@uvc@ntr@l\et@tfigptorthoprojplane%
    \setc@ntr@l{2}\figvectPTD-1[#3,#4]\vecunit@TD{-2}{#5}%
    \c@lproscalTD\v@leur[-1,-2]\edef\v@lcoef{\repdecn@mb{\v@leur}}%
    \figpttraTD#1:{#2}=#3/\v@lcoef,-2/\resetc@ntr@l\et@tfigptorthoprojplane}\ignorespaces\fi}
\ctr@ld@f\def\figpthom#1:#2=#3/#4,#5/{\ifps@cri{\s@uvc@ntr@l\et@tfigpthom%
    \setc@ntr@l{2}\figvectP-1[#4,#3]\figpttra#1:{#2}=#4/#5,-1/%
    \resetc@ntr@l\et@tfigpthom}\ignorespaces\fi}
\ctr@ln@m\figptrot
\ctr@ld@f\def\figptrotDD#1:#2=#3/#4,#5/{\ifps@cri{\s@uvc@ntr@l\et@tfigptrotDD%
    \c@ssin{\C@}{\S@}{#5}\setc@ntr@l{2}\figvectPDD-1[#4,#3]\Figg@tXY{-1}%
    \v@lXa=\C@\v@lX\advance\v@lXa-\S@\v@lY%
    \v@lYa=\S@\v@lX\advance\v@lYa\C@\v@lY%
    \Figv@ctCreg-1(\v@lXa,\v@lYa)\figpttraDD#1:{#2}=#4/1,-1/%
    \resetc@ntr@l\et@tfigptrotDD}\ignorespaces\fi}
\ctr@ld@f\def\figptrotTD#1:#2=#3/#4,#5,#6/{\ifps@cri{\s@uvc@ntr@l\et@tfigptrotTD%
    \c@ssin{\C@}{\S@}{#5}%
    \setc@ntr@l{2}\figptorthoprojplaneTD-3:=#4/#3,#6/\figvectPTD-2[-3,#3]%
    \n@rmeucTD\v@leur{-2}\ifdim\v@leur<\Cepsil@n\Figg@tXYa{#3}\else%
    \edef\v@lcoef{\repdecn@mb{\v@leur}}\figvectNVTD-1[#6,-2]%
    \Figg@tXYa{-1}\v@lXa=\v@lcoef\v@lXa\v@lYa=\v@lcoef\v@lYa\v@lZa=\v@lcoef\v@lZa%
    \v@lXa=\S@\v@lXa\v@lYa=\S@\v@lYa\v@lZa=\S@\v@lZa\Figg@tXY{-2}%
    \advance\v@lXa\C@\v@lX\advance\v@lYa\C@\v@lY\advance\v@lZa\C@\v@lZ%
    \Figg@tXY{-3}\advance\v@lXa\v@lX\advance\v@lYa\v@lY\advance\v@lZa\v@lZ\fi%
    \Figp@intregTD#1:{#2}(\v@lXa,\v@lYa,\v@lZa)\resetc@ntr@l\et@tfigptrotTD}\ignorespaces\fi}
\ctr@ln@m\figptsym
\ctr@ld@f\def\figptsymDD#1:#2=#3/#4,#5/{\ifps@cri{\s@uvc@ntr@l\et@tfigptsymDD%
    \resetc@ntr@l{2}\figptorthoprojlineDD-5:=#3/#4,#5/\figvectPDD-2[#3,-5]%
    \figpttraDD#1:{#2}=#3/2,-2/\resetc@ntr@l\et@tfigptsymDD}\ignorespaces\fi}
\ctr@ld@f\def\figptsymTD#1:#2=#3/#4,#5/{\ifps@cri{\s@uvc@ntr@l\et@tfigptsymTD%
    \resetc@ntr@l{2}\figptorthoprojplaneTD-3:=#3/#4,#5/\figvectPTD-2[#3,-3]%
    \figpttraTD#1:{#2}=#3/2,-2/\resetc@ntr@l\et@tfigptsymTD}\ignorespaces\fi}
\ctr@ln@m\figpttra
\ctr@ld@f\def\figpttraDD#1:#2=#3/#4,#5/{\ifps@cri{\Figg@tXYa{#5}\v@lXa=#4\v@lXa\v@lYa=#4\v@lYa%
    \Figg@tXY{#3}\advance\v@lX\v@lXa\advance\v@lY\v@lYa%
    \Figp@intregDD#1:{#2}(\v@lX,\v@lY)}\ignorespaces\fi}
\ctr@ld@f\def\figpttraTD#1:#2=#3/#4,#5/{\ifps@cri{\Figg@tXYa{#5}\v@lXa=#4\v@lXa\v@lYa=#4\v@lYa%
    \v@lZa=#4\v@lZa\Figg@tXY{#3}\advance\v@lX\v@lXa\advance\v@lY\v@lYa%
    \advance\v@lZ\v@lZa\Figp@intregTD#1:{#2}(\v@lX,\v@lY,\v@lZ)}\ignorespaces\fi}
\ctr@ln@m\figpttraC
\ctr@ld@f\def\figpttraCDD#1:#2=#3/#4,#5/{\ifps@cri{\v@lXa=#4\unit@\v@lYa=#5\unit@%
    \Figg@tXY{#3}\advance\v@lX\v@lXa\advance\v@lY\v@lYa%
    \Figp@intregDD#1:{#2}(\v@lX,\v@lY)}\ignorespaces\fi}
\ctr@ld@f\def\figpttraCTD#1:#2=#3/#4,#5,#6/{\ifps@cri{\v@lXa=#4\unit@\v@lYa=#5\unit@\v@lZa=#6\unit@%
    \Figg@tXY{#3}\advance\v@lX\v@lXa\advance\v@lY\v@lYa\advance\v@lZ\v@lZa%
    \Figp@intregTD#1:{#2}(\v@lX,\v@lY,\v@lZ)}\ignorespaces\fi}
\ctr@ld@f\def\figptsaxes#1:#2(#3){\ifps@cri{\an@lys@xes#3,:\ifx\t@xt@\empty%
    \ifTr@isDim\Figpts@xes#1:#2(0,#3,0,#3,0,#3)\else\Figpts@xes#1:#2(0,#3,0,#3)\fi%
    \else\Figpts@xes#1:#2(#3)\fi}\ignorespaces\fi}
\ctr@ln@m\Figpts@xes
\ctr@ld@f\def\Figpts@xesDD#1:#2(#3,#4,#5,#6){%
    \s@mme=#1\figpttraC\the\s@mme:$x$=#2/#4,0/%
    \advance\s@mme\@ne\figpttraC\the\s@mme:$y$=#2/0,#6/}
\ctr@ld@f\def\Figpts@xesTD#1:#2(#3,#4,#5,#6,#7,#8){%
    \s@mme=#1\figpttraC\the\s@mme:$x$=#2/#4,0,0/%
    \advance\s@mme\@ne\figpttraC\the\s@mme:$y$=#2/0,#6,0/%
    \advance\s@mme\@ne\figpttraC\the\s@mme:$z$=#2/0,0,#8/}
\ctr@ld@f\def\figptsmap#1=#2/#3/#4/{\ifps@cri{\s@uvc@ntr@l\et@tfigptsmap%
    \setc@ntr@l{2}\def\list@num{#2}\s@mme=#1%
    \@ecfor\p@int:=\list@num\do{\figvectP-1[#3,\p@int]\Figg@tXY{-1}%
    \pr@dMatV/#4/\figpttra\the\s@mme:=#3/1,-1/\advance\s@mme\@ne}%
    \resetc@ntr@l\et@tfigptsmap}\ignorespaces\fi}
\ctr@ln@m\figptscontrol
\ctr@ld@f\def\figptscontrolDD#1[#2,#3,#4,#5]{\ifps@cri{\s@uvc@ntr@l\et@tfigptscontrolDD\setc@ntr@l{2}%
    \v@lX=\z@\v@lY=\z@\Figtr@nptDD{-5}{#2}\Figtr@nptDD{2}{#5}%
    \divide\v@lX\@vi\divide\v@lY\@vi%
    \Figtr@nptDD{3}{#3}\Figtr@nptDD{-1.5}{#4}\Figp@intregDD-1:(\v@lX,\v@lY)%
    \v@lX=\z@\v@lY=\z@\Figtr@nptDD{2}{#2}\Figtr@nptDD{-5}{#5}%
    \divide\v@lX\@vi\divide\v@lY\@vi\Figtr@nptDD{-1.5}{#3}\Figtr@nptDD{3}{#4}%
    \s@mme=#1\advance\s@mme\@ne\Figp@intregDD\the\s@mme:(\v@lX,\v@lY)%
    \figptcopyDD#1:/-1/\resetc@ntr@l\et@tfigptscontrolDD}\ignorespaces\fi}
\ctr@ld@f\def\figptscontrolTD#1[#2,#3,#4,#5]{\ifps@cri{\s@uvc@ntr@l\et@tfigptscontrolTD\setc@ntr@l{2}%
    \v@lX=\z@\v@lY=\z@\v@lZ=\z@\Figtr@nptTD{-5}{#2}\Figtr@nptTD{2}{#5}%
    \divide\v@lX\@vi\divide\v@lY\@vi\divide\v@lZ\@vi%
    \Figtr@nptTD{3}{#3}\Figtr@nptTD{-1.5}{#4}\Figp@intregTD-1:(\v@lX,\v@lY,\v@lZ)%
    \v@lX=\z@\v@lY=\z@\v@lZ=\z@\Figtr@nptTD{2}{#2}\Figtr@nptTD{-5}{#5}%
    \divide\v@lX\@vi\divide\v@lY\@vi\divide\v@lZ\@vi\Figtr@nptTD{-1.5}{#3}\Figtr@nptTD{3}{#4}%
    \s@mme=#1\advance\s@mme\@ne\Figp@intregTD\the\s@mme:(\v@lX,\v@lY,\v@lZ)%
    \figptcopyTD#1:/-1/\resetc@ntr@l\et@tfigptscontrolTD}\ignorespaces\fi}
\ctr@ld@f\def\Figtr@nptDD#1#2{\Figg@tXYa{#2}\v@lXa=#1\v@lXa\v@lYa=#1\v@lYa%
    \advance\v@lX\v@lXa\advance\v@lY\v@lYa}
\ctr@ld@f\def\Figtr@nptTD#1#2{\Figg@tXYa{#2}\v@lXa=#1\v@lXa\v@lYa=#1\v@lYa\v@lZa=#1\v@lZa%
    \advance\v@lX\v@lXa\advance\v@lY\v@lYa\advance\v@lZ\v@lZa}
\ctr@ld@f\def\figptscontrolcurve#1,#2[#3]{\ifps@cri{\s@uvc@ntr@l\et@tfigptscontrolcurve%
    \def\list@num{#3}\extrairelepremi@r\Ak@\de\list@num%
    \extrairelepremi@r\Ai@\de\list@num\extrairelepremi@r\Aj@\de\list@num%
    \s@mme=#1\figptcopy\the\s@mme:/\Ai@/%
    \setc@ntr@l{2}\figvectP -1[\Ak@,\Aj@]%
    \@ecfor\Ak@:=\list@num\do{\advance\s@mme\@ne\figpttra\the\s@mme:=\Ai@/\curv@roundness,-1/%
       \figvectP -1[\Ai@,\Ak@]\advance\s@mme\@ne\figpttra\the\s@mme:=\Aj@/-\curv@roundness,-1/%
       \advance\s@mme\@ne\figptcopy\the\s@mme:/\Aj@/%
       \edef\Ai@{\Aj@}\edef\Aj@{\Ak@}}\advance\s@mme-#1\divide\s@mme\thr@@%
       \xdef#2{\the\s@mme}%
    \resetc@ntr@l\et@tfigptscontrolcurve}\ignorespaces\fi}
\ctr@ln@m\figptsintercirc
\ctr@ld@f\def\figptsintercircDD#1[#2,#3;#4,#5]{\ifps@cri{\s@uvc@ntr@l\et@tfigptsintercircDD%
    \setc@ntr@l{2}\let\c@lNVintc=\c@lNVintcDD\Figptsintercirc@#1[#2,#3;#4,#5]%
    \resetc@ntr@l\et@tfigptsintercircDD}\ignorespaces\fi}
\ctr@ld@f\def\figptsintercircTD#1[#2,#3;#4,#5;#6]{\ifps@cri{\s@uvc@ntr@l\et@tfigptsintercircTD%
    \setc@ntr@l{2}\let\c@lNVintc=\c@lNVintcTD\vecunitC@TD[#2,#6]%
    \Figv@ctCreg-3(\v@lX,\v@lY,\v@lZ)\Figptsintercirc@#1[#2,#3;#4,#5]%
    \resetc@ntr@l\et@tfigptsintercircTD}\ignorespaces\fi}
\ctr@ld@f\def\Figptsintercirc@#1[#2,#3;#4,#5]{\figvectP-1[#2,#4]%
    \vecunit@{-1}{-1}\delt@=\result@t\f@ctech=\result@tent%
    \s@mme=#1\advance\s@mme\@ne\figptcopy#1:/#2/\figptcopy\the\s@mme:/#4/%
    \ifdim\delt@=\z@\else%
    \v@lmin=#3\unit@\v@lmax=#5\unit@\v@leur=\v@lmin\advance\v@leur\v@lmax%
    \ifdim\v@leur>\delt@%
    \v@leur=\v@lmin\advance\v@leur-\v@lmax\maxim@m{\v@leur}{\v@leur}{-\v@leur}%
    \ifdim\v@leur<\delt@%
    \divide\v@lmin\f@ctech\divide\v@lmax\f@ctech\divide\delt@\f@ctech%
    \v@lmin=\repdecn@mb{\v@lmin}\v@lmin\v@lmax=\repdecn@mb{\v@lmax}\v@lmax%
    \invers@{\v@leur}{\delt@}\advance\v@lmax-\v@lmin%
    \v@lmax=-\repdecn@mb{\v@leur}\v@lmax\advance\delt@\v@lmax\delt@=.5\delt@%
    \v@lmax=\delt@\multiply\v@lmax\f@ctech%
    \edef\t@ille{\repdecn@mb{\v@lmax}}\figpttra-2:=#2/\t@ille,-1/%
    \delt@=\repdecn@mb{\delt@}\delt@\advance\v@lmin-\delt@%
    \sqrt@{\v@leur}{\v@lmin}\multiply\v@leur\f@ctech\edef\t@ille{\repdecn@mb{\v@leur}}%
    \c@lNVintc\figpttra#1:=-2/-\t@ille,-1/\figpttra\the\s@mme:=-2/\t@ille,-1/\fi\fi\fi}
\ctr@ld@f\def\c@lNVintcDD{\Figg@tXY{-1}\Figv@ctCreg-1(-\v@lY,\v@lX)} 
\ctr@ld@f\def\c@lNVintcTD{{\Figg@tXY{-3}\v@lmin=\v@lX\v@lmax=\v@lY\v@leur=\v@lZ%
    \Figg@tXY{-1}\c@lprovec{-3}\vecunit@{-3}{-3}
    \Figg@tXY{-1}\v@lmin=\v@lX\v@lmax=\v@lY%
    \v@leur=\v@lZ\Figg@tXY{-3}\c@lprovec{-1}}} 
\ctr@ln@m\figptsinterlinell
\ctr@ld@f\def\figptsinterlinellDD#1[#2,#3,#4,#5;#6,#7]{\ifps@cri{\s@uvc@ntr@l\et@tfigptsinterlinellDD%
    \figptcopy#1:/#6/\s@mme=#1\advance\s@mme\@ne\figptcopy\the\s@mme:/#7/%
    \v@lmin=#3\unit@\v@lmax=#4\unit@
    \setc@ntr@l{2}\figptbaryDD-4:[#6,#7;1,1]\figptsrotDD-3=-4,#7/#2,-#5/
    \Figg@tXY{-3}\Figg@tXYa{#2}\advance\v@lX-\v@lXa\advance\v@lY-\v@lYa
    \figvectP-1[-3,-2]\Figg@tXYa{-1}\figvectP-3[-4,#7]\Figptsint@rLE{#1}
    \resetc@ntr@l\et@tfigptsinterlinellDD}\ignorespaces\fi}
\ctr@ld@f\def\figptsinterlinellP#1[#2,#3,#4;#5,#6]{\ifps@cri{\s@uvc@ntr@l\et@tfigptsinterlinellP%
    \figptcopy#1:/#5/\s@mme=#1\advance\s@mme\@ne\figptcopy\the\s@mme:/#6/\setc@ntr@l{2}%
    \figvectP-1[#2,#3]\vecunit@{-1}{-1}\v@lmin=\result@t
    \figvectP-2[#2,#4]\vecunit@{-2}{-2}\v@lmax=\result@t
    \figptbary-4:[#5,#6;1,1]
    \figvectP-3[#2,-4]\c@lproscal\v@lX[-3,-1]\c@lproscal\v@lY[-3,-2]
    \figvectP-3[-4,#6]\c@lproscal\v@lXa[-3,-1]\c@lproscal\v@lYa[-3,-2]
    \Figptsint@rLE{#1}\resetc@ntr@l\et@tfigptsinterlinellP}\ignorespaces\fi}
\ctr@ld@f\def\Figptsint@rLE#1{%
    \getredf@ctDD\f@ctech(\v@lmin,\v@lmax)%
    \getredf@ctDD\p@rtent(\v@lX,\v@lY)\ifnum\p@rtent>\f@ctech\f@ctech=\p@rtent\fi%
    \getredf@ctDD\p@rtent(\v@lXa,\v@lYa)\ifnum\p@rtent>\f@ctech\f@ctech=\p@rtent\fi%
    \divide\v@lmin\f@ctech\divide\v@lmax\f@ctech\divide\v@lX\f@ctech\divide\v@lY\f@ctech%
    \divide\v@lXa\f@ctech\divide\v@lYa\f@ctech%
    \c@rre=\repdecn@mb\v@lXa\v@lmax\mili@u=\repdecn@mb\v@lYa\v@lmin%
    \getredf@ctDD\f@ctech(\c@rre,\mili@u)%
    \c@rre=\repdecn@mb\v@lX\v@lmax\mili@u=\repdecn@mb\v@lY\v@lmin%
    \getredf@ctDD\p@rtent(\c@rre,\mili@u)\ifnum\p@rtent>\f@ctech\f@ctech=\p@rtent\fi%
    \divide\v@lmin\f@ctech\divide\v@lmax\f@ctech\divide\v@lX\f@ctech\divide\v@lY\f@ctech%
    \divide\v@lXa\f@ctech\divide\v@lYa\f@ctech%
    \v@lmin=\repdecn@mb{\v@lmin}\v@lmin\v@lmax=\repdecn@mb{\v@lmax}\v@lmax%
    \edef\G@xde{\repdecn@mb\v@lmin}\edef\P@xde{\repdecn@mb\v@lmax}%
    \c@rre=-\v@lmax\v@leur=\repdecn@mb\v@lY\v@lY\advance\c@rre\v@leur\c@rre=\G@xde\c@rre%
    \v@leur=\repdecn@mb\v@lX\v@lX\v@leur=\P@xde\v@leur\advance\c@rre\v@leur
    \v@lmin=\repdecn@mb\v@lYa\v@lmin\v@lmax=\repdecn@mb\v@lXa\v@lmax%
    \mili@u=\repdecn@mb\v@lX\v@lmax\advance\mili@u\repdecn@mb\v@lY\v@lmin
    \v@lmax=\repdecn@mb\v@lXa\v@lmax\advance\v@lmax\repdecn@mb\v@lYa\v@lmin
    \ifdim\v@lmax>\epsil@n%
    \maxim@m{\v@leur}{\c@rre}{-\c@rre}\maxim@m{\v@lmin}{\mili@u}{-\mili@u}%
    \maxim@m{\v@leur}{\v@leur}{\v@lmin}\maxim@m{\v@lmin}{\v@lmax}{-\v@lmax}%
    \maxim@m{\v@leur}{\v@leur}{\v@lmin}\p@rtentiere{\p@rtent}{\v@leur}\advance\p@rtent\@ne%
    \divide\c@rre\p@rtent\divide\mili@u\p@rtent\divide\v@lmax\p@rtent%
    \delt@=\repdecn@mb{\mili@u}\mili@u\v@leur=\repdecn@mb{\v@lmax}\c@rre%
    \advance\delt@-\v@leur\ifdim\delt@<\z@\else\sqrt@\delt@\delt@%
    \invers@\v@lmax\v@lmax\edef\Uns@rAp{\repdecn@mb\v@lmax}%
    \v@leur=-\mili@u\advance\v@leur-\delt@\v@leur=\Uns@rAp\v@leur%
    \edef\t@ille{\repdecn@mb\v@leur}\figpttra#1:=-4/\t@ille,-3/\s@mme=#1\advance\s@mme\@ne%
    \v@leur=-\mili@u\advance\v@leur\delt@\v@leur=\Uns@rAp\v@leur%
    \edef\t@ille{\repdecn@mb\v@leur}\figpttra\the\s@mme:=-4/\t@ille,-3/\fi\fi}
\ctr@ln@m\figptsorthoprojline
\ctr@ld@f\def\figptsorthoprojlineDD#1=#2/#3,#4/{\ifps@cri{\s@uvc@ntr@l\et@tfigptsorthoprojlineDD%
    \setc@ntr@l{2}\figvectPDD-3[#3,#4]\figvectNVDD-4[-3]\resetc@ntr@l{2}%
    \def\list@num{#2}\s@mme=#1\@ecfor\p@int:=\list@num\do{%
    \inters@cDD\the\s@mme:[\p@int,-4;#3,-3]\advance\s@mme\@ne}%
    \resetc@ntr@l\et@tfigptsorthoprojlineDD}\ignorespaces\fi}
\ctr@ld@f\def\figptsorthoprojlineTD#1=#2/#3,#4/{\ifps@cri{\s@uvc@ntr@l\et@tfigptsorthoprojlineTD%
    \setc@ntr@l{2}\figvectPTD-2[#3,#4]\vecunit@TD{-2}{-2}%
    \def\list@num{#2}\s@mme=#1\@ecfor\p@int:=\list@num\do{%
    \figvectPTD-1[#3,\p@int]\c@lproscalTD\v@leur[-1,-2]%
    \edef\v@lcoef{\repdecn@mb{\v@leur}}\figpttraTD\the\s@mme:=#3/\v@lcoef,-2/%
    \advance\s@mme\@ne}\resetc@ntr@l\et@tfigptsorthoprojlineTD}\ignorespaces\fi}
\ctr@ln@m\figptsorthoprojplane
\ctr@ld@f\def\figptsorthoprojplaneDD{\un@v@ilable{figptsorthoprojplane}}
\ctr@ld@f\def\figptsorthoprojplaneTD#1=#2/#3,#4/{\ifps@cri{\s@uvc@ntr@l\et@tfigptsorthoprojplane%
    \setc@ntr@l{2}\vecunit@TD{-2}{#4}%
    \def\list@num{#2}\s@mme=#1\@ecfor\p@int:=\list@num\do{\figvectPTD-1[\p@int,#3]%
    \c@lproscalTD\v@leur[-1,-2]\edef\v@lcoef{\repdecn@mb{\v@leur}}%
    \figpttraTD\the\s@mme:=\p@int/\v@lcoef,-2/\advance\s@mme\@ne}%
    \resetc@ntr@l\et@tfigptsorthoprojplane}\ignorespaces\fi}
\ctr@ld@f\def\figptshom#1=#2/#3,#4/{\ifps@cri{\s@uvc@ntr@l\et@tfigptshom%
    \setc@ntr@l{2}\def\list@num{#2}\s@mme=#1%
    \@ecfor\p@int:=\list@num\do{\figvectP-1[#3,\p@int]%
    \figpttra\the\s@mme:=#3/#4,-1/\advance\s@mme\@ne}%
    \resetc@ntr@l\et@tfigptshom}\ignorespaces\fi}
\ctr@ln@m\figptsrot
\ctr@ld@f\def\figptsrotDD#1=#2/#3,#4/{\ifps@cri{\s@uvc@ntr@l\et@tfigptsrotDD%
    \c@ssin{\C@}{\S@}{#4}\setc@ntr@l{2}\def\list@num{#2}\s@mme=#1%
    \@ecfor\p@int:=\list@num\do{\figvectPDD-1[#3,\p@int]\Figg@tXY{-1}%
    \v@lXa=\C@\v@lX\advance\v@lXa-\S@\v@lY%
    \v@lYa=\S@\v@lX\advance\v@lYa\C@\v@lY%
    \Figv@ctCreg-1(\v@lXa,\v@lYa)\figpttraDD\the\s@mme:=#3/1,-1/\advance\s@mme\@ne}%
    \resetc@ntr@l\et@tfigptsrotDD}\ignorespaces\fi}
\ctr@ld@f\def\figptsrotTD#1=#2/#3,#4,#5/{\ifps@cri{\s@uvc@ntr@l\et@tfigptsrotTD%
    \c@ssin{\C@}{\S@}{#4}%
    \setc@ntr@l{2}\def\list@num{#2}\s@mme=#1%
    \@ecfor\p@int:=\list@num\do{\figptorthoprojplaneTD-3:=#3/\p@int,#5/%
    \figvectPTD-2[-3,\p@int]%
    \figvectNVTD-1[#5,-2]\n@rmeucTD\v@leur{-2}\edef\v@lcoef{\repdecn@mb{\v@leur}}%
    \Figg@tXYa{-1}\v@lXa=\v@lcoef\v@lXa\v@lYa=\v@lcoef\v@lYa\v@lZa=\v@lcoef\v@lZa%
    \v@lXa=\S@\v@lXa\v@lYa=\S@\v@lYa\v@lZa=\S@\v@lZa\Figg@tXY{-2}%
    \advance\v@lXa\C@\v@lX\advance\v@lYa\C@\v@lY\advance\v@lZa\C@\v@lZ%
    \Figg@tXY{-3}\advance\v@lXa\v@lX\advance\v@lYa\v@lY\advance\v@lZa\v@lZ%
    \Figp@intregTD\the\s@mme:(\v@lXa,\v@lYa,\v@lZa)\advance\s@mme\@ne}%
    \resetc@ntr@l\et@tfigptsrotTD}\ignorespaces\fi}
\ctr@ln@m\figptssym
\ctr@ld@f\def\figptssymDD#1=#2/#3,#4/{\ifps@cri{\s@uvc@ntr@l\et@tfigptssymDD%
    \setc@ntr@l{2}\figvectPDD-3[#3,#4]\Figg@tXY{-3}\Figv@ctCreg-4(-\v@lY,\v@lX)%
    \resetc@ntr@l{2}\def\list@num{#2}\s@mme=#1%
    \@ecfor\p@int:=\list@num\do{\inters@cDD-5:[#3,-3;\p@int,-4]\figvectPDD-2[\p@int,-5]%
    \figpttraDD\the\s@mme:=\p@int/2,-2/\advance\s@mme\@ne}%
    \resetc@ntr@l\et@tfigptssymDD}\ignorespaces\fi}
\ctr@ld@f\def\figptssymTD#1=#2/#3,#4/{\ifps@cri{\s@uvc@ntr@l\et@tfigptssymTD%
    \setc@ntr@l{2}\vecunit@TD{-2}{#4}\def\list@num{#2}\s@mme=#1%
    \@ecfor\p@int:=\list@num\do{\figvectPTD-1[\p@int,#3]%
    \c@lproscalTD\v@leur[-1,-2]\v@leur=2\v@leur\edef\v@lcoef{\repdecn@mb{\v@leur}}%
    \figpttraTD\the\s@mme:=\p@int/\v@lcoef,-2/\advance\s@mme\@ne}%
    \resetc@ntr@l\et@tfigptssymTD}\ignorespaces\fi}
\ctr@ln@m\figptstra
\ctr@ld@f\def\figptstraDD#1=#2/#3,#4/{\ifps@cri{\Figg@tXYa{#4}\v@lXa=#3\v@lXa\v@lYa=#3\v@lYa%
    \def\list@num{#2}\s@mme=#1\@ecfor\p@int:=\list@num\do{\Figg@tXY{\p@int}%
    \advance\v@lX\v@lXa\advance\v@lY\v@lYa%
    \Figp@intregDD\the\s@mme:(\v@lX,\v@lY)\advance\s@mme\@ne}}\ignorespaces\fi}
\ctr@ld@f\def\figptstraTD#1=#2/#3,#4/{\ifps@cri{\Figg@tXYa{#4}\v@lXa=#3\v@lXa\v@lYa=#3\v@lYa%
    \v@lZa=#3\v@lZa\def\list@num{#2}\s@mme=#1\@ecfor\p@int:=\list@num\do{\Figg@tXY{\p@int}%
    \advance\v@lX\v@lXa\advance\v@lY\v@lYa\advance\v@lZ\v@lZa%
    \Figp@intregTD\the\s@mme:(\v@lX,\v@lY,\v@lZ)\advance\s@mme\@ne}}\ignorespaces\fi}
\ctr@ln@m\figptvisilimSL
\ctr@ld@f\def\figptvisilimSLDD{\un@v@ilable{figptvisilimSL}}
\ctr@ld@f\def\figptvisilimSLTD#1:#2[#3,#4;#5,#6]{\ifps@cri{\s@uvc@ntr@l\et@tfigptvisilimSLTD%
    \setc@ntr@l{2}\figvectP-1[#3,#4]\n@rminf{\delt@}{-1}%
    \ifcase\curr@ntproj\v@lX=\cxa@\p@\v@lY=-\p@\v@lZ=\cxb@\p@
    \Figv@ctCreg-2(\v@lX,\v@lY,\v@lZ)\figvectP-3[#5,#6]\figvectNV-1[-2,-3]%
    \or\figvectP-1[#5,#6]\vecunitCV@TD{-1}\v@lmin=\v@lX\v@lmax=\v@lY
    \v@leur=\v@lZ\v@lX=\cza@\p@\v@lY=\czb@\p@\v@lZ=\czc@\p@\c@lprovec{-1}%
    \or\c@ley@pt{-2}\figvectN-1[#5,#6,-2]\fi
    \edef\Ai@{#3}\edef\Aj@{#4}\figvectP-2[#5,\Ai@]\c@lproscal\v@leur[-1,-2]%
    \ifdim\v@leur>\z@\p@rtent=\@ne\else\p@rtent=\m@ne\fi%
    \figvectP-2[#5,\Aj@]\c@lproscal\v@leur[-1,-2]%
    \ifdim\p@rtent\v@leur>\z@\figptcopy#1:#2/#3/%
    \message{*** \BS@ figptvisilimSL: points are on the same side.}\else%
    \figptcopy-3:/#3/\figptcopy-4:/#4/%
    \loop\figptbary-5:[-3,-4;1,1]\figvectP-2[#5,-5]\c@lproscal\v@leur[-1,-2]%
    \ifdim\p@rtent\v@leur>\z@\figptcopy-3:/-5/\else\figptcopy-4:/-5/\fi%
    \divide\delt@\tw@\ifdim\delt@>\epsil@n\repeat%
    \figptbary#1:#2[-3,-4;1,1]\fi\resetc@ntr@l\et@tfigptvisilimSLTD}\ignorespaces\fi}
\ctr@ld@f\def\c@ley@pt#1{\t@stp@r\ifitis@K\v@lX=\cza@\p@\v@lY=\czb@\p@\v@lZ=\czc@\p@%
    \Figv@ctCreg-1(\v@lX,\v@lY,\v@lZ)\Figp@intreg-2:(\wd\Bt@rget,\ht\Bt@rget,\dp\Bt@rget)%
    \figpttra#1:=-2/-\disob@intern,-1/\else\end\fi}
\ctr@ld@f\def\t@stp@r{\itis@Ktrue\ifnewt@rgetpt\else\itis@Kfalse%
    \message{*** \BS@ figptvisilimXX: target point undefined.}\fi\ifnewdis@b\else%
    \itis@Kfalse\message{*** \BS@ figptvisilimXX: observation distance undefined.}\fi%
    \ifitis@K\else\message{*** This macro must be called after \BS@ psbeginfig or after
    having set the missing parameter(s) with \BS@ figset proj()}\fi}
\ctr@ld@f\def\figscan#1(#2,#3){{\s@uvc@ntr@l\et@tfigscan\@psfgetbb{#1}\if@psfbbfound\else%
    \def\@psfllx{0}\def\@psflly{20}\def\@psfurx{540}\def\@psfury{640}\fi\figscan@{#2}{#3}%
    \resetc@ntr@l\et@tfigscan}\ignorespaces}
\ctr@ld@f\def\figscan@#1#2{%
    \unit@=\@ne bp\setc@ntr@l{2}\figsetmark{}%
    \def\minst@p{20pt}%
    \v@lX=\@psfllx\p@\v@lX=\Sc@leFact\v@lX\r@undint\v@lX\v@lX%
    \v@lY=\@psflly\p@\v@lY=\Sc@leFact\v@lY\ifdim\v@lY>\z@\r@undint\v@lY\v@lY\fi%
    \delt@=\@psfury\p@\delt@=\Sc@leFact\delt@%
    \advance\delt@-\v@lY\v@lXa=\@psfurx\p@\v@lXa=\Sc@leFact\v@lXa\v@leur=\minst@p%
    \edef\valv@lY{\repdecn@mb{\v@lY}}\edef\LgTr@it{\the\delt@}%
    \loop\ifdim\v@lX<\v@lXa\edef\valv@lX{\repdecn@mb{\v@lX}}%
    \figptDD -1:(\valv@lX,\valv@lY)\figwriten -1:\hbox{\vrule height\LgTr@it}(0)%
    \ifdim\v@leur<\minst@p\else\figsetmark{\raise-8bp\hbox{$\scriptscriptstyle\triangle$}}%
    \figwrites -1:\@ffichnb{0}{\valv@lX}(6)\v@leur=\z@\figsetmark{}\fi%
    \advance\v@leur#1pt\advance\v@lX#1pt\repeat%
    \def\minst@p{10pt}%
    \v@lX=\@psfllx\p@\v@lX=\Sc@leFact\v@lX\ifdim\v@lX>\z@\r@undint\v@lX\v@lX\fi%
    \v@lY=\@psflly\p@\v@lY=\Sc@leFact\v@lY\r@undint\v@lY\v@lY%
    \delt@=\@psfurx\p@\delt@=\Sc@leFact\delt@%
    \advance\delt@-\v@lX\v@lYa=\@psfury\p@\v@lYa=\Sc@leFact\v@lYa\v@leur=\minst@p%
    \edef\valv@lX{\repdecn@mb{\v@lX}}\edef\LgTr@it{\the\delt@}%
    \loop\ifdim\v@lY<\v@lYa\edef\valv@lY{\repdecn@mb{\v@lY}}%
    \figptDD -1:(\valv@lX,\valv@lY)\figwritee -1:\vbox{\hrule width\LgTr@it}(0)%
    \ifdim\v@leur<\minst@p\else\figsetmark{$\triangleright$\kern4bp}%
    \figwritew -1:\@ffichnb{0}{\valv@lY}(6)\v@leur=\z@\figsetmark{}\fi%
    \advance\v@leur#2pt\advance\v@lY#2pt\repeat}
\ctr@ld@f\let\figscanI=\figscan
\ctr@ld@f\def\figscan@E#1(#2,#3){{\s@uvc@ntr@l\et@tfigscan@E%
    \Figdisc@rdLTS{#1}{\t@xt@}\pdfximage{\t@xt@}%
    \setbox\Gb@x=\hbox{\pdfrefximage\pdflastximage}%
    \edef\@psfllx{0}\v@lY=-\dp\Gb@x\edef\@psflly{\repdecn@mb{\v@lY}}%
    \edef\@psfurx{\repdecn@mb{\wd\Gb@x}}%
    \v@lY=\dp\Gb@x\advance\v@lY\ht\Gb@x\edef\@psfury{\repdecn@mb{\v@lY}}%
    \figscan@{#2}{#3}\resetc@ntr@l\et@tfigscan@E}\ignorespaces}
\ctr@ld@f\def\figshowpts[#1,#2]{{\figsetmark{$\bullet$}\figsetptname{\bf ##1}%
    \p@rtent=#2\relax\ifnum\p@rtent<\z@\p@rtent=\z@\fi%
    \s@mme=#1\relax\ifnum\s@mme<\z@\s@mme=\z@\fi%
    \loop\ifnum\s@mme<\p@rtent\pt@rvect{\s@mme}%
    \ifitis@K\figwriten{\the\s@mme}:(4pt)\fi\advance\s@mme\@ne\repeat%
    \pt@rvect{\s@mme}\ifitis@K\figwriten{\the\s@mme}:(4pt)\fi}\ignorespaces}
\ctr@ld@f\def\pt@rvect#1{\set@bjc@de{#1}%
    \expandafter\expandafter\expandafter\inqpt@rvec\csname\objc@de\endcsname:}
\ctr@ld@f\def\inqpt@rvec#1#2:{\if#1\C@dCl@spt\itis@Ktrue\else\itis@Kfalse\fi}
\ctr@ld@f\def\figshowsettings{{%
    \immediate\write16{====================================================================}%
    \immediate\write16{ Current settings about:}%
    \immediate\write16{ --- GENERAL ---}%
    \immediate\write16{Scale factor and Unit = \unit@util\space (\the\unit@)
     \space -> \BS@ figinit{ScaleFactorUnit}}%
    \immediate\write16{Update mode = \ifpsupdatem@de yes\else no\fi
     \space-> \BS@ psset(update=yes/no) or \BS@ pssetdefault(update=yes/no)}%
    \immediate\write16{ --- PRINTING ---}%
    \immediate\write16{Implicit point name = \ptn@me{i} \space-> \BS@ figsetptname{Name}}%
    \immediate\write16{Point marker = \the\c@nsymb \space -> \BS@ figsetmark{Mark}}%
    \immediate\write16{Print rounded coordinates = \ifr@undcoord yes\else no\fi
     \space-> \BS@ figsetroundcoord{yes/no}}%
    \immediate\write16{ --- GRAPHICAL (general) ---}%
    \immediate\write16{First-level (or primary) settings:}%
    \immediate\write16{ Color = \curr@ntcolor \space-> \BS@ psset(color=ColorDefinition)}%
    \immediate\write16{ Filling mode = \iffillm@de yes\else no\fi
     \space-> \BS@ psset(fillmode=yes/no)}%
    \immediate\write16{ Line join = \curr@ntjoin \space-> \BS@ psset(join=miter/round/bevel)}%
    \immediate\write16{ Line style = \curr@ntdash \space-> \BS@ psset(dash=Index/Pattern)}%
    \immediate\write16{ Line width = \curr@ntwidth
     \space-> \BS@ psset(width=real in PostScript units)}%
    \immediate\write16{Second-level (or secondary) settings:}%
    \immediate\write16{ Color = \sec@ndcolor \space-> \BS@ psset second(color=ColorDefinition)}%
    \immediate\write16{ Line style = \curr@ntseconddash
     \space-> \BS@ psset second(dash=Index/Pattern)}%
    \immediate\write16{ Line width = \curr@ntsecondwidth
     \space-> \BS@ psset second(width=real in PostScript units)}%
    \immediate\write16{Third-level (or ternary) settings:}%
    \immediate\write16{ Color = \th@rdcolor \space-> \BS@ psset third(color=ColorDefinition)}%
    \immediate\write16{ --- GRAPHICAL (specific) ---}%
    \immediate\write16{Arrow-head:}%
    \immediate\write16{ (half-)Angle = \@rrowheadangle
     \space-> \BS@ psset arrowhead(angle=real in degrees)}%
    \immediate\write16{ Filling mode = \if@rrowhfill yes\else no\fi
     \space-> \BS@ psset arrowhead(fillmode=yes/no)}%
    \immediate\write16{ "Outside" = \if@rrowhout yes\else no\fi
     \space-> \BS@ psset arrowhead(out=yes/no)}%
    \immediate\write16{ Length = \@rrowheadlength
     \if@rrowratio\space(not active)\else\space(active)\fi
     \space-> \BS@ psset arrowhead(length=real in user coord.)}%
    \immediate\write16{ Ratio = \@rrowheadratio
     \if@rrowratio\space(active)\else\space(not active)\fi
     \space-> \BS@ psset arrowhead(ratio=real in [0,1])}%
    \immediate\write16{Curve: Roundness = \curv@roundness
     \space-> \BS@ psset curve(roundness=real in [0,0.5])}%
    \immediate\write16{Mesh: Diagonal = \c@ntrolmesh
     \space-> \BS@ psset mesh(diag=integer in {-1,0,1})}%
    \immediate\write16{Flow chart:}%
    \immediate\write16{ Arrow position = \@rrowp@s
     \space-> \BS@ psset flowchart(arrowposition=real in [0,1])}%
    \immediate\write16{ Arrow reference point = \ifcase\@rrowr@fpt start\else end\fi
     \space-> \BS@ psset flowchart(arrowrefpt = start/end)}%
    \immediate\write16{ Line type = \ifcase\fclin@typ@ curve\else polygon\fi
     \space-> \BS@ psset flowchart(line=polygon/curve)}%
    \immediate\write16{ Padding = (\Xp@dd, \Yp@dd)
     \space-> \BS@ psset flowchart(padding = real in user coord.)}%
    \immediate\write16{\space\space\space\space(or
     \BS@ psset flowchart(xpadding=real, ypadding=real) )}%
    \immediate\write16{ Radius = \fclin@r@d
     \space-> \BS@ psset flowchart(radius=positive real in user coord.)}%
    \immediate\write16{ Shape = \fcsh@pe
     \space-> \BS@ psset flowchart(shape = rectangle, ellipse or lozenge)}%
    \immediate\write16{ Thickness = \thickn@ss
     \space-> \BS@ psset flowchart(thickness = real in user coord.)}%
    \ifTr@isDim%
    \immediate\write16{ --- 3D to 2D PROJECTION ---}%
    \immediate\write16{Projection : \typ@proj \space-> \BS@ figinit{ScaleFactorUnit, ProjType}}%
    \immediate\write16{Longitude (psi) = \v@lPsi \space-> \BS@ figset proj(psi=real in degrees)}%
    \ifcase\curr@ntproj\immediate\write16{Depth coeff. (Lambda)
     \space = \v@lTheta \space-> \BS@ figset proj(lambda=real in [0,1])}%
    \else\immediate\write16{Latitude (theta)
     \space = \v@lTheta \space-> \BS@ figset proj(theta=real in degrees)}%
    \fi%
    \ifnum\curr@ntproj=\tw@%
    \immediate\write16{Observation distance = \disob@unit
     \space-> \BS@ figset proj(dist=real in user coord.)}%
    \immediate\write16{Target point = \t@rgetpt \space-> \BS@ figset proj(targetpt=pt number)}%
     \v@lX=\ptT@unit@\wd\Bt@rget\v@lY=\ptT@unit@\ht\Bt@rget\v@lZ=\ptT@unit@\dp\Bt@rget%
    \immediate\write16{ Its coordinates are
     (\repdecn@mb{\v@lX}, \repdecn@mb{\v@lY}, \repdecn@mb{\v@lZ})}%
    \fi%
    \fi%
    \immediate\write16{====================================================================}%
    \ignorespaces}}
\ctr@ln@w{newif}\ifitis@vect@r
\ctr@ld@f\def\figvectC#1(#2,#3){{\itis@vect@rtrue\figpt#1:(#2,#3)}\ignorespaces}
\ctr@ld@f\def\Figv@ctCreg#1(#2,#3){{\itis@vect@rtrue\Figp@intreg#1:(#2,#3)}\ignorespaces}
\ctr@ln@m\figvectDBezier
\ctr@ld@f\def\figvectDBezierDD#1:#2,#3[#4,#5,#6,#7]{\ifps@cri{\s@uvc@ntr@l\et@tfigvectDBezierDD%
    \FigvectDBezier@#2,#3[#4,#5,#6,#7]\v@lX=\c@ef\v@lX\v@lY=\c@ef\v@lY%
    \Figv@ctCreg#1(\v@lX,\v@lY)\resetc@ntr@l\et@tfigvectDBezierDD}\ignorespaces\fi}
\ctr@ld@f\def\figvectDBezierTD#1:#2,#3[#4,#5,#6,#7]{\ifps@cri{\s@uvc@ntr@l\et@tfigvectDBezierTD%
    \FigvectDBezier@#2,#3[#4,#5,#6,#7]\v@lX=\c@ef\v@lX\v@lY=\c@ef\v@lY\v@lZ=\c@ef\v@lZ%
    \Figv@ctCreg#1(\v@lX,\v@lY,\v@lZ)\resetc@ntr@l\et@tfigvectDBezierTD}\ignorespaces\fi}
\ctr@ld@f\def\FigvectDBezier@#1,#2[#3,#4,#5,#6]{\setc@ntr@l{2}%
    \edef\T@{#2}\v@leur=\p@\advance\v@leur-#2pt\edef\UNmT@{\repdecn@mb{\v@leur}}%
    \ifnum#1=\tw@\def\c@ef{6}\else\def\c@ef{3}\fi%
    \figptcopy-4:/#3/\figptcopy-3:/#4/\figptcopy-2:/#5/\figptcopy-1:/#6/%
    \l@mbd@un=-4 \l@mbd@de=-\thr@@\p@rtent=\m@ne\c@lDecast%
    \ifnum#1=\tw@\c@lDCDeux{-4}{-3}\c@lDCDeux{-3}{-2}\c@lDCDeux{-4}{-3}\else%
    \l@mbd@un=-4 \l@mbd@de=-\thr@@\p@rtent=-\tw@\c@lDecast%
    \c@lDCDeux{-4}{-3}\fi\Figg@tXY{-4}}
\ctr@ln@m\c@lDCDeux
\ctr@ld@f\def\c@lDCDeuxDD#1#2{\Figg@tXY{#2}\Figg@tXYa{#1}%
    \advance\v@lX-\v@lXa\advance\v@lY-\v@lYa\Figp@intregDD#1:(\v@lX,\v@lY)}
\ctr@ld@f\def\c@lDCDeuxTD#1#2{\Figg@tXY{#2}\Figg@tXYa{#1}\advance\v@lX-\v@lXa%
    \advance\v@lY-\v@lYa\advance\v@lZ-\v@lZa\Figp@intregTD#1:(\v@lX,\v@lY,\v@lZ)}
\ctr@ln@m\figvectN
\ctr@ld@f\def\figvectNDD#1[#2,#3]{\ifps@cri{\Figg@tXYa{#2}\Figg@tXY{#3}%
    \advance\v@lX-\v@lXa\advance\v@lY-\v@lYa%
    \Figv@ctCreg#1(-\v@lY,\v@lX)}\ignorespaces\fi}
\ctr@ld@f\def\figvectNTD#1[#2,#3,#4]{\ifps@cri{\vecunitC@TD[#2,#4]\v@lmin=\v@lX\v@lmax=\v@lY%
    \v@leur=\v@lZ\vecunitC@TD[#2,#3]\c@lprovec{#1}}\ignorespaces\fi}
\ctr@ln@m\figvectNV
\ctr@ld@f\def\figvectNVDD#1[#2]{\ifps@cri{\Figg@tXY{#2}\Figv@ctCreg#1(-\v@lY,\v@lX)}\ignorespaces\fi}
\ctr@ld@f\def\figvectNVTD#1[#2,#3]{\ifps@cri{\vecunitCV@TD{#3}\v@lmin=\v@lX\v@lmax=\v@lY%
    \v@leur=\v@lZ\vecunitCV@TD{#2}\c@lprovec{#1}}\ignorespaces\fi}
\ctr@ln@m\figvectP
\ctr@ld@f\def\figvectPDD#1[#2,#3]{\ifps@cri{\Figg@tXYa{#2}\Figg@tXY{#3}%
    \advance\v@lX-\v@lXa\advance\v@lY-\v@lYa%
    \Figv@ctCreg#1(\v@lX,\v@lY)}\ignorespaces\fi}
\ctr@ld@f\def\figvectPTD#1[#2,#3]{\ifps@cri{\Figg@tXYa{#2}\Figg@tXY{#3}%
    \advance\v@lX-\v@lXa\advance\v@lY-\v@lYa\advance\v@lZ-\v@lZa%
    \Figv@ctCreg#1(\v@lX,\v@lY,\v@lZ)}\ignorespaces\fi}
\ctr@ln@m\figvectU
\ctr@ld@f\def\figvectUDD#1[#2]{\ifps@cri{\n@rmeuc\v@leur{#2}\invers@\v@leur\v@leur%
    \delt@=\repdecn@mb{\v@leur}\unit@\edef\v@ldelt@{\repdecn@mb{\delt@}}%
    \Figg@tXY{#2}\v@lX=\v@ldelt@\v@lX\v@lY=\v@ldelt@\v@lY%
    \Figv@ctCreg#1(\v@lX,\v@lY)}\ignorespaces\fi}
\ctr@ld@f\def\figvectUTD#1[#2]{\ifps@cri{\n@rmeuc\v@leur{#2}\invers@\v@leur\v@leur%
    \delt@=\repdecn@mb{\v@leur}\unit@\edef\v@ldelt@{\repdecn@mb{\delt@}}%
    \Figg@tXY{#2}\v@lX=\v@ldelt@\v@lX\v@lY=\v@ldelt@\v@lY\v@lZ=\v@ldelt@\v@lZ%
    \Figv@ctCreg#1(\v@lX,\v@lY,\v@lZ)}\ignorespaces\fi}
\ctr@ld@f\def\figvisu#1#2#3{\c@ldefproj\initb@undb@x\xdef\figforTeXFigno{\figforTeXnextFigno}%
    \s@mme=\figforTeXnextFigno\advance\s@mme\@ne\xdef\figforTeXnextFigno{\number\s@mme}%
    \setbox\b@xvisu=\hbox{\ifnum\@utoFN>\z@\figinsert{}\gdef\@utoFInDone{0}\fi\ignorespaces#3}%
    \gdef\@utoFInDone{1}\gdef\@utoFN{0}%
    \v@lXa=-\c@@rdYmin\v@lYa=\c@@rdYmax\advance\v@lYa-\c@@rdYmin%
    \v@lX=\c@@rdXmax\advance\v@lX-\c@@rdXmin%
    \setbox#1=\hbox{#2}\v@lY=-\v@lX\maxim@m{\v@lX}{\v@lX}{\wd#1}%
    \advance\v@lY\v@lX\divide\v@lY\tw@\advance\v@lY-\c@@rdXmin%
    \setbox#1=\vbox{\parindent0mm\hsize=\v@lX\vskip\v@lYa%
    \rlap{\hskip\v@lY\smash{\raise\v@lXa\box\b@xvisu}}%
    \def\t@xt@{#2}\ifx\t@xt@\empty\else\medskip\centerline{#2}\fi}\wd#1=\v@lX}
\ctr@ld@f\def\figDecrementFigno{{\xdef\figforTeXnextFigno{\figforTeXFigno}%
    \s@mme=\figforTeXFigno\advance\s@mme\m@ne\xdef\figforTeXFigno{\number\s@mme}}}
\ctr@ln@w{newbox}\Bt@rget\setbox\Bt@rget=\null
\ctr@ln@w{newbox}\BminTD@\setbox\BminTD@=\null
\ctr@ln@w{newbox}\BmaxTD@\setbox\BmaxTD@=\null
\ctr@ln@w{newif}\ifnewt@rgetpt\ctr@ln@w{newif}\ifnewdis@b
\ctr@ld@f\def\b@undb@xTD#1#2#3{%
    \relax\ifdim#1<\wd\BminTD@\global\wd\BminTD@=#1\fi%
    \relax\ifdim#2<\ht\BminTD@\global\ht\BminTD@=#2\fi%
    \relax\ifdim#3<\dp\BminTD@\global\dp\BminTD@=#3\fi%
    \relax\ifdim#1>\wd\BmaxTD@\global\wd\BmaxTD@=#1\fi%
    \relax\ifdim#2>\ht\BmaxTD@\global\ht\BmaxTD@=#2\fi%
    \relax\ifdim#3>\dp\BmaxTD@\global\dp\BmaxTD@=#3\fi}
\ctr@ld@f\def\c@ldefdisob{{\ifdim\wd\BminTD@<\maxdimen\v@leur=\wd\BmaxTD@\advance\v@leur-\wd\BminTD@%
    \delt@=\ht\BmaxTD@\advance\delt@-\ht\BminTD@\maxim@m{\v@leur}{\v@leur}{\delt@}%
    \delt@=\dp\BmaxTD@\advance\delt@-\dp\BminTD@\maxim@m{\v@leur}{\v@leur}{\delt@}%
    \v@leur=5\v@leur\else\v@leur=800pt\fi\c@ldefdisob@{\v@leur}}}
\ctr@ln@m\disob@intern
\ctr@ln@m\disob@
\ctr@ln@m\divf@ctproj
\ctr@ld@f\def\c@ldefdisob@#1{{\v@leur=#1\ifdim\v@leur<\p@\v@leur=800pt\fi%
    \xdef\disob@intern{\repdecn@mb{\v@leur}}%
    \delt@=\ptT@unit@\v@leur\xdef\disob@unit{\repdecn@mb{\delt@}}%
    \f@ctech=\@ne\loop\ifdim\v@leur>\t@n pt\divide\v@leur\t@n\multiply\f@ctech\t@n\repeat%
    \xdef\disob@{\repdecn@mb{\v@leur}}\xdef\divf@ctproj{\the\f@ctech}}%
    \global\newdis@btrue}
\ctr@ln@m\t@rgetpt
\ctr@ld@f\def\c@ldeft@rgetpt{\newt@rgetpttrue\def\t@rgetpt{CenterBoundBox}{%
    \delt@=\wd\BmaxTD@\advance\delt@-\wd\BminTD@\divide\delt@\tw@%
    \v@leur=\wd\BminTD@\advance\v@leur\delt@\global\wd\Bt@rget=\v@leur%
    \delt@=\ht\BmaxTD@\advance\delt@-\ht\BminTD@\divide\delt@\tw@%
    \v@leur=\ht\BminTD@\advance\v@leur\delt@\global\ht\Bt@rget=\v@leur%
    \delt@=\dp\BmaxTD@\advance\delt@-\dp\BminTD@\divide\delt@\tw@%
    \v@leur=\dp\BminTD@\advance\v@leur\delt@\global\dp\Bt@rget=\v@leur}}
\ctr@ln@m\c@ldefproj
\ctr@ld@f\def\c@ldefprojTD{\ifnewt@rgetpt\else\c@ldeft@rgetpt\fi\ifnewdis@b\else\c@ldefdisob\fi}
\ctr@ld@f\def\c@lprojcav{
    \v@lZa=\cxa@\v@lY\advance\v@lX\v@lZa%
    \v@lZa=\cxb@\v@lY\v@lY=\v@lZ\advance\v@lY\v@lZa\ignorespaces}
\ctr@ln@m\v@lcoef
\ctr@ld@f\def\c@lprojrea{
    \advance\v@lX-\wd\Bt@rget\advance\v@lY-\ht\Bt@rget\advance\v@lZ-\dp\Bt@rget%
    \v@lZa=\cza@\v@lX\advance\v@lZa\czb@\v@lY\advance\v@lZa\czc@\v@lZ%
    \divide\v@lZa\divf@ctproj\advance\v@lZa\disob@ pt\invers@{\v@lZa}{\v@lZa}%
    \v@lZa=\disob@\v@lZa\edef\v@lcoef{\repdecn@mb{\v@lZa}}%
    \v@lXa=\cxa@\v@lX\advance\v@lXa\cxb@\v@lY\v@lXa=\v@lcoef\v@lXa%
    \v@lY=\cyb@\v@lY\advance\v@lY\cya@\v@lX\advance\v@lY\cyc@\v@lZ%
    \v@lY=\v@lcoef\v@lY\v@lX=\v@lXa\ignorespaces}
\ctr@ld@f\def\c@lprojort{
    \v@lXa=\cxa@\v@lX\advance\v@lXa\cxb@\v@lY%
    \v@lY=\cyb@\v@lY\advance\v@lY\cya@\v@lX\advance\v@lY\cyc@\v@lZ%
    \v@lX=\v@lXa\ignorespaces}
\ctr@ld@f\def\Figptpr@j#1:#2/#3/{{\Figg@tXY{#3}\superc@lprojSP%
    \Figp@intregDD#1:{#2}(\v@lX,\v@lY)}\ignorespaces}
\ctr@ln@m\figsetobdist
\ctr@ld@f\def\figsetobdistDD{\un@v@ilable{figsetobdist}}
\ctr@ld@f\def\figsetobdistTD(#1){{\ifcurr@ntPS%
    \immediate\write16{*** \BS@ figsetobdist is ignored inside a
     \BS@ psbeginfig-\BS@ psendfig block.}%
    \else\v@leur=#1\unit@\c@ldefdisob@{\v@leur}\fi}\ignorespaces}
\ctr@ln@m\c@lprojSP
\ctr@ln@m\curr@ntproj
\ctr@ln@m\typ@proj
\ctr@ln@m\superc@lprojSP
\ctr@ld@f\def\Figs@tproj#1{%
    \if#13 \d@faultproj\else\if#1c\d@faultproj%
    \else\if#1o\xdef\curr@ntproj{1}\xdef\typ@proj{orthogonal}%
         \figsetviewTD(\def@ultpsi,\def@ulttheta)%
         \global\let\c@lprojSP=\c@lprojort\global\let\superc@lprojSP=\c@lprojort%
    \else\if#1r\xdef\curr@ntproj{2}\xdef\typ@proj{realistic}%
         \figsetviewTD(\def@ultpsi,\def@ulttheta)%
         \global\let\c@lprojSP=\c@lprojrea\global\let\superc@lprojSP=\c@lprojrea%
    \else\d@faultproj\message{*** Unknown projection. Cavalier projection assumed.}%
    \fi\fi\fi\fi}
\ctr@ld@f\def\d@faultproj{\xdef\curr@ntproj{0}\xdef\typ@proj{cavalier}\figsetviewTD(\def@ultpsi,0.5)%
         \global\let\c@lprojSP=\c@lprojcav\global\let\superc@lprojSP=\c@lprojcav}
\ctr@ln@m\figsettarget
\ctr@ld@f\def\figsettargetDD{\un@v@ilable{figsettarget}}
\ctr@ld@f\def\figsettargetTD[#1]{{\ifcurr@ntPS%
    \immediate\write16{*** \BS@ figsettarget is ignored inside a
     \BS@ psbeginfig-\BS@ psendfig block.}%
    \else\global\newt@rgetpttrue\xdef\t@rgetpt{#1}\Figg@tXY{#1}\global\wd\Bt@rget=\v@lX%
    \global\ht\Bt@rget=\v@lY\global\dp\Bt@rget=\v@lZ\fi}\ignorespaces}
\ctr@ln@m\figsetview
\ctr@ld@f\def\figsetviewDD{\un@v@ilable{figsetview}}
\ctr@ld@f\def\figsetviewTD(#1){\ifcurr@ntPS%
     \immediate\write16{*** \BS@ figsetview is ignored inside a
     \BS@ psbeginfig-\BS@ psendfig block.}\else\Figsetview@#1,:\fi\ignorespaces}
\ctr@ld@f\def\Figsetview@#1,#2:{{\xdef\v@lPsi{#1}\def\t@xt@{#2}%
    \ifx\t@xt@\empty\def\@rgdeux{\v@lTheta}\else\X@rgdeux@#2\fi%
    \c@ssin{\costhet@}{\sinthet@}{#1}\v@lmin=\costhet@ pt\v@lmax=\sinthet@ pt%
    \ifcase\curr@ntproj%
    \v@leur=\@rgdeux\v@lmin\xdef\cxa@{\repdecn@mb{\v@leur}}%
    \v@leur=\@rgdeux\v@lmax\xdef\cxb@{\repdecn@mb{\v@leur}}\v@leur=\@rgdeux pt%
    \relax\ifdim\v@leur>\p@\message{*** Lambda too large ! See \BS@ figset proj() !}\fi%
    \else%
    \v@lmax=-\v@lmax\xdef\cxa@{\repdecn@mb{\v@lmax}}\xdef\cxb@{\costhet@}%
    \ifx\t@xt@\empty\edef\@rgdeux{\def@ulttheta}\fi\c@ssin{\C@}{\S@}{\@rgdeux}%
    \v@lmax=-\S@ pt%
    \v@leur=\v@lmax\v@leur=\costhet@\v@leur\xdef\cya@{\repdecn@mb{\v@leur}}%
    \v@leur=\v@lmax\v@leur=\sinthet@\v@leur\xdef\cyb@{\repdecn@mb{\v@leur}}%
    \xdef\cyc@{\C@}\v@lmin=-\C@ pt%
    \v@leur=\v@lmin\v@leur=\costhet@\v@leur\xdef\cza@{\repdecn@mb{\v@leur}}%
    \v@leur=\v@lmin\v@leur=\sinthet@\v@leur\xdef\czb@{\repdecn@mb{\v@leur}}%
    \xdef\czc@{\repdecn@mb{\v@lmax}}\fi%
    \xdef\v@lTheta{\@rgdeux}}}
\ctr@ld@f\def\def@ultpsi{40}
\ctr@ld@f\def\def@ulttheta{25}
\ctr@ln@m\l@debut
\ctr@ln@m\n@mref
\ctr@ld@f\def\figset#1(#2){\def\t@xt@{#1}\ifx\t@xt@\empty\trtlis@rg{#2}{\Figsetwr@te}
    \else\keln@mde#1|%
    \def\n@mref{pr}\ifx\l@debut\n@mref\ifcurr@ntPS
     \immediate\write16{*** \BS@ figset proj(...) is ignored inside a
     \BS@ psbeginfig-\BS@ psendfig block.}\else\trtlis@rg{#2}{\Figsetpr@j}\fi\else%
    \def\n@mref{wr}\ifx\l@debut\n@mref\trtlis@rg{#2}{\Figsetwr@te}\else
    \immediate\write16{*** Unknown keyword: \BS@ figset #1(...)}%
    \fi\fi\fi\ignorespaces}
\ctr@ld@f\def\Figsetpr@j#1=#2|{\keln@mtr#1|%
    \def\n@mref{dep}\ifx\l@debut\n@mref\Figsetd@p{#2}\else
    \def\n@mref{dis}\ifx\l@debut\n@mref%
     \ifnum\curr@ntproj=\tw@\figsetobdist(#2)\else\Figset@rr\fi\else
    \def\n@mref{lam}\ifx\l@debut\n@mref\Figsetd@p{#2}\else
    \def\n@mref{lat}\ifx\l@debut\n@mref\Figsetth@{#2}\else
    \def\n@mref{lon}\ifx\l@debut\n@mref\figsetview(#2)\else
    \def\n@mref{psi}\ifx\l@debut\n@mref\figsetview(#2)\else
    \def\n@mref{tar}\ifx\l@debut\n@mref%
     \ifnum\curr@ntproj=\tw@\figsettarget[#2]\else\Figset@rr\fi\else
    \def\n@mref{the}\ifx\l@debut\n@mref\Figsetth@{#2}\else
    \immediate\write16{*** Unknown attribute: \BS@ figset proj(..., #1=...).}%
    \fi\fi\fi\fi\fi\fi\fi\fi}
\ctr@ld@f\def\Figsetd@p#1{\ifnum\curr@ntproj=\z@\figsetview(\v@lPsi,#1)\else\Figset@rr\fi}
\ctr@ld@f\def\Figsetth@#1{\ifnum\curr@ntproj=\z@\Figset@rr\else\figsetview(\v@lPsi,#1)\fi}
\ctr@ld@f\def\Figset@rr{\message{*** \BS@ figset proj(): Attribute "\n@mref" ignored, incompatible
    with current projection}}
\ctr@ld@f\def\initb@undb@xTD{\wd\BminTD@=\maxdimen\ht\BminTD@=\maxdimen\dp\BminTD@=\maxdimen%
    \wd\BmaxTD@=-\maxdimen\ht\BmaxTD@=-\maxdimen\dp\BmaxTD@=-\maxdimen}
\ctr@ln@w{newbox}\Gb@x      
\ctr@ln@w{newbox}\Gb@xSC    
\ctr@ln@w{newtoks}\c@nsymb  
\ctr@ln@w{newif}\ifr@undcoord\ctr@ln@w{newif}\ifunitpr@sent
\ctr@ld@f\def\unssqrttw@{0.707106 }
\ctr@ld@f\def\figAst{\raise-1.15ex\hbox{$\ast$}}
\ctr@ld@f\def\figBullet{\raise-1.15ex\hbox{$\bullet$}}
\ctr@ld@f\def\figCirc{\raise-1.15ex\hbox{$\circ$}}
\ctr@ld@f\def\figDiamond{\raise-1.15ex\hbox{$\diamond$}}%
\ctr@ld@f\def\boxit#1#2{\leavevmode\hbox{\vrule\vbox{\hrule\vglue#1%
    \vtop{\hbox{\kern#1{#2}\kern#1}\vglue#1\hrule}}\vrule}}
\ctr@ld@f\def\centertext#1#2{\vbox{\hsize#1\parindent0cm%
    \leftskip=0pt plus 1fil\rightskip=0pt plus 1fil\parfillskip=0pt{#2}}}
\ctr@ld@f\def\lefttext#1#2{\vbox{\hsize#1\parindent0cm\rightskip=0pt plus 1fil#2}}
\ctr@ld@f\def\c@nterpt{\ignorespaces%
    \kern-.5\wd\Gb@xSC%
    \raise-.5\ht\Gb@xSC\rlap{\hbox{\raise.5\dp\Gb@xSC\hbox{\copy\Gb@xSC}}}%
    \kern .5\wd\Gb@xSC\ignorespaces}
\ctr@ld@f\def\b@undb@xSC#1#2{{\v@lXa=#1\v@lYa=#2%
    \v@leur=\ht\Gb@xSC\advance\v@leur\dp\Gb@xSC%
    \advance\v@lXa-.5\wd\Gb@xSC\advance\v@lYa-.5\v@leur\b@undb@x{\v@lXa}{\v@lYa}%
    \advance\v@lXa\wd\Gb@xSC\advance\v@lYa\v@leur\b@undb@x{\v@lXa}{\v@lYa}}}
\ctr@ln@m\Dist@n
\ctr@ln@m\l@suite
\ctr@ld@f\def\@keldist#1#2{\edef\Dist@n{#2}\y@tiunit{\Dist@n}%
    \ifunitpr@sent#1=\Dist@n\else#1=\Dist@n\unit@\fi}
\ctr@ld@f\def\y@tiunit#1{\unitpr@sentfalse\expandafter\y@tiunit@#1:}
\ctr@ld@f\def\y@tiunit@#1#2:{\ifcat#1a\unitpr@senttrue\else\def\l@suite{#2}%
    \ifx\l@suite\empty\else\y@tiunit@#2:\fi\fi}
\ctr@ln@m\figcoord
\ctr@ld@f\def\figcoordDD#1{{\v@lX=\ptT@unit@\v@lX\v@lY=\ptT@unit@\v@lY%
    \ifr@undcoord\ifcase#1\v@leur=0.5pt\or\v@leur=0.05pt\or\v@leur=0.005pt%
    \or\v@leur=0.0005pt\else\v@leur=\z@\fi%
    \ifdim\v@lX<\z@\advance\v@lX-\v@leur\else\advance\v@lX\v@leur\fi%
    \ifdim\v@lY<\z@\advance\v@lY-\v@leur\else\advance\v@lY\v@leur\fi\fi%
    (\@ffichnb{#1}{\repdecn@mb{\v@lX}},\ifmmode\else\thinspace\fi%
    \@ffichnb{#1}{\repdecn@mb{\v@lY}})}}
\ctr@ld@f\def\@ffichnb#1#2{{\def\@@ffich{\@ffich#1(}\edef\n@mbre{#2}%
    \expandafter\@@ffich\n@mbre)}}
\ctr@ld@f\def\@ffich#1(#2.#3){{#2\ifnum#1>\z@.\fi\def\dig@ts{#3}\s@mme=\z@%
    \loop\ifnum\s@mme<#1\expandafter\@ffichdec\dig@ts:\advance\s@mme\@ne\repeat}}
\ctr@ld@f\def\@ffichdec#1#2:{\relax#1\def\dig@ts{#20}}
\ctr@ld@f\def\figcoordTD#1{{\v@lX=\ptT@unit@\v@lX\v@lY=\ptT@unit@\v@lY\v@lZ=\ptT@unit@\v@lZ%
    \ifr@undcoord\ifcase#1\v@leur=0.5pt\or\v@leur=0.05pt\or\v@leur=0.005pt%
    \or\v@leur=0.0005pt\else\v@leur=\z@\fi%
    \ifdim\v@lX<\z@\advance\v@lX-\v@leur\else\advance\v@lX\v@leur\fi%
    \ifdim\v@lY<\z@\advance\v@lY-\v@leur\else\advance\v@lY\v@leur\fi%
    \ifdim\v@lZ<\z@\advance\v@lZ-\v@leur\else\advance\v@lZ\v@leur\fi\fi%
    (\@ffichnb{#1}{\repdecn@mb{\v@lX}},\ifmmode\else\thinspace\fi%
     \@ffichnb{#1}{\repdecn@mb{\v@lY}},\ifmmode\else\thinspace\fi%
     \@ffichnb{#1}{\repdecn@mb{\v@lZ}})}}
\ctr@ld@f\def\figsetroundcoord#1{\expandafter\Figsetr@undcoord#1:\ignorespaces}
\ctr@ld@f\def\Figsetr@undcoord#1#2:{\if#1n\r@undcoordfalse\else\r@undcoordtrue\fi}
\ctr@ld@f\def\Figsetwr@te#1=#2|{\keln@mun#1|%
    \def\n@mref{m}\ifx\l@debut\n@mref\figsetmark{#2}\else
    \immediate\write16{*** Unknown attribute: \BS@ figset (..., #1=...)}%
    \fi}
\ctr@ld@f\def\figsetmark#1{\c@nsymb={#1}\setbox\Gb@xSC=\hbox{\the\c@nsymb}\ignorespaces}
\ctr@ln@m\ptn@me
\ctr@ld@f\def\figsetptname#1{\def\ptn@me##1{#1}\ignorespaces}
\ctr@ld@f\def\FigWrit@L#1:#2(#3,#4){\ignorespaces\@keldist\v@leur{#3}\@keldist\delt@{#4}%
    \C@rp@r@m\def\list@num{#1}\@ecfor\p@int:=\list@num\do{\FigWrit@pt{\p@int}{#2}}}
\ctr@ld@f\def\FigWrit@pt#1#2{\FigWp@r@m{#1}{#2}\Vc@rrect\figWp@si%
    \ifdim\wd\Gb@xSC>\z@\b@undb@xSC{\v@lX}{\v@lY}\fi\figWBB@x}
\ctr@ld@f\def\FigWp@r@m#1#2{\Figg@tXY{#1}%
    \setbox\Gb@x=\hbox{\def\t@xt@{#2}\ifx\t@xt@\empty\Figg@tT{#1}\else#2\fi}\c@lprojSP}
\ctr@ld@f\let\Vc@rrect=\relax
\ctr@ld@f\let\C@rp@r@m=\relax
\ctr@ld@f\def\figwrite[#1]#2{{\ignorespaces\def\list@num{#1}\@ecfor\p@int:=\list@num\do{%
    \setbox\Gb@x=\hbox{\def\t@xt@{#2}\ifx\t@xt@\empty\Figg@tT{\p@int}\else#2\fi}%
    \Figwrit@{\p@int}}}\ignorespaces}
\ctr@ld@f\def\Figwrit@#1{\Figg@tXY{#1}\c@lprojSP%
    \rlap{\kern\v@lX\raise\v@lY\hbox{\unhcopy\Gb@x}}\v@leur=\v@lY%
    \advance\v@lY\ht\Gb@x\b@undb@x{\v@lX}{\v@lY}\advance\v@lX\wd\Gb@x%
    \v@lY=\v@leur\advance\v@lY-\dp\Gb@x\b@undb@x{\v@lX}{\v@lY}}
\ctr@ld@f\def\figwritec[#1]#2{{\ignorespaces\def\list@num{#1}%
    \@ecfor\p@int:=\list@num\do{\Figwrit@c{\p@int}{#2}}}\ignorespaces}
\ctr@ld@f\def\Figwrit@c#1#2{\FigWp@r@m{#1}{#2}%
    \rlap{\kern\v@lX\raise\v@lY\hbox{\rlap{\kern-.5\wd\Gb@x%
    \raise-.5\ht\Gb@x\hbox{\raise.5\dp\Gb@x\hbox{\unhcopy\Gb@x}}}}}%
    \v@leur=\ht\Gb@x\advance\v@leur\dp\Gb@x%
    \advance\v@lX-.5\wd\Gb@x\advance\v@lY-.5\v@leur\b@undb@x{\v@lX}{\v@lY}%
    \advance\v@lX\wd\Gb@x\advance\v@lY\v@leur\b@undb@x{\v@lX}{\v@lY}}
\ctr@ld@f\def\figwritep[#1]{{\ignorespaces\def\list@num{#1}\setbox\Gb@x=\hbox{\c@nterpt}%
    \@ecfor\p@int:=\list@num\do{\Figwrit@{\p@int}}}\ignorespaces}
\ctr@ld@f\def\figwritew#1:#2(#3){\figwritegcw#1:{#2}(#3,0pt)}
\ctr@ld@f\def\figwritee#1:#2(#3){\figwritegce#1:{#2}(#3,0pt)}
\ctr@ld@f\def\figwriten#1:#2(#3){{\def\Vc@rrect{\v@lZ=\v@leur\advance\v@lZ\dp\Gb@x}%
    \Figwrit@NS#1:{#2}(#3)}\ignorespaces}
\ctr@ld@f\def\figwrites#1:#2(#3){{\def\Vc@rrect{\v@lZ=-\v@leur\advance\v@lZ-\ht\Gb@x}%
    \Figwrit@NS#1:{#2}(#3)}\ignorespaces}
\ctr@ld@f\def\Figwrit@NS#1:#2(#3){\let\figWp@si=\FigWp@siNS\let\figWBB@x=\FigWBB@xNS%
    \FigWrit@L#1:{#2}(#3,0pt)}
\ctr@ld@f\def\FigWp@siNS{\rlap{\kern\v@lX\raise\v@lY\hbox{\rlap{\kern-.5\wd\Gb@x%
    \raise\v@lZ\hbox{\unhcopy\Gb@x}}\c@nterpt}}}
\ctr@ld@f\def\FigWBB@xNS{\advance\v@lY\v@lZ%
    \advance\v@lY-\dp\Gb@x\advance\v@lX-.5\wd\Gb@x\b@undb@x{\v@lX}{\v@lY}%
    \advance\v@lY\ht\Gb@x\advance\v@lY\dp\Gb@x%
    \advance\v@lX\wd\Gb@x\b@undb@x{\v@lX}{\v@lY}}
\ctr@ld@f\def\figwritenw#1:#2(#3){{\let\figWp@si=\FigWp@sigW\let\figWBB@x=\FigWBB@xgWE%
    \def\C@rp@r@m{\v@leur=\unssqrttw@\v@leur\delt@=\v@leur%
    \ifdim\delt@=\z@\delt@=\epsil@n\fi}\let@xte={-}\FigWrit@L#1:{#2}(#3,0pt)}\ignorespaces}
\ctr@ld@f\def\figwritesw#1:#2(#3){{\let\figWp@si=\FigWp@sigW\let\figWBB@x=\FigWBB@xgWE%
    \def\C@rp@r@m{\v@leur=\unssqrttw@\v@leur\delt@=-\v@leur%
    \ifdim\delt@=\z@\delt@=-\epsil@n\fi}\let@xte={-}\FigWrit@L#1:{#2}(#3,0pt)}\ignorespaces}
\ctr@ld@f\def\figwritene#1:#2(#3){{\let\figWp@si=\FigWp@sigE\let\figWBB@x=\FigWBB@xgWE%
    \def\C@rp@r@m{\v@leur=\unssqrttw@\v@leur\delt@=\v@leur%
    \ifdim\delt@=\z@\delt@=\epsil@n\fi}\let@xte={}\FigWrit@L#1:{#2}(#3,0pt)}\ignorespaces}
\ctr@ld@f\def\figwritese#1:#2(#3){{\let\figWp@si=\FigWp@sigE\let\figWBB@x=\FigWBB@xgWE%
    \def\C@rp@r@m{\v@leur=\unssqrttw@\v@leur\delt@=-\v@leur%
    \ifdim\delt@=\z@\delt@=-\epsil@n\fi}\let@xte={}\FigWrit@L#1:{#2}(#3,0pt)}\ignorespaces}
\ctr@ld@f\def\figwritegw#1:#2(#3,#4){{\let\figWp@si=\FigWp@sigW\let\figWBB@x=\FigWBB@xgWE%
    \let@xte={-}\FigWrit@L#1:{#2}(#3,#4)}\ignorespaces}
\ctr@ld@f\def\figwritege#1:#2(#3,#4){{\let\figWp@si=\FigWp@sigE\let\figWBB@x=\FigWBB@xgWE%
    \let@xte={}\FigWrit@L#1:{#2}(#3,#4)}\ignorespaces}
\ctr@ld@f\def\FigWp@sigW{\v@lXa=\z@\v@lYa=\ht\Gb@x\advance\v@lYa\dp\Gb@x%
    \ifdim\delt@>\z@\relax%
    \rlap{\kern\v@lX\raise\v@lY\hbox{\rlap{\kern-\wd\Gb@x\kern-\v@leur%
          \raise\delt@\hbox{\raise\dp\Gb@x\hbox{\unhcopy\Gb@x}}}\c@nterpt}}%
    \else\ifdim\delt@<\z@\relax\v@lYa=-\v@lYa%
    \rlap{\kern\v@lX\raise\v@lY\hbox{\rlap{\kern-\wd\Gb@x\kern-\v@leur%
          \raise\delt@\hbox{\raise-\ht\Gb@x\hbox{\unhcopy\Gb@x}}}\c@nterpt}}%
    \else\v@lXa=-.5\v@lYa%
    \rlap{\kern\v@lX\raise\v@lY\hbox{\rlap{\kern-\wd\Gb@x\kern-\v@leur%
          \raise-.5\ht\Gb@x\hbox{\raise.5\dp\Gb@x\hbox{\unhcopy\Gb@x}}}\c@nterpt}}%
    \fi\fi}
\ctr@ld@f\def\FigWp@sigE{\v@lXa=\z@\v@lYa=\ht\Gb@x\advance\v@lYa\dp\Gb@x%
    \ifdim\delt@>\z@\relax%
    \rlap{\kern\v@lX\raise\v@lY\hbox{\c@nterpt\kern\v@leur%
          \raise\delt@\hbox{\raise\dp\Gb@x\hbox{\unhcopy\Gb@x}}}}%
    \else\ifdim\delt@<\z@\relax\v@lYa=-\v@lYa%
    \rlap{\kern\v@lX\raise\v@lY\hbox{\c@nterpt\kern\v@leur%
          \raise\delt@\hbox{\raise-\ht\Gb@x\hbox{\unhcopy\Gb@x}}}}%
    \else\v@lXa=-.5\v@lYa%
    \rlap{\kern\v@lX\raise\v@lY\hbox{\c@nterpt\kern\v@leur%
          \raise-.5\ht\Gb@x\hbox{\raise.5\dp\Gb@x\hbox{\unhcopy\Gb@x}}}}%
    \fi\fi}
\ctr@ld@f\def\FigWBB@xgWE{\advance\v@lY\delt@%
    \advance\v@lX\the\let@xte\v@leur\advance\v@lY\v@lXa\b@undb@x{\v@lX}{\v@lY}%
    \advance\v@lX\the\let@xte\wd\Gb@x\advance\v@lY\v@lYa\b@undb@x{\v@lX}{\v@lY}}
\ctr@ld@f\def\figwritegcw#1:#2(#3,#4){{\let\figWp@si=\FigWp@sigcW\let\figWBB@x=\FigWBB@xgcWE%
    \let@xte={-}\FigWrit@L#1:{#2}(#3,#4)}\ignorespaces}
\ctr@ld@f\def\figwritegce#1:#2(#3,#4){{\let\figWp@si=\FigWp@sigcE\let\figWBB@x=\FigWBB@xgcWE%
    \let@xte={}\FigWrit@L#1:{#2}(#3,#4)}\ignorespaces}
\ctr@ld@f\def\FigWp@sigcW{\rlap{\kern\v@lX\raise\v@lY\hbox{\rlap{\kern-\wd\Gb@x\kern-\v@leur%
     \raise-.5\ht\Gb@x\hbox{\raise\delt@\hbox{\raise.5\dp\Gb@x\hbox{\unhcopy\Gb@x}}}}%
     \c@nterpt}}}
\ctr@ld@f\def\FigWp@sigcE{\rlap{\kern\v@lX\raise\v@lY\hbox{\c@nterpt\kern\v@leur%
    \raise-.5\ht\Gb@x\hbox{\raise\delt@\hbox{\raise.5\dp\Gb@x\hbox{\unhcopy\Gb@x}}}}}}
\ctr@ld@f\def\FigWBB@xgcWE{\v@lZ=\ht\Gb@x\advance\v@lZ\dp\Gb@x%
    \advance\v@lX\the\let@xte\v@leur\advance\v@lY\delt@\advance\v@lY.5\v@lZ%
    \b@undb@x{\v@lX}{\v@lY}%
    \advance\v@lX\the\let@xte\wd\Gb@x\advance\v@lY-\v@lZ\b@undb@x{\v@lX}{\v@lY}}
\ctr@ld@f\def\figwritebn#1:#2(#3){{\def\Vc@rrect{\v@lZ=\v@leur}\Figwrit@NS#1:{#2}(#3)}\ignorespaces}
\ctr@ld@f\def\figwritebs#1:#2(#3){{\def\Vc@rrect{\v@lZ=-\v@leur}\Figwrit@NS#1:{#2}(#3)}\ignorespaces}
\ctr@ld@f\def\figwritebw#1:#2(#3){{\let\figWp@si=\FigWp@sibW\let\figWBB@x=\FigWBB@xbWE%
    \let@xte={-}\FigWrit@L#1:{#2}(#3,0pt)}\ignorespaces}
\ctr@ld@f\def\figwritebe#1:#2(#3){{\let\figWp@si=\FigWp@sibE\let\figWBB@x=\FigWBB@xbWE%
    \let@xte={}\FigWrit@L#1:{#2}(#3,0pt)}\ignorespaces}
\ctr@ld@f\def\FigWp@sibW{\rlap{\kern\v@lX\raise\v@lY\hbox{\rlap{\kern-\wd\Gb@x\kern-\v@leur%
          \hbox{\unhcopy\Gb@x}}\c@nterpt}}}
\ctr@ld@f\def\FigWp@sibE{\rlap{\kern\v@lX\raise\v@lY\hbox{\c@nterpt\kern\v@leur%
          \hbox{\unhcopy\Gb@x}}}}
\ctr@ld@f\def\FigWBB@xbWE{\v@lZ=\ht\Gb@x\advance\v@lZ\dp\Gb@x%
    \advance\v@lX\the\let@xte\v@leur\advance\v@lY\ht\Gb@x\b@undb@x{\v@lX}{\v@lY}%
    \advance\v@lX\the\let@xte\wd\Gb@x\advance\v@lY-\v@lZ\b@undb@x{\v@lX}{\v@lY}}
\ctr@ln@w{newread}\frf@g  \ctr@ln@w{newwrite}\fwf@g
\ctr@ln@w{newif}\ifcurr@ntPS
\ctr@ln@w{newif}\ifps@cri
\ctr@ln@w{newif}\ifUse@llipse
\ctr@ln@w{newif}\ifpsdebugmode \psdebugmodefalse
\ctr@ln@w{newif}\ifPDFm@ke
\ifx\pdfliteral\undefined\else\ifnum\pdfoutput>\z@\PDFm@ketrue\fi\fi
\ctr@ld@f\def\initPDF@rDVI{%
\ifPDFm@ke
 \let\figscan=\figscan@E
 \let\newGr@FN=\newGr@FNPDF
 \ctr@ld@f\def\c@mcurveto{c}
 \ctr@ld@f\def\c@mfill{f}
 \ctr@ld@f\def\c@mgsave{q}
 \ctr@ld@f\def\c@mgrestore{Q}
 \ctr@ld@f\def\c@mlineto{l}
 \ctr@ld@f\def\c@mmoveto{m}
 \ctr@ld@f\def\c@msetgray{g}     \ctr@ld@f\def\c@msetgrayStroke{G}
 \ctr@ld@f\def\c@msetcmykcolor{k}\ctr@ld@f\def\c@msetcmykcolorStroke{K}
 \ctr@ld@f\def\c@msetrgbcolor{rg}\ctr@ld@f\def\c@msetrgbcolorStroke{RG}
 \ctr@ld@f\def\d@fprimarC@lor{\curr@ntcolor\space\curr@ntcolorc@md%
               \space\curr@ntcolor\space\curr@ntcolorc@mdStroke}
 \ctr@ld@f\def\d@fsecondC@lor{\sec@ndcolor\space\sec@ndcolorc@md%
               \space\sec@ndcolor\space\sec@ndcolorc@mdStroke}
 \ctr@ld@f\def\d@fthirdC@lor{\th@rdcolor\space\th@rdcolorc@md%
              \space\th@rdcolor\space\th@rdcolorc@mdStroke}
 \ctr@ld@f\def\c@msetdash{d}
 \ctr@ld@f\def\c@msetlinejoin{j}
 \ctr@ld@f\def\c@msetlinewidth{w}
 \ctr@ld@f\def\f@gclosestroke{\immediate\write\fwf@g{s}}
 \ctr@ld@f\def\f@gfill{\immediate\write\fwf@g{\fillc@md}}
 \ctr@ld@f\def\f@gnewpath{}
 \ctr@ld@f\def\f@gstroke{\immediate\write\fwf@g{S}}
\else
 \let\figinsertE=\figinsert
 \let\newGr@FN=\newGr@FNDVI
 \ctr@ld@f\def\c@mcurveto{curveto}
 \ctr@ld@f\def\c@mfill{fill}
 \ctr@ld@f\def\c@mgsave{gsave}
 \ctr@ld@f\def\c@mgrestore{grestore}
 \ctr@ld@f\def\c@mlineto{lineto}
 \ctr@ld@f\def\c@mmoveto{moveto}
 \ctr@ld@f\def\c@msetgray{setgray}          \ctr@ld@f\def\c@msetgrayStroke{}
 \ctr@ld@f\def\c@msetcmykcolor{setcmykcolor}\ctr@ld@f\def\c@msetcmykcolorStroke{}
 \ctr@ld@f\def\c@msetrgbcolor{setrgbcolor}  \ctr@ld@f\def\c@msetrgbcolorStroke{}
 \ctr@ld@f\def\d@fprimarC@lor{\curr@ntcolor\space\curr@ntcolorc@md}
 \ctr@ld@f\def\d@fsecondC@lor{\sec@ndcolor\space\sec@ndcolorc@md}
 \ctr@ld@f\def\d@fthirdC@lor{\th@rdcolor\space\th@rdcolorc@md}
 \ctr@ld@f\def\c@msetdash{setdash}
 \ctr@ld@f\def\c@msetlinejoin{setlinejoin}
 \ctr@ld@f\def\c@msetlinewidth{setlinewidth}
 \ctr@ld@f\def\f@gclosestroke{\immediate\write\fwf@g{closepath\space stroke}}
 \ctr@ld@f\def\f@gfill{\immediate\write\fwf@g{\fillc@md}}
 \ctr@ld@f\def\f@gnewpath{\immediate\write\fwf@g{newpath}}
 \ctr@ld@f\def\f@gstroke{\immediate\write\fwf@g{stroke}}
\fi}
\ctr@ld@f\def\c@pypsfile#1#2{\c@pyfil@{\immediate\write#1}{#2}}
\ctr@ld@f\def\Figinclud@PDF#1#2{\openin\frf@g=#1\pdfliteral{q #2 0 0 #2 0 0 cm}%
    \c@pyfil@{\pdfliteral}{\frf@g}\pdfliteral{Q}\closein\frf@g}
\ctr@ln@w{newif}\ifmored@ta
\ctr@ln@m\bl@nkline
\ctr@ld@f\def\c@pyfil@#1#2{\def\bl@nkline{\par}{\catcode`\%=12
    \loop\ifeof#2\mored@tafalse\else\mored@tatrue\immediate\read#2 to\tr@c
    \ifx\tr@c\bl@nkline\else#1{\tr@c}\fi\fi\ifmored@ta\repeat}}
\ctr@ld@f\def\keln@mun#1#2|{\def\l@debut{#1}\def\l@suite{#2}}
\ctr@ld@f\def\keln@mde#1#2#3|{\def\l@debut{#1#2}\def\l@suite{#3}}
\ctr@ld@f\def\keln@mtr#1#2#3#4|{\def\l@debut{#1#2#3}\def\l@suite{#4}}
\ctr@ld@f\def\keln@mqu#1#2#3#4#5|{\def\l@debut{#1#2#3#4}\def\l@suite{#5}}
\ctr@ld@f\let\@psffilein=\frf@g 
\ctr@ln@w{newif}\if@psffileok    
\ctr@ln@w{newif}\if@psfbbfound   
\ctr@ln@w{newif}\if@psfverbose   
\@psfverbosetrue
\ctr@ln@m\@psfllx \ctr@ln@m\@psflly
\ctr@ln@m\@psfurx \ctr@ln@m\@psfury
\ctr@ln@m\resetcolonc@tcode
\ctr@ld@f\def\@psfgetbb#1{\global\@psfbbfoundfalse%
\global\def\@psfllx{0}\global\def\@psflly{0}%
\global\def\@psfurx{30}\global\def\@psfury{30}%
\openin\@psffilein=#1\relax
\ifeof\@psffilein\errmessage{I couldn't open #1, will ignore it}\else
   \edef\resetcolonc@tcode{\catcode`\noexpand\:\the\catcode`\:\relax}%
   {\@psffileoktrue \chardef\other=12
    \def\do##1{\catcode`##1=\other}\dospecials \catcode`\ =10 \resetcolonc@tcode
    \loop
       \read\@psffilein to \@psffileline
       \ifeof\@psffilein\@psffileokfalse\else
          \expandafter\@psfaux\@psffileline:. \\%
       \fi
   \if@psffileok\repeat
   \if@psfbbfound\else
    \if@psfverbose\message{No bounding box comment in #1; using defaults}\fi\fi
   }\closein\@psffilein\fi}%
\ctr@ln@m\@psfbblit
\ctr@ln@m\@psfpercent
{\catcode`\%=12 \global\let\@psfpercent=
\ctr@ln@m\@psfaux
\long\def\@psfaux#1#2:#3\\{\ifx#1\@psfpercent
   \def\testit{#2}\ifx\testit\@psfbblit
      \@psfgrab #3 . . . \\%
      \@psffileokfalse
      \global\@psfbbfoundtrue
   \fi\else\ifx#1\par\else\@psffileokfalse\fi\fi}%
\ctr@ld@f\def\@psfempty{}%
\ctr@ld@f\def\@psfgrab #1 #2 #3 #4 #5\\{%
\global\def\@psfllx{#1}\ifx\@psfllx\@psfempty
      \@psfgrab #2 #3 #4 #5 .\\\else
   \global\def\@psflly{#2}%
   \global\def\@psfurx{#3}\global\def\@psfury{#4}\fi}%
\ctr@ld@f\def\PSwrit@cmd#1#2#3{{\Figg@tXY{#1}\c@lprojSP\b@undb@x{\v@lX}{\v@lY}%
    \v@lX=\ptT@ptps\v@lX\v@lY=\ptT@ptps\v@lY%
    \immediate\write#3{\repdecn@mb{\v@lX}\space\repdecn@mb{\v@lY}\space#2}}}
\ctr@ld@f\def\PSwrit@cmdS#1#2#3#4#5{{\Figg@tXY{#1}\c@lprojSP\b@undb@x{\v@lX}{\v@lY}%
    \global\result@t=\v@lX\global\result@@t=\v@lY%
    \v@lX=\ptT@ptps\v@lX\v@lY=\ptT@ptps\v@lY%
    \immediate\write#3{\repdecn@mb{\v@lX}\space\repdecn@mb{\v@lY}\space#2}}%
    \edef#4{\the\result@t}\edef#5{\the\result@@t}}
\ctr@ld@f\def\psaltitude#1[#2,#3,#4]{{\ifcurr@ntPS\ifps@cri%
    \PSc@mment{psaltitude Square Dim=#1, Triangle=[#2 / #3,#4]}%
    \s@uvc@ntr@l\et@tpsaltitude\resetc@ntr@l{2}\figptorthoprojline-5:=#2/#3,#4/%
    \figvectP -1[#3,#4]\n@rminf{\v@leur}{-1}\vecunit@{-3}{-1}%
    \figvectP -1[-5,#3]\n@rminf{\v@lmin}{-1}\figvectP -2[-5,#4]\n@rminf{\v@lmax}{-2}%
    \ifdim\v@lmin<\v@lmax\s@mme=#3\else\v@lmax=\v@lmin\s@mme=#4\fi%
    \figvectP -4[-5,#2]\vecunit@{-4}{-4}\delt@=#1\unit@%
    \edef\t@ille{\repdecn@mb{\delt@}}\figpttra-1:=-5/\t@ille,-3/%
    \figptstra-3=-5,-1/\t@ille,-4/\psline[#2,-5]\s@uvdash{\typ@dash}%
    \pssetdash{\defaultdash}\psline[-1,-2,-3]\pssetdash{\typ@dash}%
    \ifdim\v@leur<\v@lmax\Pss@tsecondSt\psline[-5,\the\s@mme]\Psrest@reSt\fi%
    \PSc@mment{End psaltitude}\resetc@ntr@l\et@tpsaltitude\fi\fi}}
\ctr@ld@f\def\Ps@rcerc#1;#2(#3,#4){\ellBB@x#1;#2,#2(#3,#4,0)%
    \f@gnewpath{\delt@=#2\unit@\delt@=\ptT@ptps\delt@%
    \BdingB@xfalse%
    \PSwrit@cmd{#1}{\repdecn@mb{\delt@}\space #3\space #4\space arc}{\fwf@g}}}
\ctr@ln@m\psarccirc
\ctr@ld@f\def\psarccircDD#1;#2(#3,#4){\ifcurr@ntPS\ifps@cri%
    \PSc@mment{psarccircDD Center=#1 ; Radius=#2 (Ang1=#3, Ang2=#4)}%
    \iffillm@de\Ps@rcerc#1;#2(#3,#4)%
    \f@gfill%
    \else\Ps@rcerc#1;#2(#3,#4)\f@gstroke\fi%
    \PSc@mment{End psarccircDD}\fi\fi}
\ctr@ld@f\def\psarccircTD#1,#2,#3;#4(#5,#6){{\ifcurr@ntPS\ifps@cri\s@uvc@ntr@l\et@tpsarccircTD%
    \PSc@mment{psarccircTD Center=#1,P1=#2,P2=#3 ; Radius=#4 (Ang1=#5, Ang2=#6)}%
    \setc@ntr@l{2}\c@lExtAxes#1,#2,#3(#4)\psarcellPATD#1,-4,-5(#5,#6)%
    \PSc@mment{End psarccircTD}\resetc@ntr@l\et@tpsarccircTD\fi\fi}}
\ctr@ld@f\def\c@lExtAxes#1,#2,#3(#4){%
    \figvectPTD-5[#1,#2]\vecunit@{-5}{-5}\figvectNTD-4[#1,#2,#3]\vecunit@{-4}{-4}%
    \figvectNVTD-3[-4,-5]\delt@=#4\unit@\edef\r@yon{\repdecn@mb{\delt@}}%
    \figpttra-4:=#1/\r@yon,-5/\figpttra-5:=#1/\r@yon,-3/}
\ctr@ln@m\psarccircP
\ctr@ld@f\def\psarccircPDD#1;#2[#3,#4]{{\ifcurr@ntPS\ifps@cri\s@uvc@ntr@l\et@tpsarccircPDD%
    \PSc@mment{psarccircPDD Center=#1; Radius=#2, [P1=#3, P2=#4]}%
    \Ps@ngleparam#1;#2[#3,#4]\ifdim\v@lmin>\v@lmax\advance\v@lmax\DePI@deg\fi%
    \edef\@ngdeb{\repdecn@mb{\v@lmin}}\edef\@ngfin{\repdecn@mb{\v@lmax}}%
    \psarccirc#1;\r@dius(\@ngdeb,\@ngfin)%
    \PSc@mment{End psarccircPDD}\resetc@ntr@l\et@tpsarccircPDD\fi\fi}}
\ctr@ld@f\def\psarccircPTD#1;#2[#3,#4,#5]{{\ifcurr@ntPS\ifps@cri\s@uvc@ntr@l\et@tpsarccircPTD%
    \PSc@mment{psarccircPTD Center=#1; Radius=#2, [P1=#3, P2=#4, P3=#5]}%
    \setc@ntr@l{2}\c@lExtAxes#1,#3,#5(#2)\psarcellPP#1,-4,-5[#3,#4]%
    \PSc@mment{End psarccircPTD}\resetc@ntr@l\et@tpsarccircPTD\fi\fi}}
\ctr@ld@f\def\Ps@ngleparam#1;#2[#3,#4]{\setc@ntr@l{2}%
    \figvectPDD-1[#1,#3]\vecunit@{-1}{-1}\Figg@tXY{-1}\arct@n\v@lmin(\v@lX,\v@lY)%
    \figvectPDD-2[#1,#4]\vecunit@{-2}{-2}\Figg@tXY{-2}\arct@n\v@lmax(\v@lX,\v@lY)%
    \v@lmin=\rdT@deg\v@lmin\v@lmax=\rdT@deg\v@lmax%
    \v@leur=#2pt\maxim@m{\mili@u}{-\v@leur}{\v@leur}%
    \edef\r@dius{\repdecn@mb{\mili@u}}}
\ctr@ld@f\def\Ps@rcercBz#1;#2(#3,#4){\Ps@rellBz#1;#2,#2(#3,#4,0)}
\ctr@ld@f\def\Ps@rellBz#1;#2,#3(#4,#5,#6){%
    \ellBB@x#1;#2,#3(#4,#5,#6)\BdingB@xfalse%
    \c@lNbarcs{#4}{#5}\v@leur=#4pt\setc@ntr@l{2}\figptell-13::#1;#2,#3(#4,#6)%
    \f@gnewpath\PSwrit@cmd{-13}{\c@mmoveto}{\fwf@g}%
    \s@mme=\z@\bcl@rellBz#1;#2,#3(#6)\BdingB@xtrue}
\ctr@ld@f\def\bcl@rellBz#1;#2,#3(#4){\relax%
    \ifnum\s@mme<\p@rtent\advance\s@mme\@ne%
    \advance\v@leur\delt@\edef\@ngle{\repdecn@mb\v@leur}\figptell-14::#1;#2,#3(\@ngle,#4)%
    \advance\v@leur\delt@\edef\@ngle{\repdecn@mb\v@leur}\figptell-15::#1;#2,#3(\@ngle,#4)%
    \advance\v@leur\delt@\edef\@ngle{\repdecn@mb\v@leur}\figptell-16::#1;#2,#3(\@ngle,#4)%
    \figptscontrolDD-18[-13,-14,-15,-16]%
    \PSwrit@cmd{-18}{}{\fwf@g}\PSwrit@cmd{-17}{}{\fwf@g}%
    \PSwrit@cmd{-16}{\c@mcurveto}{\fwf@g}%
    \figptcopyDD-13:/-16/\bcl@rellBz#1;#2,#3(#4)\fi}
\ctr@ld@f\def\Ps@rell#1;#2,#3(#4,#5,#6){\ellBB@x#1;#2,#3(#4,#5,#6)%
    \f@gnewpath{\v@lmin=#2\unit@\v@lmin=\ptT@ptps\v@lmin%
    \v@lmax=#3\unit@\v@lmax=\ptT@ptps\v@lmax\BdingB@xfalse%
    \PSwrit@cmd{#1}%
    {#6\space\repdecn@mb{\v@lmin}\space\repdecn@mb{\v@lmax}\space #4\space #5\space ellipse}{\fwf@g}}%
    \global\Use@llipsetrue}
\ctr@ln@m\psarcell
\ctr@ld@f\def\psarcellDD#1;#2,#3(#4,#5,#6){{\ifcurr@ntPS\ifps@cri%
    \PSc@mment{psarcellDD Center=#1 ; XRad=#2, YRad=#3 (Ang1=#4, Ang2=#5, Inclination=#6)}%
    \iffillm@de\Ps@rell#1;#2,#3(#4,#5,#6)%
    \f@gfill%
    \else\Ps@rell#1;#2,#3(#4,#5,#6)\f@gstroke\fi%
    \PSc@mment{End psarcellDD}\fi\fi}}
\ctr@ld@f\def\psarcellTD#1;#2,#3(#4,#5,#6){{\ifcurr@ntPS\ifps@cri\s@uvc@ntr@l\et@tpsarcellTD%
    \PSc@mment{psarcellTD Center=#1 ; XRad=#2, YRad=#3 (Ang1=#4, Ang2=#5, Inclination=#6)}%
    \setc@ntr@l{2}\figpttraC -8:=#1/#2,0,0/\figpttraC -7:=#1/0,#3,0/%
    \figvectC -4(0,0,1)\figptsrot -8=-8,-7/#1,#6,-4/\psarcellPATD#1,-8,-7(#4,#5)%
    \PSc@mment{End psarcellTD}\resetc@ntr@l\et@tpsarcellTD\fi\fi}}
\ctr@ln@m\psarcellPA
\ctr@ld@f\def\psarcellPADD#1,#2,#3(#4,#5){{\ifcurr@ntPS\ifps@cri\s@uvc@ntr@l\et@tpsarcellPADD%
    \PSc@mment{psarcellPADD Center=#1,PtAxis1=#2,PtAxis2=#3 (Ang1=#4, Ang2=#5)}%
    \setc@ntr@l{2}\figvectPDD-1[#1,#2]\vecunit@DD{-1}{-1}\v@lX=\ptT@unit@\result@t%
    \edef\XR@d{\repdecn@mb{\v@lX}}\Figg@tXY{-1}\arct@n\v@lmin(\v@lX,\v@lY)%
    \v@lmin=\rdT@deg\v@lmin\edef\Inclin@{\repdecn@mb{\v@lmin}}%
    \figgetdist\YR@d[#1,#3]\psarcellDD#1;\XR@d,\YR@d(#4,#5,\Inclin@)%
    \PSc@mment{End psarcellPADD}\resetc@ntr@l\et@tpsarcellPADD\fi\fi}}
\ctr@ld@f\def\psarcellPATD#1,#2,#3(#4,#5){{\ifcurr@ntPS\ifps@cri\s@uvc@ntr@l\et@tpsarcellPATD%
    \PSc@mment{psarcellPATD Center=#1,PtAxis1=#2,PtAxis2=#3 (Ang1=#4, Ang2=#5)}%
    \iffillm@de\Ps@rellPATD#1,#2,#3(#4,#5)%
    \f@gfill%
    \else\Ps@rellPATD#1,#2,#3(#4,#5)\f@gstroke\fi%
    \PSc@mment{End psarcellPATD}\resetc@ntr@l\et@tpsarcellPATD\fi\fi}}
\ctr@ld@f\def\Ps@rellPATD#1,#2,#3(#4,#5){\let\c@lprojSP=\relax%
    \setc@ntr@l{2}\figvectPTD-1[#1,#2]\figvectPTD-2[#1,#3]\c@lNbarcs{#4}{#5}%
    \v@leur=#4pt\c@lptellP{#1}{-1}{-2}\Figptpr@j-5:/-3/%
    \f@gnewpath\PSwrit@cmdS{-5}{\c@mmoveto}{\fwf@g}{\X@un}{\Y@un}%
    \edef\C@nt@r{#1}\s@mme=\z@\bcl@rellPATD}
\ctr@ld@f\def\bcl@rellPATD{\relax%
    \ifnum\s@mme<\p@rtent\advance\s@mme\@ne%
    \advance\v@leur\delt@\c@lptellP{\C@nt@r}{-1}{-2}\Figptpr@j-4:/-3/%
    \advance\v@leur\delt@\c@lptellP{\C@nt@r}{-1}{-2}\Figptpr@j-6:/-3/%
    \advance\v@leur\delt@\c@lptellP{\C@nt@r}{-1}{-2}\Figptpr@j-3:/-3/%
    \v@lX=\z@\v@lY=\z@\Figtr@nptDD{-5}{-5}\Figtr@nptDD{2}{-3}%
    \divide\v@lX\@vi\divide\v@lY\@vi%
    \Figtr@nptDD{3}{-4}\Figtr@nptDD{-1.5}{-6}\v@lmin=\v@lX\v@lmax=\v@lY%
    \v@lX=\z@\v@lY=\z@\Figtr@nptDD{2}{-5}\Figtr@nptDD{-5}{-3}%
    \divide\v@lX\@vi\divide\v@lY\@vi\Figtr@nptDD{-1.5}{-4}\Figtr@nptDD{3}{-6}%
    \BdingB@xfalse%
    \Figp@intregDD-4:(\v@lmin,\v@lmax)\PSwrit@cmdS{-4}{}{\fwf@g}{\X@de}{\Y@de}%
    \Figp@intregDD-4:(\v@lX,\v@lY)\PSwrit@cmdS{-4}{}{\fwf@g}{\X@tr}{\Y@tr}%
    \BdingB@xtrue\PSwrit@cmdS{-3}{\c@mcurveto}{\fwf@g}{\X@qu}{\Y@qu}%
    \B@zierBB@x{1}{\Y@un}(\X@un,\X@de,\X@tr,\X@qu)%
    \B@zierBB@x{2}{\X@un}(\Y@un,\Y@de,\Y@tr,\Y@qu)%
    \edef\X@un{\X@qu}\edef\Y@un{\Y@qu}\figptcopyDD-5:/-3/\bcl@rellPATD\fi}
\ctr@ld@f\def\c@lNbarcs#1#2{%
    \delt@=#2pt\advance\delt@-#1pt\maxim@m{\v@lmax}{\delt@}{-\delt@}%
    \v@leur=\v@lmax\divide\v@leur45 \p@rtentiere{\p@rtent}{\v@leur}\advance\p@rtent\@ne%
    \s@mme=\p@rtent\multiply\s@mme\thr@@\divide\delt@\s@mme}
\ctr@ld@f\def\psarcellPP#1,#2,#3[#4,#5]{{\ifcurr@ntPS\ifps@cri\s@uvc@ntr@l\et@tpsarcellPP%
    \PSc@mment{psarcellPP Center=#1,PtAxis1=#2,PtAxis2=#3 [Point1=#4, Point2=#5]}%
    \setc@ntr@l{2}\figvectP-2[#1,#3]\vecunit@{-2}{-2}\v@lmin=\result@t%
    \invers@{\v@lmax}{\v@lmin}%
    \figvectP-1[#1,#2]\vecunit@{-1}{-1}\v@leur=\result@t%
    \v@leur=\repdecn@mb{\v@lmax}\v@leur\edef\AsB@{\repdecn@mb{\v@leur}}
    \c@lAngle{#1}{#4}{\v@lmin}\edef\@ngdeb{\repdecn@mb{\v@lmin}}%
    \c@lAngle{#1}{#5}{\v@lmax}\ifdim\v@lmin>\v@lmax\advance\v@lmax\DePI@deg\fi%
    \edef\@ngfin{\repdecn@mb{\v@lmax}}\psarcellPA#1,#2,#3(\@ngdeb,\@ngfin)%
    \PSc@mment{End psarcellPP}\resetc@ntr@l\et@tpsarcellPP\fi\fi}}
\ctr@ld@f\def\c@lAngle#1#2#3{\figvectP-3[#1,#2]%
    \c@lproscal\delt@[-3,-1]\c@lproscal\v@leur[-3,-2]%
    \v@leur=\AsB@\v@leur\arct@n#3(\delt@,\v@leur)#3=\rdT@deg#3}
\ctr@ln@w{newif}\if@rrowratio\@rrowratiotrue
\ctr@ln@w{newif}\if@rrowhfill
\ctr@ln@w{newif}\if@rrowhout
\ctr@ld@f\def\Psset@rrowhe@d#1=#2|{\keln@mun#1|%
    \def\n@mref{a}\ifx\l@debut\n@mref\pssetarrowheadangle{#2}\else
    \def\n@mref{f}\ifx\l@debut\n@mref\pssetarrowheadfill{#2}\else
    \def\n@mref{l}\ifx\l@debut\n@mref\pssetarrowheadlength{#2}\else
    \def\n@mref{o}\ifx\l@debut\n@mref\pssetarrowheadout{#2}\else
    \def\n@mref{r}\ifx\l@debut\n@mref\pssetarrowheadratio{#2}\else
    \immediate\write16{*** Unknown attribute: \BS@ psset arrowhead(..., #1=...)}%
    \fi\fi\fi\fi\fi}
\ctr@ln@m\@rrowheadangle
\ctr@ln@m\C@AHANG \ctr@ln@m\S@AHANG \ctr@ln@m\UNSS@N
\ctr@ld@f\def\pssetarrowheadangle#1{\edef\@rrowheadangle{#1}{\c@ssin{\C@}{\S@}{#1}%
    \xdef\C@AHANG{\C@}\xdef\S@AHANG{\S@}\v@lmax=\S@ pt%
    \invers@{\v@leur}{\v@lmax}\maxim@m{\v@leur}{\v@leur}{-\v@leur}%
    \xdef\UNSS@N{\the\v@leur}}}
\ctr@ld@f\def\pssetarrowheadfill#1{\expandafter\set@rrowhfill#1:}
\ctr@ld@f\def\set@rrowhfill#1#2:{\if#1n\@rrowhfillfalse\else\@rrowhfilltrue\fi}
\ctr@ld@f\def\pssetarrowheadout#1{\expandafter\set@rrowhout#1:}
\ctr@ld@f\def\set@rrowhout#1#2:{\if#1n\@rrowhoutfalse\else\@rrowhouttrue\fi}
\ctr@ln@m\@rrowheadlength
\ctr@ld@f\def\pssetarrowheadlength#1{\edef\@rrowheadlength{#1}\@rrowratiofalse}
\ctr@ln@m\@rrowheadratio
\ctr@ld@f\def\pssetarrowheadratio#1{\edef\@rrowheadratio{#1}\@rrowratiotrue}
\ctr@ln@m\defaultarrowheadlength
\ctr@ld@f\def\psresetarrowhead{%
    \pssetarrowheadangle{\defaultarrowheadangle}%
    \pssetarrowheadfill{\defaultarrowheadfill}%
    \pssetarrowheadout{\defaultarrowheadout}%
    \pssetarrowheadratio{\defaultarrowheadratio}%
    \d@fm@cdim\defaultarrowheadlength{\defaulth@rdahlength}
    \pssetarrowheadlength{\defaultarrowheadlength}}
\ctr@ld@f\def\defaultarrowheadratio{0.1}
\ctr@ld@f\def\defaultarrowheadangle{20}
\ctr@ld@f\def\defaultarrowheadfill{no}
\ctr@ld@f\def\defaultarrowheadout{no}
\ctr@ld@f\def\defaulth@rdahlength{8pt}
\ctr@ln@m\psarrow
\ctr@ld@f\def\psarrowDD[#1,#2]{{\ifcurr@ntPS\ifps@cri\s@uvc@ntr@l\et@tpsarrow%
    \PSc@mment{psarrowDD [Pt1,Pt2]=[#1,#2]}\pssetfillmode{no}%
    \psarrowheadDD[#1,#2]\setc@ntr@l{2}\psline[#1,-3]%
    \PSc@mment{End psarrowDD}\resetc@ntr@l\et@tpsarrow\fi\fi}}
\ctr@ld@f\def\psarrowTD[#1,#2]{{\ifcurr@ntPS\ifps@cri\s@uvc@ntr@l\et@tpsarrowTD%
    \PSc@mment{psarrowTD [Pt1,Pt2]=[#1,#2]}\resetc@ntr@l{2}%
    \Figptpr@j-5:/#1/\Figptpr@j-6:/#2/\let\c@lprojSP=\relax\psarrowDD[-5,-6]%
    \PSc@mment{End psarrowTD}\resetc@ntr@l\et@tpsarrowTD\fi\fi}}
\ctr@ln@m\psarrowhead
\ctr@ld@f\def\psarrowheadDD[#1,#2]{{\ifcurr@ntPS\ifps@cri\s@uvc@ntr@l\et@tpsarrowheadDD%
    \if@rrowhfill\def\@hangle{-\@rrowheadangle}\else\def\@hangle{\@rrowheadangle}\fi%
    \if@rrowratio%
    \if@rrowhout\def\@hratio{-\@rrowheadratio}\else\def\@hratio{\@rrowheadratio}\fi%
    \PSc@mment{psarrowheadDD Ratio=\@hratio, Angle=\@hangle, [Pt1,Pt2]=[#1,#2]}%
    \Ps@rrowhead\@hratio,\@hangle[#1,#2]%
    \else%
    \if@rrowhout\def\@hlength{-\@rrowheadlength}\else\def\@hlength{\@rrowheadlength}\fi%
    \PSc@mment{psarrowheadDD Length=\@hlength, Angle=\@hangle, [Pt1,Pt2]=[#1,#2]}%
    \Ps@rrowheadfd\@hlength,\@hangle[#1,#2]%
    \fi%
    \PSc@mment{End psarrowheadDD}\resetc@ntr@l\et@tpsarrowheadDD\fi\fi}}
\ctr@ld@f\def\psarrowheadTD[#1,#2]{{\ifcurr@ntPS\ifps@cri\s@uvc@ntr@l\et@tpsarrowheadTD%
    \PSc@mment{psarrowheadTD [Pt1,Pt2]=[#1,#2]}\resetc@ntr@l{2}%
    \Figptpr@j-5:/#1/\Figptpr@j-6:/#2/\let\c@lprojSP=\relax\psarrowheadDD[-5,-6]%
    \PSc@mment{End psarrowheadTD}\resetc@ntr@l\et@tpsarrowheadTD\fi\fi}}
\ctr@ld@f\def\Ps@rrowhead#1,#2[#3,#4]{\v@leur=#1\p@\maxim@m{\v@leur}{\v@leur}{-\v@leur}%
    \ifdim\v@leur>\Cepsil@n{
    \PSc@mment{ps@rrowhead Ratio=#1, Angle=#2, [Pt1,Pt2]=[#3,#4]}\v@leur=\UNSS@N%
    \v@leur=\curr@ntwidth\v@leur\v@leur=\ptpsT@pt\v@leur\delt@=.5\v@leur
    \setc@ntr@l{2}\figvectPDD-3[#4,#3]%
    \Figg@tXY{-3}\v@lX=#1\v@lX\v@lY=#1\v@lY\Figv@ctCreg-3(\v@lX,\v@lY)%
    \vecunit@{-4}{-3}\mili@u=\result@t%
    \ifdim#2pt>\z@\v@lXa=-\C@AHANG\delt@%
     \edef\c@ef{\repdecn@mb{\v@lXa}}\figpttraDD-3:=-3/\c@ef,-4/\fi%
    \edef\c@ef{\repdecn@mb{\delt@}}%
    \v@lXa=\mili@u\v@lXa=\C@AHANG\v@lXa%
    \v@lYa=\ptpsT@pt\p@\v@lYa=\curr@ntwidth\v@lYa\v@lYa=\sDcc@ngle\v@lYa%
    \advance\v@lXa-\v@lYa\gdef\sDcc@ngle{0}%
    \ifdim\v@lXa>\v@leur\edef\c@efendpt{\repdecn@mb{\v@leur}}%
    \else\edef\c@efendpt{\repdecn@mb{\v@lXa}}\fi%
    \Figg@tXY{-3}\v@lmin=\v@lX\v@lmax=\v@lY%
    \v@lXa=\C@AHANG\v@lmin\v@lYa=\S@AHANG\v@lmax\advance\v@lXa\v@lYa%
    \v@lYa=-\S@AHANG\v@lmin\v@lX=\C@AHANG\v@lmax\advance\v@lYa\v@lX%
    \setc@ntr@l{1}\Figg@tXY{#4}\advance\v@lX\v@lXa\advance\v@lY\v@lYa%
    \setc@ntr@l{2}\Figp@intregDD-2:(\v@lX,\v@lY)%
    \v@lXa=\C@AHANG\v@lmin\v@lYa=-\S@AHANG\v@lmax\advance\v@lXa\v@lYa%
    \v@lYa=\S@AHANG\v@lmin\v@lX=\C@AHANG\v@lmax\advance\v@lYa\v@lX%
    \setc@ntr@l{1}\Figg@tXY{#4}\advance\v@lX\v@lXa\advance\v@lY\v@lYa%
    \setc@ntr@l{2}\Figp@intregDD-1:(\v@lX,\v@lY)%
    \ifdim#2pt<\z@\fillm@detrue\psline[-2,#4,-1]
    \else\figptstraDD-3=#4,-2,-1/\c@ef,-4/\psline[-2,-3,-1]\fi
    \ifdim#1pt>\z@\figpttraDD-3:=#4/\c@efendpt,-4/\else\figptcopyDD-3:/#4/\fi%
    \PSc@mment{End ps@rrowhead}}\fi}
\ctr@ld@f\def\sDcc@ngle{0}
\ctr@ld@f\def\Ps@rrowheadfd#1,#2[#3,#4]{{%
    \PSc@mment{ps@rrowheadfd Length=#1, Angle=#2, [Pt1,Pt2]=[#3,#4]}%
    \setc@ntr@l{2}\figvectPDD-1[#3,#4]\n@rmeucDD{\v@leur}{-1}\v@leur=\ptT@unit@\v@leur%
    \invers@{\v@leur}{\v@leur}\v@leur=#1\v@leur\edef\R@tio{\repdecn@mb{\v@leur}}%
    \Ps@rrowhead\R@tio,#2[#3,#4]\PSc@mment{End ps@rrowheadfd}}}
\ctr@ln@m\psarrowBezier
\ctr@ld@f\def\psarrowBezierDD[#1,#2,#3,#4]{{\ifcurr@ntPS\ifps@cri\s@uvc@ntr@l\et@tpsarrowBezierDD%
    \PSc@mment{psarrowBezierDD Control points=#1,#2,#3,#4}\setc@ntr@l{2}%
    \if@rrowratio\c@larclengthDD\v@leur,10[#1,#2,#3,#4]\else\v@leur=\z@\fi%
    \Ps@rrowB@zDD\v@leur[#1,#2,#3,#4]%
    \PSc@mment{End psarrowBezierDD}\resetc@ntr@l\et@tpsarrowBezierDD\fi\fi}}
\ctr@ld@f\def\psarrowBezierTD[#1,#2,#3,#4]{{\ifcurr@ntPS\ifps@cri\s@uvc@ntr@l\et@tpsarrowBezierTD%
    \PSc@mment{psarrowBezierTD Control points=#1,#2,#3,#4}\resetc@ntr@l{2}%
    \Figptpr@j-7:/#1/\Figptpr@j-8:/#2/\Figptpr@j-9:/#3/\Figptpr@j-10:/#4/%
    \let\c@lprojSP=\relax\ifnum\curr@ntproj<\tw@\psarrowBezierDD[-7,-8,-9,-10]%
    \else\f@gnewpath\PSwrit@cmd{-7}{\c@mmoveto}{\fwf@g}%
    \if@rrowratio\c@larclengthDD\mili@u,10[-7,-8,-9,-10]\else\mili@u=\z@\fi%
    \p@rtent=\NBz@rcs\advance\p@rtent\m@ne\subB@zierTD\p@rtent[#1,#2,#3,#4]%
    \f@gstroke%
    \advance\v@lmin\p@rtent\delt@
    \v@leur=\v@lmin\advance\v@leur0.33333 \delt@\edef\unti@rs{\repdecn@mb{\v@leur}}%
    \v@leur=\v@lmin\advance\v@leur0.66666 \delt@\edef\deti@rs{\repdecn@mb{\v@leur}}%
    \figptcopyDD-8:/-10/\c@lsubBzarc\unti@rs,\deti@rs[#1,#2,#3,#4]%
    \figptcopyDD-8:/-4/\figptcopyDD-9:/-3/\Ps@rrowB@zDD\mili@u[-7,-8,-9,-10]\fi%
    \PSc@mment{End psarrowBezierTD}\resetc@ntr@l\et@tpsarrowBezierTD\fi\fi}}
\ctr@ld@f\def\c@larclengthDD#1,#2[#3,#4,#5,#6]{{\p@rtent=#2\figptcopyDD-5:/#3/%
    \delt@=\p@\divide\delt@\p@rtent\c@rre=\z@\v@leur=\z@\s@mme=\z@%
    \loop\ifnum\s@mme<\p@rtent\advance\s@mme\@ne\advance\v@leur\delt@%
    \edef\T@{\repdecn@mb{\v@leur}}\figptBezierDD-6::\T@[#3,#4,#5,#6]%
    \figvectPDD-1[-5,-6]\n@rmeucDD{\mili@u}{-1}\advance\c@rre\mili@u%
    \figptcopyDD-5:/-6/\repeat\global\result@t=\ptT@unit@\c@rre}#1=\result@t}
\ctr@ld@f\def\Ps@rrowB@zDD#1[#2,#3,#4,#5]{{\pssetfillmode{no}%
    \if@rrowratio\delt@=\@rrowheadratio#1\else\delt@=\@rrowheadlength pt\fi%
    \v@leur=\C@AHANG\delt@\edef\R@dius{\repdecn@mb{\v@leur}}%
    \FigptintercircB@zDD-5::0,\R@dius[#5,#4,#3,#2]%
    \pssetarrowheadlength{\repdecn@mb{\delt@}}\psarrowheadDD[-5,#5]%
    \let\n@rmeuc=\n@rmeucDD\figgetdist\R@dius[#5,-3]%
    \FigptintercircB@zDD-6::0,\R@dius[#5,#4,#3,#2]%
    \figptBezierDD-5::0.33333[#5,#4,#3,#2]\figptBezierDD-3::0.66666[#5,#4,#3,#2]%
    \figptscontrolDD-5[-6,-5,-3,#2]\psBezierDD1[-6,-5,-4,#2]}}
\ctr@ln@m\psarrowcirc
\ctr@ld@f\def\psarrowcircDD#1;#2(#3,#4){{\ifcurr@ntPS\ifps@cri\s@uvc@ntr@l\et@tpsarrowcircDD%
    \PSc@mment{psarrowcircDD Center=#1 ; Radius=#2 (Ang1=#3,Ang2=#4)}%
    \pssetfillmode{no}\Pscirc@rrowhead#1;#2(#3,#4)%
    \setc@ntr@l{2}\figvectPDD -4[#1,-3]\vecunit@{-4}{-4}%
    \Figg@tXY{-4}\arct@n\v@lmin(\v@lX,\v@lY)%
    \v@lmin=\rdT@deg\v@lmin\v@leur=#4pt\advance\v@leur-\v@lmin%
    \maxim@m{\v@leur}{\v@leur}{-\v@leur}%
    \ifdim\v@leur>\DemiPI@deg\relax\ifdim\v@lmin<#4pt\advance\v@lmin\DePI@deg%
    \else\advance\v@lmin-\DePI@deg\fi\fi\edef\ar@ngle{\repdecn@mb{\v@lmin}}%
    \ifdim#3pt<#4pt\psarccirc#1;#2(#3,\ar@ngle)\else\psarccirc#1;#2(\ar@ngle,#3)\fi%
    \PSc@mment{End psarrowcircDD}\resetc@ntr@l\et@tpsarrowcircDD\fi\fi}}
\ctr@ld@f\def\psarrowcircTD#1,#2,#3;#4(#5,#6){{\ifcurr@ntPS\ifps@cri\s@uvc@ntr@l\et@tpsarrowcircTD%
    \PSc@mment{psarrowcircTD Center=#1,P1=#2,P2=#3 ; Radius=#4 (Ang1=#5, Ang2=#6)}%
    \resetc@ntr@l{2}\c@lExtAxes#1,#2,#3(#4)\let\c@lprojSP=\relax%
    \figvectPTD-11[#1,-4]\figvectPTD-12[#1,-5]\c@lNbarcs{#5}{#6}%
    \if@rrowratio\v@lmax=\degT@rd\v@lmax\edef\D@lpha{\repdecn@mb{\v@lmax}}\fi%
    \advance\p@rtent\m@ne\mili@u=\z@%
    \v@leur=#5pt\c@lptellP{#1}{-11}{-12}\Figptpr@j-9:/-3/%
    \f@gnewpath\PSwrit@cmdS{-9}{\c@mmoveto}{\fwf@g}{\X@un}{\Y@un}%
    \edef\C@nt@r{#1}\s@mme=\z@\bcl@rcircTD\f@gstroke%
    \advance\v@leur\delt@\c@lptellP{#1}{-11}{-12}\Figptpr@j-5:/-3/%
    \advance\v@leur\delt@\c@lptellP{#1}{-11}{-12}\Figptpr@j-6:/-3/%
    \advance\v@leur\delt@\c@lptellP{#1}{-11}{-12}\Figptpr@j-10:/-3/%
    \figptscontrolDD-8[-9,-5,-6,-10]%
    \if@rrowratio\c@lcurvradDD0.5[-9,-8,-7,-10]\advance\mili@u\result@t%
    \maxim@m{\mili@u}{\mili@u}{-\mili@u}\mili@u=\ptT@unit@\mili@u%
    \mili@u=\D@lpha\mili@u\advance\p@rtent\@ne\divide\mili@u\p@rtent\fi%
    \Ps@rrowB@zDD\mili@u[-9,-8,-7,-10]%
    \PSc@mment{End psarrowcircTD}\resetc@ntr@l\et@tpsarrowcircTD\fi\fi}}
\ctr@ld@f\def\bcl@rcircTD{\relax%
    \ifnum\s@mme<\p@rtent\advance\s@mme\@ne%
    \advance\v@leur\delt@\c@lptellP{\C@nt@r}{-11}{-12}\Figptpr@j-5:/-3/%
    \advance\v@leur\delt@\c@lptellP{\C@nt@r}{-11}{-12}\Figptpr@j-6:/-3/%
    \advance\v@leur\delt@\c@lptellP{\C@nt@r}{-11}{-12}\Figptpr@j-10:/-3/%
    \figptscontrolDD-8[-9,-5,-6,-10]\BdingB@xfalse%
    \PSwrit@cmdS{-8}{}{\fwf@g}{\X@de}{\Y@de}\PSwrit@cmdS{-7}{}{\fwf@g}{\X@tr}{\Y@tr}%
    \BdingB@xtrue\PSwrit@cmdS{-10}{\c@mcurveto}{\fwf@g}{\X@qu}{\Y@qu}%
    \if@rrowratio\c@lcurvradDD0.5[-9,-8,-7,-10]\advance\mili@u\result@t\fi%
    \B@zierBB@x{1}{\Y@un}(\X@un,\X@de,\X@tr,\X@qu)%
    \B@zierBB@x{2}{\X@un}(\Y@un,\Y@de,\Y@tr,\Y@qu)%
    \edef\X@un{\X@qu}\edef\Y@un{\Y@qu}\figptcopyDD-9:/-10/\bcl@rcircTD\fi}
\ctr@ld@f\def\Pscirc@rrowhead#1;#2(#3,#4){{%
    \PSc@mment{pscirc@rrowhead Center=#1 ; Radius=#2 (Ang1=#3,Ang2=#4)}%
    \v@leur=#2\unit@\edef\s@glen{\repdecn@mb{\v@leur}}\v@lY=\z@\v@lX=\v@leur%
    \resetc@ntr@l{2}\Figv@ctCreg-3(\v@lX,\v@lY)\figpttraDD-5:=#1/1,-3/%
    \figptrotDD-5:=-5/#1,#4/%
    \figvectPDD-3[#1,-5]\Figg@tXY{-3}\v@leur=\v@lX%
    \ifdim#3pt<#4pt\v@lX=\v@lY\v@lY=-\v@leur\else\v@lX=-\v@lY\v@lY=\v@leur\fi%
    \Figv@ctCreg-3(\v@lX,\v@lY)\vecunit@{-3}{-3}%
    \if@rrowratio\v@leur=#4pt\advance\v@leur-#3pt\maxim@m{\mili@u}{-\v@leur}{\v@leur}%
    \mili@u=\degT@rd\mili@u\v@leur=\s@glen\mili@u\edef\s@glen{\repdecn@mb{\v@leur}}%
    \mili@u=#2\mili@u\mili@u=\@rrowheadratio\mili@u\else\mili@u=\@rrowheadlength pt\fi%
    \figpttraDD-6:=-5/\s@glen,-3/\v@leur=#2pt\v@leur=2\v@leur%
    \invers@{\v@leur}{\v@leur}\c@rre=\repdecn@mb{\v@leur}\mili@u
    \mili@u=\c@rre\mili@u=\repdecn@mb{\c@rre}\mili@u%
    \v@leur=\p@\advance\v@leur-\mili@u
    \invers@{\mili@u}{2\v@leur}\delt@=\c@rre\delt@=\repdecn@mb{\mili@u}\delt@%
    \xdef\sDcc@ngle{\repdecn@mb{\delt@}}
    \sqrt@{\mili@u}{\v@leur}\arct@n\v@leur(\mili@u,\c@rre)%
    \v@leur=\rdT@deg\v@leur
    \ifdim#3pt<#4pt\v@leur=-\v@leur\fi%
    \if@rrowhout\v@leur=-\v@leur\fi\edef\cor@ngle{\repdecn@mb{\v@leur}}%
    \figptrotDD-6:=-6/-5,\cor@ngle/\psarrowheadDD[-6,-5]%
    \PSc@mment{End pscirc@rrowhead}}}
\ctr@ln@m\psarrowcircP
\ctr@ld@f\def\psarrowcircPDD#1;#2[#3,#4]{{\ifcurr@ntPS\ifps@cri%
    \PSc@mment{psarrowcircPDD Center=#1; Radius=#2, [P1=#3,P2=#4]}%
    \s@uvc@ntr@l\et@tpsarrowcircPDD\Ps@ngleparam#1;#2[#3,#4]%
    \ifdim\v@leur>\z@\ifdim\v@lmin>\v@lmax\advance\v@lmax\DePI@deg\fi%
    \else\ifdim\v@lmin<\v@lmax\advance\v@lmin\DePI@deg\fi\fi%
    \edef\@ngdeb{\repdecn@mb{\v@lmin}}\edef\@ngfin{\repdecn@mb{\v@lmax}}%
    \psarrowcirc#1;\r@dius(\@ngdeb,\@ngfin)%
    \PSc@mment{End psarrowcircPDD}\resetc@ntr@l\et@tpsarrowcircPDD\fi\fi}}
\ctr@ld@f\def\psarrowcircPTD#1;#2[#3,#4,#5]{{\ifcurr@ntPS\ifps@cri\s@uvc@ntr@l\et@tpsarrowcircPTD%
    \PSc@mment{psarrowcircPTD Center=#1; Radius=#2, [P1=#3,P2=#4,P3=#5]}%
    \figgetangleTD\@ngfin[#1,#3,#4,#5]\v@leur=#2pt%
    \maxim@m{\mili@u}{-\v@leur}{\v@leur}\edef\r@dius{\repdecn@mb{\mili@u}}%
    \ifdim\v@leur<\z@\v@lmax=\@ngfin pt\advance\v@lmax-\DePI@deg%
    \edef\@ngfin{\repdecn@mb{\v@lmax}}\fi\psarrowcircTD#1,#3,#5;\r@dius(0,\@ngfin)%
    \PSc@mment{End psarrowcircPTD}\resetc@ntr@l\et@tpsarrowcircPTD\fi\fi}}
\ctr@ld@f\def\psaxes#1(#2){{\ifcurr@ntPS\ifps@cri\s@uvc@ntr@l\et@tpsaxes%
    \PSc@mment{psaxes Origin=#1 Range=(#2)}\an@lys@xes#2,:\resetc@ntr@l{2}%
    \ifx\t@xt@\empty\ifTr@isDim\ps@xes#1(0,#2,0,#2,0,#2)\else\ps@xes#1(0,#2,0,#2)\fi%
    \else\ps@xes#1(#2)\fi\PSc@mment{End psaxes}\resetc@ntr@l\et@tpsaxes\fi\fi}}
\ctr@ld@f\def\an@lys@xes#1,#2:{\def\t@xt@{#2}}
\ctr@ln@m\ps@xes
\ctr@ld@f\def\ps@xesDD#1(#2,#3,#4,#5){%
    \figpttraC-5:=#1/#2,0/\figpttraC-6:=#1/#3,0/\psarrowDD[-5,-6]%
    \figpttraC-5:=#1/0,#4/\figpttraC-6:=#1/0,#5/\psarrowDD[-5,-6]}
\ctr@ld@f\def\ps@xesTD#1(#2,#3,#4,#5,#6,#7){%
    \figpttraC-7:=#1/#2,0,0/\figpttraC-8:=#1/#3,0,0/\psarrowTD[-7,-8]%
    \figpttraC-7:=#1/0,#4,0/\figpttraC-8:=#1/0,#5,0/\psarrowTD[-7,-8]%
    \figpttraC-7:=#1/0,0,#6/\figpttraC-8:=#1/0,0,#7/\psarrowTD[-7,-8]}
\ctr@ln@m\newGr@FN
\ctr@ld@f\def\newGr@FNPDF#1{\s@mme=\Gr@FNb\advance\s@mme\@ne\xdef\Gr@FNb{\number\s@mme}}
\ctr@ld@f\def\newGr@FNDVI#1{\newGr@FNPDF{}\xdef#1{\jobname GI\Gr@FNb.anx}}
\ctr@ld@f\def\psbeginfig#1{\newGr@FN\DefGIfilen@me\gdef\@utoFN{0}%
    \def\t@xt@{#1}\relax\ifx\t@xt@\empty\psupdatem@detrue%
    \gdef\@utoFN{1}\Psb@ginfig\DefGIfilen@me\else\expandafter\Psb@ginfigNu@#1 :\fi}
\ctr@ld@f\def\Psb@ginfigNu@#1 #2:{\def\t@xt@{#1}\relax\ifx\t@xt@\empty\def\t@xt@{#2}%
    \ifx\t@xt@\empty\psupdatem@detrue\gdef\@utoFN{1}\Psb@ginfig\DefGIfilen@me%
    \else\Psb@ginfigNu@#2:\fi\else\Psb@ginfig{#1}\fi}
\ctr@ln@m\PSfilen@me \ctr@ln@m\auxfilen@me
\ctr@ld@f\def\Psb@ginfig#1{\ifcurr@ntPS\else%
    \edef\PSfilen@me{#1}\edef\auxfilen@me{\jobname.anx}%
    \ifpsupdatem@de\ps@critrue\else\openin\frf@g=\PSfilen@me\relax%
    \ifeof\frf@g\ps@critrue\else\ps@crifalse\fi\closein\frf@g\fi%
    \curr@ntPStrue\c@ldefproj\expandafter\setupd@te\defaultupdate:%
    \ifps@cri\initb@undb@x%
    \immediate\openout\fwf@g=\auxfilen@me\initpss@ttings\fi%
    \fi}
\ctr@ld@f\def\Gr@FNb{0}
\ctr@ld@f\def\figforTeXFileno{\Gr@FNb}
\ctr@ld@f\def\figforTeXFigno{0 }
\ctr@ld@f\def\figforTeXnextFigno{1 }
\ctr@ld@f\edef\DefGIfilen@me{\jobname GI.anx}
\ctr@ld@f\def\initpss@ttings{\psreset{arrowhead,curve,first,flowchart,mesh,second,third}%
    \Use@llipsefalse}
\ctr@ld@f\def\B@zierBB@x#1#2(#3,#4,#5,#6){{\c@rre=\t@n\epsil@n
    \v@lmax=#4\advance\v@lmax-#5\v@lmax=\thr@@\v@lmax\advance\v@lmax#6\advance\v@lmax-#3%
    \mili@u=#4\mili@u=-\tw@\mili@u\advance\mili@u#3\advance\mili@u#5%
    \v@lmin=#4\advance\v@lmin-#3\maxim@m{\v@leur}{-\v@lmax}{\v@lmax}%
    \maxim@m{\delt@}{-\mili@u}{\mili@u}\maxim@m{\v@leur}{\v@leur}{\delt@}%
    \maxim@m{\delt@}{-\v@lmin}{\v@lmin}\maxim@m{\v@leur}{\v@leur}{\delt@}%
    \ifdim\v@leur>\c@rre\invers@{\v@leur}{\v@leur}\edef\Uns@rM@x{\repdecn@mb{\v@leur}}%
    \v@lmax=\Uns@rM@x\v@lmax\mili@u=\Uns@rM@x\mili@u\v@lmin=\Uns@rM@x\v@lmin%
    \maxim@m{\v@leur}{-\v@lmax}{\v@lmax}\ifdim\v@leur<\c@rre%
    \maxim@m{\v@leur}{-\mili@u}{\mili@u}\ifdim\v@leur<\c@rre\else%
    \invers@{\mili@u}{\mili@u}\v@leur=-0.5\v@lmin%
    \v@leur=\repdecn@mb{\mili@u}\v@leur\m@jBBB@x{\v@leur}{#1}{#2}(#3,#4,#5,#6)\fi%
    \else\delt@=\repdecn@mb{\mili@u}\mili@u\v@leur=\repdecn@mb{\v@lmax}\v@lmin%
    \advance\delt@-\v@leur\ifdim\delt@<\z@\else\invers@{\v@lmax}{\v@lmax}%
    \edef\Uns@rAp{\repdecn@mb{\v@lmax}}\sqrt@{\delt@}{\delt@}%
    \v@leur=-\mili@u\advance\v@leur\delt@\v@leur=\Uns@rAp\v@leur%
    \m@jBBB@x{\v@leur}{#1}{#2}(#3,#4,#5,#6)%
    \v@leur=-\mili@u\advance\v@leur-\delt@\v@leur=\Uns@rAp\v@leur%
    \m@jBBB@x{\v@leur}{#1}{#2}(#3,#4,#5,#6)\fi\fi\fi}}
\ctr@ld@f\def\m@jBBB@x#1#2#3(#4,#5,#6,#7){{\relax\ifdim#1>\z@\ifdim#1<\p@%
    \edef\T@{\repdecn@mb{#1}}\v@lX=\p@\advance\v@lX-#1\edef\UNmT@{\repdecn@mb{\v@lX}}%
    \v@lX=#4\v@lY=#5\v@lZ=#6\v@lXa=#7\v@lX=\UNmT@\v@lX\advance\v@lX\T@\v@lY%
    \v@lY=\UNmT@\v@lY\advance\v@lY\T@\v@lZ\v@lZ=\UNmT@\v@lZ\advance\v@lZ\T@\v@lXa%
    \v@lX=\UNmT@\v@lX\advance\v@lX\T@\v@lY\v@lY=\UNmT@\v@lY\advance\v@lY\T@\v@lZ%
    \v@lX=\UNmT@\v@lX\advance\v@lX\T@\v@lY%
    \ifcase#2\or\v@lY=#3\or\v@lY=\v@lX\v@lX=#3\fi\b@undb@x{\v@lX}{\v@lY}\fi\fi}}
\ctr@ld@f\def\PsB@zier#1[#2]{{\f@gnewpath%
    \s@mme=\z@\def\list@num{#2,0}\extrairelepremi@r\p@int\de\list@num%
    \PSwrit@cmdS{\p@int}{\c@mmoveto}{\fwf@g}{\X@un}{\Y@un}\p@rtent=#1\bclB@zier}}
\ctr@ld@f\def\bclB@zier{\relax%
    \ifnum\s@mme<\p@rtent\advance\s@mme\@ne\BdingB@xfalse%
    \extrairelepremi@r\p@int\de\list@num\PSwrit@cmdS{\p@int}{}{\fwf@g}{\X@de}{\Y@de}%
    \extrairelepremi@r\p@int\de\list@num\PSwrit@cmdS{\p@int}{}{\fwf@g}{\X@tr}{\Y@tr}%
    \BdingB@xtrue%
    \extrairelepremi@r\p@int\de\list@num\PSwrit@cmdS{\p@int}{\c@mcurveto}{\fwf@g}{\X@qu}{\Y@qu}%
    \B@zierBB@x{1}{\Y@un}(\X@un,\X@de,\X@tr,\X@qu)%
    \B@zierBB@x{2}{\X@un}(\Y@un,\Y@de,\Y@tr,\Y@qu)%
    \edef\X@un{\X@qu}\edef\Y@un{\Y@qu}\bclB@zier\fi}
\ctr@ln@m\psBezier
\ctr@ld@f\def\psBezierDD#1[#2]{\ifcurr@ntPS\ifps@cri%
    \PSc@mment{psBezierDD N arcs=#1, Control points=#2}%
    \iffillm@de\PsB@zier#1[#2]%
    \f@gfill%
    \else\PsB@zier#1[#2]\f@gstroke\fi%
    \PSc@mment{End psBezierDD}\fi\fi}
\ctr@ln@m\et@tpsBezierTD
\ctr@ld@f\def\psBezierTD#1[#2]{\ifcurr@ntPS\ifps@cri\s@uvc@ntr@l\et@tpsBezierTD%
    \PSc@mment{psBezierTD N arcs=#1, Control points=#2}%
    \iffillm@de\PsB@zierTD#1[#2]%
    \f@gfill%
    \else\PsB@zierTD#1[#2]\f@gstroke\fi%
    \PSc@mment{End psBezierTD}\resetc@ntr@l\et@tpsBezierTD\fi\fi}
\ctr@ld@f\def\PsB@zierTD#1[#2]{\ifnum\curr@ntproj<\tw@\PsB@zier#1[#2]\else\PsB@zier@TD#1[#2]\fi}
\ctr@ld@f\def\PsB@zier@TD#1[#2]{{\f@gnewpath%
    \s@mme=\z@\def\list@num{#2,0}\extrairelepremi@r\p@int\de\list@num%
    \let\c@lprojSP=\relax\setc@ntr@l{2}\Figptpr@j-7:/\p@int/%
    \PSwrit@cmd{-7}{\c@mmoveto}{\fwf@g}%
    \loop\ifnum\s@mme<#1\advance\s@mme\@ne\extrairelepremi@r\p@intun\de\list@num%
    \extrairelepremi@r\p@intde\de\list@num\extrairelepremi@r\p@inttr\de\list@num%
    \subB@zierTD\NBz@rcs[\p@int,\p@intun,\p@intde,\p@inttr]\edef\p@int{\p@inttr}\repeat}}
\ctr@ld@f\def\subB@zierTD#1[#2,#3,#4,#5]{\delt@=\p@\divide\delt@\NBz@rcs\v@lmin=\z@%
    {\Figg@tXY{-7}\edef\X@un{\the\v@lX}\edef\Y@un{\the\v@lY}%
    \s@mme=\z@\loop\ifnum\s@mme<#1\advance\s@mme\@ne%
    \v@leur=\v@lmin\advance\v@leur0.33333 \delt@\edef\unti@rs{\repdecn@mb{\v@leur}}%
    \v@leur=\v@lmin\advance\v@leur0.66666 \delt@\edef\deti@rs{\repdecn@mb{\v@leur}}%
    \advance\v@lmin\delt@\edef\trti@rs{\repdecn@mb{\v@lmin}}%
    \figptBezierTD-8::\trti@rs[#2,#3,#4,#5]\Figptpr@j-8:/-8/%
    \c@lsubBzarc\unti@rs,\deti@rs[#2,#3,#4,#5]\BdingB@xfalse%
    \PSwrit@cmdS{-4}{}{\fwf@g}{\X@de}{\Y@de}\PSwrit@cmdS{-3}{}{\fwf@g}{\X@tr}{\Y@tr}%
    \BdingB@xtrue\PSwrit@cmdS{-8}{\c@mcurveto}{\fwf@g}{\X@qu}{\Y@qu}%
    \B@zierBB@x{1}{\Y@un}(\X@un,\X@de,\X@tr,\X@qu)%
    \B@zierBB@x{2}{\X@un}(\Y@un,\Y@de,\Y@tr,\Y@qu)%
    \edef\X@un{\X@qu}\edef\Y@un{\Y@qu}\figptcopyDD-7:/-8/\repeat}}
\ctr@ld@f\def\NBz@rcs{2}
\ctr@ld@f\def\c@lsubBzarc#1,#2[#3,#4,#5,#6]{\figptBezierTD-5::#1[#3,#4,#5,#6]%
    \figptBezierTD-6::#2[#3,#4,#5,#6]\Figptpr@j-4:/-5/\Figptpr@j-5:/-6/%
    \figptscontrolDD-4[-7,-4,-5,-8]}
\ctr@ln@m\pscirc
\ctr@ld@f\def\pscircDD#1(#2){\ifcurr@ntPS\ifps@cri\PSc@mment{pscircDD Center=#1 (Radius=#2)}%
    \psarccircDD#1;#2(0,360)\PSc@mment{End pscircDD}\fi\fi}
\ctr@ld@f\def\pscircTD#1,#2,#3(#4){\ifcurr@ntPS\ifps@cri%
    \PSc@mment{pscircTD Center=#1,P1=#2,P2=#3 (Radius=#4)}%
    \psarccircTD#1,#2,#3;#4(0,360)\PSc@mment{End pscircTD}\fi\fi}
\ctr@ln@m\p@urcent
{\catcode`\%=12\gdef\p@urcent{
\ctr@ld@f\def\PSc@mment#1{\ifpsdebugmode\immediate\write\fwf@g{\p@urcent\space#1}\fi}
\ctr@ln@m\acc@louv \ctr@ln@m\acc@lfer
{\catcode`\[=1\catcode`\{=12\gdef\acc@louv[{}}
{\catcode`\]=2\catcode`\}=12\gdef\acc@lfer{}]]
\ctr@ld@f\def\PSdict@{\ifUse@llipse%
    \immediate\write\fwf@g{/ellipsedict 9 dict def ellipsedict /mtrx matrix put}%
    \immediate\write\fwf@g{/ellipse \acc@louv ellipsedict begin}%
    \immediate\write\fwf@g{ /endangle exch def /startangle exch def}%
    \immediate\write\fwf@g{ /yrad exch def /xrad exch def}%
    \immediate\write\fwf@g{ /rotangle exch def /y exch def /x exch def}%
    \immediate\write\fwf@g{ /savematrix mtrx currentmatrix def}%
    \immediate\write\fwf@g{ x y translate rotangle rotate xrad yrad scale}%
    \immediate\write\fwf@g{ 0 0 1 startangle endangle arc}%
    \immediate\write\fwf@g{ savematrix setmatrix end\acc@lfer def}%
    \fi\PShe@der{EndProlog}}
\ctr@ld@f\def\Pssetc@rve#1=#2|{\keln@mun#1|%
    \def\n@mref{r}\ifx\l@debut\n@mref\pssetroundness{#2}\else
    \immediate\write16{*** Unknown attribute: \BS@ psset curve(..., #1=...)}%
    \fi}
\ctr@ln@m\curv@roundness
\ctr@ld@f\def\pssetroundness#1{\edef\curv@roundness{#1}}
\ctr@ld@f\def\defaultroundness{0.2} 
\ctr@ln@m\pscurve
\ctr@ld@f\def\pscurveDD[#1]{{\ifcurr@ntPS\ifps@cri\PSc@mment{pscurveDD Points=#1}%
    \s@uvc@ntr@l\et@tpscurveDD%
    \iffillm@de\Psc@rveDD\curv@roundness[#1]%
    \f@gfill%
    \else\Psc@rveDD\curv@roundness[#1]\f@gstroke\fi%
    \PSc@mment{End pscurveDD}\resetc@ntr@l\et@tpscurveDD\fi\fi}}
\ctr@ld@f\def\pscurveTD[#1]{{\ifcurr@ntPS\ifps@cri%
    \PSc@mment{pscurveTD Points=#1}\s@uvc@ntr@l\et@tpscurveTD\let\c@lprojSP=\relax%
    \iffillm@de\Psc@rveTD\curv@roundness[#1]%
    \f@gfill%
    \else\Psc@rveTD\curv@roundness[#1]\f@gstroke\fi%
    \PSc@mment{End pscurveTD}\resetc@ntr@l\et@tpscurveTD\fi\fi}}
\ctr@ld@f\def\Psc@rveDD#1[#2]{%
    \def\list@num{#2}\extrairelepremi@r\Ak@\de\list@num%
    \extrairelepremi@r\Ai@\de\list@num\extrairelepremi@r\Aj@\de\list@num%
    \f@gnewpath\PSwrit@cmdS{\Ai@}{\c@mmoveto}{\fwf@g}{\X@un}{\Y@un}%
    \setc@ntr@l{2}\figvectPDD -1[\Ak@,\Aj@]%
    \@ecfor\Ak@:=\list@num\do{\figpttraDD-2:=\Ai@/#1,-1/\BdingB@xfalse%
       \PSwrit@cmdS{-2}{}{\fwf@g}{\X@de}{\Y@de}%
       \figvectPDD -1[\Ai@,\Ak@]\figpttraDD-2:=\Aj@/-#1,-1/%
       \PSwrit@cmdS{-2}{}{\fwf@g}{\X@tr}{\Y@tr}\BdingB@xtrue%
       \PSwrit@cmdS{\Aj@}{\c@mcurveto}{\fwf@g}{\X@qu}{\Y@qu}%
       \B@zierBB@x{1}{\Y@un}(\X@un,\X@de,\X@tr,\X@qu)%
       \B@zierBB@x{2}{\X@un}(\Y@un,\Y@de,\Y@tr,\Y@qu)%
       \edef\X@un{\X@qu}\edef\Y@un{\Y@qu}\edef\Ai@{\Aj@}\edef\Aj@{\Ak@}}}
\ctr@ld@f\def\Psc@rveTD#1[#2]{\ifnum\curr@ntproj<\tw@\Psc@rvePPTD#1[#2]\else\Psc@rveCPTD#1[#2]\fi}
\ctr@ld@f\def\Psc@rvePPTD#1[#2]{\setc@ntr@l{2}%
    \def\list@num{#2}\extrairelepremi@r\Ak@\de\list@num\Figptpr@j-5:/\Ak@/%
    \extrairelepremi@r\Ai@\de\list@num\Figptpr@j-3:/\Ai@/%
    \extrairelepremi@r\Aj@\de\list@num\Figptpr@j-4:/\Aj@/%
    \f@gnewpath\PSwrit@cmdS{-3}{\c@mmoveto}{\fwf@g}{\X@un}{\Y@un}%
    \figvectPDD -1[-5,-4]%
    \@ecfor\Ak@:=\list@num\do{\Figptpr@j-5:/\Ak@/\figpttraDD-2:=-3/#1,-1/%
       \BdingB@xfalse\PSwrit@cmdS{-2}{}{\fwf@g}{\X@de}{\Y@de}%
       \figvectPDD -1[-3,-5]\figpttraDD-2:=-4/-#1,-1/%
       \PSwrit@cmdS{-2}{}{\fwf@g}{\X@tr}{\Y@tr}\BdingB@xtrue%
       \PSwrit@cmdS{-4}{\c@mcurveto}{\fwf@g}{\X@qu}{\Y@qu}%
       \B@zierBB@x{1}{\Y@un}(\X@un,\X@de,\X@tr,\X@qu)%
       \B@zierBB@x{2}{\X@un}(\Y@un,\Y@de,\Y@tr,\Y@qu)%
       \edef\X@un{\X@qu}\edef\Y@un{\Y@qu}\figptcopyDD-3:/-4/\figptcopyDD-4:/-5/}}
\ctr@ld@f\def\Psc@rveCPTD#1[#2]{\setc@ntr@l{2}%
    \def\list@num{#2}\extrairelepremi@r\Ak@\de\list@num%
    \extrairelepremi@r\Ai@\de\list@num\extrairelepremi@r\Aj@\de\list@num%
    \Figptpr@j-7:/\Ai@/%
    \f@gnewpath\PSwrit@cmd{-7}{\c@mmoveto}{\fwf@g}%
    \figvectPTD -9[\Ak@,\Aj@]%
    \@ecfor\Ak@:=\list@num\do{\figpttraTD-10:=\Ai@/#1,-9/%
       \figvectPTD -9[\Ai@,\Ak@]\figpttraTD-11:=\Aj@/-#1,-9/%
       \subB@zierTD\NBz@rcs[\Ai@,-10,-11,\Aj@]\edef\Ai@{\Aj@}\edef\Aj@{\Ak@}}}
\ctr@ld@f\def\psendfig{\ifcurr@ntPS\ifps@cri\immediate\closeout\fwf@g%
    \immediate\openout\fwf@g=\PSfilen@me\relax%
    \ifPDFm@ke\PSBdingB@x\else%
    \immediate\write\fwf@g{\p@urcent\string!PS-Adobe-2.0 EPSF-2.0}%
    \PShe@der{Creator\string: TeX (fig4tex.tex)}%
    \PShe@der{Title\string: \PSfilen@me}%
    \PShe@der{CreationDate\string: \the\day/\the\month/\the\year}%
    \PSBdingB@x%
    \PShe@der{EndComments}\PSdict@\fi%
    \immediate\write\fwf@g{\c@mgsave}%
    \openin\frf@g=\auxfilen@me\c@pypsfile\fwf@g\frf@g\closein\frf@g%
    \immediate\write\fwf@g{\c@mgrestore}%
    \PSc@mment{End of file.}\immediate\closeout\fwf@g%
    \immediate\openout\fwf@g=\auxfilen@me\immediate\closeout\fwf@g%
    \immediate\write16{File \PSfilen@me\space created.}\fi\fi\curr@ntPSfalse\ps@critrue}
\ctr@ld@f\def\PShe@der#1{\immediate\write\fwf@g{\p@urcent\p@urcent#1}}
\ctr@ld@f\def\PSBdingB@x{{\v@lX=\ptT@ptps\c@@rdXmin\v@lY=\ptT@ptps\c@@rdYmin%
     \v@lXa=\ptT@ptps\c@@rdXmax\v@lYa=\ptT@ptps\c@@rdYmax%
     \PShe@der{BoundingBox\string: \repdecn@mb{\v@lX}\space\repdecn@mb{\v@lY}%
     \space\repdecn@mb{\v@lXa}\space\repdecn@mb{\v@lYa}}}}
\ctr@ld@f\def\psfcconnect[#1]{{\ifcurr@ntPS\ifps@cri\PSc@mment{psfcconnect Points=#1}%
    \pssetfillmode{no}\s@uvc@ntr@l\et@tpsfcconnect\resetc@ntr@l{2}%
    \fcc@nnect@[#1]\resetc@ntr@l\et@tpsfcconnect\PSc@mment{End psfcconnect}\fi\fi}}
\ctr@ld@f\def\fcc@nnect@[#1]{\let\N@rm=\n@rmeucDD\def\list@num{#1}%
    \extrairelepremi@r\Ai@\de\list@num\edef\pr@m{\Ai@}\v@leur=\z@\p@rtent=\@ne\c@llgtot%
    \ifcase\fclin@typ@\edef\list@num{[\pr@m,#1,\Ai@}\expandafter\pscurve\list@num]%
    \else\ifdim\fclin@r@d\p@>\z@\Pslin@conge[#1]\else\psline[#1]\fi\fi%
    \v@leur=\@rrowp@s\v@leur\edef\list@num{#1,\Ai@,0}%
    \extrairelepremi@r\Ai@\de\list@num\mili@u=\epsil@n\c@llgpart%
    \advance\mili@u-\epsil@n\advance\mili@u-\delt@\advance\v@leur-\mili@u%
    \ifcase\fclin@typ@\invers@\mili@u\delt@%
    \ifnum\@rrowr@fpt>\z@\advance\delt@-\v@leur\v@leur=\delt@\fi%
    \v@leur=\repdecn@mb\v@leur\mili@u\edef\v@lt{\repdecn@mb\v@leur}%
    \extrairelepremi@r\Ak@\de\list@num%
    \figvectPDD-1[\pr@m,\Aj@]\figpttraDD-6:=\Ai@/\curv@roundness,-1/%
    \figvectPDD-1[\Ak@,\Ai@]\figpttraDD-7:=\Aj@/\curv@roundness,-1/%
    \delt@=\@rrowheadlength\p@\delt@=\C@AHANG\delt@\edef\R@dius{\repdecn@mb{\delt@}}%
    \ifcase\@rrowr@fpt%
    \FigptintercircB@zDD-8::\v@lt,\R@dius[\Ai@,-6,-7,\Aj@]\psarrowheadDD[-5,-8]\else%
    \FigptintercircB@zDD-8::\v@lt,\R@dius[\Aj@,-7,-6,\Ai@]\psarrowheadDD[-8,-5]\fi%
    \else\advance\delt@-\v@leur%
    \p@rtentiere{\p@rtent}{\delt@}\edef\C@efun{\the\p@rtent}%
    \p@rtentiere{\p@rtent}{\v@leur}\edef\C@efde{\the\p@rtent}%
    \figptbaryDD-5:[\Ai@,\Aj@;\C@efun,\C@efde]\ifcase\@rrowr@fpt%
    \delt@=\@rrowheadlength\unit@\delt@=\C@AHANG\delt@\edef\t@ille{\repdecn@mb{\delt@}}%
    \figvectPDD-2[\Ai@,\Aj@]\vecunit@{-2}{-2}\figpttraDD-5:=-5/\t@ille,-2/\fi%
    \psarrowheadDD[\Ai@,-5]\fi}
\ctr@ld@f\def\c@llgtot{\@ecfor\Aj@:=\list@num\do{\figvectP-1[\Ai@,\Aj@]\N@rm\delt@{-1}%
    \advance\v@leur\delt@\advance\p@rtent\@ne\edef\Ai@{\Aj@}}}
\ctr@ld@f\def\c@llgpart{\extrairelepremi@r\Aj@\de\list@num\figvectP-1[\Ai@,\Aj@]\N@rm\delt@{-1}%
    \advance\mili@u\delt@\ifdim\mili@u<\v@leur\edef\pr@m{\Ai@}\edef\Ai@{\Aj@}\c@llgpart\fi}
\ctr@ld@f\def\Pslin@conge[#1]{\ifnum\p@rtent>\tw@{\def\list@num{#1}%
    \extrairelepremi@r\Ai@\de\list@num\extrairelepremi@r\Aj@\de\list@num%
    \figptcopy-6:/\Ai@/\figvectP-3[\Ai@,\Aj@]\vecunit@{-3}{-3}\v@lmax=\result@t%
    \@ecfor\Ak@:=\list@num\do{\figvectP-4[\Aj@,\Ak@]\vecunit@{-4}{-4}%
    \minim@m\v@lmin\v@lmax\result@t\v@lmax=\result@t%
    \det@rm\delt@[-3,-4]\maxim@m\mili@u{\delt@}{-\delt@}\ifdim\mili@u>\Cepsil@n%
    \ifdim\delt@>\z@\figgetangleDD\Angl@[\Aj@,\Ak@,\Ai@]\else%
    \figgetangleDD\Angl@[\Aj@,\Ai@,\Ak@]\fi%
    \v@leur=\PI@deg\advance\v@leur-\Angl@\p@\divide\v@leur\tw@%
    \edef\Angl@{\repdecn@mb\v@leur}\c@ssin{\C@}{\S@}{\Angl@}\v@leur=\fclin@r@d\unit@%
    \v@leur=\S@\v@leur\mili@u=\C@\p@\invers@\mili@u\mili@u%
    \v@leur=\repdecn@mb{\mili@u}\v@leur%
    \minim@m\v@leur\v@leur\v@lmin\edef\t@ille{\repdecn@mb{\v@leur}}%
    \figpttra-5:=\Aj@/-\t@ille,-3/\psline[-6,-5]\figpttra-6:=\Aj@/\t@ille,-4/%
    \figvectNVDD-3[-3]\figvectNVDD-8[-4]\inters@cDD-7:[-5,-3;-6,-8]%
    \ifdim\delt@>\z@\psarccircP-7;\fclin@r@d[-5,-6]\else\psarccircP-7;\fclin@r@d[-6,-5]\fi%
    \else\psline[-6,\Aj@]\figptcopy-6:/\Aj@/\fi
    \edef\Ai@{\Aj@}\edef\Aj@{\Ak@}\figptcopy-3:/-4/}\psline[-6,\Aj@]}\else\psline[#1]\fi}
\ctr@ld@f\def\psfcnode[#1]#2{{\ifcurr@ntPS\ifps@cri\PSc@mment{psfcnode Points=#1}%
    \s@uvc@ntr@l\et@tpsfcnode\resetc@ntr@l{2}%
    \def\t@xt@{#2}\ifx\t@xt@\empty\def\g@tt@xt{\setbox\Gb@x=\hbox{\Figg@tT{\p@int}}}%
    \else\def\g@tt@xt{\setbox\Gb@x=\hbox{#2}}\fi%
    \v@lmin=\h@rdfcXp@dd\advance\v@lmin\Xp@dd\unit@\multiply\v@lmin\tw@%
    \v@lmax=\h@rdfcYp@dd\advance\v@lmax\Yp@dd\unit@\multiply\v@lmax\tw@%
    \Figv@ctCreg-8(\unit@,-\unit@)\def\list@num{#1}%
    \delt@=\curr@ntwidth bp\divide\delt@\tw@%
    \fcn@de\PSc@mment{End psfcnode}\resetc@ntr@l\et@tpsfcnode\fi\fi}}
\ctr@ld@f\def\d@butn@de{\g@tt@xt\v@lX=\wd\Gb@x%
    \v@lY=\ht\Gb@x\advance\v@lY\dp\Gb@x\advance\v@lX\v@lmin\advance\v@lY\v@lmax}
\ctr@ld@f\def\fcn@deE{%
    \@ecfor\p@int:=\list@num\do{\d@butn@de\v@lX=\unssqrttw@\v@lX\v@lY=\unssqrttw@\v@lY%
    \ifdim\thickn@ss\p@>\z@
    \v@lXa=\v@lX\advance\v@lXa\delt@\v@lXa=\ptT@unit@\v@lXa\edef\XR@d{\repdecn@mb\v@lXa}%
    \v@lYa=\v@lY\advance\v@lYa\delt@\v@lYa=\ptT@unit@\v@lYa\edef\YR@d{\repdecn@mb\v@lYa}%
    \arct@n\v@leur(\v@lXa,\v@lYa)\v@leur=\rdT@deg\v@leur\edef\@nglde{\repdecn@mb\v@leur}%
    {\c@lptellDD-2::\p@int;\XR@d,\YR@d(\@nglde)}
    \advance\v@leur-\PI@deg\edef\@nglun{\repdecn@mb\v@leur}%
    {\c@lptellDD-3::\p@int;\XR@d,\YR@d(\@nglun)}%
    \figptstra-6=-3,-2,\p@int/\thickn@ss,-8/\pssetfillmode{yes}\us@secondC@lor%
    \psline[-2,-3,-6,-5]\psarcell-4;\XR@d,\YR@d(\@nglun,\@nglde,0)\fi
    \v@lX=\ptT@unit@\v@lX\v@lY=\ptT@unit@\v@lY%
    \edef\XR@d{\repdecn@mb\v@lX}\edef\YR@d{\repdecn@mb\v@lY}%
    \pssetfillmode{yes}\us@thirdC@lor\psarcell\p@int;\XR@d,\YR@d(0,360,0)%
    \pssetfillmode{no}\us@primarC@lor\psarcell\p@int;\XR@d,\YR@d(0,360,0)}}
\ctr@ld@f\def\fcn@deL{\delt@=\ptT@unit@\delt@\edef\t@ille{\repdecn@mb\delt@}%
    \@ecfor\p@int:=\list@num\do{\Figg@tXYa{\p@int}\d@butn@de%
    \ifdim\v@lX>\v@lY\itis@Ktrue\else\itis@Kfalse\fi%
    \advance\v@lXa-\v@lX\Figp@intreg-1:(\v@lXa,\v@lYa)%
    \advance\v@lXa\v@lX\advance\v@lYa-\v@lY\Figp@intreg-2:(\v@lXa,\v@lYa)%
    \advance\v@lXa\v@lX\advance\v@lYa\v@lY\Figp@intreg-3:(\v@lXa,\v@lYa)%
    \advance\v@lXa-\v@lX\advance\v@lYa\v@lY\Figp@intreg-4:(\v@lXa,\v@lYa)%
    \ifdim\thickn@ss\p@>\z@\Figg@tXYa{\p@int}\pssetfillmode{yes}\us@secondC@lor
    \c@lpt@xt{-1}{-4}\c@lpt@xt@\v@lXa\v@lYa\v@lX\v@lY\c@rre\delt@%
    \Figp@intregDD-9:(\v@lZ,\v@lYa)\Figp@intregDD-11:(\v@lZa,\v@lYa)%
    \c@lpt@xt{-4}{-3}\c@lpt@xt@\v@lYa\v@lXa\v@lY\v@lX\delt@\c@rre%
    \Figp@intregDD-12:(\v@lXa,\v@lZ)\Figp@intregDD-10:(\v@lXa,\v@lZa)%
    \ifitis@K\figptstra-7=-9,-10,-11/\thickn@ss,-8/\psline[-9,-11,-5,-6,-7]\else%
    \figptstra-7=-10,-11,-12/\thickn@ss,-8/\psline[-10,-12,-5,-6,-7]\fi\fi
    \pssetfillmode{yes}\us@thirdC@lor\psline[-1,-2,-3,-4]%
    \pssetfillmode{no}\us@primarC@lor\psline[-1,-2,-3,-4,-1]}}
\ctr@ld@f\def\c@lpt@xt#1#2{\figvectN-7[#1,#2]\vecunit@{-7}{-7}\figpttra-5:=#1/\t@ille,-7/%
    \figvectP-7[#1,#2]\Figg@tXY{-7}\c@rre=\v@lX\delt@=\v@lY\Figg@tXY{-5}}
\ctr@ld@f\def\c@lpt@xt@#1#2#3#4#5#6{\v@lZ=#6\invers@{\v@lZ}{\v@lZ}\v@leur=\repdecn@mb{#5}\v@lZ%
    \v@lZ=#2\advance\v@lZ-#4\mili@u=\repdecn@mb{\v@leur}\v@lZ%
    \v@lZ=#3\advance\v@lZ\mili@u\v@lZa=-\v@lZ\advance\v@lZa\tw@#1}
\ctr@ld@f\def\fcn@deR{\@ecfor\p@int:=\list@num\do{\Figg@tXYa{\p@int}\d@butn@de%
    \advance\v@lXa-0.5\v@lX\advance\v@lYa-0.5\v@lY\Figp@intreg-1:(\v@lXa,\v@lYa)%
    \advance\v@lXa\v@lX\Figp@intreg-2:(\v@lXa,\v@lYa)%
    \advance\v@lYa\v@lY\Figp@intreg-3:(\v@lXa,\v@lYa)%
    \advance\v@lXa-\v@lX\Figp@intreg-4:(\v@lXa,\v@lYa)%
    \ifdim\thickn@ss\p@>\z@\pssetfillmode{yes}\us@secondC@lor
    \Figv@ctCreg-5(-\delt@,-\delt@)\figpttra-9:=-1/1,-5/%
    \Figv@ctCreg-5(\delt@,-\delt@)\figpttra-10:=-2/1,-5/%
    \Figv@ctCreg-5(\delt@,\delt@)\figpttra-11:=-3/1,-5/%
    \figptstra-7=-9,-10,-11/\thickn@ss,-8/\psline[-9,-11,-5,-6,-7]\fi
    \pssetfillmode{yes}\us@thirdC@lor\psline[-1,-2,-3,-4]%
    \pssetfillmode{no}\us@primarC@lor\psline[-1,-2,-3,-4,-1]}}
\ctr@ln@m\@rrowp@s
\ctr@ln@m\Xp@dd     \ctr@ln@m\Yp@dd
\ctr@ln@m\fclin@r@d \ctr@ln@m\thickn@ss
\ctr@ld@f\def\Pssetfl@wchart#1=#2|{\keln@mtr#1|%
    \def\n@mref{arr}\ifx\l@debut\n@mref\expandafter\keln@mtr\l@suite|%
     \def\n@mref{owp}\ifx\l@debut\n@mref\edef\@rrowp@s{#2}\else
     \def\n@mref{owr}\ifx\l@debut\n@mref\setfcr@fpt#2|\else
     \immediate\write16{*** Unknown attribute: \BS@ psset flowchart(..., #1=...)}%
     \fi\fi\else%
    \def\n@mref{lin}\ifx\l@debut\n@mref\setfccurv@#2|\else
    \def\n@mref{pad}\ifx\l@debut\n@mref\edef\Xp@dd{#2}\edef\Yp@dd{#2}\else
    \def\n@mref{rad}\ifx\l@debut\n@mref\edef\fclin@r@d{#2}\else
    \def\n@mref{sha}\ifx\l@debut\n@mref\setfcshap@#2|\else
    \def\n@mref{thi}\ifx\l@debut\n@mref\edef\thickn@ss{#2}\else
    \def\n@mref{xpa}\ifx\l@debut\n@mref\edef\Xp@dd{#2}\else
    \def\n@mref{ypa}\ifx\l@debut\n@mref\edef\Yp@dd{#2}\else
    \immediate\write16{*** Unknown attribute: \BS@ psset flowchart(..., #1=...)}%
    \fi\fi\fi\fi\fi\fi\fi\fi}
\ctr@ln@m\@rrowr@fpt \ctr@ln@m\fclin@typ@
\ctr@ld@f\def\setfcr@fpt#1#2|{\if#1e\def\@rrowr@fpt{1}\else\def\@rrowr@fpt{0}\fi}
\ctr@ld@f\def\setfccurv@#1#2|{\if#1c\def\fclin@typ@{0}\else\def\fclin@typ@{1}\fi}
\ctr@ln@m\h@rdfcXp@dd \ctr@ln@m\h@rdfcYp@dd
\ctr@ln@m\fcn@de \ctr@ln@m\fcsh@pe
\ctr@ld@f\def\setfcshap@#1#2|{%
    \if#1e\let\fcn@de=\fcn@deE\def\h@rdfcXp@dd{4pt}\def\h@rdfcYp@dd{4pt}%
     \edef\fcsh@pe{ellipse}\else%
    \if#1l\let\fcn@de=\fcn@deL\def\h@rdfcXp@dd{4pt}\def\h@rdfcYp@dd{4pt}%
     \edef\fcsh@pe{lozenge}\else%
          \let\fcn@de=\fcn@deR\def\h@rdfcXp@dd{6pt}\def\h@rdfcYp@dd{6pt}%
     \edef\fcsh@pe{rectangle}\fi\fi}
\ctr@ld@f\def\psline[#1]{{\ifcurr@ntPS\ifps@cri\PSc@mment{psline Points=#1}%
    \let\pslign@=\Pslign@P\Pslin@{#1}\PSc@mment{End psline}\fi\fi}}
\ctr@ld@f\def\pslineF#1{{\ifcurr@ntPS\ifps@cri\PSc@mment{pslineF Filename=#1}%
    \let\pslign@=\Pslign@F\Pslin@{#1}\PSc@mment{End pslineF}\fi\fi}}
\ctr@ld@f\def\pslineC(#1){{\ifcurr@ntPS\ifps@cri\PSc@mment{pslineC}%
    \let\pslign@=\Pslign@C\Pslin@{#1}\PSc@mment{End pslineC}\fi\fi}}
\ctr@ld@f\def\Pslin@#1{\iffillm@de\pslign@{#1}%
    \f@gfill%
    \else\pslign@{#1}\ifx\derp@int\premp@int%
    \f@gclosestroke%
    \else\f@gstroke\fi\fi}
\ctr@ld@f\def\Pslign@P#1{\def\list@num{#1}\extrairelepremi@r\p@int\de\list@num%
    \edef\premp@int{\p@int}\f@gnewpath%
    \PSwrit@cmd{\p@int}{\c@mmoveto}{\fwf@g}%
    \@ecfor\p@int:=\list@num\do{\PSwrit@cmd{\p@int}{\c@mlineto}{\fwf@g}%
    \edef\derp@int{\p@int}}}
\ctr@ld@f\def\Pslign@F#1{\s@uvc@ntr@l\et@tPslign@F\setc@ntr@l{2}\openin\frf@g=#1\relax%
    \ifeof\frf@g\message{*** File #1 not found !}\end\else%
    \read\frf@g to\tr@c\edef\premp@int{\tr@c}\expandafter\extr@ctCF\tr@c:%
    \f@gnewpath\PSwrit@cmd{-1}{\c@mmoveto}{\fwf@g}%
    \loop\read\frf@g to\tr@c\ifeof\frf@g\mored@tafalse\else\mored@tatrue\fi%
    \ifmored@ta\expandafter\extr@ctCF\tr@c:\PSwrit@cmd{-1}{\c@mlineto}{\fwf@g}%
    \edef\derp@int{\tr@c}\repeat\fi\closein\frf@g\resetc@ntr@l\et@tPslign@F}
\ctr@ln@m\extr@ctCF
\ctr@ld@f\def\extr@ctCFDD#1 #2:{\v@lX=#1\unit@\v@lY=#2\unit@\Figp@intregDD-1:(\v@lX,\v@lY)}
\ctr@ld@f\def\extr@ctCFTD#1 #2 #3:{\v@lX=#1\unit@\v@lY=#2\unit@\v@lZ=#3\unit@%
    \Figp@intregTD-1:(\v@lX,\v@lY,\v@lZ)}
\ctr@ld@f\def\Pslign@C#1{\s@uvc@ntr@l\et@tPslign@C\setc@ntr@l{2}%
    \def\list@num{#1}\extrairelepremi@r\p@int\de\list@num%
    \edef\premp@int{\p@int}\f@gnewpath%
    \expandafter\Pslign@C@\p@int:\PSwrit@cmd{-1}{\c@mmoveto}{\fwf@g}%
    \@ecfor\p@int:=\list@num\do{\expandafter\Pslign@C@\p@int:%
    \PSwrit@cmd{-1}{\c@mlineto}{\fwf@g}\edef\derp@int{\p@int}}%
    \resetc@ntr@l\et@tPslign@C}
\ctr@ld@f\def\Pslign@C@#1 #2:{{\def\t@xt@{#1}\ifx\t@xt@\empty\Pslign@C@#2:
    \else\extr@ctCF#1 #2:\fi}}
\ctr@ln@m\c@ntrolmesh
\ctr@ld@f\def\Pssetm@sh#1=#2|{\keln@mun#1|%
    \def\n@mref{d}\ifx\l@debut\n@mref\pssetmeshdiag{#2}\else
    \immediate\write16{*** Unknown attribute: \BS@ psset mesh(..., #1=...)}%
    \fi}
\ctr@ld@f\def\pssetmeshdiag#1{\edef\c@ntrolmesh{#1}}
\ctr@ld@f\def\defaultmeshdiag{0}    
\ctr@ld@f\def\psmesh#1,#2[#3,#4,#5,#6]{{\ifcurr@ntPS\ifps@cri%
    \PSc@mment{psmesh N1=#1, N2=#2, Quadrangle=[#3,#4,#5,#6]}%
    \s@uvc@ntr@l\et@tpsmesh\Pss@tsecondSt\setc@ntr@l{2}%
    \ifnum#1>\@ne\Psmeshp@rt#1[#3,#4,#5,#6]\fi%
    \ifnum#2>\@ne\Psmeshp@rt#2[#4,#5,#6,#3]\fi%
    \ifnum\c@ntrolmesh>\z@\Psmeshdi@g#1,#2[#3,#4,#5,#6]\fi%
    \ifnum\c@ntrolmesh<\z@\Psmeshdi@g#2,#1[#4,#5,#6,#3]\fi\Psrest@reSt%
    \psline[#3,#4,#5,#6,#3]\PSc@mment{End psmesh}\resetc@ntr@l\et@tpsmesh\fi\fi}}
\ctr@ld@f\def\Psmeshp@rt#1[#2,#3,#4,#5]{{\l@mbd@un=\@ne\l@mbd@de=#1\loop%
    \ifnum\l@mbd@un<#1\advance\l@mbd@de\m@ne\figptbary-1:[#2,#3;\l@mbd@de,\l@mbd@un]%
    \figptbary-2:[#5,#4;\l@mbd@de,\l@mbd@un]\psline[-1,-2]\advance\l@mbd@un\@ne\repeat}}
\ctr@ld@f\def\Psmeshdi@g#1,#2[#3,#4,#5,#6]{\figptcopy-2:/#3/\figptcopy-3:/#6/%
    \l@mbd@un=\z@\l@mbd@de=#1\loop\ifnum\l@mbd@un<#1%
    \advance\l@mbd@un\@ne\advance\l@mbd@de\m@ne\figptcopy-1:/-2/\figptcopy-4:/-3/%
    \figptbary-2:[#3,#4;\l@mbd@de,\l@mbd@un]%
    \figptbary-3:[#6,#5;\l@mbd@de,\l@mbd@un]\Psmeshdi@gp@rt#2[-1,-2,-3,-4]\repeat}
\ctr@ld@f\def\Psmeshdi@gp@rt#1[#2,#3,#4,#5]{{\l@mbd@un=\z@\l@mbd@de=#1\loop%
    \ifnum\l@mbd@un<#1\figptbary-5:[#2,#5;\l@mbd@de,\l@mbd@un]%
    \advance\l@mbd@de\m@ne\advance\l@mbd@un\@ne%
    \figptbary-6:[#3,#4;\l@mbd@de,\l@mbd@un]\psline[-5,-6]\repeat}}
\ctr@ln@m\psnormal
\ctr@ld@f\def\psnormalDD#1,#2[#3,#4]{{\ifcurr@ntPS\ifps@cri%
    \PSc@mment{psnormal Length=#1, Lambda=#2 [Pt1,Pt2]=[#3,#4]}%
    \s@uvc@ntr@l\et@tpsnormal\resetc@ntr@l{2}\figptendnormal-6::#1,#2[#3,#4]%
    \figptcopyDD-5:/-1/\psarrow[-5,-6]%
    \PSc@mment{End psnormal}\resetc@ntr@l\et@tpsnormal\fi\fi}}
\ctr@ld@f\def\psreset#1{\trtlis@rg{#1}{\Psreset@}}
\ctr@ld@f\def\Psreset@#1|{\keln@mde#1|%
    \def\n@mref{ar}\ifx\l@debut\n@mref\psresetarrowhead\else
    \def\n@mref{cu}\ifx\l@debut\n@mref\psset curve(roundness=\defaultroundness)\else
    \def\n@mref{fi}\ifx\l@debut\n@mref\psset (color=\defaultcolor,dash=\defaultdash,%
         fill=\defaultfill,join=\defaultjoin,width=\defaultwidth)\else
    \def\n@mref{fl}\ifx\l@debut\n@mref\psset flowchart(arrowp=\defaultfcarrowposition,%
 arrowr=\defaultfcarrowrefpt,line=\defaultfcline,xpadd=\defaultfcxpadding,%
 ypadd=\defaultfcypadding,radius=\defaultfcradius,shape=\defaultfcshape,%
 thick=\defaultfcthickness)\else
    \def\n@mref{me}\ifx\l@debut\n@mref\psset mesh(diag=\defaultmeshdiag)\else
    \def\n@mref{se}\ifx\l@debut\n@mref\psresetsecondsettings\else
    \def\n@mref{th}\ifx\l@debut\n@mref\psset third(color=\defaultthirdcolor)\else
    \immediate\write16{*** Unknown keyword #1 (\BS@ psreset).}%
    \fi\fi\fi\fi\fi\fi\fi}
\ctr@ld@f\def\psset#1(#2){\def\t@xt@{#1}\ifx\t@xt@\empty\trtlis@rg{#2}{\Pssetf@rst}
    \else\keln@mde#1|%
    \def\n@mref{ar}\ifx\l@debut\n@mref\trtlis@rg{#2}{\Psset@rrowhe@d}\else
    \def\n@mref{cu}\ifx\l@debut\n@mref\trtlis@rg{#2}{\Pssetc@rve}\else
    \def\n@mref{fi}\ifx\l@debut\n@mref\trtlis@rg{#2}{\Pssetf@rst}\else
    \def\n@mref{fl}\ifx\l@debut\n@mref\trtlis@rg{#2}{\Pssetfl@wchart}\else
    \def\n@mref{me}\ifx\l@debut\n@mref\trtlis@rg{#2}{\Pssetm@sh}\else
    \def\n@mref{se}\ifx\l@debut\n@mref\trtlis@rg{#2}{\Pssets@cond}\else
    \def\n@mref{th}\ifx\l@debut\n@mref\trtlis@rg{#2}{\Pssetth@rd}\else
    \immediate\write16{*** Unknown keyword: \BS@ psset #1(...)}%
    \fi\fi\fi\fi\fi\fi\fi\fi}
\ctr@ld@f\def\pssetdefault#1(#2){\ifcurr@ntPS\immediate\write16{*** \BS@ pssetdefault is ignored
    inside a \BS@ psbeginfig-\BS@ psendfig block.}%
    \immediate\write16{*** It must be called before \BS@ psbeginfig.}\else%
    \def\t@xt@{#1}\ifx\t@xt@\empty\trtlis@rg{#2}{\Pssd@f@rst}\else\keln@mde#1|%
    \def\n@mref{ar}\ifx\l@debut\n@mref\trtlis@rg{#2}{\Pssd@@rrowhe@d}\else
    \def\n@mref{cu}\ifx\l@debut\n@mref\trtlis@rg{#2}{\Pssd@c@rve}\else
    \def\n@mref{fi}\ifx\l@debut\n@mref\trtlis@rg{#2}{\Pssd@f@rst}\else
    \def\n@mref{fl}\ifx\l@debut\n@mref\trtlis@rg{#2}{\Pssd@fl@wchart}\else
    \def\n@mref{me}\ifx\l@debut\n@mref\trtlis@rg{#2}{\Pssd@m@sh}\else
    \def\n@mref{se}\ifx\l@debut\n@mref\trtlis@rg{#2}{\Pssd@s@cond}\else
    \def\n@mref{th}\ifx\l@debut\n@mref\trtlis@rg{#2}{\Pssd@th@rd}\else
    \immediate\write16{*** Unknown keyword: \BS@ pssetdefault #1(...)}%
    \fi\fi\fi\fi\fi\fi\fi\fi\initpss@ttings\fi}
\ctr@ld@f\def\Pssd@f@rst#1=#2|{\keln@mun#1|%
    \def\n@mref{c}\ifx\l@debut\n@mref\edef\defaultcolor{#2}\else
    \def\n@mref{d}\ifx\l@debut\n@mref\edef\defaultdash{#2}\else
    \def\n@mref{f}\ifx\l@debut\n@mref\edef\defaultfill{#2}\else
    \def\n@mref{j}\ifx\l@debut\n@mref\edef\defaultjoin{#2}\else
    \def\n@mref{u}\ifx\l@debut\n@mref\edef\defaultupdate{#2}\pssetupdate{#2}\else
    \def\n@mref{w}\ifx\l@debut\n@mref\edef\defaultwidth{#2}\else
    \immediate\write16{*** Unknown attribute: \BS@ pssetdefault (..., #1=...)}%
    \fi\fi\fi\fi\fi\fi}
\ctr@ld@f\def\Pssd@@rrowhe@d#1=#2|{\keln@mun#1|%
    \def\n@mref{a}\ifx\l@debut\n@mref\edef\defaultarrowheadangle{#2}\else
    \def\n@mref{f}\ifx\l@debut\n@mref\edef\defaultarrowheadangle{#2}\else
    \def\n@mref{l}\ifx\l@debut\n@mref\y@tiunit{#2}\ifunitpr@sent%
     \edef\defaulth@rdahlength{#2}\else\edef\defaulth@rdahlength{#2pt}%
     \message{*** \BS@ pssetdefault (..., #1=#2, ...) : unit is missing, pt is assumed.}%
     \fi\else
    \def\n@mref{o}\ifx\l@debut\n@mref\edef\defaultarrowheadout{#2}\else
    \def\n@mref{r}\ifx\l@debut\n@mref\edef\defaultarrowheadratio{#2}\else
    \immediate\write16{*** Unknown attribute: \BS@ pssetdefault arrowhead(..., #1=...)}%
    \fi\fi\fi\fi\fi}
\ctr@ld@f\def\Pssd@c@rve#1=#2|{\keln@mun#1|%
    \def\n@mref{r}\ifx\l@debut\n@mref\edef\defaultroundness{#2}\else%
    \immediate\write16{*** Unknown attribute: \BS@ pssetdefault curve(..., #1=...)}%
    \fi}
\ctr@ld@f\def\Pssd@fl@wchart#1=#2|{\keln@mtr#1|%
    \def\n@mref{arr}\ifx\l@debut\n@mref\expandafter\keln@mtr\l@suite|%
     \def\n@mref{owp}\ifx\l@debut\n@mref\edef\defaultfcarrowposition{#2}\else
     \def\n@mref{owr}\ifx\l@debut\n@mref\edef\defaultfcarrowrefpt{#2}\else
     \immediate\write16{*** Unknown attribute: \BS@ pssetdefault flowchart(..., #1=...)}%
     \fi\fi\else%
    \def\n@mref{lin}\ifx\l@debut\n@mref\edef\defaultfcline{#2}\else
    \def\n@mref{pad}\ifx\l@debut\n@mref\edef\defaultfcxpadding{#2}%
                    \edef\defaultfcypadding{#2}\else
    \def\n@mref{rad}\ifx\l@debut\n@mref\edef\defaultfcradius{#2}\else
    \def\n@mref{sha}\ifx\l@debut\n@mref\edef\defaultfcshape{#2}\else
    \def\n@mref{thi}\ifx\l@debut\n@mref\edef\defaultfcthickness{#2}\else
    \def\n@mref{xpa}\ifx\l@debut\n@mref\edef\defaultfcxpadding{#2}\else
    \def\n@mref{ypa}\ifx\l@debut\n@mref\edef\defaultfcypadding{#2}\else
    \immediate\write16{*** Unknown attribute: \BS@ pssetdefault flowchart(..., #1=...)}%
    \fi\fi\fi\fi\fi\fi\fi\fi}
\ctr@ld@f\def\defaultfcarrowposition{0.5}
\ctr@ld@f\def\defaultfcarrowrefpt{start}
\ctr@ld@f\def\defaultfcline{polygon}
\ctr@ld@f\def\defaultfcradius{0}
\ctr@ld@f\def\defaultfcshape{rectangle}
\ctr@ld@f\def\defaultfcthickness{0}
\ctr@ld@f\def\defaultfcxpadding{0}
\ctr@ld@f\def\defaultfcypadding{0}
\ctr@ld@f\def\Pssd@m@sh#1=#2|{\keln@mun#1|%
    \def\n@mref{d}\ifx\l@debut\n@mref\edef\defaultmeshdiag{#2}\else%
    \immediate\write16{*** Unknown attribute: \BS@ pssetdefault mesh(..., #1=...)}%
    \fi}
\ctr@ld@f\def\Pssd@s@cond#1=#2|{\keln@mun#1|%
    \def\n@mref{c}\ifx\l@debut\n@mref\edef\defaultsecondcolor{#2}\else%
    \def\n@mref{d}\ifx\l@debut\n@mref\edef\defaultseconddash{#2}\else%
    \def\n@mref{w}\ifx\l@debut\n@mref\edef\defaultsecondwidth{#2}\else%
    \immediate\write16{*** Unknown attribute: \BS@ pssetdefault second(..., #1=...)}%
    \fi\fi\fi}
\ctr@ld@f\def\Pssd@th@rd#1=#2|{\keln@mun#1|%
    \def\n@mref{c}\ifx\l@debut\n@mref\edef\defaultthirdcolor{#2}\else%
    \immediate\write16{*** Unknown attribute: \BS@ pssetdefault third(..., #1=...)}%
    \fi}
\ctr@ln@w{newif}\iffillm@de
\ctr@ld@f\def\pssetfillmode#1{\expandafter\setfillm@de#1:}
\ctr@ld@f\def\setfillm@de#1#2:{\if#1n\fillm@defalse\else\fillm@detrue\fi}
\ctr@ld@f\def\defaultfill{no}     
\ctr@ln@w{newif}\ifpsupdatem@de
\ctr@ld@f\def\pssetupdate#1{\ifcurr@ntPS\immediate\write16{*** \BS@ pssetupdate is ignored inside a
     \BS@ psbeginfig-\BS@ psendfig block.}%
    \immediate\write16{*** It must be called before \BS@ psbeginfig.}%
    \else\expandafter\setupd@te#1:\fi}
\ctr@ld@f\def\setupd@te#1#2:{\if#1n\psupdatem@defalse\else\psupdatem@detrue\fi}
\ctr@ld@f\def\defaultupdate{no}     
\ctr@ln@m\curr@ntcolor \ctr@ln@m\curr@ntcolorc@md
\ctr@ld@f\def\Pssetc@lor#1{\ifps@cri\result@tent=\@ne\expandafter\c@lnbV@l#1 :%
    \def\curr@ntcolor{}\def\curr@ntcolorc@md{}%
    \ifcase\result@tent\or\pssetgray{#1}\or\or\pssetrgb{#1}\or\pssetcmyk{#1}\fi\fi}
\ctr@ln@m\curr@ntcolorc@mdStroke
\ctr@ld@f\def\pssetcmyk#1{\ifps@cri\def\curr@ntcolor{#1}\def\curr@ntcolorc@md{\c@msetcmykcolor}%
    \def\curr@ntcolorc@mdStroke{\c@msetcmykcolorStroke}%
    \ifcurr@ntPS\PSc@mment{pssetcmyk Color=#1}\us@primarC@lor\fi\fi}
\ctr@ld@f\def\pssetrgb#1{\ifps@cri\def\curr@ntcolor{#1}\def\curr@ntcolorc@md{\c@msetrgbcolor}%
    \def\curr@ntcolorc@mdStroke{\c@msetrgbcolorStroke}%
    \ifcurr@ntPS\PSc@mment{pssetrgb Color=#1}\us@primarC@lor\fi\fi}
\ctr@ld@f\def\pssetgray#1{\ifps@cri\def\curr@ntcolor{#1}\def\curr@ntcolorc@md{\c@msetgray}%
    \def\curr@ntcolorc@mdStroke{\c@msetgrayStroke}%
    \ifcurr@ntPS\PSc@mment{pssetgray Gray level=#1}\us@primarC@lor\fi\fi}
\ctr@ln@m\fillc@md
\ctr@ld@f\def\us@primarC@lor{\immediate\write\fwf@g{\d@fprimarC@lor}%
    \let\fillc@md=\prfillc@md}
\ctr@ld@f\def\prfillc@md{\d@fprimarC@lor\space\c@mfill}
\ctr@ld@f\def\defaultcolor{0}       
\ctr@ld@f\def\c@lnbV@l#1 #2:{\def\t@xt@{#1}\relax\ifx\t@xt@\empty\c@lnbV@l#2:
    \else\c@lnbV@l@#1 #2:\fi}
\ctr@ld@f\def\c@lnbV@l@#1 #2:{\def\t@xt@{#2}\ifx\t@xt@\empty%
    \def\t@xt@{#1}\ifx\t@xt@\empty\advance\result@tent\m@ne\fi
    \else\advance\result@tent\@ne\c@lnbV@l@#2:\fi}
\ctr@ld@f\def\Blackcmyk{0 0 0 1}
\ctr@ld@f\def\Whitecmyk{0 0 0 0}
\ctr@ld@f\def\Cyancmyk{1 0 0 0}
\ctr@ld@f\def\Magentacmyk{0 1 0 0}
\ctr@ld@f\def\Yellowcmyk{0 0 1 0}
\ctr@ld@f\def\Redcmyk{0 1 1 0}
\ctr@ld@f\def\Greencmyk{1 0 1 0}
\ctr@ld@f\def\Bluecmyk{1 1 0 0}
\ctr@ld@f\def\Graycmyk{0 0 0 0.50}
\ctr@ld@f\def\BrickRedcmyk{0 0.89 0.94 0.28} 
\ctr@ld@f\def\Browncmyk{0 0.81 1 0.60} 
\ctr@ld@f\def\ForestGreencmyk{0.91 0 0.88 0.12} 
\ctr@ld@f\def\Goldenrodcmyk{ 0 0.10 0.84 0} 
\ctr@ld@f\def\Marooncmyk{0 0.87 0.68 0.32} 
\ctr@ld@f\def\Orangecmyk{0 0.61 0.87 0} 
\ctr@ld@f\def\Purplecmyk{0.45 0.86 0 0} 
\ctr@ld@f\def\RoyalBluecmyk{1. 0.50 0 0} 
\ctr@ld@f\def\Violetcmyk{0.79 0.88 0 0} 
\ctr@ld@f\def\Blackrgb{0 0 0}
\ctr@ld@f\def\Whitergb{1 1 1}
\ctr@ld@f\def\Redrgb{1 0 0}
\ctr@ld@f\def\Greenrgb{0 1 0}
\ctr@ld@f\def\Bluergb{0 0 1}
\ctr@ld@f\def\Cyanrgb{0 1 1}
\ctr@ld@f\def\Magentargb{1 0 1}
\ctr@ld@f\def\Yellowrgb{1 1 0}
\ctr@ld@f\def\Grayrgb{0.5 0.5 0.5}
\ctr@ld@f\def\Chocolatergb{0.824 0.412 0.118}
\ctr@ld@f\def\DarkGoldenrodrgb{0.722 0.525 0.043}
\ctr@ld@f\def\DarkOrangergb{1 0.549 0}
\ctr@ld@f\def\Firebrickrgb{0.698 0.133 0.133}
\ctr@ld@f\def\ForestGreenrgb{0.133 0.545 0.133}
\ctr@ld@f\def\Goldrgb{1 0.843 0}
\ctr@ld@f\def\HotPinkrgb{1 0.412 0.706}
\ctr@ld@f\def\Maroonrgb{0.690 0.188 0.376}
\ctr@ld@f\def\Pinkrgb{1 0.753 0.796}
\ctr@ld@f\def\RoyalBluergb{0.255 0.412 0.882}
\ctr@ld@f\def\Pssetf@rst#1=#2|{\keln@mun#1|%
    \def\n@mref{c}\ifx\l@debut\n@mref\Pssetc@lor{#2}\else
    \def\n@mref{d}\ifx\l@debut\n@mref\pssetdash{#2}\else
    \def\n@mref{f}\ifx\l@debut\n@mref\pssetfillmode{#2}\else
    \def\n@mref{j}\ifx\l@debut\n@mref\pssetjoin{#2}\else
    \def\n@mref{u}\ifx\l@debut\n@mref\pssetupdate{#2}\else
    \def\n@mref{w}\ifx\l@debut\n@mref\pssetwidth{#2}\else
    \immediate\write16{*** Unknown attribute: \BS@ psset (..., #1=...)}%
    \fi\fi\fi\fi\fi\fi}
\ctr@ln@m\curr@ntdash
\ctr@ld@f\def\s@uvdash#1{\edef#1{\curr@ntdash}}
\ctr@ld@f\def\defaultdash{1}        
\ctr@ld@f\def\pssetdash#1{\ifps@cri\edef\curr@ntdash{#1}\ifcurr@ntPS\expandafter\Pssetd@sh#1 :\fi\fi}
\ctr@ld@f\def\Pssetd@shI#1{\PSc@mment{pssetdash Index=#1}\ifcase#1%
    \or\immediate\write\fwf@g{[] 0 \c@msetdash}
    \or\immediate\write\fwf@g{[6 2] 0 \c@msetdash}
    \or\immediate\write\fwf@g{[4 2] 0 \c@msetdash}
    \or\immediate\write\fwf@g{[2 2] 0 \c@msetdash}
    \or\immediate\write\fwf@g{[1 2] 0 \c@msetdash}
    \or\immediate\write\fwf@g{[2 4] 0 \c@msetdash}
    \or\immediate\write\fwf@g{[3 5] 0 \c@msetdash}
    \or\immediate\write\fwf@g{[3 3] 0 \c@msetdash}
    \or\immediate\write\fwf@g{[3 5 1 5] 0 \c@msetdash}
    \or\immediate\write\fwf@g{[6 4 2 4] 0 \c@msetdash}
    \fi}
\ctr@ld@f\def\Pssetd@sh#1 #2:{{\def\t@xt@{#1}\ifx\t@xt@\empty\Pssetd@sh#2:
    \else\def\t@xt@{#2}\ifx\t@xt@\empty\Pssetd@shI{#1}\else\s@mme=\@ne\def\debutp@t{#1}%
    \an@lysd@sh#2:\ifodd\s@mme\edef\debutp@t{\debutp@t\space\finp@t}\def\finp@t{0}\fi%
    \PSc@mment{pssetdash Pattern=#1 #2}%
    \immediate\write\fwf@g{[\debutp@t] \finp@t\space\c@msetdash}\fi\fi}}
\ctr@ld@f\def\an@lysd@sh#1 #2:{\def\t@xt@{#2}\ifx\t@xt@\empty\def\finp@t{#1}\else%
    \edef\debutp@t{\debutp@t\space#1}\advance\s@mme\@ne\an@lysd@sh#2:\fi}
\ctr@ln@m\curr@ntwidth
\ctr@ld@f\def\s@uvwidth#1{\edef#1{\curr@ntwidth}}
\ctr@ld@f\def\defaultwidth{0.4}     
\ctr@ld@f\def\pssetwidth#1{\ifps@cri\edef\curr@ntwidth{#1}\ifcurr@ntPS%
    \PSc@mment{pssetwidth Width=#1}\immediate\write\fwf@g{#1 \c@msetlinewidth}\fi\fi}
\ctr@ln@m\curr@ntjoin
\ctr@ld@f\def\pssetjoin#1{\ifps@cri\edef\curr@ntjoin{#1}\ifcurr@ntPS\expandafter\Pssetj@in#1:\fi\fi}
\ctr@ld@f\def\Pssetj@in#1#2:{\PSc@mment{pssetjoin join=#1}%
    \if#1r\def\t@xt@{1}\else\if#1b\def\t@xt@{2}\else\def\t@xt@{0}\fi\fi%
    \immediate\write\fwf@g{\t@xt@\space\c@msetlinejoin}}
\ctr@ld@f\def\defaultjoin{miter}   
\ctr@ld@f\def\Pssets@cond#1=#2|{\keln@mun#1|%
    \def\n@mref{c}\ifx\l@debut\n@mref\Pssets@condcolor{#2}\else%
    \def\n@mref{d}\ifx\l@debut\n@mref\pssetseconddash{#2}\else%
    \def\n@mref{w}\ifx\l@debut\n@mref\pssetsecondwidth{#2}\else%
    \immediate\write16{*** Unknown attribute: \BS@ psset second(..., #1=...)}%
    \fi\fi\fi}
\ctr@ln@m\curr@ntseconddash
\ctr@ld@f\def\pssetseconddash#1{\edef\curr@ntseconddash{#1}}
\ctr@ld@f\def\defaultseconddash{4}  
\ctr@ln@m\curr@ntsecondwidth
\ctr@ld@f\def\pssetsecondwidth#1{\edef\curr@ntsecondwidth{#1}}
\ctr@ld@f\edef\defaultsecondwidth{\defaultwidth} 
\ctr@ld@f\def\psresetsecondsettings{%
    \pssetseconddash{\defaultseconddash}\pssetsecondwidth{\defaultsecondwidth}%
    \Pssets@condcolor{\defaultsecondcolor}}
\ctr@ln@m\sec@ndcolor \ctr@ln@m\sec@ndcolorc@md
\ctr@ld@f\def\Pssets@condcolor#1{\ifps@cri\result@tent=\@ne\expandafter\c@lnbV@l#1 :%
    \def\sec@ndcolor{}\def\sec@ndcolorc@md{}%
    \ifcase\result@tent\or\pssetsecondgray{#1}\or\or\pssetsecondrgb{#1}%
    \or\pssetsecondcmyk{#1}\fi\fi}
\ctr@ln@m\sec@ndcolorc@mdStroke
\ctr@ld@f\def\pssetsecondcmyk#1{\def\sec@ndcolor{#1}\def\sec@ndcolorc@md{\c@msetcmykcolor}%
    \def\sec@ndcolorc@mdStroke{\c@msetcmykcolorStroke}}
\ctr@ld@f\def\pssetsecondrgb#1{\def\sec@ndcolor{#1}\def\sec@ndcolorc@md{\c@msetrgbcolor}%
    \def\sec@ndcolorc@mdStroke{\c@msetrgbcolorStroke}}
\ctr@ld@f\def\pssetsecondgray#1{\def\sec@ndcolor{#1}\def\sec@ndcolorc@md{\c@msetgray}%
    \def\sec@ndcolorc@mdStroke{\c@msetgrayStroke}}
\ctr@ld@f\def\us@secondC@lor{\immediate\write\fwf@g{\d@fsecondC@lor}%
    \let\fillc@md=\sdfillc@md}
\ctr@ld@f\def\sdfillc@md{\d@fsecondC@lor\space\c@mfill}
\ctr@ld@f\edef\defaultsecondcolor{\defaultcolor} 
\ctr@ld@f\def\Pss@tsecondSt{%
    \s@uvdash{\typ@dash}\pssetdash{\curr@ntseconddash}%
    \s@uvwidth{\typ@width}\pssetwidth{\curr@ntsecondwidth}\us@secondC@lor}
\ctr@ld@f\def\Psrest@reSt{\pssetwidth{\typ@width}\pssetdash{\typ@dash}\us@primarC@lor}
\ctr@ld@f\def\Pssetth@rd#1=#2|{\keln@mun#1|%
    \def\n@mref{c}\ifx\l@debut\n@mref\Pssetth@rdcolor{#2}\else%
    \immediate\write16{*** Unknown attribute: \BS@ psset third(..., #1=...)}%
    \fi}
\ctr@ln@m\th@rdcolor \ctr@ln@m\th@rdcolorc@md
\ctr@ld@f\def\Pssetth@rdcolor#1{\ifps@cri\result@tent=\@ne\expandafter\c@lnbV@l#1 :%
    \def\th@rdcolor{}\def\th@rdcolorc@md{}%
    \ifcase\result@tent\or\Pssetth@rdgray{#1}\or\or\Pssetth@rdrgb{#1}%
    \or\Pssetth@rdcmyk{#1}\fi\fi}
\ctr@ln@m\th@rdcolorc@mdStroke
\ctr@ld@f\def\Pssetth@rdcmyk#1{\def\th@rdcolor{#1}\def\th@rdcolorc@md{\c@msetcmykcolor}%
    \def\th@rdcolorc@mdStroke{\c@msetcmykcolorStroke}}
\ctr@ld@f\def\Pssetth@rdrgb#1{\def\th@rdcolor{#1}\def\th@rdcolorc@md{\c@msetrgbcolor}%
    \def\th@rdcolorc@mdStroke{\c@msetrgbcolorStroke}}
\ctr@ld@f\def\Pssetth@rdgray#1{\def\th@rdcolor{#1}\def\th@rdcolorc@md{\c@msetgray}%
    \def\th@rdcolorc@mdStroke{\c@msetgrayStroke}}
\ctr@ld@f\def\us@thirdC@lor{\immediate\write\fwf@g{\d@fthirdC@lor}%
    \let\fillc@md=\thfillc@md}
\ctr@ld@f\def\thfillc@md{\d@fthirdC@lor\space\c@mfill}
\ctr@ld@f\def\defaultthirdcolor{1}  
\ctr@ld@f\def\pstrimesh#1[#2,#3,#4]{{\ifcurr@ntPS\ifps@cri%
    \PSc@mment{pstrimesh Type=#1, Triangle=[#2,#3,#4]}%
    \s@uvc@ntr@l\et@tpstrimesh\ifnum#1>\@ne\Pss@tsecondSt\setc@ntr@l{2}%
    \Pstrimeshp@rt#1[#2,#3,#4]\Pstrimeshp@rt#1[#3,#4,#2]%
    \Pstrimeshp@rt#1[#4,#2,#3]\Psrest@reSt\fi\psline[#2,#3,#4,#2]%
    \PSc@mment{End pstrimesh}\resetc@ntr@l\et@tpstrimesh\fi\fi}}
\ctr@ld@f\def\Pstrimeshp@rt#1[#2,#3,#4]{{\l@mbd@un=\@ne\l@mbd@de=#1\loop\ifnum\l@mbd@de>\@ne%
    \advance\l@mbd@de\m@ne\figptbary-1:[#2,#3;\l@mbd@de,\l@mbd@un]%
    \figptbary-2:[#2,#4;\l@mbd@de,\l@mbd@un]\psline[-1,-2]%
    \advance\l@mbd@un\@ne\repeat}}
\initpr@lim\initpss@ttings\initPDF@rDVI
\ctr@ln@w{newbox}\figBoxA
\ctr@ln@w{newbox}\figBoxB
\ctr@ln@w{newbox}\figBoxC
\catcode`\@=12